\documentclass[openany]{memo-l}

\usepackage{amsfonts, amssymb}
\usepackage{graphicx}

\usepackage{color} 
   \definecolor{cites}{rgb}{0.50 , 0.00 , 0.00}  
   \definecolor{urls} {rgb}{0.00 , 0.00 , 0.50}  
   \definecolor{links}{rgb}{0.00 , 0.00 , 0.50}   

\usepackage{hyperref} 
\hypersetup{  
      colorlinks=true,   
      citecolor=cites,   
      urlcolor=urls,     
      linkcolor=links,   
      pdftitle={Limit Operators, Collective Compactness, and the Spectral Theory of Infinite Matrices},  
      pdfauthor={Simon N. Chandler-Wilde and Marko Lindner}, 
}

\def\birthday{67} 

\newcommand\C{{\mathbb C}}
\newcommand\Z{{\mathbb Z}}
\newcommand\N{{\mathbb N}}
\newcommand\R{{\mathbb R}}

\newcommand\B{{\mathcal B}}

\newcommand\J{{\mathcal J}}
\newcommand\K{{\mathcal K}}

\renewcommand\P{{\mathcal P}}
\newcommand\Sh{{\mathcal T}}

\newcommand\V{{\mathcal V}}
\newcommand\Vt{\tilde{\mathcal V}}
\newcommand\W{{\mathcal W}}

\newcommand\cA{{\mathcal A}}
\newcommand\cAAP{{\cA_{\sf AP}}}
\newcommand\JAP{{\J_{\sf AP}}}
\newcommand\M{{\sf M}}    
\newcommand\Co{{\sf C}}   
\newcommand\Ko{{\sf K}}   
\newcommand\Lo{{\sf L}}   
\newcommand\Io{{\sf I}}   

\renewcommand\l{\ell}
\newcommand\eps{\varepsilon}
\newcommand\ph{\varphi}
\newcommand\dist{{\rm dist\,}}

\newcommand\toto{\rightrightarrows}
\newcommand\pto{\stackrel\P\to}
\newcommand\sto{\stackrel s\to}
\newcommand\Sto{\stackrel S\to}

\newcommand\NR{{L^{\sf n\$}}}          
\newcommand\AP{{Y^\infty_{\sf AP}}}   
\newcommand\LAP{{L^\infty_{\sf AP}}}   
\newcommand\UM{{UM}}                  
\newcommand\rich{{L^\infty_\$}}        
\newcommand\opsp{\sigma^{\sf op}}
\newcommand\plim{\P\hspace{-1mm}-\hspace{-1mm}\lim}
\newcommand\clos{{\rm clos}}

\renewcommand\mod{{\rm mod}}
\newcommand\ind{{\rm ind}}
\renewcommand\sp{{\rm spec}}
\newcommand\speps{{\rm spec}_\eps}
\newcommand\spess{{\rm spec}_{\rm ess}}
\newcommand\spp{{\rm spec}_{\rm point}}
\newcommand\Lim{{\rm Lim}\,}
\newcommand\pa{\triangleleft}

\newcommand{\refeq}[1]{(\ref{#1})}       
\newcommand{\esssup}{\mathop{\rm ess\,sup}}
\newcommand{\qtext}[1]{\quad\textrm{#1}\quad}
\newcommand{\qqtext}[1]{\qquad\textrm{#1}\qquad}

\newcommand\Proof{\noindent {\it Proof. }}
\newcommand\Proofend{\rule{2mm}{2mm}}
\newcommand\Remarkend{\raisebox{1mm}{\framebox{}}}

\newtheorem{Theorem}{Theorem}[chapter]
\newtheorem{Lemma}[Theorem]{Lemma}
\newtheorem{Corollary}[Theorem]{Corollary}

\newtheorem{Definition}[Theorem]{Definition}

\newenvironment{Remark}
 {\par\noindent\refstepcounter{Theorem}{\indent \sc Remark \theTheorem}\ }
 {\Remarkend\pagebreak[2]}

\newenvironment{Example}
 {\par\noindent\refstepcounter{Theorem}{\indent \sc Example \theTheorem}\ }
 {\Remarkend\pagebreak[2]}

\numberwithin{section}{chapter}
\numberwithin{equation}{chapter}
\numberwithin{figure}{chapter}

\makeindex

\begin{document}

\frontmatter

\title{Limit Operators, Collective Compactness,\\
and the Spectral Theory of Infinite Matrices}


\author{Simon N.\ Chandler-Wilde}
\address{Department of Mathematics, University of Reading, Whiteknights, \indent PO Box
220, Reading RG6 6AX, United Kingdom} \curraddr{}
\email{S.N.Chandler-Wilde@reading.ac.uk}

\author{Marko Lindner}
\address{Fakult\"at Mathematik, TU Chemnitz, D-09107 Chemnitz, Germany} \curraddr{}
\email{Marko.Lindner@mathematik.tu-chemnitz.de}

\date{}




\begin{abstract}
In the first half of this text we explore the interrelationships
between the abstract theory of limit operators (see e.g.\ the recent
monographs of Rabinovich, Roch and Silbermann \cite{RaRoSiBook} and
Lindner \cite{LiBook}) and the concepts and results of the
generalised collectively compact operator theory introduced by
Chandler-Wilde and Zhang \cite{ChanZh02}. We build up to results
obtained by applying this generalised collectively compact operator
theory to the set of limit operators of an operator $A$ (its
operator spectrum).  In the second half of this text we study
bounded linear operators on the generalised sequence space
$\ell^p(\Z^N,U)$, where $p\in [1,\infty]$ and $U$ is some complex
Banach space. We make what seems to be a more complete study than
hitherto of the connections between Fredholmness, invertibility,
invertibility at infinity, and invertibility or injectivity of the
set of limit operators, with some emphasis on the case when the
operator $A$ is a locally compact perturbation of the identity.
Especially, we obtain stronger results than previously known for the
subtle limiting cases of $p=1$ and $\infty$. Our tools in this study
are the results from the first half of the text and an
exploitation of the partial duality between $\ell^1$ and
$\ell^\infty$ and its implications for bounded linear operators
which are also continuous with respect to the weaker topology (the
strict topology) introduced in the first half of the text. Results
in this second half of the text include a new proof that
injectivity of all limit operators (the classic Favard condition)
implies invertibility for a general class of almost periodic
operators, and characterisations of invertibility at infinity and
Fredholmness for operators in the so-called Wiener algebra. In two
final chapters our results are illustrated by and applied to
concrete examples. Firstly, we study the spectra and essential
spectra of discrete Schr\"odinger operators (both self-adjoint and
non-self-adjoint), including operators with almost periodic and
random potentials. In the final chapter we apply our results to
integral operators on $\R^N$.
\end{abstract}

\maketitle

\cleardoublepage \thispagestyle{empty} \vspace*{13.5pc}
\begin{center}
{\Large \sl Dedicated to Professor Bernd Silbermann\\[1ex]
on the occasion of his \birthday th birthday.}
\end{center}
\cleardoublepage

\setcounter{page}{7}

\tableofcontents


\mainmatter


\chapter{Introduction} \label{chap_intro}

\section{Overview} \label{sec_overview}
This text develops an abstract theory of limit operators and a
generalised collectively compact operator theory which can be used
separately or together to obtain information on the location in the
complex plane of the spectrum, essential spectrum, and
pseudospectrum for large classes of linear operators arising in
applications. We have in mind here differential, integral,
pseudo-differential, difference, and pseudo-difference operators, in
particular operators of all these types on unbounded domains. This
text also illustrates this general theory by developing, in a more
complete form than hitherto, a theory of the limit operator method
in one of its most concrete forms, as it applies to bounded linear
operators on spaces of sequences, where each component of the
sequence takes values in some Banach space. Finally, we apply this
concrete form of the theory to the analysis of lattice Schr\"odinger
operators and to the study of integral operators on $\R^N$.

Let us give an idea of the methods and results that we will develop
and the problems that they enable us to study. Let
$Y=\ell^p=\ell^p(\Z,\C)$, for $1\leq p\leq \infty$, denote the usual
Banach space of complex-valued bilateral sequences
$x=(x(m))_{m\in\Z}$ for which the norm $\|x\|$ is finite; here
$\|x\|:= \sup_m |x(m)|$, in the case $p=\infty$, while
$\|x\|:=(\sum_{m\in\Z}|x(m)|^p)^{1/p}$ for $1\leq p<\infty$. Let
$L(Y)$ \index{$L(Y)$} denote the space of bounded linear operators
on $Y$, and suppose $A\in L(Y)$ is given by the rule
\begin{equation} \label{eqn_first}
Ax(m) = \sum_{n\in\Z} a_{mn}x(n), \quad m\in\Z,
\end{equation}
for some coefficients $a_{mn}\in\C$ which we think of as elements of
an infinite matrix $[A]=[a_{mn}]_{m,n\in\Z}$ \index{$A_n$@$[A]$,
matrix representation of $A$} associated with the operator $A$. Of
course $A$, given by (\ref{eqn_first}), is only a bounded operator
on $Y$ under certain constraints on the entries $a_{mn}$. Simple
conditions that are sufficient to guarantee that $A\in L(Y)$, for
$1\leq p\leq \infty$, are to require that the entries are uniformly
bounded, i.e.\
\begin{equation} \label{eqn_second}
\sup_{m,n} |a_{mn}| <\infty,
\end{equation}
and to require that, for some $w\geq 0$, $a_{mn}=0$ if $|m-n|>w$. If
these  conditions hold we say that $[A]$ is a band matrix\index{band
matrix} with band-width $w$\index{band width} and that $A$ is a band
operator\index{band operator}\index{operator!band}. We note that the
tri-diagonal case $w=1$, when $A$ is termed a Jacobi
operator\index{Jacobi operator}\index{operator!Jacobi}, is
much-studied in the mathematical physics literature (e.g.\
\cite{Teschl99,LastSimon06}). This class includes, in particular,
the one-dimensional discrete Schr\"odinger
operator\index{Schr\"odinger operator}\index{operator!Schr\"odinger}
for which $a_{mn}=1$ for $|m-n|=1$.

It is well known (see Lemma \ref{Wiener1} below and the surrounding
remarks) that, under these conditions on $[A]$ (that $[A]$ is a band
matrix and (\ref{eqn_second}) holds), the spectrum\index{spectrum}
of $A$, i.e.\ the set of $\lambda\in\C$ for which $\lambda I-A$ is
not invertible as a member of the algebra $L(Y)$, is independent of
$p$. One of our main results in Section \ref{subsection_wiener}
implies that also the essential spectrum\index{essential
spectrum}\index{spectrum!essential} of $A$ (by which we mean the set
of $\lambda$ for which $\lambda I-A$ is not a Fredholm
operator\footnote{Throughout we will say that a bounded linear
operator $C$ from Banach space $X$ to Banach space $Y$ is: {\em
normally solvable} if its range $C(X)$ is closed; {\em
semi-Fredholm}\index{semi-Fredholm
operator}\index{operator!Fredholm!semi-}\index{Fredholm
operator!semi-} if, additionally, either $\alpha(C):=\dim(\ker C)$
or $\beta(C):=\dim(Y/C(X))$ are finite; a $\Phi_+$
operator\index{Fz+@$\Phi_+$ operator}\index{operator!$\Phi_+$} if it
is a semi-Fredholm operator with $\alpha<\infty$, and a $\Phi_-$
operator\index{Fz-@$\Phi_-$ operator}\index{operator!$\Phi_-$} if it
is semi-Fredholm with $\beta<\infty$; {\em Fredholm}\index{Fredholm
operator}\index{operator!Fredholm} if it is semi-Fredholm and both
$\alpha$ and $\beta$ are finite. If $C$ is semi-Fredholm then
$\alpha(C)-\beta(C)$ is called the {\em
index}\index{index}\index{Fredholm index} of
$C$.})\index{$\alpha(\cdot)$}\index{$\beta(\cdot)$} is independent
of $p$. Moreover, we prove that the essential spectrum is determined
by the behaviour of $A$ at infinity in the following precise sense.

Let $h=(h(j))_{j\in\N}\subset \Z$ be a sequence tending to infinity
for which it holds that $a_{m+h(j),n+h(j)}$ approaches a limit
$\tilde a_{m,n}$ for every $m,n\in\Z$. (The existence of many such
sequences is ensured by the theorem of Bolzano-Weierstrass and a
diagonal argument.) Then we call the operator $A_h$\index{$A_h$,
limit operator of $A$}, with matrix $[A_h]=[\tilde a_{mn}]$, a {\em
limit operator}\index{limit operator}\index{operator!limit} of the
operator $A$. Moreover, following e.g.\ \cite{RaRoSiBook}, we call
the set of all limit operators of $A$ the {\em operator
spectrum}\index{operator spectrum}\index{spectrum!operator} of $A$,
which we denote by $\opsp(A)$\index{$\opsp(A)$}. In terms of these
definitions our results imply that the essential spectrum of $A$
(which is independent of $p\in[1,\infty]$) is the union of the
spectra of the elements $A_h$ of the operator spectrum of $A$
(again, each of these spectra is independent of $p$). Moreover, this
is also precisely the union of the point spectra\index{point
spectrum}\index{spectrum!point} (sets of eigenvalues) of the limit
operators $A_h$ in the case $p=\infty$, in symbols
\begin{equation} \label{eqn_ess}
\spess(A) = \cup_{A_h\in\opsp(A)} \{\lambda: A_h x = \lambda x
\mbox{ has a bounded solution } x\ne 0\}.
\end{equation}
This formula and other related results have implications for the
spectrum of $A$. In particular, if it happens that $A\in \opsp(A)$
(we call $A$ {\em self-similar}\index{self-similar
operator}\index{operator!self-similar} in that case), then it holds
that
\begin{equation} \label{eqn_ess2}
\sp(A) = \spess(A) = \cup_{A_h\in\opsp(A)} \{\lambda: A_h x =
\lambda x \mbox{ has a bounded solution}\}.
\end{equation}
In the case $A\not\in \opsp(A)$ we do not have such a precise
characterisation, but if we construct $B\in L(Y)$ such that
$A\in\opsp(B)$ (see e.g.\ \cite[\S3.8.2]{LiBook} for how to do
this), then it holds that
\begin{equation} \label{eqn_ess3}
\sp(A) \subset \spess(B) = \cup_{B_h\in\opsp(B)} \{\lambda: B_h x =
\lambda x \mbox{ has a bounded solution}\}.
\end{equation}

A main aim of this text is to prove results of the above type
which apply in the simple setting just outlined, but also in the
more general setting where
$Y=\ell^p(\Z^N,U)$\index{$\ell^p(\Z^N,U)$} is a space of generalised
sequences $x=(x(m))_{m\in\Z^N}$, for some $N\in\N$, taking values in
some Banach space $U$. In this general setting the definition
(\ref{eqn_first}) makes sense if we replace $\Z$ by $\Z^N$ and
understand each matrix entry $a_{mn}$ as an element of $L(U)$. Such
results are the concern of Chapter \ref{sec_ellp}, and are applied
to discrete Schr\"odinger operators and to integral operators on
$\R^N$ in the final two chapters.

This integral operator application in Chapter \ref{sec_intop} will
illustrate how operators on $\R^N$ can be studied via {\em
discretisation}\index{discretisation}. To see how this simple idea
works in the case $N=1$,  let $G$ denote the isometric isomorphism
which sends $f\in L^p(\R)$ to the sequence $x=(x(m))_{m\in\Z}\in
\ell^p(\Z,L^p[0,1])$, where $x(m)\in L^p[0,1]$ is given by
$$
(x(m))(t) = f(m+t), \quad m\in\Z, \;0<t<1.
$$
Then the spectral properties of an integral operator\index{integral
operator}\index{operator!integral} $\Ko$ on $L^p(\R)$, whose action
is given by
$$
\Ko f(t) = \int_\R k(s,t) f(s)\,ds, \quad t\in\R,
$$
for some kernel function $k$, can be studied by considering its {\em
discretisation}\index{discretisation} $K := G\Ko G^{-1}$. In  turn
$K$ is determined by its matrix $[K]=[\kappa_{mn}]_{m,n\in\Z}$, with
$\kappa_{mn}\in L(L^p[0,1])$ the integral operator given by
$$
\kappa_{mn} g(t) = \int_0^1 k(m+s,n+t)g(s)ds, \quad 0\leq t\leq 1.
$$

Let us also indicate how the results we will develop are relevant to
differential operators\index{differential
operator}\index{operator!differential} (and other non-zero order
pseudo-differential operators\index{pseudo-differential
operator}\index{differential
operator!pseudo-}\index{operator!differential!pseudo-}). Consider
the first order linear differential operator $\Lo$, which we can
think of as an operator from $BC^1(\R)$ to $BC(\R)$, defined by
$$
\Lo y(t) = y^\prime(t) + a(t) y(t), \quad t\in\R,
$$
for some $a\in BC(\R)$. (Here $BC(\R)\subset
L^\infty(\R)$\index{$BC(\R^N)$} denotes the space of bounded
continuous functions on $\R$ and $BC^1(\R):=\{x\in
BC(\R):x^\prime\in BC(\R)\}$\index{$BC^1(\R^N)$}.) In the case when
$a(s) \equiv 1$ it is easy to see that $\Lo$ is invertible.
Specifically, denoting $\Lo$ by $\Lo_1$ in this case and defining
$\Co_1:BC(\R)\to BC^1(\R)$ by
$$
\Co_1 y(t) = \int_\R \kappa(s-t) y(s)\, ds,
$$
where
$$
\kappa(s) := \left\{\begin{array}{cc}
                    e^{s}, & s < 0, \\
                    0, & \mbox{otherwise}, \\
                  \end{array}\right.
$$
it is easy to check by explicit calculation  that
$\Lo_1\Co_1=\Co_1\Lo_1 = \Io$ (the identity operator). Thus the
study of spectral properties of the differential operator $\Lo$ is
reduced, through the identity
\begin{equation} \label{eqn_ds}
\Lo = \Lo_1 + \M_{a-1} = \Lo_1(\Io + \Ko),
\end{equation}
where $\M_{a-1}$ denotes the operator of multiplication by the
function $a-1$, to the study of spectral properties of the integral
operator $\Ko = \Co_1 \M_{a-1}$.

This procedure of reduction of a differential equation to an
integral equation applies much more generally; indeed the above
example can be viewed as a special case of a general reduction of
study of a pseudo-differential operator of non-zero order to one of
zero order (see e.g.\ \cite[\S4.4.4]{RaRoSiBook}). One interesting
and simple generalisation is to the case where $\underline{\Lo}$ is
a matrix differential operator\index{matrix differential
operator}\index{differential
operator!matrix}\index{operator!differential!matrix}, a bounded
operator from $(BC^1(\R))^M$ to $(BC(\R))^M$ given by
$$
\underline{\Lo} x(t) = x^\prime(t) + A(t) x(t), \quad t\in\R,
$$
where $A$ is an $M\times M$ matrix whose entries are in $BC(\R)$.
Then, modifying the above argument, the study of $\underline{\Lo}$
can be reduced to the study of the matrix integral
operator\index{matrix integral operator}\index{integral
operator!matrix}\index{operator!integral!matrix} $\underline{\Ko} =
\underline{\Co}\,\underline{\M}_{A-I}$. Here $\underline{\M}_{A-I}$
is the operator of multiplication by the matrix $A-I$ ($I$ the
identity matrix) and $\underline{\Co}$ is the diagonal matrix whose
entries are the (scalar) integral operator $\Co_1$.

Large parts of the generalisation to the case when the Banach space
$U$ is infinite-dimensional apply only in the case when $A=I+K$,
where $I$ is the identity operator and the entries of
$[K]=[\kappa_{mn}]_{m,n\in\Z}$ are collectively compact. (Where
$\mathcal{I}$ is some index set, a family $\{A_i:i\in \mathcal{I}\}$
of linear operators on a Banach space $U$ is said to be {\em
collectively compact}\index{collectively
compact}\index{compact!collectively} if $\{A_i x:
i\in\mathcal{I},\,x\in U,\, \|x\|\leq 1\}$ is relatively compact in
$U$.) The first half of this text (Chapters
\ref{sec_strict}-\ref{sec_LimOps}) is devoted to developing an
abstract theory of limit operators, in which $Y$ is a general Banach
space and in which the role of compactness and collective
compactness ideas (in an appropriate weak sense) play a prominent
role. Specifically we combine the abstract theory of limit operators
as expounded recently in \cite[Chapter 1]{RaRoSiBook} with the
generalised collectively compact operator theory developed in
\cite{ChanZh02}, building up in Chapter \ref{sec_LimOps} to general
results in the theory of limit operators whose power we illustrate
in the second half of the text, deriving results of the type
(\ref{eqn_ess}).

Let us give a flavour of the general theory we expound in the first
half of the text. To do this it is helpful to first put the
example we have introduced above in more abstract notation. In the
case $Y=\ell^p=\ell^p(\Z,\C)$, let $V_k\in L(Y)$, for $k\in \Z$,
denote the translation operator\index{translation
operator}\index{shift
operator}\index{operator!translation}\index{operator!shift} defined
by
\begin{equation} \label{eq_shift1}
V_kx(m) = x(m-k),\quad m\in\Z.\index{$V_k$}
\end{equation}
Then it follows from our definition above that $A_h$ is a limit
operator of the operator $A$ defined by (\ref{eqn_first}) if
$[V_{-h(j)} A V_{h(j)}]$ (the matrix representation of $V_{-h(j)} A
V_{h(j)}$) converges elementwise to $[A_h]$ as $j\to\infty$. Let us
introduce, moreover, $P_n\in L(Y)$\index{$P_n$} defined by
$$
P_n x(m) = \left\{\begin{array}{cc}
                    x(m), & |m|\leq n, \\
                    0, & \mbox{otherwise}. \\
                  \end{array}\right.
$$
Given sequences $(y_n)\subset Y$ and $(B_n)\subset L(Y)$ and
elements $y\in Y$ and $B\in L(Y)$ let us write $y_n\sto
y$\index{$\sto$} and say that $(y_n)$ converges
strictly\index{strict convergence}\index{convergence!strict} to $y$
if the sequence $(y_n)$ is bounded and
\begin{equation} \label{eqn_str}
\|P_m(x_n-x)\|\to 0 \mbox{ as } n\to\infty,
\end{equation}
for every $m$, and
write
$B_n\pto
B$\index{$\pto$}\index{P-conv@$\P$-convergence}\index{convergence!$\P$-}
if the sequence $(B_n)$ is bounded and
\begin{equation} \label{eqn_pto}
\|P_m(B_n-B)\| \to 0 \mbox{ and } \|(B_n-B)P_m\| \to 0 \mbox{ as }
n\to\infty,
\end{equation}
for every $m$. Then $A_h$ is a limit operator\index{limit
operator}\index{operator!limit} of $A$ if
\begin{equation} \label{eqn_lo}
V_{-h(n)} A V_{h(n)} \pto A_h.
\end{equation}
Defining, moreover, for $b=(b(m))_{m\in Z}\in \ell^\infty$, the
multiplication operator\index{multiplication
operator}\index{operator!multiplication} $M_b\in L(Y)$\index{$M_b$}
by
\begin{equation} \label{eq_mult1}
M_b x (m) = b(m) x(m), \quad m\in\Z,
\end{equation}
we note that $A$ is a band operator\index{band
operator}\index{operator!band} with band width\index{band width} $w$
if and only if $A$ has a representation in the form
\begin{equation} \label{eq_BO1}
A = \sum_{|k|\le w} M_{b_k}V_k,
\end{equation}
for some $b_k\in \ell^\infty$. The set $BO(Y)$\index{$BO(Y)$} of
band operators on $Y$ is an algebra. The Banach subalgebra of $L(Y)$
that is the closure of $BO(Y)$ in operator norm will be called the
algebra of {\em band-dominated operators}\index{band-dominated
operator}\index{operator!band-dominated}, will be denoted by
$BDO(Y)$\index{$BDO(Y)$}, and will play a main role in the second
half of the text, from Chapter \ref{sec_ellp} onwards.

In the general theory we present in the first five chapters,
following \cite{RaRoSiBook} and \cite{ChanZh02}, $Y$ becomes an
arbitrary Banach space, the specific operators $P_n$ are replaced by
a a sequence $\P=(P_n)_{n=0}^\infty$\index{$P$@$\P$} of bounded
linear operators on $Y$, satisfying constraints specified at the
beginning of Chapter \ref{sec_strict}, the specific translation
operators $V_n$ are replaced by a more general discrete group of
isometric isomorphisms, and then the definitions (\ref{eqn_str}),
(\ref{eqn_pto}), and (\ref{eqn_lo}) are retained in essentially the
same form. The notion of compactness that proves important is with
respect to what we term (adapting the definition of Buck
\cite{Buck58}) the {\em  strict topology}\index{strict
topology}\index{topology!strict} on $Y$, a topology in which $\sto$
is the sequential convergence. Moreover, when we study operators of
the form $A=I+K$ it is not compactness of $K$ with respect to the
strict topology that we require (that $K$ maps a neighbourhood of
zero to a relatively compact set)\index{compact
operator}\index{operator!compact}, but a weaker notion, that $K$
maps bounded sets to relatively compact sets, operators having this
property sometimes denoted {\em Montel}\index{Montel
operator}\index{operator!Montel} in the topological vector space
literature. (The notions `compact' and `Montel' coincide in normed
spaces; indeed this is also the case in metrisable topological
vector spaces.)

In the remainder of this introductory chapter, building on the above
short overview and flavour of the text, we detail a history of the
limit operator method and compactness ideas applied in this context,
with the aim of putting the current text in the context of
extensive previous developments in the study of differential and
pseudo-differential equations on unbounded domains; in this history,
as we shall see, a prominent role and motivating force has been the
development of theories for operators with almost periodic
coefficients. In the last section we make a short, but slightly more
detailed summary of the contents of the chapters to come.

\section{A Brief History} \label{sec_history}

The work reported in this text has a number of historical roots.
One we have already mentioned is the paper by Buck \cite{Buck58}
whose strict topology we adapt and use throughout this text. A
main thread is the development of limit operator ideas. The
historical development of this thread of research, which commences
with the study of differential equations with almost periodic
solutions, can be traced through the papers of Favard \cite{Favard},
Muhamadiev \cite{Muh1971,Muh1972,Muh1981,Muh1985}, Lange and
Rabinovich \cite{LaRa1985a,LaRa1985b,LaRa1987}, culminating in more
recent work of Rabinovich, Roch and Silbermann
\cite{RaRoSi1998,RaRoSi2001,RaRoSiBook}. The other main historical
thread, which has developed rather independently but overlaps
strongly, is the development of collectively compact operator theory
and generalisations of this theory, and its use to study
well-posedness and stability of approximation methods for integral
and other operator equations.

{\bf Limit Operators. } To our knowledge, the first appearance of
limit operator ideas is in a 1927 paper of Favard \cite{Favard}, who
studied linear ordinary differential equations with almost periodic
coefficients. His paper deals with systems of ODEs on the real line
with almost periodic coefficients, taking the form
\begin{equation} \label{eqn_fav}
x^{\prime}(t) + A(t)x(t) = f(t),
\quad t \in\R,
\end{equation}
where the $M\times M$ matrix $A(t)$ has entries that are almost
periodic functions of $t$ and the function $f$ is almost periodic. A
standard characterisation of almost periodicity is the following.
Let $\Sh(A):=\{V_s A:s\in \R\}$ denote the set of translates of $A$
(here $(V_s A)(t) = A(t-s)$). Then the coefficients of $A$ are
almost periodic\index{almost periodic matrix function} if and only
if $\mathcal{T}(A)$ is relatively compact in the norm topology on
$BC(\R)$. If $A$ is almost periodic, the compact set that is the
closure of $\mathcal{T}(A)$ is often denoted $\mathcal{H}(A)$ and
called the {\em hull}\index{hull} of $A$. A main result in
\cite{Favard} is the following: if
\begin{equation} \label{eqn_FC}
x^{\prime}(t) + \tilde A(t)x(t) = 0, \quad t \in\R,
\end{equation}
has only the trivial solution in $BC^1(\R)$, for all $\tilde A\in
\mathcal{H}(A)$, and $(\ref{eqn_fav})$ has a solution in $BC^1(\R)$,
then (\ref{eqn_fav}) has a solution that is almost periodic. (Since
$A\in\mathcal{H}(A)$, this is the unique solution in $BC^1(\R)$.)

Certain of the ideas and concepts that we use in this text are
present already in this first paper, for example the role in this
concrete setting of the strict convergence $\sto$ and of compactness
arguments. In particular, conditions analogous to the requirement
that (\ref{eqn_FC}) have no non-trivial bounded solutions for all
$A\in \mathcal{H}(A)$ will play a strong role in this text.
Conditions of this type are sometimes referred to as {\em Favard
conditions}\index{Favard condition} (e.g.\ Shubin
\cite{Shubin75,Shubin78}, Kurbatov \cite{Kurbatov1989,Kurbatov},
Chandler-Wilde \& Lindner \cite{CWLi:JFA2007}).

The first appearance of limit operators per se would seem to be in
the work of Muhamadiev \cite{Muh1971,Muh1972}. In \cite{Muh1971}
Muhamadiev develops Favard's theory as follows. In terms of the
differential operator $\underline{\Lo}:(BC^1(\R))^M\to(BC(\R))^M$
given by (\ref{eqn_ds}), equation (\ref{eqn_fav}) is
$$
\underline{\Lo} x = f.
$$
Under the same assumptions as Favard (that $A$ is almost periodic
and the Favard condition holds) Muhamadiev proves that
$\underline{\Lo}:(BC^1(\R))^M\to(BC(\R))^M$ is a bijection.
Combining this  result with that of Favard, it follows that
$\underline{\Lo}$ is also a bijection from $(AP^1(\R))^M$ to
$(AP(\R))^M$. (Here $AP(\R)\subset BC(\R)$\index{$A_p$@$AP(\R^N)$}
is the set of almost periodic functions and $AP^1(\R)=AP(\R)\cap
BC^1(\R)$\index{$A_p1$@$AP^1(\R^N)$}.) New ideas which play an
important role in the proof of these results include a method of
approximating almost periodic by periodic functions and the fact
that, if $A$ is a periodic function, then injectivity of
$\underline{\Lo}$ implies invertibility. (These ideas are taken up
in the proofs of Theorems \ref{lem_A-inj-inv} and \ref{prop_NRBDO}
in Chapter \ref{sec_ellp}.)

Muhamadiev also considers in the same paper the more general
situation when the entries of $A$ are in the much larger set
$BUC(\R)\subset BC(\R)$ of bounded uniformly continuous functions. A
key property here (which follows from the Arzela-Ascoli theorem and
a diagonal argument) is that, if the sequence $(t_n)\subset \R$
tends to infinity, then $A(\cdot-t_n)$ has a subsequence which is
convergent to a limit $\tilde A$, uniformly on every finite
interval. (Cf.\ the concept of a {\em rich} operator introduced in
\S\ref{subsec_limops}.) Denoting by $\Lim(A)$ the set of limit
functions $\tilde A$ obtained in this way, the following theorem is
stated: if (\ref{eqn_FC}) only has the trivial solution in
$BC^1(\R)$ for every $\tilde A\in \Lim(A)$ then $\underline{\tilde
\Lo}:(BC^1(\R))^M\to(BC(\R))^M$ is a bijection for every $\tilde
A\in \Lim(A)$ (here $\underline{\tilde \Lo}$ denotes the operator
defined by (\ref{eqn_ds}) with $A$ replaced by $\tilde A$).

This is a key result in the development of limit operator theory and
it is a shame that \cite{Muh1971} does not sketch what must be an
interesting proof (we are told only that it `is complicated').
Denoting by $\M_A$ the operator of multiplication by $A$, the set
$\{\M_{\tilde A}:\tilde A\in\Lim(A)\}$ is a set of limit operators
of the operator $\M_A$, and so the set $\{\underline{\tilde
\Lo}:\tilde A\in\Lim(A)\}$ is a set of limit operators of the
operator $\underline{\Lo}$. Thus this result takes the form: if each
limit operator $\underline{\tilde\Lo}$ is injective, specifically
$\underline{\tilde\Lo}x=0$ has no non-trivial bounded solution, then
each limit operator is invertible. A result of this form is a
component in the proof of (\ref{eqn_ess}) and similar results in
this text (and see \cite{CWLi:JFA2007}). In the case that $A$ is
almost periodic it is an easy exercise to show that $\mathcal{H}(A)
= \Lim(A)$, i.e.\ the hull of $A$ coincides with the set of limit
functions of $A$ (cf. Theorem \ref{prop_NR1}). Thus this second
theorem of Muhamadiev includes his result for the case when $A$ is
almost periodic.

The first extension of results of this type to multidimensional
problems is the study of systems of partial differential equations
in $\R^N$ in \cite{Muh1972}. Muhamadiev studies differential
operators elliptic in the sense of Petrovskii with bounded uniformly
H\"older continuous coefficients, specifically those operators
$\underline{\Lo}$ that are what he terms {\em
recurrent}\index{recurrent operator}\index{operator!recurrent}, by
which he means that
$\opsp(\underline{\Lo})=\opsp(\underline{\tilde\Lo})$, for all
$\underline{\tilde\Lo}\in\opsp(\underline{\Lo})$. Here
$\opsp(\underline{\Lo})$ is an appropriate version of the operator
spectrum of $\underline{\Lo}$. Precisely, where $A_p(t)$, for
$t\in\R^N$ and for multi-indices $p$ with $|p|\leq r$, is the family
of coefficients of the operator $\underline{\Lo}$ (here $r$ is the
order of the operator), the differential operator of the same form
$\underline{\tilde\Lo}$ with coefficients $\tilde A_p(t)$ is a
member of $\opsp(\underline{\Lo})$ if there exists a sequence
$t_k\to\infty$ such that, for every $p$,
\begin{equation} \label{eqn_Ap}
A_p(t-t_k)\to \tilde A_p(t)
\end{equation}
uniformly on compact subsets of $\R^N$ as $k\to\infty$.

The main result he states is for the case where $\underline{\Lo}$ is
recurrent and is also, roughly speaking, almost periodic with
respect to the first $N-1$ variables. His result takes the form that
if a Favard condition is satisfied ($\underline{\tilde\Lo}x=0$ has
no non-trivial bounded solutions for all
$\underline{\tilde\Lo}\in\opsp(\underline{\Lo})$) and if
supplementary conditions are satisfied which ensure that
approximations to $\underline{\Lo}$ with periodic coefficients have
index zero as a mapping between appropriate spaces of periodic
functions, then $\underline{\Lo}$ is invertible as an operator
between appropriate spaces of bounded H\"older continuous functions.

Muhamadiev's results apply in particular in the case when the
coefficients of the differential operator are almost periodic (an
almost periodic function is recurrent and its set of limit functions
is its hull). Shubin, as part of a review of differential (and
pseudo-differential) operators with almost periodic solutions
\cite{Shubin78}, gives a detailed account of Muhamadiev's theory, in
the almost periodic scalar case (one case where Muhamadiev's
supplementary conditions are satisfied), and of results which relate
invertibility in spaces of bounded functions to invertibility in
$L^2(\R^N)$. Specifically, his paper includes a proof, for a scalar
elliptic differential operator $\Lo$ with $C^\infty$ almost periodic
coefficients, that the following are equivalent: (i) that the Favard
condition holds; (ii) that $\Lo$ is invertible as an operator on
$BC^\infty(\R^N)$; (iii) that $\Lo$ is invertible as an operator on
$L^2(\R^N)$ in an appropriate sense.

In \cite{Muh1981} Muhamadiev continues the study of the same class
of differential operators $\Lo$ on $\R^N$, elliptic in the sense of
Petrovskii, but now, for some of his results, with no constraints on
behaviour of coefficients at infinity beyond boundedness, though his
main results require also uniform H\"older continuity of all his
coefficients. With this constraint (which, inter alia, is a {\em
richness} requirement in the sense of \S\ref{subsec_limops}), he
studies Fredholmness (or Noethericity) of $\underline{\Lo}$
considered as a bounded operator between appropriate spaces of
bounded H\"older continuous functions. It is in this paper that a
connection is first made between Fredholmness of an operator and
invertibility of its limit operators. The identical Favard condition
to that in \cite{Muh1972} plays a key role. His main results are the
following: (i) that $\underline{\Lo}$ is $\Phi_+$ iff the Favard
condition holds; (ii) that if $\underline{\Lo}$ is $\Phi_-$ then all
the limit operators of $\underline{\Lo}$ are surjective; (iii) (his
Theorem 2.5 and his remark on p.\ 899) that $\underline{\Lo}$ is
Fredholm iff all the limit operators of $\underline{\Lo}$ are
invertible. We note further that his methods of argument in the
proof of his Theorem 2.1 show moreover that if $\underline{\Lo}$ is
Fredholm then the limit operators of $\underline{\Lo}$ are not only
invertible but the inverses are also uniformly bounded, i.e.\
$$
\sup_{\underline{\tilde\Lo}\in\opsp(\underline{\Lo})}
\|{\underline{\tilde\Lo}}^{-1}\|<\infty.
$$
Extensions of these results to give criteria for normal solvability
and Fredholmness of $\underline{\Lo}$ as an operator on Sobolev
spaces are made in \cite{Muh1985}.

In \cite{Muh1981} Muhamadiev also, briefly, introduces what we can term a
{\em weak limit operator}. Uniform continuity of the coefficients $A_p(t)$
is required to ensure that every sequence $t_k\to\infty$ has a
subsequence, which we denote again by $t_k$, such that the limits
(\ref{eqn_Ap}) exist uniformly on compact subsets (cf.\ the definition of
{\em richness} in \S\ref{subsec_limops}). The set of all limit operators
defined by (\ref{eqn_Ap}) where the convergence is uniform on compact sets
we have denoted by $\opsp(\underline{\Lo})$. Muhamadiev notes that it is
enough to require that the coefficients $A_p$ be bounded (and measurable)
for the same {\em richness} property to hold but with convergence
uniformly on compact sets replaced\footnote{We note that, since the
coefficients $A_p$ are bounded so that the sequence $A_p(\cdot-t_k)$ is
bounded, requiring that the limits (\ref{eqn_Ap}) exist uniformly on
compact subsets is equivalent to requiring convergence $\sto$ in the
strict topology, while weak convergence\index{convergence!weak} in
$L^2(\R^N)$ is equivalent to weak$*$
convergence\index{convergence!weak$^*$} in $L^\infty(\R^N)$} by weak
convergence in $L^2(\R^N)$. In the case when the coefficients $A_p$ are
bounded, the set of limit operators defined by (\ref{eqn_Ap}) where the
convergence is weak convergence in $L^2(\R^N)$ we will term the set of
{\em weak limit operators}\label{weaklim}\index{weak limit
operator}\index{limit operator!weak}\index{operator!limit!weak} of
$\underline{\Lo}$. We note that this set coincides with
$\opsp(\underline{\Lo})$ in the case when each $A_p$ is uniformly
continuous. In \cite{Muh1985} Muhamadiev gives criteria for Fredholmness
of $\underline{\Lo}$ on certain function spaces in terms of invertibility
of each of the weak limit operators of $\underline{\Lo}$.


Muhamadiev's work has been a source of inspiration for the decades
that followed. For example, similar to his main results in
\cite{Muh1981} but much more recently, A. and V. Volpert show that,
for a rather general class of scalar elliptic partial differential
operators $\Lo$ on rather general unbounded domains and also for systems of such,
a Favard condition is equivalent to the $\Phi_+$ property of $\Lo$
on appropriate H\"older \cite{Volpert2003,Volpert2005,Volpert2006}
or Sobolev \cite{Volpert2002,Volpert2005,Volpert2006} spaces.

Lange and Rabinovich \cite{LaRa1985a}, inspired by and building on
Muhamadiev's paper \cite{Muh1981}, carry the idea of (semi-)Fredholm
studies by means of limit operators over to the setting of operators
on the discrete domain $\Z^N$. They give sufficient and necessary
Fredholm criteria for the class $BDO(Y)$ of band-dominated operators
(as defined after (\ref{eq_BO1}) and studied in more detail below in
\S\ref{sec_BDO}) acting on $Y=\ell^p(\Z^N,\C)$ spaces. For
$1<p<\infty$, they show that such an operator is Fredholm iff all
its limit operators are invertible and if their inverses are
uniformly bounded. Their proof combines the limit operator arguments
of Muhamadiev \cite{Muh1981} with ideas of Simonenko and Kozak
\cite{KozakSimonenko,Simonenko1965,Simonenko1968} for the
construction of a Fredholm regulariser of $A$ by a clever assembly
of local regularisers. Lange and Rabinovich are thereby the first to
completely characterise Fredholmness in terms of invertibility of
limit operators for the general class of band-dominated operators on
$\ell^p(\Z^N,\C)$. Before, Simonenko
\cite{Simonenko1965,Simonenko1968} was able to deal with the
subclass of those operators whose coefficients (i.e. matrix
diagonals) converge along rays at infinity; later Shteinberg
\cite{Shteinberg78} was able to relax this requirement to a
condition of slow oscillation at infinity. Lange and Rabinovich
require nothing but boundedness of the operator coefficients.

The final section of \cite{LaRa1985a} studies (semi-)Fredholmness of
operators in the so-called {\em Wiener algebra}\index{Wiener
algebra} $\W$\index{$W$@$\W$} (see our \S\ref{subsection_wiener})
consisting of all operators
\begin{equation} \label{eq_Wintro}
A = \sum_{k\in\Z^N} M_{b_k}V_k\qqtext{with}
\sum_k\|b_k\|_\infty<\infty,
\end{equation}
where $b_k\in\ell^\infty(\Z^N,\C)$ for every $k\in\Z^N$ are the
coefficients (or diagonals) of the operator $A$ and $V_k$ and
$M_{b_k}$ are the shift and multiplication operators defined in
(\ref{eq_shift1}) and (\ref{eq_mult1}). Operators $A\in\W$ belong to
$BDO(Y)$ for all spaces $Y=\ell^p(\Z^N,\C)$, $p\in[1,\infty]$. For
$p=\infty$, an analogue of the main result of \cite{Muh1981} is
formulated (in fact, the proof in \cite{LaRa1985a} literally
consists of the sentence `The proofs of Theorems 4.1 and 4.2 repeat
the proofs of Theorems 2.1 and 2.2 in \cite{Muh1981}, with obvious
amendments.'): $A$ is $\Phi_+$ iff all its limit operators are
injective, i.e. Favard's condition holds; if $A$ is $\Phi_-$ then
all its limit operators are surjective. The paper concludes with a
first, simplified version of our Theorem \ref{Wiener2} below, with a
somewhat abbreviated proof: that $A\in\W$ is either Fredholm on all
spaces $Y=\ell^p(\Z^N,\C)$, $p\in[1,\infty]$, or on none of them.
Moreover, the uniform boundedness condition on the inverses of its
limit operators is redundant. The latter implies that
\begin{equation} \label{eq:spessW}
\spess(A)\ =\ \bigcup_{A_h\in\opsp(A)} \sp(A_h)
\end{equation}
if $A\in\W$, with all expressions independent of $p\in[1,\infty]$.

From here on we mainly follow the discrete branch of the limit
operator story since this is the focus of our text, noting that
the further generalisation from scalar-valued to vector-valued
$\ell^p$ spaces $Y=\ell^p(\Z^N,U)$ with an arbitrary complex Banach
space $U$ enables us to emulate differential, integral and
pseudo-differential operators on $L^p(\R^N)$ (e.g.\
\cite{LaRa1985b}) by operators on the discrete space $Y$ with
$U=L^p([0,1]^N)$ (see e.g.\ \cite{Kurbatov,RaRoSi2001}, the
discussion in the paragraphs after equation \eqref{eqn_ess3} above,
and Chapter \ref{sec_intop} below).

In the last 10 years, the limit operators of band-dominated
operators on the discrete spaces $Y=\ell^p(\Z^N,U)$ with
$p\in(1,\infty)$ have been extensively studied by Rabinovich, Roch,
Silbermann and a small number of their coauthors. The first work of
this troika was \cite{RaRoSi1998}, where the results of
\cite{LaRa1985a} for $p\in (1,\infty)$ are picked up, this time with
full proofs, and are extended, utilised and illuminated in
connection with other problems and concepts such as the
applicability of the so-called finite section method (a truncation
method for the approximate solution of corresponding operator
equations) and the idea of two different symbol calculi in the
factor algebra of $BDO(Y)$ modulus compact operators. Another
important result of \cite{RaRoSi1998} is the observation that the
limit operator idea is compatible with the local principle of Allan
\cite{Allan} and Douglas \cite{Douglas} for the study of
invertibility in non-commutative Banach algebras. The latter result
was used to slightly relax the uniform boundedness condition on the
inverses of the limit operators in the general Fredholm criterion
\cite[Corollary 5]{RaRoSi1998} and to completely remove this
condition in the case of slowly oscillating coefficients
\cite[Theorem 9]{RaRoSi1998}.

In \cite{RaRoSi2001}, building on results of \cite{Ra2,RaRoSi1998},
the same authors tackle the case when $U$ is an arbitrary Hilbert
space under the additional condition that $p=2$ so that
$Y=\ell^2(\Z^N,U)$ is a Hilbert space too and the set of
band-dominated operators on it is a C$^*$-algebra. In this C$^*$
setting, which makes life slightly easier than the more general case
when $BDO(Y)$ is merely a Banach algebra, the serious obstacle of an
infinite dimensional space $U$ is overcome. The matrix $[A]$ that
corresponds to an operator $A\in BDO(Y)$ now has operator entries
$a_{ij}\in L(U)$ which are infinite dimensional operators
themselves. This changes the Fredholm theory completely: An operator
$A$ with only finitely many nonzero entries $a_{ij}$ is in general
no longer of finite rank -- not even compact. That is why
Rabinovich, Roch and Silbermann replace the ideal
$K(Y)$\index{$K(Y)$} of compact operators by another set, later on
denoted by $K(Y,\P)$\index{$K(Y,\P)$}, which is the norm closure of
the set of all operators $A$ with finitely many nonzero matrix
entries. Also this set is contained in $BDO(Y)$, it is an ideal
there and it is shown that if for $A\in BDO(Y)$ there exists a
$K(Y,\P)$-regulariser $B\in L(Y)$, i.e. $AB-I$ and $BA-I$ are in
$K(Y,\P)$, then automatically $B\in BDO(Y)$. If $U$ is finite
dimensional and $p\in (1,\infty)$, which was the setting of
\cite{RaRoSi1998}, then $K(Y,\P)$ is the same as $K(Y)$ and
invertibility modulo $K(Y,\P)$, termed {\em invertibility at
infinity}\index{invertible at infinity} in \cite{RaRoSi2001},
coincides with invertibility modulo $K(Y)$ alias Fredholmness. So
one could argue that in \cite{RaRoSi1998} the subject already was
invertibility at infinity which, as a coincidence, turned out to be
Fredholmness too. In fact, the major milestone in \cite{RaRoSi2001}
was to understand that the limit operator method studies
invertibility at infinity and not Fredholmness, and therefore the
new ideal $K(Y,\P)$ was the right one to work with. Fortunately,
invertibility at infinity and Fredholmness are closely related
properties so that knowledge about one of them already says a lot
about the other and so the limit operator method can still be used
to make statements about Fredholmness -- via invertibility at
infinity.

Another problem that occurs when passing to an infinite dimensional
space $U$ is that the simple Bolzano-Weierstrass argument (coupled
with a diagonal construction) previously showing that, for $A\in
BDO(Y)$, every sequence $h=(h(k))_{k\in\N}\subset \Z^N$ with
$|h(k)|\to\infty$ has a subsequence $g$ such that the matrix of the
translates $[V_{-g(k)}AV_{g(k)}]=[a_{i+g(k),j+g(k)}]_{i,j\in\Z^N}$
converges entry-wise as $k\to\infty$, is no longer applicable as the
matrix diagonals are bounded sequences in the infinite dimensional
space $L(U)$ now. So the class of all operators $A\in BDO(Y)$ for
which every such sequence $h$ has a subsequence $g$ with this
convergence property (the limiting operator being the limit operator
$A_g$) had to be singled out in \cite{RaRoSi2001}. Operators of this
class were later on termed {\em rich operators}\index{rich
operator}\index{operator!rich}.

There is one more technical subtlety when passing to an infinite
dimensional space $U$: The so-called $\P-$convergence
(\ref{eqn_pto}) that is used in (\ref{eqn_lo}) is equivalent to
strong convergence $B_n\to B$ and $B_n^*\to B^*$ if $p\in
(1,\infty)$ and $U$ is finite dimensional; in fact, this is how it
was treated in \cite{RaRoSi1998}. So this was another difference to
\cite{RaRoSi1998} although nothing new since $\P-$convergence was de
facto introduced for exactly this purpose by Muhamadiev
\cite{Muh1981} already.

The next two works in this story were the very comprehensive
monograph \cite{RaRoSiBook} by the troika Rabinovich, Roch and
Silbermann, which summarised the state of the art to which it
largely contributed itself, and the PhD thesis \cite{LiDiss} of the
second author of this text. Both grew at roughly the same time and
under mutual inspiration and support. In \cite{RaRoSiBook}, besides
many other things that cannot be discussed here, the case
$Y=\ell^p(\Z^N,U)$ was successfully treated for arbitrary Banach
spaces $U$ and $p\in (1,\infty)$. The gaps at $p\in\{1,\infty\}$ are
filled in \cite{LiSi03} and finally in \cite{LiDiss}. The challenge
about $p=\infty$ is that duality, which is a frequent instrument in
the arguments of \cite{RaRoSi1998,RaRoSi2001,RaRoSiBook}, is more
problematic since the dual space of $Y$ is no longer one of the
$Y$-spaces at hand. Instead one works with the predual and imposes
the existence of a preadjoint operator acting on it. Note that some
of these ideas have been picked up and are significantly extended
and improved in Section \ref{sec_dsa} below.

Another important thread that should be mentioned here is the
determination not only of Fredholmness but also of the Fredholm
index by means of limit operators. The key paper in this respect is
\cite{RaRoRoe} by Rabinovich, Roch and Roe, where the case $N=1$,
$p=2$, $U=\C$ has been studied using $C^*-$algebra techniques
combined with $K-$theory. The idea is to decompose $Y=\ell^2(\Z,\C)$
into the subspaces $Y_-$ and $Y_+$ that correspond to the negative
and the non-negative half axis, respectively, thereby splitting the
twosided infinite matrix $[A]$ of $A\in BDO(Y)$ into the four
onesided infinite submatrices $[A_{--}],[A_{-+}],[A_{+-}]$ and
$[A_{++}]$. Since $[A_{-+}]$ and $[A_{+-}]$ are compact (note that
$U$ has finite dimension), these two blocks can be removed without
changing Fredholmness or the index. By a similar argument, for every
$m\in\N$, the first $m$ rows and columns of both $[A_{--}]$ and
$[A_{++}]$ can be removed without losing any information about
Fredholmness and the index. So it is not really surprising that also
the index of $A$ is exclusively stored in the asymptotic behaviour
of the matrix entries of $[A_{--}]$ and $[A_{++}]$ at infinity, i.e.
in the limit operators of $A$. Indeed, calling the index of
$A_{\pm\pm}$, understood as an operator on $Y_\pm$, the $\pm-$index
of $A$, respectively, it is shown in \cite{RaRoRoe} that all limit
operators of $A$ with respect to sequences tending to $\pm\infty$
have the same $\pm-$index as $A$ has, respectively. Since the index
of $A$ is the sum of its plus- and its minus-index, this gives a
formula for the index of $A$ in terms of plus- and minus-indices of
two of its limit operators. The index formula of \cite{RaRoRoe} was
later carried over to the case $N=1$, $p\in (1,\infty)$, $U=\C$ in
\cite{Roch_ellp} (where it was shown that the index of $A$ does not
depend on $p$ -- see \cite{LiWiener} for the same result in the
setting of a more general Banach space $U$ and $p\in[1,\infty]$),
re-proved by completely different techniques (using the sequence
of the finite sections of $A$) in \cite{RaRoSi:IndexFSM} and
generalised to the case of an arbitrary Banach space $U$ in case
$A=I+K$ with a locally compact operator $K$ (i.e. all entries of
$[K]$ are compact operators on $U$) in \cite{RaRoIndexLC}.

The most recent extended account of the limit operator method is the
monograph \cite{LiBook} by the second author. Besides a unification
of techniques and results of \cite{LiDiss} and \cite{RaRoSiBook}, an
exposition of the topic of infinite matrices, in particular
band-dominated operators, that is accessible for a wide audience and
a number of additions and clarifications to the theory, it also
contains the first fruits of the work with the other author of this
text. For example, it contains a treatment of boundary integral
equations\index{boundary integral equation} on unbounded surfaces
(also see \cite{CWLiAMSProcPaper,CWLi2}), their Fredholmness and
finite sections, as well as more complete results on the interplay
of Fredholmness and invertibility at infinity and on different
aspects of the finite section method.

The above is an account of the main development of limit operator
ideas and the theory of limit operators, starting with the work of
Favard \cite{Favard}. However, there are many other branches of this
story (such as the ``frequency limit operators'' discussed in
\cite{Kurbatov1996} or the ``zoom limit operators'' discussed in
\cite{BoKaRa} and briefly in Section 3.6 of \cite{LiBook}) that we
have not mentioned explicitly. We have also omitted mention of a
number of instances where limit-operator-type ideas have been
discovered and applied independently. In particular,
limit-operator-type ideas have been applied recently to great effect
in the spectral theory of discrete Schr\"odinger and Jacobi operators
as well as more general bounded linear operators on Hilbert spaces. One
instance is the recent work of Davies
\cite{Davies,Davies2001:SpecNSA}, where the spectrum of a random
Jacobi operator $A$ is studied by looking at strong limits of
sequences $U_nAU_n^*$, where $U_n$ denotes a sequence of unitary
operators. (This idea in the work of Davies dates back to an
earlier paper of Davies and Simon \cite{DaviesSimon}, where the idea
of the {\it limit class} of an operator is introduced, which has
some similarity to the idea of an operator spectrum.) In
\cite{Davies} the notion of a pseudo-ergodic operator is introduced
(we take up their study as a significant example in Chapter
\ref{sec_Schroed} below), this idea capturing many aspects of the
spectral behaviour of random operators while eliminating
probabilistic arguments. In limit operator terminology, a Jacobi
operator $A=\sum_{|k|\le 1} M_{b_k} V_k$ is pseudo-ergodic iff every
operator
\begin{equation} \label{eq:formB}
B=\sum_{|k|\le 1} M_{c_k} V_k,
\end{equation}
with $c_k(m)\in\clos\{b_k(n):n\in\Z^N\}$ for $m\in\Z^N$ and $|k|\le
1$, is a limit operator of $A$. Davies shows for a pseudo-ergodic
discrete Schr\"odinger operator $A$ that $\sp(A)=\spess(A)=\cup \,
\sp(B)$, with the union taken over all operators $B$ of the form
(\ref{eq:formB}) (i.e. taken over the set of limit operators of
$A$). This result can be viewed as a special case of
\eqref{eq:spessW}. In both \cite{Davies} and
\cite{Davies2001:SpecNSA} results are obtained by these methods
about the spectrum of the non-self-adjoint Anderson model of Hatano
and Nelson \cite{HatanoNelson}.

In another series of papers, Georgescu, Mantoiu and coworkers
\cite{Mantoiu1,Mantoiu,GeorgescuGolenia,
GeorgescuIfti2000,GeorgescuIfti2002,GeorgescuIfti2005} develop an
impressive  $C^*$-algebraic approach to the essential spectrum of
self-adjoint (possibly unbounded) operators $A$ on $L^2(X)$ with a
locally compact non-compact abelian group $X$. If $A$ is subject to
two so-called Landstad conditions \cite{GeorgescuIfti2005}, where
one of them is equivalent to the condition $A\in
BDO(\ell^2(\Z^N,\C))$ in case $X=\Z^N$, they prove formula
\eqref{eq:spessW} with the closure taken on the right-hand side, i.e.
\begin{equation} \label{eq:spessWcl}
\spess(A)\ =\ \clos \left(\bigcup_{h} \sp(A_h)\right).
\end{equation}
The
operators $A_h$ are called ``localizations at infinity'' of $A$ and
each $A_h$ is defined as strong limit of translates $U_xAU_x^*$ of
$A$, where $(U_xf)(y):=f(x\circ y)$ for $y\in X$ and $x\in X$ tends
to infinity towards an element $h$ of the Stone-\v Cech boundary of
$X$ (the non-trivial characters of $L^\infty(X)$). The set of all
localizations at infinity $A_h$ derived in that way coincides with
$\opsp(A)$, as was shown in \cite{RaRoSi1998,Roe} (also see Section
3.5.2 of \cite{LiBook}).

In \cite{LastSimon06} Last and Simon prove \eqref{eq:spessWcl}
if $A$ is a self-adjoint and bounded Jacobi
operator whose sub/super diagonal is also bounded away from zero. We
discuss another result from \cite{LastSimon06} in Example
\ref{ex:LastSimon} below. Finally, Last, Simon and Remling
\cite{LastSimon99,Remling1,Remling2} show that limit operators are
also an efficient tool in the study of the absolutely continuous
spectrum of self-adjoint Jacobi operators $A$! The reason is that,
like the essential spectrum, also the a.c. spectrum of $A$ is
independent of (self-adjoint) perturbations of $A$ in finitely many
entries of its three matrix diagonals, whence the a.c.\!\! spectrum
of $A$ can be fully described in terms of the limit operators of
$A$. A remarkable difference, discovered in \cite{LastSimon06} for
discrete Schr\"odinger operators in 1D, between the essential and
a.c.\!\! spectrum is however that, while $\spess(B)\subset\spess(A)$
holds if $B$ is a limit operator of $A$, one has the opposite
inclusion $\Sigma_{ac}(B)\supset\Sigma_{ac}(A)$ for the essential
supports of the absolutely continuous parts of the spectral measure.
As one corollary one gets that $\Sigma_{ac}(B)=\Sigma_{ac}(A)$ for
all limit operators $B$ if $A$ has an almost periodic potential.
Remling \cite{Remling1,Remling2} goes well beyond these points by
proving further connections between the set of all limit operators
of the discrete Schr\"odinger operator $A$ and the set
$\Sigma_{ac}(A)$, leading to his so-called oracle theorem that says
that in the presence of a.c.\!\! spectrum, one can predict values
of the potential of $A$ from information on earlier values of the
potential -- an indication of how restricted the class of discrete
Schr\"odinger operators with non-empty a.c.\!\! spectrum is. One
illustration of this result is that if the potential $V$ of $A$ only
assumes finitely many values and $A$ has a.c. spectrum then $V$ must
be eventually periodic. Remling also indicates how to carry his
results over to Jacobi operators.

So for different classes of self-adjoint operators $A$ on $\ell^2$,
the formula \eqref{eq:spessWcl} has been derived by different
authors. We would like to mention that, when $A$ is a normal (e.g.
self-adjoint) operator (so that all of its limit operators $A_h$ are
also normal), then this formula is equivalent to the statement from
\cite{LaRa1985a,RaRoSi1998} that $A-\lambda I$ is Fredholm iff all
limit operators of $A-\lambda I$ are invertible and their inverses
are uniformly bounded since
$\dist(\lambda,\sp(A_h))=1/\|(A_h-\lambda I)^{-1}\|$ for normal
operators $A_h$. However, if $A\in\W$ then it is known
(\cite{LaRa1985a,RaRoSi1998}, see \cite{LiDiss,LiBook} for the
vector-valued case) that the right-hand side of \eqref{eq:spessW} is
closed. For general $A\in BDO(\ell^2)$ this is an open problem (see
number 8 in our list of open problems in Chapter
\ref{chap_openprobs}).

Another fairly substantial body of work, in which limit operator
ideas are important, has grown out of the collectively compact
operator theory of Anselone and co-authors \cite{AnseloneBook}. It
seems appropriate to summarise the main historical developments in
this line of research in separate paragraphs below since this body
of work work has in common that it is characterised by collective
compactness concepts. However, in much of this work
limit-operator-type ideas also play a key role. Specifically, limit
operator ideas combined with collectively compact operator theory
are used already in the 1985 paper of Anselone and Sloan
\cite{AnseloneSloan} to show the stability in $BC[0,\infty)$ of the
finite section method for the classical Wiener-Hopf integral
equation
\begin{equation} \label{eq:wh}
y(s) = x(s) + \int_0^\infty \kappa(s-t)y(t)\, dt, \quad s\geq 0,
\end{equation}
with $\kappa\in L^1(\R)$. (The finite section method\index{finite
section method} is just the approximation of \eqref{eq:wh} by the
equation on the finite interval
$$
y(s) = x(s) + \int_0^A \kappa(s-t)y(t)\, dt, \quad 0\leq s\leq A,
$$
and the main issue is to study stability\index{stability} and
convergence as $A\to\infty$.)
 The methods and results of \cite{AnseloneSloan} are generalised in
Chandler-Wilde \cite{CW93} to obtain criteria in $BC(\R)$ for both
stability of the finite section method and solvability for the
equation
$$
y(s) = x(s) + \int_{-\infty}^\infty \kappa(s-t)z(t)y(t)\, dt, \quad s\in\R,
$$
in operator form
\begin{equation} \label{CW93eq}
y = x + K_z y,
\end{equation}
where $K_z$ is the integral operator with kernel $\kappa(s-t)z(t)$,
and it is assumed that $\kappa\in L^1(\R)$ and $z\in L^\infty(\R)$.

Limit operators do not appear explicitly in \cite{CW93}, or in
generalisations of this work to multidimensional cases
\cite{ChanZh97,ChanZh02}, to more general classes of kernels
\cite{CWZhRoss00}, to other functions spaces ($L^p(\R)$, $1\le p\le
\infty$, or weighted spaces) \cite{ArensCWHas1,ArensCWHas2}, or to
general operator equations on Banach spaces \cite{ChanZh02}. Rather,
as we discuss in the paragraphs below, the results of these papers
provide criteria for unique solvability of \eqref{CW93eq} expressed
in terms of injectivity in $BC(\R)$ (or equivalently in
$L^\infty(\R)$) of the elements of a particular family of operators.
The connection to limit operators, explored  in Section
\ref{subsec_limops} below, is that this family of operators
necessarily contains both the operator $I-K_z$ and all the weak
limit operators of $I-K_z$. (Here {\em weak limit operator} has the
same meaning as in our discussion of the paper \cite{Muh1981} on
page \pageref{weaklim} above; we call $K_{\tilde z}$ a weak limit
operator of $K_z$ if, for some unbounded sequence $(t_k)\subset \R$,
it holds that $z(\cdot-t_k)\stackrel{w*}{\to}\tilde z$ as
$k\to\infty$, where $\stackrel{w*}{\to}$ is weak$*$ convergence in
$L^\infty(\R)$.)

{\bf Collective Compactness. } In the mid 1960's Anselone and
co-workers (see \cite{AnseloneBook} and the references therein)
introduced the concept of collectively compact operators. A family
$\K$ of linear operators on a Banach space $Y$ is called {\em
collectively compact}\index{collectively
compact}\index{compact!collectively} if, for any sequences
$(K_m)\subset\K$ and $(x_m)\subset Y$ with $\|x_m\|\le 1$, there is
always a subsequence of $(K_mx_m)$ that converges in the norm of
$Y$. It is immediate that every collectively compact family $\K$ is
bounded and that all of its members are compact operators.

There are some important features of collectively compact sets of
operators. First, recall that if $K$ is a compact operator on $Y$
and a sequence $A_m$ of operators on $Y$ converges strongly (i.e.
pointwise) to $0$, then $A_mK$ converges to $0$ in the operator norm
on $Y$. But under the same assumption, even $A_mK_m$ converges to
$0$ in the norm for any sequence $(K_m)\subset\K$ provided $\K$ is
collectively compact. This fact was probably the motivation for the
introduction of this notion. It was used by Anselone for the
convergence analysis of approximation methods like the Nystr\"om
method for second kind integral equations.

Another important feature \cite[Theorem 1.6]{AnseloneBook} is that
if $\{K_m\}_{m=1}^\infty$ is collectively compact and strongly
convergent to $K$, then also $K$ is compact, and the following
holds:
\begin{eqnarray} I-K \textrm{ is invertible }
\qquad\Longleftrightarrow\nonumber \\
I-K_m \textrm{ is invertible for large $m$, say } m>m_0, \textrm{
and } \sup_{m>m_0}\|(I-K_m)^{-1}\|<\infty \nonumber
\end{eqnarray}
Since $K$ and $K_m$ are compact, the above is equivalent to the
following statement
\begin{eqnarray} I-K \textrm{ is injective }
\qquad\Longleftrightarrow\qquad \exists m_0:
\inf_{m>m_0}\nu(I-K_m)>0 \label{collcomp2}
\end{eqnarray}
where $\nu(A):=\inf\{\|Ax\|:\|x\|=1\}$\index{$\nu(A)$} is the
so-called {\em lower norm}\index{lower norm} of an operator $A$.

There are many important examples where $K$ is not compact in the
norm topology on $Y$ but does have compactness properties in a
weaker topology. To be precise, $K$, while not compact (mapping a
neighbourhood of zero to a relatively compact set) has the property
that, in the weaker topology, it maps bounded sets to sets that are
relatively compact (such operators are sometimes termed {\em
Montel}\index{Montel
operator}\index{operator!Montel})\footnote{Already in the 1970's
DePree and co-authors \cite{DePree1,DePree2,DePree3} studied
collectively compact operator theory in a topological vector space
setting, but they retained compactness of $K$ rather than studying
the weaker case where $K$ is Montel.}. In particular, this is
generically the case when $K$ is an integral operator on an
unbounded domain with a continuous or weakly singular kernel; these
properties of the kernel make $K$ a `smoothing' operator, so that
$K$ has local compactness properties, but $K$ fails to be compact
because the domain is not compact. Anselone and Sloan
\cite{AnseloneSloan} were the first to extend the arguments of
collectively compact operator theory to tackle a case of this type,
namely to study the finite section method for classical Wiener-Hopf
operators on the half-axis. As mentioned already above, the
arguments introduced were developed into a methodology for
establishing existence from uniqueness for classes of second kind
integral equations on unbounded domains and for analyzing the
convergence and stability of approximation methods in a series of
papers by the first author and collaborators
\cite{CW93,PepCW95,ChanZh97,CWZhRoss00,MeierCW01,CW_NumMath,ChanZh02,ArensCWHas1,ArensCWHas2}.
A particular motivation for this was the analysis of integral
equation methods for problems of scattering of acoustic, elastic and
electromagnetic waves by unbounded surfaces
\cite{CW_MathMeth97,ChanZh_PRS98,ZhCW98,CWRossZh_PRS99,ChanZh_SIAM_JMA99,MeierCW01,CW_NumMath,ZhChan03,Natrosh_Ar_CW03,ArensCWHas2,CWLi2}.
Other applications included the study of multidimensional
Wiener-Hopf operators and, related to the Schr\"odinger operator, a
study of Lippmann-Schwinger integral equations \cite{ChanZh97}.
Related developments of the ideas of Anselone and Sloan
\cite{AnseloneSloan} to the analysis of nonlinear integral equations
on unbounded domains are described in
\cite{Agarwal,AnseloneLee,ORegan}.

In \cite{ChanZh02} the first author and Zhang put these ideas into
the setting of an abstract Banach space $Y$, in which a key role is
played by the notion of a {\it generalised collectively
compact}\index{generalised collectively compact}\index{collectively
compact!generalised}\index{compact!collectively!generalised} family
$\K$. Now the sequence $(K_mx_m)$ has a subsequence that converges
in a topology that is weaker than the norm topology on $Y$, whenever
$(K_m)\subset\K$ and $(x_m)\subset Y$ with $\|x_m\|\le 1$. This
notion no longer requires the elements of $\K$ to be compact
operators and therefore covers a lot more operators originating from
applications. But still, the following similar result to
\refeq{collcomp2} was established in \cite{ChanZh02}. If $\K$ is
generalised collectively compact (or {\em uniformly
Montel}\index{uniformly Montel}\index{Montel!uniformly} as we shall
term the same property in Definition \ref{def_collcomp} in this
text) and some important additional assumptions hold (see the end
of our discussion of limit operators above and Theorem \ref{th4.5}
below), then
\begin{eqnarray}
I-K \textrm{ is injective for all } K\in\K
\qquad\Longleftrightarrow\qquad \inf_{K\in\K}\nu(I-K)>0.
\label{collcomp3}
\end{eqnarray}
If the family $\K$ satisfies rather strong additional constraints
(see Theorem \ref{th4.5} below for details), then also invertibility
of $I-K$ for every $K\in \K$ follows from injectivity for all
$K\in\K$.

To give a concrete flavour of these results (this was the first
concrete application of these ideas made to boundary integral
equations in wave scattering \cite{CW93,CW_MathMeth97}), one case
where they apply is to the integral equation \eqref{CW93eq}, with
the family $\K$ defined by
$$
\K:= \{K_z:z\in L^\infty(\R) \mbox{ and } z(s)\in Q, \mbox{ for almost all } s\in\R\},
$$
for some $Q\subset \C$ which is compact and convex. That is,
existence and uniform boundedness of $(I-K_z)^{-1}$ (as an operator
on $BC(\R)$) for all $K_z\in \K$, can be shown to follow from
injectivity of $I-K_z$ for all $K_z\in \K$ (see
\cite{CW93,CW_MathMeth97} or \cite{ChanZh02} for details).


{\bf Generalised Collective Compactness and Limit Operators: A Main
Aim of this Memoir. } We have briefly indicated above some
connections between the two bodies of work that we have described
under the headings `Limit Operators' and `Collective Compactness'.
It was a main aim for us in writing\footnote{It was results in
\cite{LiDiss}, in particular where it is pointed out for the first
time that the operator spectrum $\opsp(A)$ of a rich operator is
always relatively sequentially compact with respect to
$\P$-convergence -- one of the additional assumptions required for
Theorem \ref{th4.5} below, taken from \cite{ChanZh02} -- which
prompted the authors to start to investigate this symbiosis, in
discussions after the first author examined the thesis of the
second!} this text to explore these connections in a methodical
way, in particular to explore the possibilities for applying the
generalised collectively compact operator theory \cite{ChanZh02} in
a limit operator context, and for combining collective compactness
and limit operator ideas. It turns out that these ideas have a very
fruitful interplay (other recent examples of this interplay in
addition to this text are \cite{RaRoIndexLC} and
\cite{CWLi:JFA2007}, the latter paper making use of some of the
results we will present below). We finish this introduction by
summarising the main new results of the text and of this interplay
of ideas.

\section{Summary and the Main New Results} \label{sec_newresults}

In this final section of the introduction, having given an overview
of the ideas of the text and their historical development in
Sections \ref{sec_overview} and \ref{sec_history}, we briefly
summarise, chapter by chapter, the main contents and results we
obtain.

The text falls into two connected parts. The first part, chapters
2-5, is concerned with extensions to the general abstract theory of
limit operators, as expounded in Chapter 1 of
\cite{RaRoSiBook}, with exploring the connections with the abstract
generalised collectively compact operator theory of \cite{ChanZh02},
and with making applications of this theory in the limit operator
context.

The short initial Chapter \ref{sec_strict} introduces, following
\cite{RaRoSiBook}, the idea of a sequence $(P_n)$ of bounded linear
operators on a Banach space $Y$ that form an {\em approximate
identity} on $Y$, satisfying conditions (i) and (ii) at the
beginning of the chapter. In the theory of limit operators developed
in \cite{RaRoSiBook} the following notion of sequential convergence
plays a crucial role: we write that $x_n\sto x$ if $(x_n)$ is a
bounded sequence and $P_m x_n\to P_m x$ as $n\to\infty$ for every
$m$. In this text, as we have noted already, we call $\sto$ {\em
strict convergence}, by analogy with the strict topology of
\cite{Buck58}. Chapter 2 is concerned with study of a topology,
which we term the strict topology, in which $\sto$ is the sequential
convergence. We recall properties of this topology from
\cite{ChanZh02} (which derive in large part from similar results in
\cite{Buck58}) and show further properties, for example
characterising the compact and sequentially compact sets in the
strict topology, characterising when the strict and norm topologies
coincide, and introducing many examples that we build on later.


In Chapter \ref{chap_classes} we study a number of subspaces of
$L(Y)$, the space of bounded linear operators on the Banach space
$Y$, namely those subspaces that play an important role in the
abstract theory of limit operators \cite{RaRoSiBook} and in the
generalised collectively compact operator theory of \cite{ChanZh02},
and so play an important role in the rest of the text. These
subspaces include the classes $L(Y,\P)$ and $K(Y,\P)$  central to
the theory of limit operators \cite{RaRoSiBook}\footnote{Indeed, in
\cite{RaRoSiBook} the notion of a limit operator of an operator is
defined only for operators in $L(Y,\P)$. We will not be this
restrictive but extend this notion and other definitions as far as
possible to the whole class $L(Y)$.}, the class $S(Y)$ of operators
that are sequentially continuous on $(Y,s)$ ($Y$ equipped with the
strict topology), and the class $SN(Y)$ of operators that are
sequentially continuous from $(Y,s)$ to $(Y, \|\cdot\|)$ ($Y$
equipped with the norm topology). We also study the usual class
$K(Y)$ of compact operators on $(Y, \|\cdot\|)$, the class $KS(Y)$
of compact operators on $(Y,s)$, and the class $M(Y)$ of {\em
Montel} operators on $(Y,s)$ (that map bounded sets to relatively
compact sets), this latter class playing a particularly key role in
the later text. Our concern is to derive explicit
characterisations of these sets and to explore the relationships
between them. For example we see in this chapter that $S(Y)$
coincides with the set of continuous linear operators on $(Y,s)$,
that $SN(Y)$ is the set of continuous linear operators from $(Y,s)$
to $(Y,\|\cdot\|)$, and that $S(Y)$ and $SN(Y)$ are, roughly
speaking, ``one-sided'' versions of $L(Y,\P)$ and $K(Y,\P)$ (see
Lemmas \ref{lem1} and \ref{LYP_iff}, Corollary \ref{cor4} and
(\ref{KYP_def})). An easy but informative result, which makes clear
that being Montel is a much weaker property than being compact on
$(Y,s)$, is Corollary \ref{cor_MYLY},  which has application for
example in Chapter \ref{sec_Schroed}, that $M(Y)=L(Y)$ if $(Y,s)$ is
sequentially complete and each $P_n\in K(Y)$. In Section
\ref{sec:algebraic} we study algebraic properties, for example
showing that all of $S(Y)$, $SN(Y)$, $M(Y)$, and $KS(Y)$ are Banach
subalgebras of $L(Y)$, and that $SN(Y)$, $KS(Y)$, and $M(Y)\cap
S(Y)$ are all ideals of $S(Y)$.

In the short Chapter \ref{chap_converg} we introduce and contrast
various notions of convergence of sequences of operators  in $L(Y)$,
with an emphasis on those used in the abstract theory of limit
operators \cite{RaRoSiBook} and in the generalised collectively
compact operator theory of \cite{ChanZh02}. Specifically, our main
emphasis is on the $\pto$ convergence of Definition \ref{def_P_conv}
(and see \cite{RoSi,RaRoSiBook}, though the $\pto$ convergence in
\cite{RaRoSiBook} is restricted to operators in $L(Y,\P)$), and on
the weaker notions of convergence ($\sto$ and $\Sto$) important in
\cite{ChanZh02}. These three notions of convergence will play a
strong role in subsequent chapters, but we also compare these
notions of convergence to ordinary norm ($\toto$) and strong ($\to$)
convergence in $L(Y)$. Informative characterisations we derive are
that a sequence $(A_n)\subset L(Y)$ satisfies $A_n\pto0$ if and only
if $(A_n)$ is bounded and both $P_mA_n\toto0$ and $A_nP_m\toto0$ as
$n\to\infty$, for every $m$, while $A_n\Sto0$ if and only if $(A_n)$
is bounded and $P_mA_nx\sto 0$ as $n\to\infty$ for every $m$ and
every $x\in Y$. Thus $A_n\toto 0 \Rightarrow A_n\pto 0\Rightarrow
A_n\Sto 0 \Leftarrow A_n\to 0$ (cf.\ Corollary \ref{strong2}).

Chapter \ref{sec_LimOps} introduces the main abstract concepts and
results of the text. In Section \ref{sec_inv_inf_Fred} we
introduce the concept of invertibility at infinity (Definition
\ref{inv_inf}) of an operator $A\in L(Y)$ 
and show that if $A\in S(Y)$ is a Montel perturbation of an invertible
operator in $S(Y)$, then $A$ invertible at infinity implies that $A$ is
Fredholm (Theorem \ref{lem_Fredholm}). In Section \ref{subsec_collcomp} we
present a specialisation of some main results of the generalised
collectively compact operator theory on a Banach space $Y$ of
\cite{ChanZh02}, to the case where the family of seminorms on $Y$ required
in the theory of \cite{ChanZh02} is given by (\ref{eq_seminorm}). In
Section \ref{subsec_limops} we make an abstract definition of the limit
operators and the operator spectrum, $\opsp(A)$ (the set of limit
operators), of an arbitrary operator $A\in L(Y)$. We also list the main
properties of the operator (including Fredholmness, invertibility at
infinity, and injectivity and invertibility of the members of the operator
spectrum) that we seek to make connections between in Chapters
\ref{sec_LimOps} and \ref{sec_ellp}. New results in this section include:
a refinement of \cite[Proposition 1.2.9]{RaRoSiBook}, that, in the
so-called {\em perfect case}\index{perfect} (when both $P_n$ and its
adjoint converge strongly to the identity as $n\to\infty$), Fredholmness
implies invertibility of all the limit operators of $A$ for every $A\in
L(Y)$ (not just for $A\in L(Y,\P)$); the observation that, moreover,
Fredholmness is equivalent to invertibility when $A$ is {\em
self-similar}\index{self-similar operator}\index{operator!self-similar}
($A\in\opsp(A)$); and a strong existence result for self-similar
operators, that every {\em rich} operator\index{rich
operator}\index{operator!rich} (`richness' of $A$, defined in Section
\ref{subsec_limops} and characterised in Theorem \ref{lem_opsp} ,
guarantees the existence of limit operators) has a self-similar limit
operator. Section \ref{sec:coll_comp_op_sp} contains what are probably the
main new results of this text related to limit operators at an abstract
level (Theorems \ref{prop_appl1} and \ref{prop5.10}), obtained by applying
the generalised collectively compact operator theory of Section
\ref{subsec_collcomp}. These results, roughly speaking, relate
invertibility of an operator $A=I+K$, in the case when $K$ is Montel on
$(Y,s)$ (and $\opsp(K)$ is uniformly Montel in the sense of Definition
\ref{def_collcomp}), and stability of an approximating sequence
$A_n=I+K_n$ to $A$, to injectivity of the elements of $\opsp(A)$ or of the
members of a slightly larger set. One main theme of the rest of the text
is demonstrating the applicability of these results, and in fact they are
central to the proofs of Theorems \ref{prop_appl2}, \ref{prop_FC_JFA}, and
\ref{prop_NRBO}, and are important components of the proofs of Corollaries
\ref{cor_combine}, \ref{cor_combine2}, \ref{prop_W_selfadj},
\ref{cor_spessformula}, \ref{cor7.25} and \ref{cor_oned_wie}.

The largest chapter of this text is Chapter \ref{sec_ellp}, where
we apply results of Chapter \ref{sec_LimOps} to the concrete spaces
$Y=Y^p:=\ell^p(\Z^N,U)$ and $Y=Y^0:=c_0(\Z^N,U)$ with $N\in\N$,
$p\in [1,\infty]$ and $U$ a complex Banach space. This is still a
fairly general situation, as illustrated e.g. in Chapters
\ref{sec_Schroed} and \ref{sec_intop}, but it is concrete enough to
allow much more precise statements on Fredholmness and invertibility
at infinity of operators on $Y$ than was possible in the very
general setting of Chapter \ref{sec_LimOps}. We start with a short
section on the Banach algebra of almost periodic aka norm-rich
operators on $Y$ and their operator spectrum. These results are
picked up and improved in a later section on band-dominated and
norm-rich operators with the main result that $A=I+K$, with $K$
band-dominated, norm-rich and $\opsp(K)$ uniformly Montel, is
invertible on $Y^\infty$ iff all its limit operators are injective
on $Y^\infty$, i.e. Favard's condition holds. In Section
\ref{sec_dsa} we study interrelations of (semi-) Fred\-holmness and
index of an operator $A$ on $Y^\infty$ and its restriction $A_0$ to
$Y^0$. The main techniques here are to isometrically embed
$Y^\infty$ into $(Y^0)^{**}$ and $Y^1(U^*):=\ell^1(\Z^N,U^*)$ into
$(Y^\infty)^*$ and to show that, in terms of these embeddings, every
$A\in S(Y)$ is the restriction of $A_0^{**}$ to $Y^\infty$ and
$A_0^*$ is the restriction of $A^*$ to $Y^1(U^*)$. Further
connections are established in the case that $U$ has a predual
$U^\pa$ and $A$ is the adjoint of an operator $A^\pa$ on $U^\pa$. In
the following section we introduce $BDO(Y)$, the set of all
band-dominated operators on $Y$, which is the norm closure in $L(Y)$
of the set of all band operators or, equivalently, the smallest
Banach subalgebra of $L(Y)$ that contains all shift operators
(\ref{eq_shift1}) and all multiplication operators (\ref{eq_mult1})
with an operator-valued multiplier $b\in\ell^\infty(\Z^N,L(U))$.
Using results from Chapters \ref{chap_classes}, \ref{sec_LimOps} and
Section \ref{sec_dsa} we then derive criteria for and connections
between invertibility at infinity and Fredholmness of $A\in
BDO(Y^p)$ that are new or more complete than before if
$p\in\{0,1,\infty\}$. The final section of Chapter \ref{sec_ellp}
studies operators in the Wiener algebra $\W$ that we mentioned in
(\ref{eq_Wintro}) already. The attraction for studying this subset
of $BDO(Y)$ in more depth is that operators in $\W$ act boundedly on
every space $Y^p$ with $p\in\{0\}\cup[1,\infty]$ and neither their
invertibility nor invertibility at infinity, Fredholmness or the
Fredholm index depend on the choice of $p$. It results that, while,
for general operators $A\in BDO(Y^p)$ in the previous section, the
collective compactness arguments from Chapter \ref{sec_LimOps}
mainly contribute to the theory for $p=\infty$ (and, via duality, to
$p=0$ and $p=1$), now for operators $A\in\W$, the same results hold
for $p\in (1,\infty)$. For example, one gets that, for $N=1$ and
$K\in\W$ rich with $\opsp(K)$ uniformly Montel, the operator $A=I+K$
is Fredholm on any of the spaces $Y^p$ iff all of its limit
operators are injective on $Y^\infty$, i.e. Favard's condition
holds.

In Chapters \ref{sec_Schroed} and \ref{sec_intop} we illustrate our
main results for two operator classes of major interest in
applications. In Chapter \ref{sec_Schroed} we discuss discrete
Schr\"odinger operators (both self-adjoint and non-self-adjoint),
making links to other recent work (e.g Davies \cite{Davies} and Last
\& Simon \cite{LastSimon06}). In particular, we derive what appear
to be new characterisations of the spectrum and essential spectrum
(Theorem \ref{prop_schr_first} and (\ref{eqn_conj})), in terms of
the point spectrum in $\ell^\infty$ of the discrete operator and its
limit operators, and apply these results to Schr\"odinger operators
with potentials that are almost-periodic, perturbations of
almost-periodic, pseudo-ergodic in the sense of Davies
\cite{Davies}, and random, in this latter case reproducing results
of Trefethen, Contedini \& Embree \cite{TrefContEmb}. In Chapter
\ref{sec_intop} we demonstrate how the results of Chapter
\ref{sec_ellp} can be applied to continuous operators, deducing
criteria for Fredholmness and invertibility for members of a large
Banach algebra of integral operators on $\R^N$.

We conclude this text by a final chapter that contains a small
list of open problems that we consider highly interesting for the
further development of this field.
\newpage

{\bf Acknowledgements. } We would like to dedicate this text to
Professor Bernd Silbermann, who has been, and still is, an icon and
a very influential figure for both of us. When we started this
text he was still 62. But it soon turned out to be a rather long
term project, and our plans to dedicate it to his 65th birthday were
shattered when we couldn't stop ourselves from continuing to work
and write on it. So finally, being an extraordinary man, he is
getting the extraordinary honour of a dedication to his {\bf
\birthday th} birthday. {\it Happy birthday, Bernd!}

\bigskip

We would also like to acknowledge useful discussions of parts of
this work with colleagues including Professor Les Bunce (Reading),
Professor Brian Davies (King's College, London), Professor Peter
Kuchment (Texas A\&M University), Professor Steffen Roch
(Darmstadt), and of course Professor Bernd Silbermann (Chemnitz).
Finally, we gratefully acknowledge the financial support of a
Leverhulme Fellowship and a visiting Fellowship of the Isaac Newton
Institute for the first author and of a Marie Curie Fellowship
(MEIF-CT-2005-009758) and the Research Grant PERG02-GA-2007-224761
of the European Commission for the second author.\\[10mm]
Reading and Chemnitz\hfill Simon~N.~Chandler-Wilde\\[1mm]
6 April 2008\hfill Marko~Lindner\\[15mm]

To add to the above acknowledgements, we would also like to thank
Professor Barry Simon (Caltech) for helpful comments on the version
of this text that was submitted for publication in 2008. We have
taken the opportunity, before uploading a final version for
publication, to make amendments and additions in the light of these
comments, specifically to shorten  Remark \ref{rem7.4}, and to add
in Section \ref{sec_history} further historical comments on recent
developments in the application of limit-operator type ideas to
discrete Schr\"odinger and Jacobi operators. At the same time we
have taken the opportunity to add references in the bibliography to
two of our own recent papers/preprints which provide further
application of the ideas of Chapter \ref{sec_Schroed}. In
particular, \cite{Li:BiDiag} extends Example \ref{ex7.8} and the
results of \cite{TrefContEmb}, and \cite{CWRatLi10} studies a
related application, computing information about the spectrum of a
class of non-self-adjoint random operators and
matrices previously studied by Feinberg and Zee \cite{FeinZee99}. \\[10mm]
Reading and Chemnitz\hfill Simon~N.~Chandler-Wilde\\[1mm]
1 May 2010\hfill Marko~Lindner


\chapter{The Strict Topology} \label{sec_strict}
As usual, by $\N$\index{$N$@$\N$}, $\Z$\index{$Z$@$\Z$},
$\R$\index{$R$@$\R$} and $\C$\index{$C$@$\C$} we refer to the sets
of all natural, integer, real and complex numbers, respectively. We
also put $\N_0:=\{0,1,2,...\}$\index{$N0$@$\N_0$} and
$\R_+:=[0,\infty)$\index{$R+$@$\R_+$}.

Throughout, $Y$ will denote a real or complex Banach space, with
norm denoted by $\|\cdot\|$, and $L(Y)$ refers to the set of bounded
linear operators on $Y$. $\P
=(P_n)_{n=0}^\infty$\index{$P$@$\P$}\index{$P_n$} will denote a
sequence of operators in $L(Y)$ with the properties that:

(i) $\sup_n\|P_nx\|=\|x\|$, $x\in Y$;\index{condition (i)}

(ii) for every $m\in\N_0$ there exists $N(m)\ge m$ such that
\begin{equation} \label{pnpm}
P_nP_m = P_m = P_mP_n, \quad n\ge  N(m).\index{condition (ii)}
\end{equation}
Throughout we will write $n\gg m$\index{$\gg$} or $m\ll
n$\index{$\ll$} if $P_nP_k=P_k=P_kP_n$ for all $k\le m$. Note that
$l\ll m$ and $m\ll n$ imply $l\ll n$ since
$P_nP_k=P_nP_mP_k=P_mP_k=P_k=P_kP_m=P_kP_mP_n=P_kP_n$ for all $k\le
l$.

In \cite{RaRoSiBook} a bounded sequence $\P$ satisfying (ii) is
called an {\em increasing approximate projection} and an increasing
approximate projection satisfying (i) (or (i) with the `$=$'
replaced by a `$\ge$') is called an {\em approximate
identity}\index{approximate identity}. Thus $\P$ is an approximate
identity in the terminology of \cite{RaRoSiBook}.

It is easy to see 
that (i) implies that $\|P_n\|\le 1$ for all $n$. Moreover, if (i)
and (ii) hold, then $\|P_n\|=1$ for all sufficiently large $n$ and,
for all $x\in Y$, the limit $\lim_{n\to\infty}\|P_nx\|$ exists and
\begin{equation} \label{eq_lim}
\lim_{n\to\infty} \|P_nx\|=\|x\|,\quad x\in Y.
\end{equation}
For all $n\in\N_0$ let $Q_n := I-P_n$\index{$Q_n$} and note that
(ii) implies that, for all $m\in\N_0$,
\begin{equation} \label{qnqm}
Q_nQ_m = Q_n = Q_mQ_n, \quad n\ge N(m),
\end{equation}
and (i) that $\|Q_n\|\le 2$ for all $n\in\N_0$.

\begin{Remark} \label{rem21}
There are two possibilities inherent in the assumptions (i) and
(ii): either $P_n=I$ for some $n$ or $P_n\neq I$ for all $n$. In
the first case (ii) implies that $P_n=I$ for all sufficiently
large $n$. In the second case $Q_n\neq0$ for all $n$ and
(\ref{qnqm}) implies that $\|Q_n\|\ge 1$, so that $1\le \|Q_n\|\le
2$ for all $n$. Therefore, unless $P_n=I$ for all sufficiently
large $n$, it does not hold that $\|P_n-I\|\to 0$ as $n\to\infty$.
It may happen that $P_n$ converges strongly to $I$; this is the
case in Example \ref{e23} below and in Example \ref{ex_Lp}, for
$1\le p <\infty$, but not in Examples \ref{e21}, \ref{e21a}, \ref{e22} or
Example \ref{ex_Lp} for $p =\infty$.
\end{Remark}

\begin{Example} \label{e21} Let $Y=\l^\infty$, the Banach space of bounded
real- or complex-valued sequences $x = (x(m))_{m\in\Z}$, with norm
$\|x\|=\sup_m|x(m)|$. Define, for $x\in Y$, $n\in\N_0$,
$$
P_n x(m) = \left\{\begin{array}{ll}
  x(m), &|m|\le n, \\
  0, &|m|>n. \\
\end{array} \right.
$$
Then $\P=(P_n)$ satisfies (i) and (ii) with $N(m) = m$, so that
$P_n$ is a projection operator for each $n$. In this case
$\|Q_n\|=1$ for all $n$.
\end{Example}

\begin{Example} \label{e21a} Let $Y=\l^\infty$ as in Example \ref{e21}.
Define, for $n\in\N_0$, $P_n\in L(Y)$ by $P_nx(m)=x(m)$, for
$|m|\leq (3^n-1)/2$, and by the requirement that
$P_nx(m+3^n)=P_nx(m)$, for $m\in\Z$. Then $\P=(P_n)$ satisfies (i)
and (ii) with $N(m) = m$, so that $P_n$ is a projection operator for
each $n$. In this case $\|Q_n\|=2$ for all $n$.
\end{Example}

\begin{Example} \label{ex_Lp}
Let $Y=L^p(\R^N)$\index{$L^p(\R^N)$}, the Banach space of those
Lebesgue measurable real- or complex-valued functions $x$ on $\R^N$,
for which the norm $\|x\|_p$ is finite, where $\|x\|_p :=
(\int_{\R^N}|x(s)|^pds)^{1/p}$, for $1\le p<\infty$, and
$\|x\|_\infty := \esssup_{s\in\R^N}|x(s)|$. Define, for $x\in Y$,
$n\in\N_0$,
$$
P_n x(s) = \left\{\begin{array}{ll}
  x(s), &|s|\le n, \\
  0, &|s|> n. \\
\end{array} \right.
$$
Then $\P=(P_n)$ satisfies (i) and (ii) with $N(m) = m$, so that
$P_n$ is a projection operator for each $n$. In this case
$\|Q_n\|=1$ for all $n$.
\end{Example}

\begin{Example} \label{e22} Let $Y=BC(\R^N)$, the Banach space of bounded
continuous real- or complex-valued functions on $\R^N$, with norm
$\|x\| = \sup_{s\in\R^N}|x(s)|$. Choose $\chi\in BC(\R)$ with
$\|\chi\|=1$ and $\chi(t)=0$, $t\le 0$, $=1$, $t \ge 1$. Define, for
$n\in\N_0$ and $x\in Y$,
$$
P_nx(s) = \chi(n+1-|s|)x(s), \quad s\in\R^N.
$$
Then $\P=(P_n)$ satisfies (i) and (ii) with $N(m)=m+1$. In this
case $\|Q_n\|=\|1-\chi\|$.
\end{Example}

\begin{Example} \label{e23} Let $Y=C[0,1]$ with $\|x\| = \sup_{0\le s\le
1}|x(s)|$ and let $P_nx$ denote the piecewise linear function
which interpolates $x$ at $j/2^n$, $j=0,1,...,2^n$. Then
$\P=(P_n)$ satisfies (i) and (ii) with $N(m)=m$.
\end{Example}

Note that in Examples \ref{e21} and \ref{e23} $P_n$ has
finite-dimensional range for each $n$, so that $P_n\in
K(Y)$\index{$K(Y)$}, the set of compact linear operators on $Y$.
(The significance of this is discussed in Remark \ref{PnCompact}
below.)

Throughout, for $(x_n)\subset Y$, $x\in Y$, we will write $x_n\to x$
(as $n\to\infty$) if $(x_n)$ converges to $x$ in the norm topology,
i.e. $\|x_n-x\|\to0$ as $n\to\infty$. We will also be concerned with
convergence in weaker topologies, defined in terms of semi-norms on
$Y$ that are related to the sequence $\P$. For $n\in\N_0$ define
$|\cdot|_n:Y\to\R_+$\index{$\mid\cdot\mid_n$} by

\begin{equation} \label{eq_seminorm}
| x|_n := \max_{0\le m\le n}\|P_mx\|, \quad x\in Y.
\end{equation}

It is easy to check that $|\cdot|_n$ is a semi-norm on $Y$, and
(i) and (ii) imply that, for $m<n$, $x\in Y$,
$$
| x|_m\le |x|_n\le \|x\|
$$
and
$$
\|x\| = \sup_n|x|_n = \lim_{n\to\infty}|x|_n.
$$

This last equation implies that the family of semi-norms,
$\{|\cdot|_n:n\in\N_0\}$, is separating. We will term the metrisable
topology generated by this family of semi-norms the {\em local
topology}\index{local topology}\index{topology!local}: with this
topology, $Y$ is a separated locally convex topological vector space
(TVS). By definition, a sequence $(x_n)$ converges to $x$ in the
local topology if and only if $|x_n-x|_m\to0$ as $n\to\infty$ for
all $m$, i.e. if and only if
\begin{equation} \label{ltc}
P_mx_n\to P_mx \mbox{ as } n\to\infty, \mbox{ for all } m.
\end{equation}

We will see below that, unless $P_n=I$ for each $n$, $Y$ is not
complete in the local topology. Let $X$ denote the Fr\'echet space
that is the completion of $Y$ in the local topology. The linear
operators $P_n$, which are continuous as mappings from $Y$ equipped
with the local topology to $Y$ equipped with the norm topology, have
unique extensions to continuous linear operators from $X$ to $Y$
\cite{Robertson}. The equations (\ref{pnpm}) hold for the extended
operators $P_n:X\to Y$, by the continuity of the operators $P_n$ and
since $Y$ is dense in its completion $X$. Extending the definition
of the semi-norm $|\cdot|_n$ from $Y$ to $X$ using equation
(\ref{eq_seminorm}), clearly $\{|\cdot|_n:n\in\N_0\}$ generates the
extension of the local topology from $Y$ to $X$.

Let $\hat Y$\index{$Yhat$@$\hat Y$} and $\tilde
Y$\index{$Ytil$@$\tilde Y$} denote the linear subspaces of $X$,
$$
\hat Y := \{x\in X: \|x\|:= \sup_n|x|_n <\infty\},
$$
and
$$
\tilde Y := \cup_{n\in \N_0} P_n(X),
$$
and note that, for every $n$, $\tilde Y = \cup_{m\ge n}
P_m(X)\subset Y$ by (ii). Let $Y_0$\index{$Y_0$} denote the
norm-closure of $\tilde Y$. It follows from \cite[Theorem 2.1
(ii)]{ChanZh02} and the completeness of $X$ that, equipped with the
norm $\|\cdot\|$, $\hat Y$ is a Banach space with $Y$ and
$Y_0\subset Y$ as closed subspaces.

\begin{Lemma}\cite[Lemma 1.1.20]{RaRoSiBook} \label{lem_Y0}
$Y_0=\{x\in Y: Q_nx\to 0 \mbox{ as } n\to\infty\}$ so that $Y=Y_0$
iff $P_n\in L(Y)$ converges strongly to $I$ as $n\to\infty$.
\end{Lemma}

\begin{Example} \label{ex_Lp2}
Let $Y$ and $\P$ be as in Example \ref{ex_Lp}. Then $X=L^p_{\rm
loc}(\R^N)$\index{$L^p_{\rm loc}(\R^N)$}, $\tilde Y$ is the set of
compactly supported functions $x\in L^p(\R^N)$, $\hat
Y=Y=L^p(\R^N)$, $Y_0=L^p(\R^N)$, if $1<p<\infty$, and
$$
Y_0=\{x\in L^\infty(\R^N): \esssup_{|s|>a}|x(s)|\to 0 \mbox{ as }
a\to\infty\}
$$
if $p=\infty$.
\end{Example}

\begin{Example} \label{ex_BUC}
Let $BUC(\R^N)\subset BC(\R^N)$\index{$BUC(\R^N)$} denote the set of
bounded uniformly continuous functions on $\R^N$, let
$C_L(\R^N)\subset BUC(\R^N)$ denote the set of those $x\in C(\R^N)$
for which $x(s)\to x_\infty$ as $|s|\to\infty$, uniformly in $s$,
for some constant $x_\infty$. Let $C_0(\R^N)$ denote the set of
those $x\in C_L(\R^N)$ for which $x(s)\to 0$ as $|s|\to\infty$, and
let $C_C(\R^N)$ denote the set of compactly supported continuous
functions. Note that $C_0(\R^N),\, C_L(\R^N),\, BUC(\R^N)$ are all
closed subspaces of $BC(\R^N)$, equipped with the usual norm.

Let $Y$ denote one of $C_0(\R^N),\,C_L(\R^N),\,BUC(\R^N)$ or
$BC(\R^N)$ and let $\P$ be as defined in Example \ref{e22}. Then
$\hat Y=BC(\R^N),\,\tilde Y=C_C(\R^N)$ and $Y_0=C_0(\R^N)$.
\end{Example}

\begin{Example} \label{ex_ap}
For $x\in\ell^\infty$ and $k\in\Z$ let $x_k\in\ell^\infty$ be
defined by $x_k(m) = x(m-k)$, $m\in\Z$. Call $x\in\ell^\infty$ {\em
almost periodic}\index{almost periodic sequence} (e.g.\
\cite[Definition 3.58]{LiBook}) if
$\{x_k:k\in\Z\}\subset\ell^\infty$ is relatively compact. Let
$\ell^\infty_{\sf AP}\subset
\ell^\infty$\index{$\ell^pa$@$\ell^\infty_{\sf AP}$} denote the set
of almost periodic functions. Let $Y=\ell^\infty_{\sf AP}$ and
define $\P$ as in Example \ref{e21a}, noting that
$P_n(\ell^\infty_{\sf AP})\subset \ell^\infty_{\sf AP}$ for every
$n\in\N_0$. Then $\hat Y = \ell^\infty$ and $Y_0$ is a strict
subspace of $Y$. In particular, if $y(m)=\exp(2\pi iam)$, $m\in\Z$,
with $a$ irrational, then $y\in Y$ (see e.g.\ \cite[Lemma
3.64,Proposition 3.65]{LiBook}) but $\|P_ny-y\|=2$, for $n\in\N_0$,
so that $y\not\in Y_0$ by Lemma \ref{lem_Y0}.
\end{Example}

We will also be interested in a third topology on $\hat Y\supset Y$,
intermediate between the local and norm topologies. Given a positive
null-sequence $a:\N_0\to (0,\infty)$, define
$$
| x|_a := \sup_n a(n)|x|_n.\index{$\mid\cdot\mid_a$}
$$
Then $\{|\cdot|_a:a$ is a positive null sequence$\}$ is a second
separating family on $\hat Y$ and generates another separated
locally convex topology on $\hat Y$ which, by analogy with
\cite{Buck58}, we will term the {\em strict topology}\index{strict
topology}\index{topology!strict}. For $(x_n)\subset \hat Y$, $x\in
\hat Y$, we will write $x_n\sto x$\index{$\sto$} if $x_n$ converges
to $x$ in the strict topology, i.e. if $|x_n-x|_a\to 0$ as
$n\to\infty$ for every null sequence $a$.

The topology we have called the strict topology is termed the
$\beta$ topology in \cite{ChanZh02}. Various properties of the
$\beta$/strict topology are shown in \cite[Theorem 2.1]{ChanZh02},
in large part adapting arguments from \cite{Buck58}. The properties
that we need for our arguments are summarised in the next lemma. As
usual \cite{Robertson} we will call a set $S$ in a TVS {\em
bounded}\index{bounded} if it is absorbed by every neighbourhood of
zero, {\em totally bounded}\index{totally bounded} if, for every
neighbourhood of zero, $U$, there exists a finite set
$\{a_1,...,a_N\}$ such that $S\subset \cup_{1\le j\le N} (a_j+U)$,
and {\em compact}\index{compact} if every open cover of $S$ has a
finite subcover. Every totally bounded set is bounded
\cite{Robertson}.

\begin{Lemma} \label{lem_strict_prop}
\begin{itemize}
\item[(i)] In $\hat Y$ the bounded sets in the strict topology and the
norm topology are the same.

\item[(ii)] On every norm-bounded subset of $\hat Y$ the strict
topology coincides with the local topology.

\item[(iii)] A sequence $(x_n)\subset \hat Y$ is convergent in  the
strict topology iff it is convergent in the local topology and is
bounded in the norm topology, so that
\begin{equation} \label{eqstrict}
x_n\sto x\quad \Leftrightarrow\quad \sup_n\|x_n\|<\infty \mbox{ and
} P_mx_n\to P_mx \mbox{ as } n\to\infty, \mbox{ for all } m.
\end{equation}

\item[(iv)] A norm-bounded subset of $\hat Y$ is closed in the strict
topology iff it is sequentially closed.

\item[(v)] A sequence in $\hat Y$ is Cauchy in the strict topology
iff it is Cauchy in the local topology and bounded in the norm
topology.

\item[(vi)] Let $S\subset \hat Y$. Then the following statements are
equivalent:
\begin{itemize}
\item[(a)] $S$ is totally bounded in the strict topology.
\item[(b)] $S$ is norm-bounded and totally bounded in the local
topology. \item[(c)] Every sequence in $S$ has a  subsequence that
is Cauchy in the strict topology.
\end{itemize}
\end{itemize}

\end{Lemma}
\Proof For (i)-(iv) see \cite[Theorem 2.1]{ChanZh02}. Part (v)
follows from (iii) on noting that, if $f_j:\N\to\N,\,j=1,2$, are
such that $f=(f_1,f_2):\N\to\N^2$ is a bijection, then, defining
$y_n:=x_{f_1(n)}-x_{f_2(n)}$, $(x_n)$ is Cauchy iff $(y_n)$ is
convergent to zero.

To see (vi) note that if (a) holds then $S$ is bounded in the strict
topology and so in the norm topology by (ii). Also, $S$ is totally
bounded in the (coarser) local topology. Thus (b) holds. Conversely,
if (b) holds, $U$ is a neighbourhood of zero in the strict topology,
$M :=  \sup_{x\in S} \|x\|$, and $B:=\{x:\|x\|\le 2M\}$, then, by
(ii), there exists a neighbourhood of zero in the local topology,
$U^\prime$, such that $U\cap B = U^\prime \cap B$. Further, there
exists a finite set $\{s_1,...,s_N\}\subset S$ such that $ S \subset
\bigcup_{1\le j\le N} (s_j+U^\prime)$. It follows that
$$
S \subset \bigcup_{1\le j\le N} (s_j+U^\prime\cap B) =
\bigcup_{1\le j\le N} (s_j+U\cap B) \subset \bigcup_{1\le j\le N}
(s_j+U).
$$
Thus also (b) $\Rightarrow$ (a).

To see that (b) and (c) are equivalent note that, as the local
topology is metrisable, $S$ is totally bounded in the local topology
iff every sequence in $S$ has a subsequence that is Cauchy in the
local topology. Further, by (v), a sequence is Cauchy in the strict
topology iff it is Cauchy in the local topology and norm-bounded.
\Proofend

\noindent Note that it follows from (ii) that the linear operators
on $Y$ that are bounded with respect to the strict topology (map
bounded sets onto bounded sets) are precisely the members of
$L(Y)$.

Let $E$ denote one of $\hat Y,\,Y$ and $Y_0$. When it is necessary
to make a clear distinction we will denote the TVS consisting of $E$
(considered as a linear space) equipped with the strict topology by
$(E,s)$\index{$Ya$@$(Y,s)$} and will denote the TVS (and Banach
space) consisting of $E$ with the norm topology as
$(E,\|\cdot\|)$\index{$Yaa$@$(Y,\parallel\cdot\parallel)$}.

\begin{Lemma} \label{lem_topo}
If $P_n=I$ for some $n$, then the local, strict, and norm topologies
coincide on $\hat Y$. If $P_n\neq I$ for all $n$, then:

(a) on $\hat Y$ the local topology is strictly coarser than the
strict topology which is strictly coarser than the norm topology;
and

(b) $\hat Y$, equipped with the local topology, is not complete,
while $\hat Y$ equipped with the strict topology is complete and
non-metrisable.
\end{Lemma}
\Proof It is easy to see that any set open in the local topology is
open in the strict topology and that any set open in the strict
topology is open in the norm topology. If $P_n=I$ for some $n$ then
the converse statements clearly hold, as at least one of the
semi-norms defining each topology coincides with the norm. Thus the
topologies coincide.

If $P_n\neq I$ for any $n$ then there exists $(x_n)$ such that
$\|Q_nx_n\|=1$ for all $n$. For all $m$, $P_mQ_nx_n = 0$ for all
sufficiently large $n$, by (ii). Clearly $Q_nx_n\not\to0$, but it
follows from (\ref{eqstrict}) that $Q_nx_n\sto 0$ as $n\to\infty$.
Thus the strict and norm topologies are distinct. To see that the
local and strict topologies are distinct, note that $nQ_nx_n$
converges to zero in the local topology but
$\|nQ_nx_n\|=n\to\infty$ so that, by (\ref{eqstrict}),
$nQ_nx_n\not\sto0$.

If $\hat Y$ equipped with the local topology were complete it would
be a Fr\'echet space and it would follow from the open mapping
theorem \cite{Rudin91} applied to the identity operator that the
local and norm topologies coincide.

Let $Y^*$ denote the completion of $\hat Y$ in the strict topology.
Then $Y^*\subset X$, since $\hat Y\subset X$ and $X$ is complete in
the coarser local topology. Suppose $Y^*\ne\hat Y$. Then there
exists $x\in Y^*$ with $|x|_n\to\infty$ as $n\to\infty$. Let
$b_n:=1/2\,\max(1,|x|_n)$, $a_n:=1/b_n$, and $a=(a_0,a_1,...)$. Then
$y\in Y^*$ and $|x-y|_a<1$ imply that $|y|_n>|x|_n/2$ for all
sufficiently large $n$, so that $\{y\in\hat Y:
| x-y|_a<1\}=\varnothing$. This is a contradiction, for $\hat Y$ is
dense in its completion. \Proofend

By definition, $Y_0$ is the completion of $\tilde Y$ in the norm
topology and we have seen that $Q_nx\to 0$ if $x\in Y_0$ in Lemma
\ref{lem_Y0}. The next lemma states corresponding results for the
strict topology.

\begin{Lemma} \label{lem_Qny}
For $y\in \hat Y$, $Q_ny\sto 0$ as $n\to\infty$. Further $\hat Y$ is
the completion of $\tilde Y$ in the strict topology and $Y=\hat Y$
iff $Y$ is sequentially complete in the strict topology.
\end{Lemma}
\Proof If $y\in\hat Y\subset X$ then $P_ny\in\tilde Y\subset Y$ and
$P_mP_ny=P_my$ for all sufficiently large $n$. Further, by (i),
$\|P_ny\|\le\|y\|$. Thus, by Lemma \ref{lem_strict_prop} (iii),
$P_ny\sto y$. Thus the completion of $\tilde Y$ contains $\hat Y$
and in fact coincides with $\hat Y$ since $\hat Y$ is complete by
Lemma \ref{lem_topo} (ii). Since $\hat Y$ is complete, $Y\subset\hat
Y$ is sequentially complete iff it is sequentially closed. But,
since $P_ny\sto y$ for every $y\in\hat Y$, this holds iff $Y=\hat
Y$. \Proofend

As usual, we will call a subset of a topological space {\em
relatively compact}\index{relatively
compact}\index{compact!relatively} if its closure is compact. We
will call a subset of a topological space {\em relatively
sequentially compact}\index{relatively sequentially
compact}\index{sequentially
compact!relatively}\index{compact!sequentially!relatively} if every
sequence in the subset has a subsequence converging to a point in
the topological space.

\begin{Lemma} \label{lem_compact}
Let $S\subset Y$. Then $S$ is compact in $(Y,s)$ iff it is
sequentially compact. Further,
$$
\mbox{(a)} \Leftrightarrow \mbox{(b)} \Rightarrow \mbox{(c)}
\Leftrightarrow \mbox{(d)}
$$
where (a)-(d) are the statements:
\begin{itemize}
\item[(a)] $S$ is relatively compact in the strict topology.

\item[(b)] $S$ is relatively sequentially compact in the strict
topology.

\item[(c)] $S$ is totally bounded in the strict topology.

\item[(d)] $S$ is norm-bounded and $P_n(S)$ is relatively compact
in the norm topology for each $n$.
\end{itemize}
If $(Y,s)$ is sequentially complete (i.e. $Y=\hat Y$) then (a)--(d)
are equivalent.
\end{Lemma}
\Proof To show that compactness (relative compactness) of $S$ is
equivalent to sequential compactness (relative sequential compactness) it
is enough to show this in the strict topology restricted to $\bar S$, the
closure in $(Y,s)$ of $S$. But, if $S$ is relatively sequentially compact
or relatively compact then it is bounded and so $\bar S$ is bounded. But,
by (ii) of Lemma \ref{lem_strict_prop}, the strict topology coincides with
the metrisable local topology on bounded sets, and in metric spaces
compactness and sequential compactness coincide. Thus the first statement
of the theorem holds and also (a) $\Leftrightarrow$ (b). That (b) implies
(c), and the converse if $Y$ is sequentially complete, is immediate from
(vi) of Lemma \ref{lem_strict_prop}. If (c) holds then, also by (vi) of
Lemma \ref{lem_strict_prop}, $S$ is norm-bounded and every sequence in $S$
has a subsequence that is Cauchy in the strict topology. Since $P_n$ is
continuous from $(\hat Y,s)$ to $(\hat Y,\|\cdot\|)$ and $(\hat
Y,\|\cdot\|)$ is complete, this implies that $P_n(S)$ is relatively
compact in the norm topology. Finally, suppose (d) holds and take an
arbitrary bounded sequence $(x_n)\subset\hat Y$. Choose a subsequence
$(x_n^{(1)})$ such that $P_1x_n^{(1)}$ norm-converges as $n\to\infty$.
From $(x_n^{(1)})$ choose a subsequence $(x_n^{(2)})$ such that
$P_2x_n^{(2)}$ norm-converges, and so on. Then $(y_n)$, with
$y_n:=x_n^{(n)}$, which is bounded and Cauchy in the local topology is
Cauchy in the strict topology by Lemma \ref{lem_strict_prop} (iv). Thus
every sequence in $S$ has a subsequence that is Cauchy in the strict
topology, so that, by Lemma \ref{lem_strict_prop} (vi), (c) holds.
\Proofend

\begin{Example} \label{ex_arzasc}As an important example of relative
compactness in the strict topology, consider the situation of
Examples \ref{e22} and \ref{ex_BUC}, where $\hat Y=Y=BC(\R^N)$ and
$(Y,s)$ is just $Y$ with the standard strict topology of
\cite{Buck58}. It follows from the equivalence of (a) and (d) in the
above lemma and the Arzela-Ascoli theorem that $S\subset Y$ is
relatively compact in $(Y,s)$ iff $S$ is bounded and equicontinuous.
Recall that $S\subset Y$ is equicontinuous if
$$
\sup_{x\in S} |x(s)-x(t)| \to 0 \mbox{ as } t\to s,
$$
for all $s\in\R^N$.
\end{Example}

\begin{Remark} \label{PnCompact} As the following corollary of the above lemma
already indicates,
many of the results we obtain in the text
will simplify and become more complete in the case that $P_n\in
K(Y)$ for all $n$. We note that, by (ii), $P_n$ is compact for all
$n$ if, for every $N$, $P_n$ is compact for some $n>N$.
\end{Remark}

\begin{Corollary} \label{lem_finite}
If $Y$ is sequentially complete in the strict topology (i.e. $Y=\hat
Y$) and $P_n\in K(Y)$ for all $n$, then $S\subset Y$ is relatively
compact in the strict topology iff it is norm-bounded.
\end{Corollary}


\chapter{Classes of Operators} \label{chap_classes}

We have introduced already $L(Y)$\index{$L(Y)$} and
$K(Y)$\index{$K(Y)$}, the sets of linear operators that are,
respectively, bounded and compact on $(Y,\|\cdot\|)$. We have noted
that $L(Y)$ coincides with the set of linear operators that are
bounded on $(Y,s)$. Let $C(Y)$\index{$C(Y)$} and
$S(Y)$\index{$S(Y)$} denote the sets of those linear operators that
are, respectively, continuous and sequentially continuous on
$(Y,s)$. Thus $A\in S(Y)$ if and only if, for every sequence
$(x_n)\subset Y$ and $x\in Y$,
\begin{equation} \label{scont}
x_n\sto x\quad \Rightarrow\quad Ax_n\sto Ax.
\end{equation}
Let $SN(Y)$\index{$SN(Y)$} denote the set of those linear operators
that are sequentially continuous from $(Y,s)$ to $(Y,\|\cdot \|)$,
so that $A\in SN(Y)$ iff
\begin{equation} \label{sncont}
x_n\sto x\quad \Rightarrow\quad Ax_n\to Ax.
\end{equation}
We remark that the operators in $S(Y)$ and $SN(Y)$ are precisely
those termed $s-$continuous and $sn-$continuous, respectively, in
\cite{ArensCWHas1}.

From standard properties of topological vector spaces \cite[Theorems
A6 and 1.30]{Rudin91}, and Lemma \ref{lem_topo}, it follows that
$C(Y)\subset S(Y) \subset L(Y)$. In fact we have the following
stronger result.

\begin{Lemma}
$C(Y)=S(Y)$.
\end{Lemma}
\Proof
Let $C(\hat Y)$, $S(\hat Y)$ denote the sets of linear
operators on $\hat Y$ that are, respectively, continuous and
sequentially continuous. For $n\in \N_0$ let $Y_n$ denote the linear
subspace of $Y$,
\begin{equation} \label{Ymdef}
Y_n := \{x\in Y:|x|_n=0\} = \{x\in Y: P_mx=0,\;\;0\le m\le n\}.
\end{equation}
Note that, by (ii), for every $m\in \N_0$, $Q_n(\hat Y)\subset Y_m$
for all sufficiently large $n$, and, for all $x\in Y$, $\|x-Q_n x\|
= \|P_n x\| \le |x|_n$. Thus Assumption A$^\prime$ of
\cite{ChanZh02} holds and it follows from \cite[Theorem
3.7]{ChanZh02} that $C(\hat Y)=S(\hat Y)$.

By Lemma \ref{lem_Qny}, the sequential closure of $Y\subset \hat Y$
in the strict topology is $\hat Y$. In Lemma \ref{lem3.13} we will
show that every $A\in S(Y)$ has an extension $\hat A\in S(\hat Y)$
defined by $\hat A x = \lim_{n\to\infty} A P_n x$, where the limit
exists in the strict topology. Then $\hat A\in C(\hat Y)$ and $A =
\hat A|_Y\in C(Y)$. \Proofend

\noindent In view of this lemma it holds that
\begin{equation} \label{eqinclus1}
SN(Y) \subset C(Y)=S(Y)\subset L(Y).
\end{equation}

As Lemmas \ref{lem1}-\ref{lem3} below clarify, in general $SN(Y)$ is
a strict subset of $S(Y)$. The following example shows that
$S(Y)\neq L(Y)$ in general, indeed that $A$ may be compact on
$(Y,\|\cdot\|)$ but not sequentially continuous on $(Y,s)$.

\begin{Example} \label{exnotinS}
Let $Y=\l^\infty$ and $P_n$ be as in Example \ref{e21}. Let $c_\l^+$
denote the set of those $x\in \l^\infty$ for which
$\lim_{m\to+\infty}x(m)$ exists. By the Hahn-Banach theorem a
bounded linear functional $\l_+:Y\to\C$ exists such that $\l_+(x) =
\lim_{m\to+\infty}x(m)$, $x\in c_\l^+$. Define $A\in L(Y)$ by $Ax =
\l_+(x)y$, $x\in Y$, where $y\in Y$ is non-zero and fixed. Then the
range of $A$ is one-dimensional so that $A\in K(Y)\subset L(Y)$.
However, defining $x= (...,1,1,1,...)$ and $x_n=Q_nx$, $Q_nx\sto0$
as $n\to\infty$ but $AQ_n x=1$ for all $n$. Thus $A\not\in S(Y)$.
\end{Example}

The following lemmas provide alternative characterisations of the
classes $SN(Y)$ and $S(Y)$ and shed some light on the relationship
with $K(Y)$.

\begin{Lemma} \label{lem1}
$A\in SN(Y)$ iff $A\in L(Y)$ and $\|AQ_n\|\to 0$ as $n\to\infty$.
\end{Lemma}
\Proof Suppose $A\in SN(Y)$. Then $A\in L(Y)$.  To see that also
$\|AQ_n\|\to 0$ as $n\to\infty$, suppose that this does not hold.
Then there is a bounded sequence $(x_n)\subset Y$ such that
$AQ_nx_n\not\to0$. But this is impossible as $Q_nx_n\sto 0$, and
hence $\|AQ_nx_n\|\to 0$ as $n\to\infty$, which is a
contradiction.

For the reverse implication, take an arbitrary sequence
$(x_n)\subset Y$ with $x_n\sto 0$ as $n\to\infty$. Then $\|x_n\|$
is bounded and $\|P_mx_n\|\to 0$ as $n\to\infty$ for every $m$.
Now, for every $m$ and $n$,
\begin{eqnarray*}
\|Ax_n\|&\le&\|AP_mx_n\|\ +\ \|AQ_mx_n\|\\
&\le&\|A\|\|P_mx_n\|\ +\ \|AQ_m\|\sup_n\|x_n\|
\end{eqnarray*}
holds, where $\|AQ_m\|$ can be made as small as desired by
choosing $m$ large enough, and $\|P_mx_n\|$ tends to zero as
$n\to\infty$. \Proofend

\begin{Lemma} \label{lem3}
$A\in S(Y)$ iff $A\in L(Y)$ and $P_mA\in SN(Y)$ for every $m$.
\end{Lemma}
\Proof If $A\in S(Y)$ then $A\in L(Y)$. The rest trivially follows
from
\[
Ax_n\sto 0\ \textrm{ as }\ n\to\infty\qquad\iff\qquad
\|P_mAx_n\|\to 0\ \textrm{ as }\ n\to\infty\ \forall m
\]
for every bounded operator $A$ and every bounded sequence
$(x_n)\subset Y$. \Proofend

\begin{Corollary} \label{cor4}
$A\in S(Y)$ iff $A\in L(Y)$ and $\|P_mAQ_n\|\to 0$ as $n\to\infty$
for every $m$.
\end{Corollary}

\begin{Example} \label{ex_conv} Let $Y=L^p(\R^N)$ and $P_n$ be defined as in
Example \ref{ex_Lp}. Let $\kappa\in L^1(\R^N)$ and let $A$ be the
linear integral operator, a so-called {\em convolution
operator}\index{convolution operator}\index{operator!convolution},
defined by
$$
A x(s) = \int_{\R^N} \kappa(s-t) x(t) dt, \quad s\in\R^N.
$$
Then, by Young's inequality, $A\in L(Y)$ with $\|A\|\le
\|\kappa\|_1$. Further, $A\in S(Y)$, since $\|P_mAQ_n\|\to 0$ as
$n\to\infty$, for all $m$. To see this, note that, for every $m$,
$P_mAQ_n=0$ for all sufficiently large $n$ if $\kappa\in
C^\infty_0(\R^N)$, and then use the density of $C^\infty_0(\R^N)$ in
$L^1(\R^N)$.
\end{Example}

\begin{Lemma} \label{lem2}
$S(Y)\cap K(Y) \subseteq SN(Y)$, with equality if and only if
$P_n\in K(Y)$ for all $n$.
\end{Lemma}
\Proof Suppose $A\in S(Y)\cap K(Y)$. Take an arbitrary sequence
$(x_n)\subset Y$ with $x_n\sto 0$ as $n\to\infty$. From $A\in
S(Y)$ we conclude that $Ax_n\sto 0$ as $n\to\infty$. Since
$\{x_n\}$ is bounded and $A$ is compact, we know that $\{Ax_n\}$
is relatively compact, so every subsequence of $(Ax_n)$ has a
norm-convergent subsequence, where the latter has limit $0$ since
$Ax_n\sto 0$ as $n\to\infty$. Of course, this property ensures
that $Ax_n$ itself norm-converges to $0$.

To see when equality holds consider that, by (ii), it holds for every $m$
that $P_mQ_n=0$ for all sufficiently large $n$. Thus, by Lemma \ref{lem1},
$P_m\in SN(Y)$ for all $m$. So clearly $SN(Y)\not\subset K(Y)$ if $P_m$ is
not compact for all $m$. If $P_m$ is compact for all $m$ and $A\in SN(Y)$
then, by Lemma \ref{lem1} again, $A$ is the norm limit $\lim_{m\to\infty}
AP_m$, with $AP_m$ compact, so that $A$ is compact and $SN(Y)\subset
K(Y)$. Thus equality holds iff $P_m\in K(Y)$ for all $m$. \Proofend

Recall that (e.g.\ \cite{AnseloneBook}) if $K\in K(Y)$ and $A_n$ converges
strongly to $A$ then, since pointwise convergence is uniform on compact sets,
$\|(A_n-A)K\|\to 0$. This and Lemma \ref{lem2} have the following implication.

\begin{Lemma} \label{lem_strong}
If $P_n$ converges strongly to $I$ then $\|Q_nK\|\to 0$ as $n\to 0$ for all
$K\in K(Y)$, while if $P_n^*$ converges strongly to $I^*$ ($P_n^*$ and $I^*$
the adjoints of $P_n$ and $I$) then $\|KQ_n\|\to 0$ as $n\to 0$ for all $K\in
K(Y)$, so that $K(Y)\subset SN(Y)$.
\end{Lemma}

\begin{Lemma} \label{lem_sny}
Let $A$ be a linear operator on $Y$. Then the following statements
are equivalent.

(a) $A\in SN(Y)$.

(b) $A\in L(Y)$ and there is a neighbourhood of zero, $U$, in
$(Y,s)$, for which $A(U)$ is norm-bounded, in fact for which
$\sup_{x\in U}\|AQ_n x\|\to 0$ as $n\to\infty$.

(c) $A$ is a continuous mapping from $(Y,s)$ to $(Y,\|\cdot\|)$.
\end{Lemma}
\Proof That continuity implies sequential continuity, so that (c)
$\Rightarrow$ (a), is standard \cite{Rudin91}.

Suppose that (a) holds. Then $A\in L(Y)\supset SN(Y)$. Moreover, by Lemma
\ref{lem1}, $\|AQ_n\|\to 0$ as $n\to\infty$. Choose positive integers
$n_1\ll n_2\ll...$ such that
$$
\|AQ_j\| \le 4^{-m}, \quad j \ge n_m,
$$
for $m\in\N$. Then, for $x\in Y$, $m\in \N$, and $n_m\le j\le
n_{m+1}$, since $P_{n_N}Q_j x\sto Q_j x$ as $N\to\infty$ and $A\in
SN(Y)$,
\begin{eqnarray*}
\|AQ_jx\| &=& \lim_{N\to\infty}\|AP_{n_N}Q_jx\|\\
& = & \lim_{N\to\infty}\|A(Q_j-Q_{n_N})x\|\\
& = & \lim_{N\to\infty}\|\sum_{i=m+1}^{N-1}A(Q_{\tilde
n_i}-Q_{\tilde n_{i+1}})x\|,
\end{eqnarray*}
where $\tilde n_i := n_i$, for $i> m+1 $, and $\tilde n_{m+1} :=
j$. Further, for $i\ge m+1$,
$$
\|A(Q_{\tilde n_i}-Q_{\tilde n_{i+1}})x\| = \|AQ_{\tilde
n_i}P_{\tilde n_{i+1}}x\| \le 4^{1-i}\|P_{\tilde n_{i+1}}x\|.
$$
Now, define $n_0:= 0$ and $a_m := 2^{-i}$, for $n_i\le m <
n_{i+1}$, $i\in\N_0$, and set $a := (a_0,a_1,...)$ and $U:=
\{x:|x|_a <1\}$. Then, from the above inequalities, we see that,
for $x\in U$, $m\in \N$, and $n_m\le j\le n_{m+1}$, it holds that
$$
\|AQ_jx\| \le \sum_{i=m+1}^{\infty} 4^{1-i}\|P_{\tilde n_{i+1}}x\|
\le \sum_{i=m+1}^{\infty} 4^{1-i}/a_{\tilde n_{i+1}} \le 2^{3-m},
$$
since $a_{\tilde n_{i+1}}\ge a_{ n_{i+1}} = 2^{-(i+1)}$. Thus
$$
\sup_{x\in U}\|AQ_jx\|\to 0 \mbox{ as } j\to\infty.
$$
Moreover, for $x\in U$, $\|P_{n_1}x\|\le 1/a_{n_1} = 2$, so that
$$
\sup_{x\in U} \|Ax\| \le \sup_{x\in U} \|AP_{n_1}x\| + \sup_{x\in
U} \|AQ_{n_1}x\| \le 2\|A\|+4.
$$
Since $U$ is a neighbourhood in $(Y,s)$, we have shown that (a)
$\Rightarrow$ (b).

Now suppose that (b) holds, so that $A\in L(Y)$ and $A(U)$ is
norm-bounded, with $U$ a neighbourhood of 0 in $(Y,s)$. To show
that (c) holds it is enough to show that $A$ is continuous at 0.
But if $V$ is a neighbourhood of 0 in $(Y,\|\cdot\|)$ then $V$
contains $\lambda A(U)$, for some $\lambda>0$, and so $A^{-1}(V)$
contains $\lambda U$. Thus (b) $\Rightarrow$ (c).
\Proofend

Recall that $K(Y)$ stands for the set of compact operators on $Y$.
Following \cite[\S 1.1.2]{RaRoSiBook}, we denote the set of all
$K\in L(Y)$ which are subject to
\begin{equation} \label{KYP_def}
\|KQ_n\|\to 0\qtext{and}\|Q_nK\|\to 0\qtext{as} n\to\infty
\end{equation}
by $K(Y,\P)$\index{$K(Y,\P)$}. Moreover, let
$L(Y,\P)$\index{$L(Y,\P)$} refer to the set of all bounded linear
operators $A$ on $Y$ such that $AK$ and $KA$ are in $K(Y,\P)$
whenever $K\in K(Y,\P)$. Both $K(Y,\P)$ and $L(Y,\P)$ are Banach
subalgebras of $L(Y)$, and $K(Y,\P)$ is a two-sided ideal in
$L(Y,\P)$. By definition, $L(Y,\P)$ is the largest subalgebra of
$L(Y)$ with that property. It is shown in \cite[Theorem
1.1.9]{RaRoSiBook} that $L(Y,\P)$ is {\em inverse
closed}\index{inverse closed}; that is, if $A\in L(Y,\P)$ is
invertible then $A^{-1}\in L(Y,\P)$.

\begin{Lemma} \label{LYP_iff}
An operator $A\in L(Y)$ is in $L(Y,\P)$ iff for every $m\in\N_0$,
\[
\|P_mAQ_n\|\to 0\qtext{and}\|Q_nAP_m\|\to 0\qtext{as} n\to\infty.
\]
\end{Lemma}
\Proof This is a straightforward computation. (See \cite[Prop.
1.1.8]{RaRoSiBook}.) \Proofend

\begin{Remark} \label{half}
Clearly the characterisations of $S(Y)$ and $SN(Y)$ in Lemma
\ref{lem1} and Corollary \ref{cor4} bear a close resemblance to the
definition of $K(Y,\P)$ and the characterisation of $L(Y,\P)$ in the
above lemma, respectively. In particular, in the case that $Y$ is a
Hilbert space and $P_n$ is self-adjoint, for each $n$, it holds that
$A\in K(Y,\P)$ ($\in L(Y,\P)$) iff $A$ and $A^*$ are in $SN(Y)$ (in
$S(Y)$), where $A^*$ denotes the adjoint of $A$.
\end{Remark}

The above characterisation also yields the following interesting
result:

\begin{Lemma} \label{lem_bothorneither}
For an operator $K\in L(Y,\P)$, either both or neither of the two
properties in \refeq{KYP_def} hold, so that $L(Y,\P)\cap SN(Y) =
K(Y,\P)$.
\end{Lemma}
\Proof Suppose $K\in L(Y,\P)$ and $\|KQ_n\|\to 0$ as $n\to\infty$.
Then for all $m,n\in\N_0$,
\[
\|Q_nK\|\ \le\ \|Q_nKP_m\|+\|Q_nKQ_m\|\ \le\ \|Q_nKP_m\|+\|KQ_m\|
\]
holds, where $\|KQ_m\|$ can be made as small as desired by choosing
$m$ large enough, and $\|Q_nKP_m\|$ tends to zero as $n\to\infty$.
Consequently, also the second property holds in \refeq{KYP_def}. By
a symmetric argument we see that the second property implies the
first. \Proofend

In analogy to Lemma \ref{lem2} we have the following result.

\begin{Lemma} \label{lem2'}
$L(Y,\P)\cap K(Y) \subseteq K(Y,\P)$, with equality if and only if
$P_n\in K(Y)$ for all $n$.
\end{Lemma}
\Proof From Corollary \ref{cor4} and Lemma \ref{LYP_iff} we know
that $L(Y,\P)\subseteq S(Y)$. Consequently,
\begin{eqnarray*}
L(Y,\P)\cap K(Y)&=&L(Y,\P)\cap L(Y,\P)\cap K(Y)\\
&\subseteq& L(Y,\P)\cap S(Y)\cap K(Y)\\
&\subseteq& L(Y,\P)\cap SN(Y)\ \ =\ \ K(Y,\P),
\end{eqnarray*}
where we used Lemmas \ref{lem2} and \ref{lem_bothorneither} for the
last two steps. Moreover, if $P_n\in K(Y)$ for all $n$ and $K\in
K(Y,\P)$ then $P_nK\in K(Y)$ for all $n$ and $K=\lim P_nK\in K(Y)$.
If $P_n\not\in K(Y)$ for some $n$ then $P_n$ is contained in the
difference of the two sets under consideration. \Proofend

The above lemma has the following refinement in the case when both $P_n$ and
its adjoint converge strongly to the identity.
\begin{Definition} \label{def_perfect} \cite{RaRoSiBook} Call $\P$ {\em
perfect}\index{perfect} if $P_n$ converges strongly to $I$ and
$P_n^*$ converges strongly to $I^*$.
\end{Definition}

\begin{Lemma} \label{per_lem} If $\P$ is perfect then $K(Y)\subset K(Y,\P)$. If
also $P_n\in K(Y)$ for every $n$, then $K(Y)=K(Y,\P)$ and $L(Y)=L(Y,\P)$.
\end{Lemma}
\Proof
That $K(Y)\subset K(Y,\P)$ follows from Lemma \ref{lem_strong}, and that
$K(Y,\P)\subset K(Y)$ if $\P\subset K(Y)$ follows from Lemma \ref{lem2'}. Then
$L(Y,\P)=L(Y)$ is immediate from the definition of $L(Y,\P)$.
\Proofend

The class of operators
\begin{equation} \label{eq_L0}
L_0(Y)\ :=\ \{A\in L(Y)\ :\ x\in Y_0\,\Rightarrow\,Ax\in
Y_0\}\index{$L(Y)a$@$L_0(Y)$}
\end{equation}
will turn out to be of particular interest to us. Recall that $Y_0$
is characterised by Lemma \ref{lem_Y0}.

\begin{Lemma} \label{lem_L0}
For $A\in L(Y)$, the condition $A\in L_0(Y)$ is equivalent to the
strong convergence $Q_nAP_m\to 0$ as $n\to\infty$ for every fixed
$m$.
\end{Lemma}
\Proof Fix an arbitrary $m\in\N$. By Lemma \ref{lem_Y0}, the strong
convergence $Q_nAP_m\to 0$ as $n\to\infty$ is equivalent to
$AP_mx\in Y_0$ for every $x\in Y$. Clearly, $A\in L_0(Y)$ implies
$AP_mx\in Y_0$ for every $x\in Y$, since $P_mx\in Y_0$. The reverse
implication follows from $P_m x\to x$ for every $x\in Y_0$, from the
continuity of $A$, and the closedness of $Y_0$. \Proofend

As an immediate consequence of Lemmas \ref{LYP_iff} and \ref{lem_L0}
we get the following.

\begin{Corollary} \label{cor_LYPinL0}
Every operator in $L(Y,\P)$ maps $Y_0$ into $Y_0$, i.e.
$L(Y,\P)\subset L_0(Y)$.
\end{Corollary}

We finish this discussion by noting that every $A\in S(Y)$ has a
unique extension to $S(\hat Y)$.

\begin{Lemma} \label{lem3.13}
Every $A\in S(Y)$ has a unique extension to an operator $\hat A\in
S(\hat Y)$, defined by
\begin{equation} \label{eq_extA}
\hat A x\ :=\ \lim_{n\to\infty} AP_nx,\qquad x\in\hat Y,
\end{equation}
where the limit is understood in the strict topology. It holds that
$\|\hat A\|=\|A\|$, and if $A\in SN(Y)$, $L(Y,\P)$ or $K(Y,\P)$,
then $\hat A\in SN(\hat Y)$, $L(\hat Y,\P)$ or $K(\hat Y,\P)$,
respectively. Conversely, if $\hat A\in S(\hat Y)$, $SN(\hat Y)$,
$L(\hat Y,\P)$ or $K(\hat Y,\P)$ and if $\hat A(Y)\subset Y$, then
$A:=\hat A|_Y\in S(Y)$, $SN(Y)$, $L(Y,\P)$ or $K(Y,\P)$,
respectively.
\end{Lemma}
\Proof It is easy to see that every sequentially continuous linear
operator on a TVS has a unique sequentially continuous extension to
the sequential completion of the TVS. The obvious construction of
this extension in our case is (\ref{eq_extA}). For $x\in Y$, we have
$AP_nx\sto Ax$ since $P_nx\sto x$ and $A\in S(Y)$, so that $\hat
A|_Y=A$. Moreover, for $x\in\hat Y$,
\[
\|\hat Ax\|\le \sup_n\|AP_nx\|\le \|A\|\sup_n\|P_n x\| = \|A\| \|x\|
\]
so that $\hat A$ is bounded and $\|\hat A\|\le\|A\|$. Together with
$\hat A|_Y=A$, this gives $\|\hat A\|=\|A\|$.

Now let us show that $\hat A\in SN(\hat Y)$ if $A\in SN(Y)$. For
every $x\in\hat Y$ and $n\in\N$, we have $P_kx\sto x$, and therefore
$AQ_nP_kx=\hat AQ_nP_kx\sto \hat AQ_nx$ as $k\to\infty$ since $\hat
A,Q_n\in S(\hat Y)$. Thus, for $x\in\hat Y$ with $\|x\|=1$,
\[
\|\hat AQ_n x\|\le
\sup_k\|AQ_nP_kx\|\le\sup_k\|AQ_nP_k\|\le\|AQ_n\|\sup_k\|P_k\|\to 0
\]
as $n\to\infty$, by Lemma \ref{lem1}, since $A\in SN(Y)$. Hence
$\hat A\in SN(\hat Y)$, by Lemma \ref{lem1} again.

From the trivial equality $\|Q_n\hat AP_m\|=\|Q_n AP_m\|$, together
with $\hat A\in S(\hat Y)$, we get that $\hat A\in L(\hat Y,\P)$ if
$A\in L(Y,\P)$, by Corollary \ref{cor4} and Lemma \ref{LYP_iff}.
Finally, it follows that $\hat A\in K(\hat Y,\P)=L(\hat Y,\P)\cap
SN(\hat Y)$ if $A\in K(Y,\P)=L(Y,\P)\cap SN(Y)$, by Lemma
\ref{lem_bothorneither}. \Proofend

\section{Compactness and Collective Compactness on $(Y,s)$}
A linear operator on a TVS is said to be {\em compact}\index{compact
operator}\index{operator!compact} if the image of some neighbourhood
of zero is relatively compact. A linear operator is often said to be
{\em Montel}\index{Montel operator}\index{operator!Montel} if it has
the weaker property that it maps bounded sets onto relatively
compact sets. These properties coincide when the TVS is a normed
space. Much of the familiar theory of compact operators on normed
spaces generalises to compact operators on locally convex separated
TVS's, for example the theory of Riesz \cite{Robertson}. In
particular, a compact operator has a discrete spectrum (as an
element of the algebra $C(Y)$), whose only accumulation point is
zero, and all non-zero points of the spectrum are eigenvalues. By
contrast, as we will see below, the spectrum of a Montel operator
may be much more complex.

Let $KS(Y)$\index{$KS(Y)$} denote the set of compact operators on
$(Y,s)$ and $M(Y)$\index{$M(Y)$} the set of Montel operators. Then
it is standard (and clear) that $KS(Y)\subset M(Y)\subset L(Y)$ and
$KS(Y)\subset C(Y)$. Also $K(Y)\subset M(Y)$, since bounded sets
coincide in the strict and norm topologies and relatively compact
sets in the norm topology are relatively compact in the strict
topology. Thus, by Example \ref{exnotinS}, it may not hold that
$M(Y)\subset S(Y)$. (A Venn diagram illustrating the various subsets
of $L(Y)$ that we have introduced in this chapter is shown in Figure
\ref{fig:venn} below.) By Lemmas \ref{lem_strict_prop} and
\ref{lem_compact}, an operator $A$ is in $M(Y)$ iff the image of
every norm-bounded set is relatively sequentially compact in the
strict topology. Operators with this property are termed {\em
sequentially compact with respect to the $(Y,s)$ topology} in
\cite{ChanZh02}. The following two lemmas are useful
characterisations of $M(Y)$ in the case when $Y$ is sequentially
complete (which is the case when $Y=\hat Y$, by Lemma
\ref{lem_Qny}).

\begin{Lemma} \label{lem_charact_mon}
If $Y=\hat Y$ then $A\in M(Y)$ iff $A\in L(Y)$ and $P_mA\in K(Y)$
for every $m$.
\end{Lemma}
\Proof The lemma follows immediately from the equivalence of (a)
and (d) in Lemma \ref{lem_compact}. This implies that $A\in M(Y)$
iff $A(S)$ is norm-bounded and $P_nA(S)$ is relatively compact in
the norm topology, for every $n$ and every norm-bounded set $S$.
\Proofend

\begin{Remark}
In the case $Y=L^p(\R^N)$ of Example \ref{ex_Lp} (in which case
$Y=\hat Y$), an operator which satisfies $P_mA\in K(Y)$ and also
$AP_m\in K(Y)$ for each $m$ is termed {\sl locally
compact}\index{locally compact operator}\index{compact
operator!locally}\index{operator!compact!locally} in
\cite{CWLi2,LiBook, RaRoIndexLC, RaRoSiBook}. In the case
$Y=BC(\R^N)$ of Example \ref{e22} (in which again $Y=\hat Y$) an
operator $A\in L(Y)$ is termed {\sl locally compact} in
\cite{Jorgens} if it holds merely that $P_mA\in K(Y)$ for every $m$,
i.e. by Lemma \ref{lem_charact_mon}, if $A\in M(Y)$.
\end{Remark}

\begin{Lemma} \label{lem_charact_mon2}
 If $A\in M(Y)$ then $AP_n\in KS(Y)$ for every $n$.
Conversely, if $Y=\hat Y$, $A\in S(Y)$ and $AP_n\in M(Y)$ for every
$n$, then $A\in M(Y)$.
\end{Lemma}
\Proof By Lemma \ref{lem1}, $P_n\in SN(Y)$, and so, by Lemma
\ref{lem_sny}, maps some neighbourhood in $(Y,s)$ to a bounded set
in $(Y,\|\cdot\|)$. (In fact every neighbourhood in $(Y,s)$ is
mapped to a bounded set.) Thus $AP_n\in KS(Y)$ if $A\in M(Y)$.

If $AP_n\in M(Y)$ for every $n$ then, by Lemma
\ref{lem_charact_mon}, $P_mAP_n\in K(Y)$ for every $m$ and $n$. If
also $A\in S(Y)$ then, by Corollary \ref{cor4}, $\|P_mA -
P_mAP_n\|\to 0$ as $n\to\infty$, so that $P_mA\in K(Y)$ for every
$m$. Thus $A\in M(Y)$ by Lemma \ref{lem_charact_mon}. \Proofend

Many of the arguments we make in this text will deal with
families of operators that have the following collective
compactness property.

\begin{Definition} \cite{ChanZh02}\label{def_collcomp}
We say that a set $\K$ of linear operators on $Y$ is {\em uniformly
Montel on $(Y,s)$}\index{uniformly Montel}\index{Montel!uniformly}
or is {\em collectively sequentially compact on
$(Y,s)$}\index{collectively sequentially compact}\index{sequentially
compact!collectively}\index{compact!sequentially!collectively} if,
for every bounded set $B$, $ \cup_{K\in\K}K(B) $ is relatively
compact in the strict topology.
\end{Definition}

\begin{Remark} \label{rem_unifmont}
Note that, by Lemma \ref{lem_compact}, $ \cup_{K\in\K}K(B) $ is
relatively compact in the strict topology iff $ \cup_{K\in\K}K(B) $
is relatively sequentially compact in the strict topology, i.e.\
iff, for every sequence $(K_n)\subset \K$ and $(x_n)\subset B$,
$(K_nx_n)$ has a strictly convergent subsequence.
\end{Remark}

That being Montel on $(Y,s)$ is significantly weaker than being
compact is very clear in the case when $P_n$ is compact for all
$n$. The next two results follow from Corollary \ref{lem_finite}
(the first is also a corollary of Lemma \ref{lem_charact_mon}).

\begin{Corollary} \label{cor_MYLY}
If $Y=\hat Y$ and $P_n\in K(Y)$ for every $n$ then $M(Y)=L(Y)$.
\end{Corollary}

\begin{Corollary} \label{cor_UMUB}
If $Y=\hat Y$ and $P_n\in K(Y)$ for every $n$ then a set $\K$ of
linear operators on $Y$ is uniformly Montel on $(Y,s)$ iff $\K$ is
uniformly bounded.
\end{Corollary}

\begin{Remark} \label{rem_A_IK}
Some of our subsequent results will only apply to operators $A$ of
the form $A=I+K$ with $K\in S(Y)\cap M(Y)$. It follows from
Corollary \ref{cor_MYLY} that, if $Y=\hat Y$ and $P_n\in K(Y)$ for
each $n$, then $A-I\in S(Y)\cap M(Y)$ whenever $A\in S(Y)$, so that
every $A\in S(Y)$ can be written in this form.
\end{Remark}

$M(Y)$ is the set of operators which map bounded sets to relatively
compact sets in $(Y,s)$, and we have seen in Lemma \ref{lem_sny}
that $SN(Y)$ is precisely the set of those operators that map some
neighbourhood in $(Y,s)$ to a bounded set. On the other hand,
$KS(Y)$ is the set of those operators that map some neighbourhood to
a relatively compact set. Clearly, if $A\in M(Y)\cap SN(Y)$ then
$A^2\in KS(Y)$. What is less clear is that $A\in KS(Y)$, which the
next lemma shows.

\begin{Lemma} \label{SN&M=KS}
It always holds that $S(Y)\cap K(Y) \subset SN(Y)\cap M(Y) = KS(Y)$.
If $P_n\in K(Y)$ for each $n$ then $S(Y)\cap K(Y) = SN(Y) = KS(Y)$.
\end{Lemma}
\Proof We have seen already that $S(Y)\cap K(Y) \subset SN(Y)$ and
$K(Y)\subset M(Y)$.

To show that $SN(Y)\cap M(Y) = KS(Y)$, suppose that $K\in KS(Y)$.
Then, as noted already, it follows that $K\in M(Y)$. Since $K$
maps some neighbourhood of zero in $(Y,s)$ to a relatively compact
set in $(Y,s)$, and relatively compact sets are bounded, it
follows from Lemma \ref{lem_sny} that $K\in SN(Y)$.

Suppose now that $K\in SN(Y)\cap M(Y)$. Then, by Lemma
\ref{lem_sny}, there exists a neighbourhood of zero, $U$, in
$(Y,s)$, for which $K(U)$ is norm bounded and $\sup_{x\in
U}\|KQ_nx\|\to 0$ as $n\to\infty$. Let $a:\N\to (0,\infty)$ be a
null sequence and let $V$ be the neighbourhood of zero, $V:=\{x\in
U:|x|_a<1\}$. Then, for every $n$, $P_n(V)\subset Y$ is bounded so
that, since $K\in M(Y)$, for every $n$, the image of every sequence
in $V$ under the mapping $KP_n$ has an $s-$convergent subsequence.
In particular, given a sequence $(x_m)\subset V$ we can construct a
chain of subsequences $(x_m)\supset
(x_m^{(1)})\supset(x_m^{(2)})\supset\cdots$ such that
$KP_jx_m^{(n)}$ is $s-$convergent as $m\to\infty$ for $j=1,...,n$.
Then $KP_jx_m^{(m)}$ is $s-$convergent for every $j$, to a limit
$y_j\in Y$. Now, for each $k,j_1,j_2$,
\[
\|P_k(y_{j_1}-y_{j_2})\|=\lim_{m\to\infty}\|P_kK(Q_{j_2}-Q_{j_1})x_m^{(m)}\|
\le \sup_{x\in V}\|K(Q_{j_2}-Q_{j_1})x\|.
\]
Thus
\[
\|y_{j_1}-y_{j_2}\|=\sup_k\|P_k(y_{j_1}-y_{j_2})\|\to 0
\]
as $j_1,j_2\to\infty$, so that, since $Y$ is a Banach space, for
some $y\in Y$, $y_j\to y$ as $j\to\infty$. It follows that
$Kx_m^{(m)}\sto y$ as $m\to\infty$, so that we have shown that
$K(V)$ is relatively sequentially compact in $(Y,s)$ and so, by
Lemma \ref{lem_compact}, relatively compact, so that $K\in KS(Y)$.
To see this last claim, note that $K(V)\subset K(U)$ is bounded and
that, for every $k$ and $j$,
\[
\|P_k(Kx_m^{(m)}-y)\|\le\|P_k(KP_jx_m^{(m)}-y_j)\|+\|P_k(y_j-y)\|+\|P_kKQ_jx_m^{(m)}\|,
\]
which can be made as small as desired by choosing first $j$ and then
$m$ sufficiently large.

The last sentence of the lemma follows immediately from Lemma
\ref{lem2} and from $SN(Y)=S(Y)\cap K(Y)\subset K(Y)\subset M(Y)$
and therefore $SN(Y)\cap M(Y)=SN(Y)$. \Proofend

\begin{figure}[h]
\begin{center}
\includegraphics[width=0.9\textwidth]{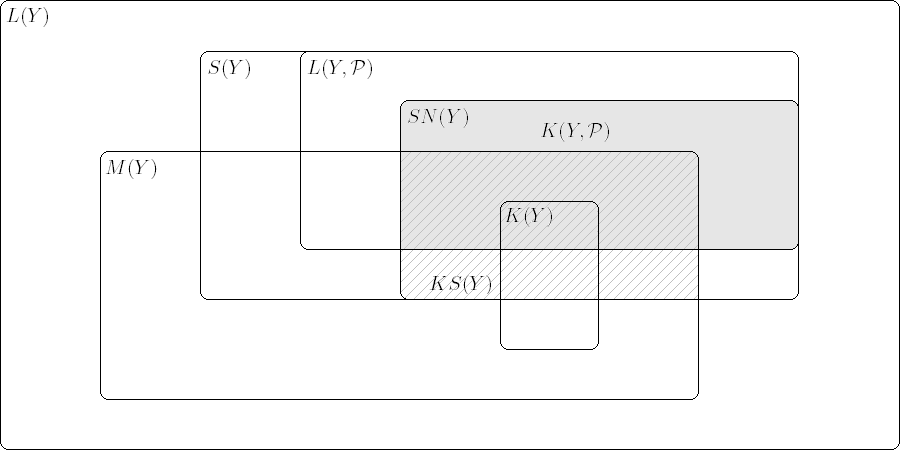}
\caption{\footnotesize Venn diagram of the operator classes studied in this chapter.
(The gray shaded area represents $K(Y,\P)$, and the hatched area is $KS(Y)$.)}
\label{fig:venn}
\end{center}
\end{figure}

We finish the section with further examples of operators in
$S(Y)$, $KS(Y)$, and $M(Y)$.
\begin{Example} \label{ex_io_SY}
Let $Y=BC(\R^N)$ and $P_n$ be defined as in Examples \ref{e22} and
\ref{ex_arzasc}, so that $Y=\hat Y$ and $Y$ is sequentially complete
in the strict topology. Suppose that $K$ is a linear integral
operator\index{integral operator}\index{operator!integral} on $Y$ of
the form
$$
Kx(s) = \int_{\R^N} k(s,t)x(t) dt, \quad s\in\R^N,
$$
with $k(s,\cdot)\in L^1(\R^N)$ for every $s\in \R^N$. Then
\cite{Jorgens} $K\in S(Y)$ iff
\begin{equation} \label{eq_kbound}
\sup_{s\in\R^N} \int_{\R^N} |k(s,t)| dt < \infty
\end{equation}
and
\begin{equation} \label{eq_kcont}
\int_{\R^N} |k(s,t)-k(s^\prime,t)| dt \to 0, \mbox{ as }
s^\prime\to s,
\end{equation}
for every $s\in\R^N$. Conditions (\ref{eq_kbound}) and
(\ref{eq_kcont}) imply that $K$ maps bounded sets to bounded
equicontinuous sets, that is to sets that are relatively compact in
the strict topology (see Example \ref{ex_arzasc}). Thus $K\in
M(Y)\cap S(Y)$ if (\ref{eq_kbound}) and (\ref{eq_kcont}) hold.
Conditions (\ref{eq_kbound}) and (\ref{eq_kcont}) hold, in
particular, if (as in Chapter \ref{sec_intop}) $K$ is in the algebra
generated by operators of multiplication by bounded continuous
functions and operators of convolution by $L^1$ functions (as in
Example \ref{ex_conv}), or in the closure of this algebra in $L(Y)$
with its norm topology. (Such operators are termed {\em
convolution-type operators}\index{convolution-type
operator}\index{operator!convolution-type} in \cite{RaRoSiBook}.)
\end{Example}

\begin{Example} As a special case of the above example, suppose that $K$ is
defined as in Example \ref{ex_conv}. Then \cite{Jorgens} the
spectrum of $K$ is $\{0\}\cup \{\hat \kappa(\xi):\xi\in\R^N\}$,
where $\hat\kappa\in BC(\R^N)$ is the Fourier transform of $\kappa$.
All non-zero points of the spectrum are eigenvalues
($\hat\kappa(\xi)$ has eigenfunction $x(s) := \exp(i\xi\cdot s)$).
Since the spectrum of $K$ is not discrete, $K\in M(Y)\cap S(Y)$ but
$K\not\in KS(Y)$.
\end{Example}

\begin{Example} (Cf. \cite{AnseloneSloan,CW_NumMath}.) As another special
case of Example \ref{ex_io_SY}, consider the one-dimensional case
$N=1$ with $K$ defined by
$$
Kx(s) = \int_0^1 \exp(ist)x(t)dt, \quad s\in\R.
$$
Then $K=KP_n$ for all $n>1$, so that $K\in KS(Y)$ by Lemma
\ref{lem_charact_mon2}. But $K\not\in K(Y)$ as, defining $x_n(s) =
\exp(-ins)$, $s\in\R$, $(Kx_n)$ has no norm-convergent subsequence
since $Kx_n(s)\to 0$ as $s\to\infty$ for every $n$ but $Kx_n(n)=1$
for each $n$.
\end{Example}

\section{Algebraic Properties}
\label{sec:algebraic} We will find the algebraic properties
collected in the following lemma useful. These are immediate from
the definitions and Lemmas \ref{lem_sny}, \ref{lem_charact_mon2} and
\ref{SN&M=KS}.

\begin{Lemma} \label{lem_alg}
Let $A$ and $B$ be linear operators on $Y$. Then
\begin{eqnarray*}
A\in M(Y), \;\;B\in L(Y) &\Rightarrow& AB \in M(Y)\\
A\in S(Y), \;\;B\in M(Y) &\Rightarrow& AB \in M(Y)\\
A\in L(Y), \;\;B\in SN(Y) &\Rightarrow& AB \in SN(Y)\\
A\in SN(Y), \;\;B\in S(Y) &\Rightarrow& AB \in SN(Y)\\
A\in SN(Y), \;\;B\in M(Y) &\Rightarrow& AB \in K(Y)\\
A\in M(Y), \;\;B\in SN(Y) &\Rightarrow& AB \in KS(Y)\\
A\in KS(Y)                &\Rightarrow& A^2 \in K(Y)
\end{eqnarray*}
\end{Lemma}

$S(Y)$, $SN(Y)$, $M(Y)$, and $KS(Y)$ are all vector subspaces of
$L(Y)$. It follows from the above lemma that they are all
subalgebras of $L(Y)$, and that $SN(Y)$, $M(Y)\cap S(Y)$, and
$KS(Y)$ are all (two-sided) ideals of $S(Y)$. Moreover, all these
subalgebras are closed when endowed with the norm topology of
$L(Y)$.

\begin{Lemma} \label{lem_banalg}
$S(Y)$, $SN(Y)$, $M(Y)$, and $KS(Y)$ are all Banach subalgebras of
$L(Y)$, with $S(Y)$ a unital subalgebra and $SN(Y)$, $M(Y)\cap
S(Y)$, and $KS(Y)$ two-sided ideals of $S(Y)$.
\end{Lemma}
\Proof It only remains to show that each subalgebra is closed. If
$A\in L(Y)$ is in the closure of $SN(Y)$ then, for every $B\in
SN(Y)$ and every $n$,
$$
\|AQ_n\| \le \|BQ_n\| + \|(B-A)Q_n\| \le \|BQ_n\| + 2\|B-A\|.
$$
Since $B$ can be chosen to make $\|B-A\|$ arbitrarily small and,
by Lemma \ref{lem1}, $\|BQ_n\|\to 0$ as $n\to\infty$ for every
$B$, it follows that $\|AQ_n\|\to 0$ as $n\to\infty$ so that $A\in
SN(Y)$. Thus $SN(Y)$ is closed.

Since $SN(Y)$ is closed it follows from Lemma \ref{lem3} that $S(Y)$
is closed.

As $K(Y)$ is closed, it follows from Lemma \ref{lem_charact_mon} in
the case $Y=\hat Y$ that $M(Y)$ is closed. In the general case, to
see that $M(Y)$ is closed, suppose that $(A_m)\subset M(Y)$ and
$A_m\toto A\in L(Y)$. Let $(x_n)$ be a bounded sequence in $Y$.
Then, by a diagonal argument as in the proof of Lemma \ref{SN&M=KS},
we can find a subsequence, denoted again by $(x_n)$, such that, for
each $m$, there exists a $y_m\in Y$ such that $A_mx_n\sto y_m$ as
$n\to\infty$. Arguing as in the proof of Lemma \ref{SN&M=KS}, we can
show that the sequence $(y_m)$ is Cauchy in $(Y,\|.\|)$ and so has a
limit $y\in Y$, and that $Ax_n\sto y$. This shows that the image of
every bounded set under $A$ is relatively sequentially compact and
so relatively compact in $(Y,s)$, by Lemma \ref{lem_compact}, i.e.
$A\in M(Y)$.

By Lemma \ref{SN&M=KS}, $KS(Y)=M(Y)\cap SN(Y)$, being the
intersection of two closed spaces, is closed itself. \Proofend

We will say that $A\in L(Y)$ is {\em invertible}\index{invertible}
if it is invertible in the algebra of linear operators on $Y$, i.e.
if it is bijective. Automatically, by the open mapping theorem, it
follows that $A^{-1}\in L(Y)$. An interesting question is whether
$S(Y)=C(Y)$ is inverse closed, i.e. whether, if $A\in S(Y)$ is
invertible, it necessarily holds that $A^{-1}\in S(Y)$. Since
$(Y,s)$ is not barrelled \cite{ChanZh02}, this question is not
settled by standard generalisations of the open mapping theorem to
non-metrisable TVS's \cite{Robertson}. Indeed, it is not clear to us
whether $S(Y)$ is inverse closed without further assumptions on $Y$.
But we do have the following result which implies that $S(Y)$ is
inverse closed in the case when $Y=\hat Y$ and $P_n\in K(Y)$ for
each $n$.

\begin{Lemma} \label{lem_inv_closed}
Suppose $A,B\in S(Y)$ are invertible and that $A^{-1}\in S(Y)$ and
$A-B\in M(Y)$. Then $B^{-1}\in S(Y)$.
\end{Lemma}
\Proof We have that $B^{-1} = D^{-1}A^{-1}$, where $D = I+C$ and
$C=A^{-1}(B-A)$. By Lemma \ref{lem_alg}, $C\in S(Y)\cap M(Y)$. To show
that $B^{-1}\in S(Y)$ we need only to show that $D^{-1}\in S(Y)$.

Suppose that $(x_n)\subset Y$, $x\in Y$, and $x_n\sto x$. Let $y_n
: = D^{-1} x_n$. By (\ref{eqstrict}), and since $D^{-1}=B^{-1}A\in
L(Y)$, $(x_n)$ and $(y_n)$ are bounded. For each $n$,
\begin{equation} \label{eq_yn}
y_n + C y_n = x_n.
\end{equation}
Since $C\in M(Y)$ there exists a subsequence $(y_{n_m})$ and $y\in
Y$ such that $x_{n_m} - C y_{n_m} \sto y$. From (\ref{eq_yn}) it
follows that $y_{n_m}\sto y$. Since $C\in S(Y)$, it follows that
$x_{n_m} - C y_{n_m} \sto x - Cy$. Thus $y = x - Cy$, i.e. $y =
D^{-1}x$. We have shown that $y_n = D^{-1}x_n$ has a subsequence
strictly converging to $y = D^{-1}x$. By the same argument, every
subsequence of $y_n$ has a subsequence strictly converging to $y$.
Thus $D^{-1}x_n \sto D^{-1} x$. So $D^{-1}\in S(Y)$. \Proofend

\begin{Corollary} \label{cor_SY_inv_clo}
If $Y=\hat Y$ and $P_n\in K(Y)$ for all $n$ then $S(Y)$ is inverse
closed.
\end{Corollary}
\Proof If $Y=\hat Y$ and $P_n\in K(Y)$ for all $n$, and $A\in S(Y)$
is invertible, then $I-A\in M(Y$) by Corollary \ref{cor_MYLY}, so
that $A^{-1}\in S(Y)$ by the above lemma. \Proofend


\chapter{Notions of Operator Convergence} \label{chap_converg}
A component in the arguments to be developed is that one needs some
notion of the convergence of a sequence of operators. For
$(A_n)\subset L(Y)$, $A\in L(Y)$, let us write $A_n\toto
A$\index{$\toto$}\index{norm convergence}\index{convergence!norm} if
$\|A_n-A\|\to 0$ (as $n\to\infty$) and $A_n\to
A$\index{$\to$}\index{strong convergence}\index{convergence!strong}
if $(A_n)$ converges strongly to $A$, in the strong operator
topology induced by the norm topology on $Y$, i.e. if $A_n x \to A
x$ for all $x\in Y$. Following \cite{LiBook,RaRoSiBook} we introduce
also the following definition.

\begin{Definition} \label{def_P_conv}
We say that a sequence $(A_n)\subset L(Y)$
$\P$-converges\index{P-conv@$\P$-convergence}\index{convergence!$\P$-}
to $A\in L(Y)$ if, for all $K\in K(Y,\P)$, both
\begin{equation} \label{eq_Pconv}
\|(A_n-A)K\| \to 0 \qtext{and} \|K(A_n-A)\| \to 0 \qtext{as}
n\to\infty.
\end{equation}
In this case we write $A_n\pto A$\index{$\pto$} or $A=\P$\!--\,$\lim
A_n$.
\end{Definition}

The following lemma is a generalisation of Proposition 1.65 from
\cite{LiBook}. It shows that every $\P$-convergent sequence is
bounded in $L(Y)$ and that, conversely, for a bounded sequence
$(A_n)$ one has to check property (\ref{eq_Pconv}) only for $K\in\P$
in order to guarantee $A_n\pto A$.

\begin{Lemma} \label{lem_Pconv}
Suppose $(A_n)\subset L(Y)$ and $A\in L(Y)$. Then $A_n\pto A$ iff
$(A_n)$ is bounded in $L(Y)$ and, for all $m$,
\begin{equation} \label{eq_Pconv2}
\|(A_n-A)P_m\| \to 0 \qtext{and} \|P_m(A_n-A)\| \to 0\qtext{as}
n\to\infty.
\end{equation}
\end{Lemma}
\Proof Suppose $(A_n)$ is bounded and \refeq{eq_Pconv2} holds.
Then, for all $m\in\N$ and all $K\in K(Y,\P)$, one has
\[
\|K(A_n-A)\|\ \leq\ \|K\|\,\|P_m(A_n-A)\|\,+\,\|KQ_m\|\,\|A_n-A\|,
\]
where the first term tends to zero as $n\to\infty$, and the second
one is as small as desired if $m$ is large enough. The first
property of \refeq{eq_Pconv} is shown absolutely analogously.

Conversely, if \refeq{eq_Pconv} holds for all $K\in K(Y,\P)$, then
\refeq{eq_Pconv2} holds for all $m\in\N$ since $\P\subset
K(Y,\P)$. It remains to show that $(A_n)$ is bounded.

Suppose the converse is true. Without loss of generality, we can
suppose that $A=0$. Now we will successively define two sequences:
$(m_i)_{i=1}^\infty\subset\N$ and $(n_k)_{k=0}^\infty\subset\N_0$.
We start with $m_2\gg m_1:=1$ and $n_0:=0$.

For every $k\in\N$, choose $n_k\in\N$ such that
\[
n_k>n_{k-1}\ ,\qquad \|A_{n_k}\|>k^2+3\qqtext{and}
\|P_{m_{3k-1}}A_{n_k}\|<1,
\]
the latter possible since $P_{m_{3k-1}}\in K(Y,\P)$ and $A_n\pto
0$. Then
\[
\|Q_{m_{3k-1}}A_{n_k}\|\ \ge\ \|A_{n_k}\|-\|P_{m_{3k-1}}A_{n_k}\|\
> \ k^2+3-1\,=\,k^2+2.
\]
Take $u_k\in Y$ with $\|u_k\|=1$ and
$\|Q_{m_{3k-1}}A_{n_k}u_k\|>k^2+1$, and choose the next three
elements of the sequence $(m_i)$ by $m_{3k+2}\gg m_{3k+1}\gg
m_{3k}\gg m_{3k-1}$ such that
$\|P_{m_{3k}}Q_{m_{3k-1}}A_{n_k}u_k\|>k^2$, which is possible by
\refeq{eq_lim}. Then
\[
\|P_{m_{3k}}Q_{m_{3k-1}}A_{n_k}\|>k^2\qtext{for all} k\in\N.
\]
Now put
\[
K:=\sum_{j=1}^\infty \frac 1{j^2}\ P_{m_{3j+1}}Q_{m_{3j-2}}.
\]
From $m_i\ll m_{i+1}$ for all $i\in\N$ we get that
$P_{m_{3k}}Q_{m_{3k-1}}P_{m_{3j+1}}Q_{m_{3j-2}}$ equals
$P_{m_{3k}}Q_{m_{3k-1}}$ if $k=j$ and $0$ otherwise. Consequently,
\[
P_{m_{3k}}Q_{m_{3k-1}}K=\frac 1{k^2}\ P_{m_{3k}}Q_{m_{3k-1}},
\]
and hence,
\[
\|KA_{n_k}\|\ \ge\ \|P_{m_{3k}}Q_{m_{3k-1}}KA_{n_k}\|\,=\,\|\frac
1{k^2}\, P_{m_{3k}}Q_{m_{3k-1}}A_{n_k}\|\ >\ \frac{k^2}{k^2}\,=\,1
\]
for every $k\in\N$. On the other hand, $K\in K(Y,\P)$, which
implies $\|KA_n\|\to 0$ as $n\to\infty$. Contradiction. \Proofend

The following is a simple but important example of $\P$-convergence
that is fundamental to the application we study in Chapter
\ref{sec_ellp}.

\begin{Example} \label{ex_mult}
Generalising Example \ref{e21}, we suppose that
$Y=\ell^p(\Z^N,U)$\index{$\ell^p(\Z^N,U)$}, for some
$p\in[1,\infty]$ and $N\in\N$ where $U$ is some Banach space. The
elements of $Y$ are of the form $x=(x(m))_{m\in\Z^N}$ with $x(m)\in
U$ for every $m=(m_1,...,m_N)\in\Z^N$ and we equip $Y$ with the
usual norm. For $m\in\Z^N$, we define $|m|:=\max(|m_1|,...,|m_N|)$
and put
\begin{equation} \label{Pndef}
P_n x(m) = \left\{\begin{array}{ll}
  x(m), &|m|\le n, \\
  0, &|m|>n, \\
\end{array} \right.\index{$P_n$}
\end{equation}
for every $x\in Y$ and $n\in\N_0$. Similarly to Examples \ref{e21}
and \ref{ex_Lp}, $\P=(P_n)$ satisfies (i) and (ii) with $N(m) = m$.

For $b=(b(m))_{m\in\Z^N}\in\ell^\infty(\Z^N,L(U))$ define the
multiplication operator\index{multiplication
operator}\index{operator!multiplication} $M_b\in L(Y)$ by
\begin{equation} \label{eq_multop}
M_b x(m) = b(m) x(m), \quad m\in \Z^N,\index{$M_b$}
\end{equation}
for $x\in Y$, and note that $\|M_b\|=\|b\|$. It is a straightforward
consequence of this equation and Lemma \ref{lem_Pconv} that, for a
sequence $b_n\in\ell^\infty(\Z^N,L(U))$,
\begin{eqnarray} \nonumber
M_{b_n} \pto 0 &\Leftrightarrow &\sup_n\|b_n\| < \infty \mbox{ and }
b_n \sto 0\\ &\Leftrightarrow &\sup_n\|b_n\|<\infty \mbox{ and }
\|b_n(m)\|\to 0, \;\forall m\in\Z^N. \label{eqn_Pconfmult}
\end{eqnarray}
In the above equation by $b_n\sto 0$ we mean that $b_n$ converges to
zero in the strict topology generated by the family $\P = (P_n)$,
where we are here using the notation $P_n$ also to denote the
operator on $\ell^\infty(\Z^N,L(U))$ defined by (\ref{Pndef}) for
$x\in \ell^\infty(\Z^N,L(U))$.
\end{Example}

We have seen already that $S(Y)$ and $L(Y,\P)$ are Banach
subalgebras of $L(Y)$. Both are also closed with respect to
$\P-$convergence.

\begin{Lemma} \label{prop_lep_pclosed}
$S(Y)$ and $L(Y,\P)$ are sequentially closed with respect to
$\P-$con\-vergence.
\end{Lemma}
\Proof First suppose $(A_n)\subset S(Y)$, $A\in L(Y)$ and $A_n\pto
A$. Then, if $x_n\sto 0$, for every $k$ and $m$, we have
\[
\|P_kAx_n\|\le \|P_k(A-A_m)\|\sup_n\|x_n\| + \|P_kA_mx_n\|.
\]
But $\|P_kA_mx_n\|\to 0$ as $n\to\infty$ since $A_m\in S(Y)$, and
$\|P_k(A-A_m)\|$ can be made as small as desired by choosing $m$
large. So we get $Ax_n\sto 0$, and therefore $A\in S(Y)$.

Similarly (see Proposition 1.1.17(a) in \cite{RaRoSiBook} for the
details) we show that also $L(Y,\P)$ is sequentially closed.
\Proofend

To make use of results from \cite{ChanZh02} we introduce also the
notions of operator convergence used there. For $(A_n)\subset L(Y)$
and $A\in L(Y)$, let us write $A_n\sto A$\index{$\sto$} if, for all
$(x_n)\subset Y$,
\begin{equation} \label{eq_sto_def}
x_n\sto x\quad \Rightarrow\quad A_n x_n\sto A x.
\end{equation}
Call $\cA\subset L(Y)$ {\it $s$-sequentially
compact}\index{s-seq@$s$-sequentially
compact}\index{compact!$s$-sequentially} if, for every sequence
$(A_n)\subset \cA$, there exists a subsequence $(A_{n_m})$ and $A\in
\cA$ such that $A_{n_m}\sto A$. Note that $A\sto A$ holds iff $A\in
S(Y)$. It follows that, if $\cA\subset L(Y)$ is $s$-sequentially
compact, then $\cA\subset S(Y)$.

A more familiar and related notion of operator convergence is that
of strong (or pointwise) convergence. For $(A_n)\subset L(Y)$, $A\in
L(Y)$, we will say that $(A_n)$ {\em converges to $A$ in the strong
operator topology on $(Y,s)$}, and write $A_n\Sto
A$,\index{$\sto2$@$\Sto$} if
\begin{equation} \label{eq_Sto}
A_n x\sto A x,\quad x\in Y.
\end{equation}
Clearly, the $S$-limit is unique, that is $A_n\Sto A$ and $A_n\Sto
B$ implies $A=B$. Hence also the $s$-limit and $\P$-limit are
unique, by Lemma \ref{rel_conv} and Corollary \ref{Pimps} below.

 Clearly,
\begin{equation} \label{strongImpS}
A_n\to A\quad\Rightarrow\quad A_n \Sto A.
\end{equation}

The following lemmas explore further properties of and relationships
between the notions of operator convergence we have introduced. We
will exhibit this relationship through Example \ref{ex_conv2}.

\begin{Lemma} \label{rel_conv}
Suppose $(A_n)\subset L(Y)$, $A\in L(Y)$. Then
\begin{equation} \label{stoSto}
A_n \sto A \quad\Rightarrow\quad A_n \Sto A \mbox{ and } A\in S(Y).
\end{equation}
Further,  $A_n \Sto A$ as $n\to\infty$ iff $(A_n)$ is bounded and
$P_m(A_n-A) \to 0$ as $n\to\infty$ for all $m\in \N$.
\end{Lemma}
\Proof It is clear from the definitions that $A_n \sto A$ implies
$A_n \Sto A$. That $A_n \sto A$ implies $A\in S(Y)$ is shown in
\cite[Lemma 3.1]{ChanZh02}. That $A_n \Sto A$ implies $P_m(A_n-A)
\to 0$ is clear from (\ref{eqstrict}), and that it also implies that
$(A_n)$ is bounded is shown in \cite[Lemma 3.3]{ChanZh02}.
Conversely, if $(A_n)$ is bounded and $P_m(A_n-A) \to 0$ for each
$m$, then, for every $x\in Y$,  $(A_n x)$ is bounded and $P_m(A_n
x-A x)\to 0$ for each $m$, so that $A_nx\sto Ax$ by
(\ref{eqstrict}). \Proofend

\begin{Example} \label{ex_conv2}
Let $Y$, $P_n$ and the multiplication operator $M_b$ be defined as
in Example \ref{ex_mult}, and suppose that $(b_n)\subset
\ell^\infty(\Z^N,L(U))$. Then, extending the results of Example
\ref{ex_mult}, we see that
\begin{eqnarray*}
M_{b_n} \toto 0 &\Leftrightarrow &\|b_n\| =
\sup_{m\in\Z^N}\|b_n(m)\|\to 0,\\
M_{b_n}\pto 0 &\Leftrightarrow &\sup_n\|b_n\|<\infty \mbox{ and }\|b_n(m)\|\to
0,\quad \forall m\in\Z^N,\\
M_{b_n}\sto 0 &\Leftrightarrow & M_{b_n}\Sto 0\\
&\Leftrightarrow & \sup_n\|b_n\|<\infty \mbox{ and
}\|b_n(m)x(m)\|\to 0,\quad \forall m\in\Z^N,\;x\in Y.
\end{eqnarray*}
Thus $M_{b_n}\pto 0$ requires that each component of $b_n$ converges
to zero in norm, while $M_{b_n} \sto 0$ requires that each component
of $b_n$ converges strongly to zero. We have (cf.\ Corollary
\ref{strong2} below) that
$$
M_{b_n}\toto 0 \ \Rightarrow\  M_{b_n}\pto 0 \ \Rightarrow\
M_{b_n}\sto 0 \ \Leftrightarrow\ M_{b_n}\Sto 0 \ \Leftarrow\
M_{b_n}\to 0.
$$
If $U$ is finite-dimensional, then $\pto$, $\Sto$ and $\sto$ all
coincide. If $p=\infty$, then $\to$ is equivalent to $\toto$. If
$1<p<\infty$ and $U$ is finite-dimensional, then $\to$ coincides
with $\pto$, $\Sto$ and $\sto$.
\end{Example}

\begin{Lemma} \label{lll}
Suppose $(A_n)\subset L(Y)$ is bounded, $A\in S(Y)$, and
$$
|| P_m(A_n-A)||\to 0 \mbox{ as } n\to\infty
$$
for each $m$. Then $A_n\sto A$.
\end{Lemma}
\Proof If the conditions of the lemma hold and $x_n\sto x$ then
$Ax_n\sto Ax$ and, by (\ref{eqstrict}), $\sup_n||x_n||<\infty$, so
that $(A_nx_n)$ is bounded, and, for each $m$,
$$
|| P_m(A_nx_n-Ax)|| \le ||P_m(A_n-A)x_n|| + ||P_mA(x_n-x)||\to 0
$$
as $n\to\infty$. Thus, by (\ref{eqstrict}), $A_nx_n\sto Ax$.
\Proofend

As a corollary of Lemmas \ref{lem_Pconv} and \ref{lll} we have
\begin{Corollary} \label{Pimps}
Suppose $(A_n)\subset L(Y)$, $A\in S(Y)$. Then
\begin{equation} \label{Ptosto}
A_n \pto A\quad \Rightarrow\quad A_n \sto A.
\end{equation}
\end{Corollary}

Let us say that $\cA\subset  L(Y)$ is {\it $s$-sequentially
equicontinuous} if
$$
(A_n)\subset \cA,\,\, x_n\sto 0\quad \Rightarrow\quad A_n x_n\sto 0.
$$
Clearly, if $\cA$ is $s$-sequentially equicontinuous, then
$\cA\subset S(Y)$. The significance here of this definition is the
following result taken from \cite{ChanZh02}.

\begin{Lemma} \label{lnc4}
Suppose $(A_n)\subset S(Y)$, $A\in S(Y)$. Then $A_n\sto A$ iff
$A_n\Sto A$ and $\{A_n:n\in \N\}$ is $s$-sequentially
equicontinuous.
\end{Lemma}

Let us say that $\cA\subset L(Y)$ is {\em sequentially compact in
the strong operator topology on $(Y,s)$} if, for every sequence
$(A_n)\subset\cA$, there exists $A\in\cA$ and a subsequence
$(A_{n_m})$ such that $A_{n_m}\Sto A$. Then Lemma \ref{lnc4} and
other observations made above imply the following corollary.

\begin{Corollary} \label{cext}
Suppose $\cA\subset  L(Y)$. Then $\cA$ is $s$-sequentially compact
iff $\cA\subset S(Y)$ and $\cA$ is $s$-sequentially equicontinuous
and sequentially compact  in the strong operator topology on
$(Y,s)$.
\end{Corollary}

In the case that the strict and norm topologies coincide, in which
case $C(Y)=S(Y)=L(Y)$, it follows from Lemma \ref{rel_conv}, i.e.
from the uniform boundedness theorem in Banach spaces, that if
$(A_n)\subset L(Y)$ and $A_n\Sto A$ then $\{A_n:n\in \N\}$ is
$s$-sequentially equicontinuous. In the case when these topologies
do not coincide, in which case, by Lemma \ref{lem_topo}, $(Y,s)$ is
not metrisable, other versions of the Banach-Steinhaus theorem would
apply \cite{Rudin91,Robertson},  if $(Y,s)$ were a Baire space or,
more generally, a barrelled TVS, to give that $\{A_n:n\in\N\}$ is
$s$-sequentially equicontinuous if $A_n\Sto A$ and $(A_n)\subset
C(Y)$. But, by \cite[Theorem 2.1]{ChanZh02}, $(Y,s)$ is not
barrelled unless the norm and strict topologies coincide. And in
fact the following example makes it clear that a version of the
Banach-Steinhaus theorem, enabling equicontinuity to be deduced from
continuity and pointwise boundedness, does not always hold for
$(Y,s)$ if the strict and norm topologies do not
coincide.\\

\begin{Example}
Let $Y$ be defined as in Example \ref{e21}. For $n\in \N$ define
$A_n\in  L(Y)$ by $A_n x (m)= x(n)$, for $x\in Y$, $m\in\Z$. It is
easy to see that $(A_n)\subset C(Y)\subset L(Y)$, and clearly
$||A_n||\le 1$ so that $(A_n)$ is bounded. But $(A_n)$ is not
$s$-sequentially equicontinuous as, defining $x_n(m)=1+\tanh(m-n)$,
$m\in \Z$, $n\in \N$, clearly $(x_n)\subset Y$, $x_n\sto 0$, but
$A_n x_n(0)= 1$, so $A_n x_n\not\sto 0$.
\end{Example}

In the case that $Y$ satisfies an additional assumption, it is shown
in \cite{ChanZh02} that a sequence $(A_n)\subset S(Y)$ that is
convergent in the strong operator topology on $(Y,s)$ is
$s$-sequentially equicontinuous.  The additional assumption is the
following one, in which $Y_m\subset Y$ is the subspace
defined by (\ref{Ymdef}):\\

\noindent {\bf Assumption A.}\index{Assumption A} For every $m\in
\N$ there exists $n>m$ and $Q:Y\to Y_m$ such that
\begin{equation} \label{assa} ||x-Qx+y||
\le \max(|x|_n,||y||),\quad x\in Y, \,y\in Y_n.
\end{equation}

\noindent That Assumption A is satisfied in some applications is illustrated by
the following example.\\

\begin{Example} \label{exAa} Suppose that $Y=BC(\R^N)$ and $P_n$ are defined as
in Example
\ref{e22}.  Then (\ref{assa}) holds with $n=m+2$ and $Q=Q_{m+1}$,
for then $||x-Qx+y|| = ||P_{n-1}x+y|| = \max(||P_{n-1}x||,||y||)\le
\max(|x|_n,||y||)$, for all $x\in Y$ and $y\in Y_n$.
\end{Example}
\begin{Lemma} \label{strong} \cite{ChanZh02}
Suppose that Assumption A holds and that $(A_n)\subset S(Y)$, $A\in
S(Y)$, and $A_n\Sto A$. Then $\{A_n:n\in\N\}$ is $s$-sequentially
equicontinuous.
\end{Lemma}
Combining Lemmas \ref{strong}, \ref{lnc4}, Corollary \ref{Pimps},
and (\ref{strongImpS}), we have the following result which shows
that, when Assumption A holds, the convergence $\sto$ is weaker than
both ordinary strong convergence and $\P$-convergence.

\begin{Corollary} \label{strong2}
Suppose that  Assumption A holds and that $(A_n)\subset S(Y)$, $A\in
S(Y)$. Then
$$
A_n \pto A\quad\Rightarrow\quad  A_n\sto A \quad\Leftrightarrow\quad
A_n\Sto A \quad\Leftarrow\quad A_n \to A.
$$
\end{Corollary}


\chapter{Key Concepts and Results}
 \label{sec_LimOps}
This chapter introduces the key concepts and develops the key
results of the text. We first recall the concepts of invertibility
at infinity and Fredholmness and start to explore their
inter-relation. Next, we summarise some main results from the
abstract generalised collectively compact operator theory developed
in \cite{ChanZh02} and from the abstract theory of limit operators
\cite{RaRoSi1998,RaRoSiBook,LiBook}. It then turns out that the
collection of all limit operators of an operator $A\in S(Y)$ is
subject to the constraints made in the operator theory of
\cite{ChanZh02}. Therefore we apply this theory and derive some new
general results which can be used to study the invertibility of
operators and their limit operators. We will illustrate the
application of the results of this chapter throughout the remainder
of the text.

\section{Invertibility at Infinity and Fredholmness} \label{sec_inv_inf_Fred}

Following \cite{RaRoSiBook,LiBook} we introduce the following
definition.

\begin{Definition} \label{inv_inf}
An operator $A\in L(Y)$ is said to be {\em invertible at
infinity}\index{invertible at infinity} if there exist operators
$B\in L(Y)$ and $T_1,T_2\in K(Y,\P)$ such that
\begin{equation} \label{eq_inv_inf}
AB = I + T_1 \qtext{and} BA = I + T_2.
\end{equation}
\end{Definition}

\begin{Remark} \label{rem_PFred}
If $A\in L(Y,\P)$ then $A+K(Y,\P)$ is invertible in the quotient
algebra $L(Y,\P)/K(Y,\P)$ iff $A$ is invertible at infinity with
$B\in L(Y,\P)$ in (\ref{eq_inv_inf}). Rabinovich et al.\
\cite{RaRoSiBook} call $A\in L(Y,\P)$ {\em
$\P$-Fredholm}\index{P-Fred@$\P$-Fredholm} when this is the case.
\end{Remark}

The following lemma gives some justification for the name
`invertible at infinity'.

\begin{Lemma} \label{lem_inv_inf}
If $A\in L(Y)$ and (\ref{eq_inv_inf}) holds with $B\in L(Y)$ and
$T_1,T_2\in SN(Y)\supset K(Y,\P)$, then there exist
$B_1^\prime,B_2^\prime\in L(Y)$ and $n\in\N_0$ such that
\begin{equation} \label{eq_inv_inf2}
Q_nAB_1^\prime = Q_n = B_2^\prime A Q_n.
\end{equation}
If $B\in S(Y)$ then we can also choose $B_1^\prime,B_2^\prime\in
S(Y)$. If $B,T_1,T_2\in L(Y,\P)$ then also $B_1^\prime,B_2^\prime$
can be chosen in $L(Y,\P)$.
\end{Lemma}
\Proof Choose $m\in\N_0$ large enough that $\|T_1Q_m\|<1$ and
$\|T_2Q_m\|<1$, and take an $n\in\N_0$ such that
$Q_nQ_m=Q_n=Q_mQ_n$. From (\ref{eq_inv_inf}) we get
$Q_nABQ_m=Q_n(I+T_1Q_m)$ and $BAQ_n=(I+T_2Q_m)Q_n$, proving
(\ref{eq_inv_inf2}) where we put $B_1':=BQ_m(I+T_1Q_m)^{-1}$ and
$B_2':=(I+T_2Q_m)^{-1}B$.

The two additional claims follow immediately from the Neumann series
formula and the fact that $S(Y)$ and $L(Y,\P)$ are Banach
subalgebras of $L(Y)$. \Proofend

Recall that $A\in L(Y)$ is called
\emph{semi-Fredholm}\index{semi-Fredholm operator}\index{Fredholm
operator!semi-}\index{operator!Fredholm!semi-} if it has a closed
range, $A(Y)$, and if one of the numbers
\begin{equation} \label{eq_alphabeta}
\alpha(A):=\dim(\ker A)\qquad\textrm{and}
\qquad\beta(A):=\dim(Y/A(Y))
\end{equation}\index{$\alpha(\cdot)$}\index{$\beta(\cdot)$}
is finite, and that $A$ is called \emph{Fredholm}\index{Fredholm
operator}\index{operator!Fredholm} if both $\alpha(A)$ and
$\beta(A)$ are finite (in which case its range is automatically
closed). In the latter case the {\em
index}\index{index}\index{Fredholm index} of $A$ is defined as
$\alpha(A)-\beta(A)$.

The following theorem shows that, in some important cases,
invertibility at infinity implies Fredholmness. Since $I\in S(Y)$,
one possible choice of $C$ in this statement is $C=I$.

\begin{Theorem} \label{lem_Fredholm}
Suppose $A = C + K$, where $C\in L(Y)$ is invertible, with
$C^{-1}\in S(Y)$, and $K\in S(Y)\cap M(Y)$. Suppose further that
(\ref{eq_inv_inf}) holds with $B\in L(Y)$ and $T_1,T_2\in SN(Y)$.
Then $A$ is Fredholm.
\end{Theorem}
\Proof We have from (\ref{eq_inv_inf}) that $CB + KB = I + T_1$, $BC
+ BK = I + T_2$, so that
$$
B = C^{-1}(I + T_1 - KB), \quad B = (I+T_2-BK)C^{-1}.
$$
Using the first of these equations we see that
$$
AC^{-1}(I-KB) = (I+KC^{-1})(I-KB) = I - KC^{-1}T_1,
$$
and note that $KC^{-1}T_1\in KS(Y)$ by Lemma \ref{lem_alg}, and thus
\begin{equation} \label{eq_reg_right}
AC^{-1}(I-KB)(I+KC^{-1}T_1) = I - (KC^{-1}T_1)^2,
\end{equation}
with $(KC^{-1}T_1)^2\in K(Y)$ by Lemma \ref{lem_alg}. Similarly,
\begin{equation} \label{eq_reg_left}
(I-BK)C^{-1}A = (I-BK)(I+C^{-1}K) = I - T_2C^{-1}K
\end{equation}
with $T_2C^{-1}K\in K(Y)$ by Lemma \ref{lem_alg}. We have
constructed right and left regularisers for $A$, so $A$ is Fredholm.
\Proofend

The following is a corollary of the above result and Lemma
\ref{lem_charact_mon}. The last sentence follows from the
observation that $A= I+(A-I)$ and that $P_n(A-I)\in K(Y)$ if $A\in
L(Y)$ and $P_n\in K(Y)$.

\begin{Corollary} \label{cor_invinf_Fred} If $A\in L(Y)$ is invertible at
infinity and $A = I + K$, with $K\in S(Y)$ and $P_nK\in K(Y)$ for
every $n$, then $A$ is Fredholm. In the case that $P_n\in K(Y)$ for
all $n$, $A\in S(Y)$ is Fredholm if it is invertible at infinity.
\end{Corollary}

In the case that $\P$ is perfect we will see in Lemma \ref{thm_per}
that, conversely, Fredholmness implies invertibility at infinity. We
will, in Chapter \ref{sec_ellp}, also establish this result for the
case $Y=\ell^p(\Z^N,U)$, for $p=1,\infty$, and a Banach space $U$,
in which case $\P$ is not perfect.

\section{A Generalised Collectively Compact Operator Theory}
\label{subsec_collcomp} Following \cite[Section 4]{ChanZh02}, we let
$iso(Y)$\index{$\easo$@$iso(Y)$} denote the set of isometric
isomorphisms on $Y$ and call a set $S\subset iso(Y)$ {\em
sufficient}\index{sufficient set} if, for some $n\in\N$ it holds
that, for every $x\in Y$ there exists $V\in S$ such that
$2\,|Vx|_n\ge\|x\|$. The following examples illustrate this
definition:

\begin{Example} \label{e51}
Let $Y=BC(\R^N)$ and $\P$ as in Example \ref{e22}. Then both
$S_1=\{V_k:k\in\R^N\}$, where $V_k$\index{$V_k$} is the translation
operator\index{translation operator}\index{shift
operator}\index{operator!translation}\index{operator!shift}
$$
V_kx(s)=x(s-k),\quad s\in\R^N,
$$
and $S_2=\{\Psi_k:k\in\N\}$, where
$$
\Psi_kx(s)=x(ks),\quad s\in\R^N,
$$
are sufficient families of isometric isomorphisms on $Y$ where we
can choose $n=1$ in both cases.
\end{Example}

\begin{Example} \label{e52}
Let $Y=L^p(\R^N)$ and $\P$ as in Example \ref{ex_Lp}. If
$p=\infty$ then both $S_1$ and $S_2$ from Example \ref{e51} are
sufficient with $n=1$ for instance.

If $p<\infty$ then neither $S_1$ nor $S_2$ is sufficient. Although
it is true that for every $x\in Y$ there are $n\in\N$ and $V_k\in
S_1$ such that $2\,|V_kx|_n\ge\|x\|$, there is no universal
$n\in\N$ which is large enough to guarantee this property for all
$x\in Y$.
\end{Example}

We say that an operator $A\in L(Y)$ is {\em bounded
below}\index{bounded below} if
$$
\nu(A)\,:=\,\inf_{\|x\|=1}\|Ax\|\,>\,0.
$$
In that case, we refer to $\nu(A)$\index{$\nu(A)$} as the {\sl lower
norm}\index{lower norm} of $A$.

$A\in L(Y)$ is bounded below iff $A$ is injective and has a closed
range. Indeed, necessity is obvious and sufficiency follows from
Banach's theorem on the inverse operator saying that $A^{-1}:A(Y)\to
Y$ acts boundedly on the range of $A$ if that is closed. Another
elementary result on the lower norm is that it depends continuously
on the operator; in particular, we have
\begin{equation} \label{eq_lonoemineq}
| \nu(A)-\nu(B)|\ \le\ \|A-B\|
\end{equation}
for all $A,B\in L(Y)$.

If $A$ is invertible then $A$ is bounded below and
$\nu(A)=1/\|A^{-1}\|$. We will say that a set $\cA\subset L(Y)$ is
{\em uniformly bounded below}\index{uniformly bounded
below}\index{bounded below!uniformly} if every $A\in\cA$ is bounded
below and if there is a $\nu>0$ such that $\nu(A)\ge\nu$ for all
$A\in\cA$, that is
$$
\|Ax\|\ge\nu\|x\|,\qquad A\in\cA,\,x\in Y.
$$

For $\K\subset L(Y)$, we abbreviate the set $\{I-K:K\in\K\}$ by
$I-\K$.

In the following theorem we use the notation $\K^\times$ to denote
the set of all subsequences of sequences
$(K_1,K_2,...)\in\K_1\times\K_2\times\cdots$ for a fixed family of
sets $\K_1,\K_2,...\subset L(Y)$. This theorem is a slight
strengthening of Theorems 4.1 and 4.4 in \cite{ChanZh02} (in
\cite{ChanZh02} the condition (\ref{inj_sur}) has `$I-K_n$ bounded
below' replaced by the weaker `$I-K_n$ injective'), but an
examination of the proof of Theorem 4.4 in \cite{ChanZh02} shows
that this slightly stronger result follows by exactly the same
argument.

\begin{Theorem} \label{th4.1} 
Suppose that $Y=\hat Y$, $S\subset iso(Y)$ is sufficient,
$\K,\K_1,\K_2,...\subset L(Y)$, and that
\begin{itemize}
\item[(i)] $\cup_{n\ge 1}\K_n$ is uniformly Montel on $(Y,s)$;
\item[(ii)] for every sequence $(K_n)\in\K^\times$, there exist
a subsequence $(K_{n(m)})$ and $K\in\K$ such that $K_{n(m)}\sto K$
as $m\to\infty$;
\item[(iii)] for all $n\in\N$, it holds that $V^{-1}KV\in\K_n$ for all
$K\in\K_n$ and $V\in S$;
\item[(iv)] $I-K$ is injective for all $K\in\K$.
\end{itemize}
Then:

{\bf a)} There is an $n_0\in\N$ such that $I-\cup_{n\ge n_0}\K_n$ is
uniformly bounded below, i.e. there is a $\nu>0$ such that
\[
\|(I-K)x\|\ge \nu \|x\|,\qquad K\in\K_n,\ n\ge n_0,\ x\in Y;
\]

{\bf b)} If, in addition, for every $K\in\K$, there exists a
sequence $(K_n)\in\K^\times$ such that $K_n\sto K$ and all operators
$I-K_n$ have the property
\begin{equation} \label{inj_sur}
I-K_n\textrm{ bounded below }\quad\Longrightarrow\quad I-K_n\textrm{
surjective},\qquad n=1,2,...
\end{equation}
then all operators in $I-\K$ are invertible, and
\[
\sup_{K\in\K}\|(I-K)^{-1}\|\le \nu^{-1}.
\]
\end{Theorem}

The following special case of the above theorem, obtained by setting
$\K_1=\K_2=\cdots=\K$ in Theorem \ref{th4.1} is worth noting. We
will say that a subset $\cA\subset\K\subset L(Y)$ is {\em
$s$-dense}\index{s-dense@$s$-dense} in $\K$ if, for every $K\in\K$,
there is a sequence $(K_n)\subset\cA$ with $K_n\sto K$.

\begin{Theorem} \label{th4.5} \cite[Theorem 4.5]{ChanZh02}
Suppose that $Y=\hat Y$, $S\subset iso(Y)$ is sufficient and
$\K\subset L(Y)$ has the following properties:
\begin{itemize}
\item[(i)] $\K$ is uniformly Montel on $(Y,s)$;
\item[(ii)] $\K$ is $s$-sequentially compact;
\item[(iii)] $V^{-1}KV\in\K$ for all $K\in\K,\,V\in S$;
\item[(iv)] $I-K$ is injective for all $K\in\K$.
\end{itemize}
Then:

{\bf a)} The set $I-\K$ is uniformly bounded below;

{\bf b)} If in $I-\K$ there is an $s$-dense subset of surjective
operators then all operators in $I-\K$ are surjective.
\end{Theorem}

Note that in statement b), as in Theorem \ref{th4.1}, all operators
in $I-\K$ are consequently invertible, and their inverses are
uniformly bounded by $1/\nu$ where $\nu>0$ is a lower bound on all
lower norms $\nu(I-K)$ with $K\in\K$ which exists by a).

\section{Limit Operators} \label{subsec_limops}
Following \cite[Section 1.2]{RaRoSiBook}, let $N\in\N$ and
$\V=\{V_k\}_{k\in\Z^N}$\index{$Va$@$\V$}\index{$V_k$} denote a group
of linear isometries on $Y$ which are subject to
\begin{equation} \label{constraint1}
V_0=I\qqtext{and} V_kV_m=V_{k+m},\quad
k,m\in\Z^N.
\end{equation}
Moreover, we impose that $\V$ is compatible with $\P$ in the
following sense:
\begin{equation} \label{constraint2}
\forall m,n\in\N\quad \exists k_0\in\N:\qquad
P_mV_kP_n=0\quad\qtext{if} |k|>k_0
\end{equation}
and
\begin{equation} \label{constraint3}
\forall m\in\N,\ k\in\Z^N\quad \exists n_0\in\N:\qquad
P_mV_kQ_n=0=Q_nV_kP_m\quad\qtext{if} n>n_0.
\end{equation}

\begin{Definition} \label{def_limop} \cite{RaRoSiBook}
We will say that a sequence $h=(h(n))_{n=1}^\infty\subset\Z^N$
tends to infinity if $|h(n)|\to\infty$ as $n\to\infty$. If $h$
tends to infinity and $A\in L(Y)$ then
$$
A_h\ :=\ \plim_{n\to\infty}\ V_{-h(n)}AV_{h(n)}\index{$A_h$, limit
operator of $A$}
$$
is called the {\em limit operator}\index{limit
operator}\index{operator!limit} of $A$ with respect to the sequence
$h$, provided the $\P-$limit exists.
\end{Definition}

\noindent From what we know about $\P-$convergence it follows that
the limit operator $A_h$ is unique if it exists.

The {\em operator spectrum}\index{operator
spectrum}\index{spectrum!operator}
$\opsp(A):=\{A_h\}$\index{$\opsp(A)$} (see e.g.
\cite{LiBook,RaRoSiBook}) is the collection of all limit operators
of $A$ where $h\subset\Z^N$ runs through all sequences tending to
infinity such that $A_h$ exists. For certain operators $A$,
invertibility at infinity, as specified in Definition \ref{inv_inf},
can be characterised in terms of properties of the operator spectrum
(see Theorem \ref{prop_rich_invinf} below). In turn, as we have seen
in Corollary \ref{cor_invinf_Fred}, for a large class of operators
invertibility at infinity implies Fredholmness, so that Fredholmness
can be determined by studying the operator spectrum; indeed, in some
cases it is known that the operator spectrum also determines the
index (\cite{RaRoIndexLC}, \cite{RaRoRoe}, \cite[Section
2.7]{RaRoSiBook}, \cite[Section 3.3.1]{LiBook}). To establish the
most complete results we need to restrict consideration to {\em rich
operators}\index{rich operator}\index{operator!rich}, where $A\in
L(Y)$ is referred to as a rich operator if every sequence
$h\subset\Z^N$ tending to infinity has an infinite subsequence $g$
such that the limit operator $A_g$ exists.

\begin{Example} \label{ex_limit_operator}
Let $Y=\ell^\infty$ and define $\P$ as in Example \ref{e21}. Define
$\V= \{V_k\}_{k\in\Z}\subset L(Y)$ by
$$
V_kx(m) = x(m-k), \quad m\in\Z.
$$
Then (\ref{constraint1})-(\ref{constraint3}) hold for every $k_0\ge
m+n$ and $n_0\ge m+k$. For $b\in \ell^\infty$ let $M_b\in L(Y)$
denote the multiplication operator defined by
$$
M_bx(m) = b(m) x(m), \quad m\in\Z.
$$
For $S\subset \Z$ let $\chi_S\in \ell^\infty$ denote the
characteristic function of $S$. Define $a\in \ell^\infty$ by
$$
a(m) := \lfloor\sqrt{|m|}\,\rfloor \mod 2, \quad m\in\Z,
$$
where $\lfloor s\rfloor\le s$ denotes the integer part of $s$, and
set $A=M_a$. Then, where $B:=
\{\chi_{\{n,...,+\infty\}},\,\chi_{\{-\infty,...,n\}}:n\in\Z\}$,
$$
\opsp(A) = \{0,I,M_b:b\in B\}.
$$
\begin{figure}[h]
\begin{center}
\includegraphics[width=\textwidth]{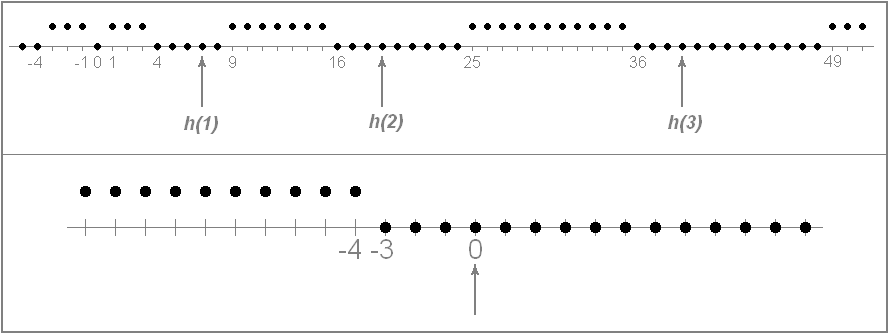}
\caption{\footnotesize Functions $a$ (top) and
$b=\chi_{\{-\infty,...,-4\}}$ (below) from Example
\ref{ex_limit_operator}.} \label{fig:ex511}
\end{center}
\end{figure}
For example, if $h(n):=4n^2+3$, then
$$
\plim_{n\to\infty}\ V_{-h(n)}AV_{h(n)} = M_b
$$
where $b:=\chi_{\{-\infty,...,-4\}}$ (see Figure \ref{fig:ex511}).
The operator $A$ is rich; this can be seen directly or by applying
Lemma \ref{lem_oprich} below.
\end{Example}

The following theorem summarises and extends known results on the
operator spectrum $\opsp(A)$ and on the relationship between $A$ and
its operator spectrum. Statements (i) and (ii) are from
\cite{RaRoSi1998}, (iii) and (iv) are from \cite{RaRoSiBook} and
statements (v)-(vii) go back to \cite[Section 3.3]{LiDiss} and can
also be found in \cite[Section 1.2]{RaRoSiBook}. (Note that the
proofs of (iii)-(vii) given in \cite{LiDiss,RaRoSiBook} work for all
$A\in L(Y)$, although the results state a requirement for $A\in
L(Y,\P)$ or make a particular choice of $Y$, and note also that (iv)
is immediate from (ii) and (iii) and that (vii) is immediate from
(ii), (v) and (vi), see \cite{CWLi:JFA2007}.) Thus we include only a
proof of (viii) and (ix), in which $Y_0\subset Y$ is as given in
Chapter \ref{sec_strict}.

For brevity, we introduce the notation
\begin{equation} \label{eq_Sh(A)}
\Sh(A)\ :=\ \{V_{-k}AV_k:k\in\Z^N\}\index{$T(A)$@$\Sh(A)$}
\end{equation}
for the set of all translates of an operator $A\in L(Y)$.

\begin{Theorem} \label{lem_opsp}
For every $A\in L(Y)$, the following statements hold.
\begin{itemize}
\item[(i)] If $B\in\opsp(A)$ then $\|B\|\le\|A\|$.
\item[(ii)] If $B\in\opsp(A)$ and $k\in\Z^N$ then also $V_{-k}BV_k\in\opsp(A)$.
\item[(iii)] $\opsp(A)$ is sequentially closed with respect to
$\P-$convergence.
\item[(iv)] If $B\in \opsp(A)$ then $\opsp(B)\subset \opsp(A)$.
\item[(v)] $A$ is rich iff $\Sh(A)$ is relatively
$\P-$sequentially compact.
\item[(vi)] If $A$ is rich then $\opsp(A)$ is $\P-$sequentially compact.
\item[(vii)] If $A$ is rich and $B\in\opsp(A)$ then $B$ is rich.
\item[(viii)] If $B\in\opsp(A)$ then $||Bx||\ge \nu(A) ||x||$ for $x\in
Y_0$, so that $\nu(B)\ge \nu(A)$ if $Y=Y_0$.
\item[(ix)] If $B\in\opsp(A)\cap L(Y,\P)$ is invertible then $\nu(B)\ge
\nu(A)$.
\end{itemize}
\end{Theorem}
\Proof (viii) If $B\in\opsp(A)$ then $B=A_h$ for some sequence
$h\subset\Z^N$. For $m\in\N$ and every $x\in Y$ we have that
$$
|| V_{-h(n)}AV_{h(n)}P_mx|| = ||AV_{h(n)}P_mx|| \ge \nu(A)
|| ||V_{h(n)}P_mx||=\nu(A) ||P_mx||.
$$
Since $V_{-h(n)}AV_{h(n)}\pto B$, taking the limit as $n\to\infty$ we get
$$
|| BP_mx|| \ge \nu(A) ||P_mx||.
$$
For $x\in Y_0$ we have, by Lemma \ref{lem_Y0}, that $P_mx\to x$ as
$m\to\infty$, so the result follows.

(ix)
For $m,n\in \N$ and $x\in Y$,
\begin{eqnarray}
\nonumber ||P_mB^{-1}x|| & \le & ||P_mB^{-1}Q_nx|| + ||P_mB^{-1}P_nx||\\
\label{eq_ix} & \le & ||P_mB^{-1}Q_nx|| + ||B^{-1}P_nx||.
\end{eqnarray}
As $L(Y,\P)$ is inverse closed \cite[Theorem 1.1.9]{RaRoSiBook}, we
have that $B^{-1}\in L(Y,\P)$, so that $||P_mB^{-1}Q_k||\to0$ and
$||Q_kB^{-1}P_n||\to0$ as $k\to\infty$, the latter implying, by
Lemma \ref{lem_Y0}, that $B^{-1}P_nx\in Y_0$. Thus from
(\ref{eq_ix}) and (viii) we have that
\begin{eqnarray*}
\nu(A)||P_mB^{-1}x||& \le & \nu(A)||P_mB^{-1}Q_nx|| + \nu(A)||B^{-1}P_nx||\\
& \le & \nu(A)||P_mB^{-1}Q_nx|| + ||P_nx||.
\end{eqnarray*}
Taking the limit first as $n\to\infty$ and then as $m\to\infty$, noting
(\ref{eq_lim}), we get that
$ \nu(A)||B^{-1}x||\le ||x||$. We have shown that $ \nu(A)||y||\le ||By||$, for
all $y\in Y$, as required. \Proofend

Within the subspace $L(Y,\P)$ of $L(Y)$, for every fixed sequence
$h$ tending to infinity, the mapping $A\mapsto A_h$ is compatible
with all of addition, composition, scalar multiplication and passing
to norm-limits \cite{RaRoSi1998}. That is, the equations
\begin{eqnarray}\nonumber
(A+B)_h=A_h+B_h,& &(AB)_h=A_hB_h,\\ \label{eq_limops} (\lambda
A)_h=\lambda A_h,& &\left(\lim_{m\to\infty} A^{(m)}\right)_h =
\lim_{m\to\infty} A^{(m)}_h
\end{eqnarray}
hold, in each case provided the limit operators on the right hand
side exist. By definition, $L(Y,\P)$ is a subalgebra of $L(Y)$. By
Lemma \ref{LYP_iff}, \refeq{constraint3} implies that $\V\subset
L(Y,\P)$. Thus, if $A\in L(Y,\P)$ then $\Sh(A)\subset L(Y,\P)$, so
that $\opsp(A)\subset L(Y,\P)$ by Lemma \ref{prop_lep_pclosed}.
Similarly, since by Lemma \ref{lem_banalg}, $S(Y)\supset L(Y,\P)$ is
a subalgebra of $L(Y)$, if $A\in S(Y)$ then $V_{-k}AV_k\in S(Y)$ for
all $k\in\Z^N$, so that $\opsp(A)\subset S(Y)$ by Lemma
\ref{prop_lep_pclosed}. As a consequence of (\ref{eq_limops})
together with a diagonal argument to see the closedness, the set of
rich operators $A\in L(Y,\P)$ is a Banach subalgebra of $L(Y,\P)$.

We have at this point introduced all the main concepts that we will use in the
rest of the text. The remainder of the text will in large part be focussed
on establishing relationships between the following properties of an operator
$A\in L(Y)$ for important operator classes:
\begin{eqnarray} \nonumber
&\mbox{(a)}& \mbox{$A$ is invertible.}\\
\nonumber
&\mbox{(b)}& \mbox{$A$ is Fredholm.}\\
\nonumber
&\mbox{(c)}& \mbox{$A$ is invertible at infinity.}\\
\label{A_props}
&\mbox{(d)}& \mbox{All limit operators of $A$ are invertible and the}\\
\nonumber
& & \mbox{inverses are uniformly bounded.}\\
\nonumber
&\mbox{(e)}& \mbox{All limit operators of $A$ are invertible.}\\
\nonumber
&\mbox{(f)}& \mbox{All limit operators of $A$ are injective.}
\end{eqnarray}
Clearly, it always holds that (a)$\Rightarrow$(b) and that
(d)$\Rightarrow$(e)$\Rightarrow$(f). In the case that $\P$ is perfect
(Definition \ref{def_perfect}) we have also the following result.

\begin{Lemma} \label{thm_per}
If $\P$ is perfect then, for all $A\in L(Y)$,
(b)$\Rightarrow$(c)$\Rightarrow$(d).
\end{Lemma}
\Proof That (b)$\Rightarrow$(c) follows from Lemma \ref{per_lem}.
That (c)$\Rightarrow$(d) follows as in the proof of
\cite[Proposition 1.2.9]{RaRoSiBook} noting that this proof applies
word for word with $L(Y,\P)$ replaced by $L(Y)$, provided that we
also replace `$\P$-Fredholm' by `invertible at infinity'. \Proofend

Thus we see that
(a)$\Rightarrow$(b)$\Rightarrow$(c)$\Rightarrow$(d)$\Rightarrow$(e)$\Rightarrow$(f)
in the case that $\P$ is perfect. Much of the remainder of the
text will be concerned with establishing these implications also
in cases when $\P$ is not perfect, and in investigating the
converse. In particular, we will see in the next section that the
generalised collectively compact operator theory introduced in the
last section sheds some light on the relationships between (d), (e)
and (f). Further, we have met in Section \ref{sec_inv_inf_Fred}
conditions on $A$ which ensure that (c)$\Rightarrow$(b).

One case in which the connection between the properties
(\ref{A_props}) is evident is the case in which $A$ is self-similar.
Here, following \cite{CWLi:JFA2007} call $A\in L(Y)$ {\em
self-similar}\index{self-similar
operator}\index{operator!self-similar} if $A\in \opsp(A)$ and,
generalising \cite{Muh1972}, call  $A\in L(Y)$ {\em recurrent}
\index{recurrent operator}\index{operator!recurrent} if, for every
$B\in\opsp(A)$, it holds that $\opsp(B)=\opsp(A)$. It is immediate
from these definitions and Theorem \ref{lem_opsp} (iv) that if $A$
is recurrent then all limit operators of $A$ are self-similar and
recurrent. In the case when $A$ is self-similar, the following
relationship between the properties (\ref{A_props}) (a)-(e) of $A$
is immediate from the definitions and Lemma \ref{thm_per}.

\begin{Corollary} \label{thm_per_ss} If $\P$ is perfect and $A\in L(Y)$ is
self-similar then all of (a)-(e) in (\ref{A_props}) are equivalent.
\end{Corollary}

That self-similar operators exist, indeed are ubiquitous, is clear
from the following lemma. For completeness, we include the proof
which is a variation on the proof of \cite[Proposition
3.7]{CWLi:JFA2007}. Proposition 3.7 in \cite{CWLi:JFA2007} is the
same statement as Lemma \ref{prop_ss}, but proved only for the
special case when $Y=\ell^p(\Z,U)$, for some $p\in[1,\infty]$ and
some Banach space $U$.

\begin{Lemma} \label{prop_ss} If $A\in L(Y)$ is rich then some
$B\in\opsp(A)$ is self-similar and recurrent.
\end{Lemma}
\Proof If $B\in\opsp(A)$ then $\opsp(B)\subset \opsp(A)$ (Theorem
\ref{lem_opsp} (iv)) and $\opsp(B)$ is non-empty since, by Theorem
\ref{lem_opsp} (vii), $B$ is rich. Further, if $B$ is recurrent then
every $C\in\opsp(B)\subset\opsp(A)$ is recurrent and self-similar.
Thus it is enough to show that there exists a $B\in\opsp(A)$ which
is recurrent.

Let
\[
\cA\ :=\ \{\ \opsp(B)\ :\ B\in\opsp(A)\ \}.
\]
By Theorem \ref{lem_opsp} (iv) this is a partially ordered set,
equipped with the order '$\supseteq$'. Note further that $\opsp(B)$
is a maximal element in this partially ordered set iff $B$ is
recurrent. So it remains to show the existence of a maximal element.
But this follows from Zorn's lemma if we can show that every totally
ordered subset of $\cA$ has an upper bound.

So let $\B$ be a totally ordered subset of $\cA$, i.e.
\[
\B\ =\ \{\ \opsp(B)\ :\ B\in\sigma\ \}
\]
where $\sigma\subseteq\opsp(A)$ is such that, for any two
$B_1,B_2\in\sigma$, we either have $\opsp(B_1)\supseteq\opsp(B_2)$
or $\opsp(B_2)\supseteq\opsp(B_1)$. On $X:=\opsp(A)$ define the following
family of seminorms. Let
\[
\varrho_{2n-1}(T)\ :=\ \|P_nT\|,\qquad \varrho_{2n}(T)\ :=\ \|TP_n\|,
\]
for $n=1,2,...$ and every $T\in X$, and denote the topology that is
generated on $X$ by $\{\varrho_1,\varrho_2,...\}$ by $\mathcal{T}$.
Note that, since $\mathcal{T}$ is generated by a countable family of
seminorms, the topological space $(X,\mathcal{T})$ is metrisable, so
that sequential compactness in $(X,\mathcal{T})$ coincides with
compactness. Further, by Lemma \ref{lem_Pconv} and since (Theorem
\ref{lem_opsp}(i)) $\|T\|\le\|A\|$ for every $T\in X$, convergence
in $(X,\mathcal{T})$ is equivalent to $\P-$convergence in $X$.
Therefore, by Theorem \ref{lem_opsp}(vi)-(vii),  $X$ itself and all
elements of $\B$ are compact sets in $(X,\mathcal{T})$.

Now put $\Sigma:=\cap_{B\in\sigma}\opsp(B)=\cap_{C\in\B} C$. It
follows, since $\B$ is totally ordered and each element of $\B$ is
compact (see the proof of \cite[Proposition 3.7]{CWLi:JFA2007} for
details) that $\Sigma$ is non-empty. But clearly, by Theorem
\ref{lem_opsp} (iv) again, for every $T\in\Sigma$, $\opsp(T)\in\cA$
is an upper bound of the chain $\B$. \Proofend

We will see an application of this lemma, taken from
\cite{CWLi:JFA2007}, at the end of Section \ref{sec_BDO}, and we
will meet concrete examples of self-similar operators in Section
\ref{sec_NR+AR} and Chapter \ref{sec_Schroed}.

\section{Collective Compactness and the Operator Spectrum} \label{sec:coll_comp_op_sp}
We continue to suppose that $\V$ satisfies the constraints
(\ref{constraint1})-(\ref{constraint3}) introduced in Section
\ref{subsec_limops}. This being the case, Theorem \ref{lem_opsp}
shows how nicely the operator spectrum fits the conditions made on
the set $\K$ in Theorem \ref{th4.5} if we put $S:=\V$. Indeed,
property (iii) of Theorem \ref{th4.5} is then guaranteed by Theorem
\ref{lem_opsp} (ii). Moreover, if the operator under consideration
is rich and in $S(Y)$ then, by Theorem \ref{lem_opsp} (vi), its
operator spectrum is $\P-$sequentially compact, and hence
$s-$sequentially compact by Corollary \ref{Pimps} (recall that we
have just seen after Theorem \ref{lem_opsp} that if $A\in S(Y)$ then
$\opsp(A)\subset S(Y)$).

Bearing in mind these observations, we now apply Theorem \ref{th4.5}
to the operator spectrum of an operator $A\in L(Y)$. We set $K:=I-A$
and apply Theorem \ref{th4.5} with
$$
\K\ :=\ \opsp(K)\ =\ I-\opsp(A)\qquad\textrm{so that}\qquad I-\K\ =\
\opsp(A),
$$
noting that $K$ is rich iff $A$ is rich.

\begin{Theorem} \label{prop_appl1}
Suppose $Y=\hat Y$, $A=I-K\in S(Y)$ is rich, $\V$ is sufficient,
$\opsp(K)$ is uniformly Montel on $(Y,s)$, and all the limit
operators of $A$ are injective. Then $\opsp(A)$ is uniformly bounded
below. If, moreover, there is an $s$-dense subset of surjective
operators in $\opsp(A)$ then all elements of $\opsp(A)$ are
invertible and their inverses are uniformly bounded.
\end{Theorem}

We can express the condition that $\opsp(K)$ be uniformly Montel on
$(Y,s)$ more directly in terms of properties of the operator $K$.
This is the content of the next lemma, for which we introduce the
following definition: call a sequence $(A_k)_{k\in\Z^N}\subset L(Y)$
{\em asymptotically Montel}\index{asymptotically
Montel}\index{Montel!asymptotically} on $(Y,s)$ if, for every
sequence $h=(h(n))_{n=1}^\infty\subset\Z^N$ tending to infinity and
every bounded sequence $(x_n)\subset Y$, it holds that $A_{h(n)}x_n$
has a strictly converging subsequence.

\begin{Lemma} \label{charact_uniformly_Montel}
If $K\in L(Y)$ and the sequence $(V_{-k}KV_k)_{k\in\Z^N}$ is
asymptotically Montel on $(Y,s)$, then $\opsp(K)$ is uniformly
Montel. Conversely, if $K$ is rich and $\opsp(K)$ is uniformly
Montel then $(V_{-k}KV_k)_{k\in\Z^N}$ is asymptotically Montel.
\end{Lemma}
\Proof Suppose $(V_{-k}KV_k)_{k\in\Z^N}$ is asymptotically Montel
and pick any sequence $(K_n)_{n\in\N}\subset \opsp(K)$. Then, by
definition of the operator spectrum and $\P-$conver\-gence, for
every $n\in\N$ we can find $h(n)\in \Z^N$ with $|h(n)|\ge n$ and
\begin{equation} \label{eq_something}
|| P_k(K_n-V_{-h(n)}K_nV_{h(n)})|| \le \frac{1}{n}, \quad 1\le k\le
n.
\end{equation}
Now choose any $(x_n)\subset Y$ with $||x_n||\le 1$. Then, as
$(V_{-k}KV_k)_{k\in\Z^N}$ is asymptotically Montel and $h$ tends to
infinity, $V_{-h(n)}K_nV_{h(n)}x_n$ has a strictly convergent
subsequence. On the other hand, by (\ref{eq_something}),
$$
P_k(K_n-V_{-h(n)}K_nV_{h(n)})x_n\to 0, \quad n\to\infty,
$$
for each $k\in\N$, so that $(K_n-V_{-h(n)}K_nV_{h(n)})x_n\sto 0$ by
(\ref{eqstrict}). Thus $K_nx_n$ has a strictly convergent
subsequence, so that, by Remark \ref{rem_unifmont}, $\opsp(K)$ is
uniformly Montel.

Conversely, suppose that $K$ is rich and $\opsp(K)$ is uniformly
Montel. Take an arbitrary sequence $h=(h(n))_{n=1}^\infty\subset
\Z^N$ which tends to infinity and an arbitrary bounded sequence
$(x_n)\subset Y$. Since $K$ is rich, $(h(n))$ and $(x_n)$ have
subsequences, denoted again by $(h(n))$ and $(x_n)$, such that
$V_{-h(n)}KV_{h(n)}\pto K_h\in \opsp(K)$. Thus
$$
P_k(V_{-h(n)}KV_{h(n)}- K_h)x_n\to 0, \quad n\to\infty,
$$
for each $k\in\N$, so that $(V_{-h(n)}KV_{h(n)}- K_h)x_n\sto 0$. On
the other hand, $K_hx_n$ has a strictly convergent subsequence since
$K_h$ is Montel. Thus $V_{-h(n)}KV_{h(n)}x_n$ has a strictly
convergent subsequence. \Proofend

\begin{Remark}
Note that, clearly, $(V_{-k}KV_k)_{k\in\Z^N}$ is asymptotically
Montel iff $(V_{-k}K)_{k\in\Z^N}$ is asymptotically Montel since
$\V=\{V_k\}_{k\in\Z^N}\subset iso(Y)$.
\end{Remark}

An extension of Theorem \ref{prop_appl1} can be derived by applying
Theorem \ref{th4.5} to
$$
\K\ :=\ \opsp(K)\,\cup\,\Sh(K),
$$
with $\Sh(K)$ defined by \refeq{eq_Sh(A)}, so that $I-\K =
\opsp(A)\cup\Sh(A)$. Properties (ii) and (iii) of Theorem
\ref{th4.5} can be checked in a similar way as before. Property (i)
of Theorem \ref{th4.5}, that $\K$ is uniformly Montel on $(Y,s)$, is
equivalently characterised by any of the properties (i)-(iii) of
Lemma \ref{lem_coll_comp_entries} below, which are equivalent even
for arbitrary $K\in L(Y)$. Note that, for a rich operator $K$, by
Lemma \ref{charact_uniformly_Montel}, any of these properties is
moreover equivalent to $\opsp(K)$ being uniformly Montel on $(Y,s)$
and $K\in M(Y)$. Then we get the following slightly enhanced version
of the first part of Theorem \ref{prop_appl1}, which in addition
allows to conclude from $A$ being injective to the closedness of the
range of $A$.

\begin{Theorem} \label{prop_appl1a}
Suppose $Y=\hat Y$, $A=I-K\in S(Y)$ is rich, $\V$ is sufficient, $K$
is subject to any of (i)-(iii) of Lemma \ref{lem_coll_comp_entries},
and $A$ as well as all its limit operators are injective. Then $A$
is bounded below and $\opsp(A)$ is uniformly bounded below.
\end{Theorem}

Note that Theorems \ref{prop_appl1} and \ref{prop_appl1a} are
applications of Theorem \ref{th4.5} which was just a special case of
Theorem \ref{th4.1}. We will now apply Theorem \ref{th4.1} directly.

\begin{Theorem} \label{prop5.10}
Suppose that $Y=\hat Y$, $\V$ is sufficient, $A=I-K\in L(Y)$,
$A_n=I-K_n\in L(Y)$ for $n\in\N$ and that:
\begin{itemize}
\item[(a)] $A_n\sto A$;
\item[(b)] $A_n$ bounded below $\ \Rightarrow\ \ A_n$ surjective, for
each $n\in\N$;
\item[(c)] $\cup_{n\in\N}\Sh(K_n)=\{V_{-k}K_nV_k:k\in\Z^N, n\in\N\}$ is
uniformly Montel
on $(Y,s)$;
\item[(d)] there exists a set $\B\subset L(Y)$, such that,
for every sequence $(k(m))\subset\Z^N$ and increasing sequence
$(n(m))\subset\N$, there exist subsequences, denoted again by
$(k(m))$ and $(n(m))$, and $B\in\B$ such that
\[
V_{-k(m)}A_{n(m)}V_{k(m)}\ \sto\ B\in\B\quad\textrm{as}\quad
m\to\infty;
\]
\item[(e)] every $B\in\B$ is injective.
\end{itemize}
Then $A$ is invertible and, for some $n_0\in\N$, $A_n$ is invertible
for all $n\ge n_0$, and
\[
\|A^{-1}\|\le\sup_{n\ge n_0}\|A_n^{-1}\| < \infty.
\]
\end{Theorem}
\Proof Let $\K_n:=\Sh(K_n)$, $n\in\N$, and set $\K:=I-\B$, $S:=\V$.
Then (c) and (d) imply that conditions (i)--(iv) of Theorem
\ref{th4.1} are satisfied, and (a) and (b) imply that the condition
in Theorem \ref{th4.1} b) is satisfied. Thus, applying Theorem
\ref{th4.1}, the result follows. \Proofend

Later on, in Chapter 6, this result will be used to derive Theorem
\ref{prop_NRBO} on the invertibility of norm-rich/almost periodic
band operators. As described at the end of Section \ref{sec_BDO},
this result has also been used to prove Propositions 3.4 and 3.8 in
\cite{CWLi:JFA2007}.

\begin{Remark} \label{rem_long}
Note that condition (a) in the above theorem implies that $A\in
S(Y)$ by Lemma \ref{rel_conv}. Moreover, from condition (d) with
$k(m)=0$ for all $m\in\N$ and condition (a) again we get that $A\in
\B$. Since $A\in S(Y)$ it holds that $\opsp(A)\subset S(Y)$ (see
discussion at the end of Section \ref{subsec_limops}); if also
$A_n\pto A$ then condition (d) also implies that
$\opsp(A)\subset\B$. To see this last claim, suppose that $\tilde
A\in\opsp(A)$. Then there exists $(k(m))\subset\Z^N$ such that
$V_{-k(m)}AV_{k(m)}\pto \tilde A$ which implies, in particular, that
$\|P_j(\tilde A - V_{-k(m)}AV_{k(m)})\|\to 0$ as $m\to \infty$, for
every $j$. Choose the sequence $(l(m))\subset\N$ such that $P_n
V_{-k(m)}P_{l(m)}=P_nV_{-k(m)}$, for $n=1,...,m$, which is possible
by (\ref{constraint3}). Then
\[
\|P_j(\tilde A-V_{-k(m)}AV_{k(m)})\|=\|P_j(\tilde
A-V_{-k(m)}P_{l(m)}AV_{k(m)})\|
\]
for all $j$ and all $m>j$. Since $A_n\pto A$, for every $m$ we can
choose $n(m)$ such that $\|P_{l(m)}(A_{n(m)}-A)\|<m^{-1}$. Then
\[
\|P_j(\tilde A-V_{-k(m)}A_{n(m)}V_{k(m)})\|=\|P_j(\tilde
A-V_{-k(m)}P_{l(m)}A_{n(m)}V_{k(m)})\|\to 0
\]
as $m\to\infty$, for every $j$, so that, by Lemma \ref{lll},
$V_{-k(m)}A_{n(m)}V_{k(m)}\sto\tilde A$, and so $\tilde A\in\B$.
\end{Remark}

In the case that $A=I-K$ with $K\in S(Y)\cap M(Y)$, one particular
choice of $A_n$ which satisfies (a) and (b) is
\[
A_n=I-KP_n.
\]
For, by Lemma \ref{lem_charact_mon2}, $KP_n\in KS(Y)$ for every $n$
so that, by a version of Riesz Fredholm theory for TVS's (see e.g.
\cite{Robertson}), assumption (b) holds. Further, by Corollary
\ref{cor4}, for every $m$,
\[
\|P_m(A_n-A)\|\ =\ \|P_mKQ_n\|\ \to\ 0,\qquad n\to\infty
\]
so that, by Lemma \ref{lll}, $A_n\sto A$. Further, if
$\Sh(K)=\{V_{-k}KV_k:k\in\Z^N\}$ is uniformly Montel on $(Y,s)$,
then so is $\{V_{-k}K:k\in\Z^N\}$ and hence also
$\{V_{-k}KP_nV_k:k\in\Z^N,n\in\N\}=\cup_{n\in\N}\Sh(KP_n)$. Thus
Theorem \ref{prop5.10} has the following corollary.

\begin{Corollary} \label{cor6.14}
Suppose that $Y=\hat Y$, $\V$ is sufficient, $A=I-K$ with $K\in
S(Y)$, and set $A_n=I-KP_n$ for $n\in\N$. Suppose that $\Sh(K)$ is
uniformly Montel on $(Y,s)$, and that there exists $\B\subset L(Y)$
such that every $B\in\B$ is injective and for every sequence
$(k(m))\subset\Z^N$ and increasing sequence $(n(m))\subset\N$, there
exist subsequences, denoted again by $k(m)$ and $n(m)$, and
$B\in\B$, such that
\[
I-V_{-k(m)}KP_{n(m)}V_{k(m)}\sto B\in\B.
\]
Then $A$ is invertible and, for some $n_0\in\N$, $A_n$ is invertible
for all $n\ge n_0$, and
\[
\|A^{-1}\|\le\sup_{n\ge n_0}\|A_n^{-1}\|<\infty.
\]
\end{Corollary}

\begin{Remark}
We note, by Remark \ref{rem_long}, that, necessarily, $A\in\B$ and
$\opsp(A)\subset\B$.
\end{Remark}

More concrete statements than Theorems \ref{prop_appl1},
\ref{prop_appl1a} and \ref{prop5.10} can be given when we pass to a
more concrete class of spaces $Y$. This is what we do in Chapter
\ref{sec_ellp}.


\chapter{Operators on $\ell^p(\Z^N,U)$} \label{sec_ellp}
In this chapter we focus on the case, introduced already briefly in
Example \ref{ex_mult}, when
$Y=\ell^p(\Z^N,U)$,\index{$\ell^p(\Z^N,U)$} where $1\le p\le\infty$,
$N\in\N$ and $U$ is an arbitrary complex Banach space. The elements
of $Y$ are of the form $x=(x(m))_{m\in\Z^N}$ with $x(m)\in U$ for
every $m=(m_1,...,m_N)\in\Z^N$. We equip $Y$ with the usual $\ell^p$
norm of the scalar sequence $(\|x(m)\|_U)$. We also consider the
case when $Y=c_0(\Z^N,U)$\index{$\c0$@$c_0(\Z^N,U)$}, the Banach
subspace of $\ell^\infty(\Z^N,U)$ consisting of the elements that
vanish at infinity, i.e. $\|x(m)\|_U\to 0$ as $m\to\infty$.

Since the parameter $N\in\N$ is of no big importance in almost all
of what follows, we will use the abbreviations
$Y^0(U):=c_0(\Z^N,U)$\index{$Y^0(U)$} and
$Y^p(U):=\ell^p(\Z^N,U)$\index{$Y^p(U)$} for $1\le p\le\infty$. If
there is no danger of confusion about what $U$ is, we will even
write $Y^0$ and $Y^p$. Some of our following statements hold for all
the spaces under consideration. In this case we will simply write
$Y$, which then can be replaced by any of $Y^0$\index{$Y^0$} and
$Y^p$\index{$Y^p$} with $1\le p\le\infty$.

In terms of dual spaces\index{dual space}, we have $(Y^0(U))^*\cong
Y^1(U^*)$, $(Y^1(U))^*\cong Y^\infty(U^*)$, and $(Y^p(U))^*\cong
Y^q(U^*)$ for $1<p<\infty$ and $1/p+1/q=1$ (see e.g. \cite{Rosier}).
To give two prominent examples, the space $Y^p(\C)$ is the usual
$\ell^p$ or $c_0$ space of complex-valued sequences over $\Z^N$, and
the space $Y^p(L^p([0,1]^N))$ is isometrically isomorphic to
$L^p(\R^N)$ for $1\le p\le\infty$; see Chapter \ref{sec_intop} or
\cite{Kurbatov,RaRoSiBook,LiBook} for a more detailed discussion of
this isomorphism.

For $m\in\Z^N$, we define
$|m|:=\max(|m_1|,...,|m_N|)$\index{$\mid\cdot\mid$} and put
$$
P_n x(m) = \left\{\begin{array}{ll}
  x(m), &|m|\le n, \\
  0, &|m|>n, \\
\end{array} \right.\index{$P_n$}
$$
for every $x\in Y$ and $n\in\N_0$. As observed already in Example
\ref{ex_mult},  $\P=(P_n)$ satisfies conditions (i) and (ii) of Chapter 2 with $N(m) = m$,
so that $P_n$ is a projection operator for each $n$. In this case
$\|Q_n\|=1$ for all $n$. For $p\in\{0\}\cup [1,\infty)$, we have
$(Y^p)_0=Y^p$, while $(Y^\infty)_0=Y^0$. Moreover, for
$p\in[1,\infty]$, we have $\widehat{Y^p}=Y^p$ (cf.\ Example
\ref{ex_Lp2}), whereas $\widehat{Y^0}=Y^\infty$.

In the setting of this chapter we have the following refinement of
Lemma \ref{lem2'}.

\begin{Lemma} \label{lem_KK}
The following two statements hold.
\begin{itemize}
\item[(i)] If $p\in\{0\}\cup(1,\infty)$ then $K(Y^p)\subset K(Y^p,\P)$.
\item[(ii)] If $U$ is finite-dimensional then $K(Y,\P)\subset K(Y)$.
\end{itemize}
\end{Lemma}
\Proof (i) Let $p\in\{0\}\cup (1,\infty)$. Then $P_n\to I$ as well
as $P_n^*\to I^*$, on $Y^p$ and $(Y^p)^*$, respectively, as
$n\to\infty$, i.e.\ $\P$ is perfect. Thus $K(Y^p)\subset K(Y^p,\P)$
by Lemma \ref{per_lem}.

(ii) If $\dim U<\infty$ then $P_n\in K(Y)$ for all $n$, and, by
Lemma \ref{lem2'}, we get $K(Y,\P)=L(Y,\P)\cap K(Y)\subset K(Y)$.
\Proofend

\begin{Remark} \label{rem_Fred_invinf}
From Lemma \ref{lem_KK} and \refeq{eq_inv_inf} we conclude another
result, in addition to Theorem \ref{lem_Fredholm} and Corollary
\ref{cor_invinf_Fred}, which relates Fredholmness to invertibility
at infinity in the setting of this chapter: For an arbitrary
operator $A\in L(Y^p(U))$, Fredholmness implies invertibility at
infinity if $p\in\{0\}\cup (1,\infty)$, and invertibility at
infinity implies Fredholmness if $\dim U<\infty$.
\end{Remark}

If we moreover put, for $x\in Y$,
\begin{equation} \label{eq_shifts}
V_kx(m)\ :=\ x(m-k),\quad m\in\Z^N\index{$V_k$}
\end{equation}
for every $k\in\Z^N$, then the family $\V=\{V_k\}_{k\in\Z^N}$
consists of isometric isomorphisms on $Y$ and satisfies conditions
\refeq{constraint1}--\refeq{constraint3} so that it is compatible
with $\P$.

In the case $p=\infty$ (but not if $1\le p<\infty$) the family
$\V=\{V_k\}_{k\in\Z^N}$ is sufficient (see Example \ref{e52}). Thus
each of Theorems \ref{prop_appl1}, \ref{prop_appl1a}, \ref{prop5.10}
and Corollary \ref{cor6.14} can be applied in the case $Y=Y^\infty$.
Note also that for $Y^\infty$ Assumption A of Chapter
\ref{chap_converg} holds with $Q=Q_m$ and $n=m+1$ so that Corollary
\ref{strong2} applies and that therefore the convergence $\sto$ can
be replaced by the formally weaker notion $\Sto$ of convergence in
the strong operator topology on $(Y^\infty,s)$, defined by
(\ref{eq_Sto}). Applying Theorem \ref{prop_appl1}, for example, we
have the following result in which we say that a subset
$A\subset\K\subset L(Y^\infty)$ is {\em
$S-$dense}\index{S-dense@$S$-dense} in $\K$ if, for every $K\in\K$,
there is a sequence $(K_n)\subset A$ with $K_n\Sto K$.

\begin{Theorem} \label{prop_appl2}
Suppose $A=I-K\in S(Y^\infty)$ is rich, $\opsp(K)$ is uniformly
Montel on $(Y^\infty,s)$, and all the limit operators of $A$ are
injective. Then $\opsp(A)$ is uniformly bounded below. If, moreover,
there is an $S$-dense subset of surjective operators in $\opsp(A)$
then all elements of $\opsp(A)$ are invertible and their inverses
are uniformly bounded.
\end{Theorem}

\begin{Remark} \label{rem_finite_dim}
We have seen in Lemma \ref{charact_uniformly_Montel} that $\opsp(K)$
is uniformly Montel on $(Y^\infty,s)$ iff the sequence
$(V_{-k}KV_k)_{k\in\Z^N}$ is asymptotically Montel. If the Banach
space $U$ is of finite dimension, then the condition that $\opsp(K)$
be uniformly Montel on $(Y^\infty,s)$ is even redundant. For $U$
finite-dimensional implies that $P_n\in K(Y^\infty)$ for all $n$, so
that $\opsp(K)$ is uniformly Montel by Corollary \ref{cor_UMUB} and
Theorem \ref{lem_opsp} (i).
\end{Remark}

\section{Periodic and Almost Periodic Operators} \label{sec_NR+AR}
With regard to Theorem \ref{prop_appl2} we remark that we know of no
examples where the requirement for an $S-$dense subset of surjective
operators is not redundant. Precisely, in all the examples we have
studied existence of an $S-$dense subset of $\opsp(A)$ of surjective
operators can be deduced from injectivity of all the elements of
$\opsp(A)$. To illustrate Theorem \ref{prop_appl2} we present next
an important example of this type in which even invertibility of $A$
(on $Y^\infty$) follows from injectivity of all limit operators of
$A$. We will present further examples in sections \ref{sec_BDO} and
\ref{sec_NRBDO}.

To introduce this example we require the following definitions.
Recall from Theorem \ref{lem_opsp} (v) that $A\in L(Y)$ is rich iff
the set $\Sh(A)$ of all translates of $A$ (recall \refeq{eq_Sh(A)})
is relatively $\P-$sequentially compact. Now call $A\in L(Y)$
\emph{norm-rich}\index{norm-rich
operator}\index{operator!norm-rich}\index{operator!rich!norm-} or
\emph{almost periodic}\index{almost periodic
operator}\index{operator!almost periodic} if the set $\Sh(A)$ is
relatively compact in the norm topology on $L(Y)$ (this is precisely
the definition of an almost periodic
 operator in Kurbatov \cite{Kurbatov1989}). Call $A\in L(Y)$ \emph{absolutely
rich}\index{absolutely rich operator}\index{operator!absolutely
rich}\index{operator!rich!absolutely} or
\emph{periodic}\index{periodic operator}\index{operator!periodic} if
every sequence in $\Sh(A)$ has a constant subsequence, i.e. iff
$\Sh(A)$ is a finite set. It is easy to establish the following
characterisation.

\begin{Lemma} \label{lem_absrich}
An operator $A\in L(Y)$ is absolutely rich/periodic iff there exist
$m_1,...,m_N\in\N$ such that
\[
VA=AV\qqtext{for all}
V\in\Vt_A:=\{V_{m_je^{(j)}}\}_{j=1}^N\index{$Vaa$@$\Vt_A$}
\]
with $e^{(1)},...,e^{(N)}$ denoting the standard unit
vectors\index{$\e$@$e^{(k)}$}\index{unit coordinate vector} in
$\R^N$, i.e. $e^{(j)}(i)=1$ if $i=j$ and $=0$ otherwise.
\end{Lemma}

\begin{Example} \label{ex_almostperiodic}
For $b=(b(m))_{m\in\Z^N}\in\ell^\infty(\Z^N,L(U))=Y^\infty(L(U))$
define the multiplication operator $M_b\in L(Y)$, as in Example
\ref{ex_mult}, by equation (\ref{eq_multop}). Then, for $k\in\Z^N$,
\begin{equation} \label{eqn_multiden}
V_{-k}M_bV_k=M_{V_{-k}b}\qqtext{and} \|M_b\|=\|b\|.
\end{equation}

From these identities and (\ref{eqn_Pconfmult}), it is immediate
that $M_b$ is rich iff the set
\begin{equation} \label{eq_shiftsb}
\{V_kb\}_{k\in\Z^N}
\end{equation}
is relatively sequentially compact in the strict topology on
$Y^\infty(L(U))$. It can be shown moreover \cite[Theorem
2.1.16]{RaRoSiBook} that this is the case iff the set
$\{b(m):m\in\Z^N\}$ is relatively compact in $L(U)$. Also, it is
clear from (\ref{eqn_multiden}) that $M_b$ is norm rich iff the set
(\ref{eq_shiftsb}) is relatively compact in the norm topology on
$Y^\infty(L(U))$; it is usual to say (cf.\ Example \ref{ex_ap} and
\cite[Definition 3.58]{LiBook}) that $b$ is \emph{almost
periodic}\index{almost periodic sequence} if this condition on $b$
holds. The set of all almost periodic sequences $b\in Y^\infty(Z)$
with a Banach space $Z$ shall be denoted by
$\AP(Z)$\index{$Y_n^\inftyAP$@$\AP$}.

Similarly, $M_b$ is absolutely rich/periodic iff every sequence in
\refeq{eq_shiftsb} has a constant subsequence, i.e. iff
\refeq{eq_shiftsb} is finite.
By Lemma \ref{lem_absrich} this is
equivalent to the requirement that there exist $m_1,...,m_N\in\N$
such that
\[
b(k+m_je^{(j)})=b(k),\qquad k\in\Z^N,\ j=1,...,N,
\]
i.e. to the requirement that $b(k)$ is periodic as a function of
each of the components of $k\in\Z^N$.
\end{Example}

Suppose that $A\in L(Y^\infty)$ is absolutely rich/periodic. For $n\in\N$,
let
\begin{equation} \label{eq_Yn}
Y_n^\infty\ =\ Y_n^\infty(U)\ :=\ \{x\in Y^\infty(U):V^nx=x\textrm{
for all } V\in\Vt_A\}\index{$Y_n^\infty$}
\end{equation}
with $\Vt_A$ as defined in Lemma \ref{lem_absrich}. Then
$Y_n^\infty$ is a closed subspace of $Y^\infty$ consisting of
periodic elements; $x\in Y_n^\infty$ iff
\[
x(k\,+\,n\,m_j\,e^{(j)})\ =\ x(k),\qquad k\in\Z^N,\ j=1,...,N,
\]
where the integers $m_1,...,m_N$ are as in the definition of $\Vt_A$
in Lemma \ref{lem_absrich}. Clearly, $x\in Y_n^\infty$ is determined
by its components in the box
\[
C_n\ :=\ \{i=(i_1,...,i_N)\in\Z^N:-n\frac{m_j}2< i_j\le
n\frac{m_j}2,\ j=1,...,N\}.
\]
Define the projection operator $\tilde P_n:Y^\infty\to Y_n^\infty$
by the requirement that
\begin{equation} \label{eq_PnYn}
\tilde P_nx(k)=x(k)\qquad\textrm{for all}\qquad k\in
C_n.\index{$P_ntil$@$\tilde P_n$}
\end{equation}
Then, clearly, for each $n$, $\tilde P_n Q_j=0$ for all sufficiently
large $j$, so that $\tilde P_n\in SN(Y^\infty)$ by Lemma \ref{lem1}.
(Note however that $\tilde P_n\not\in L(Y^\infty,\P)$.)

The last part of the following result and its proof can be seen as a
generalisation of Theorem 2.10 in \cite{ChanZh97}
\begin{Theorem} \label{lem_A-inj-inv}
If $A\in L(Y^\infty)$ is absolutely rich/periodic then $A(Y_n^\infty)\subset
Y_n^\infty$ for each $n$, and
\begin{equation} \label{eq_x}
\opsp(A)=\{V_{-i}AV_i:i\in\Z^N\}=\{V_{-i}AV_i:i\in C_1\}.
\end{equation}
If also $A=I+K$ with $K\in S(Y^\infty)\cap M(Y^\infty)$ and $A$ is
injective then $A$ is invertible.
\end{Theorem}
\Proof If $x\in Y_n^\infty$ then $Ax\in Y_n^\infty$ since
$V^n(Ax)=AV^nx=Ax$ for every $V\in\Vt_A$. From the definitions and
Lemma \ref{lem_absrich} it is clear that \refeq{eq_x} holds. Suppose
now that $K\in S(Y^\infty)\cap M(Y^\infty)$ and that $A=I+K$ is
injective. First we show that, for every $n$, $I+K\tilde P_n$ is
invertible. To see injectivity, suppose $x\in Y^\infty$ and
$(I+K\tilde P_n)x=0$. Then $x=-K\tilde P_nx\in Y^\infty_n$ since
$\tilde P_nx\in Y^\infty_n$ and $K=A-I$ is absolutely rich. Now
$x\in Y^\infty_n$ implies $\tilde P_nx=x$ and therefore
$0=(I+K\tilde P_n)x=(I+K)x=Ax$, i.e. $x=0$ by injectivity of $A$.
Now surjectivity of $I+K\tilde P_n$ follows from its injectivity by
the Riesz theory for compact operators in topological vector spaces
\cite{Robertson} since $K\tilde P_n\in KS(Y^\infty)$ by Lemma
\ref{lem_alg}.

Next, note that from \refeq{eq_x} it follows that $A\in\opsp(A)$ and
that, since $A$ is injective, all the limit operators of $A$ are
injective. Further, by \refeq{eq_x} and Remark \ref{rem_finite_dim},
it follows that $\opsp(K)$ is uniformly Montel since $K\in
M(Y^\infty)$. Applying Theorem \ref{prop_appl2} we see that the
limit operators of $A$ are uniformly bounded below, in particular
that $A$ is bounded below.

To see finally that $A$ is surjective let $y\in Y^\infty$ and set
$y_n=\tilde P_ny\in Y^\infty_n$ and $x_n=(I+K\tilde P_n)^{-1}y_n$ so
that $x_n+K\tilde P_nx_n=y_n$ which implies (as seen above) that
$x_n\in Y^\infty_n$, that $\tilde P_nx_n=x_n$, and hence that
$Ax_n=x_n+Kx_n=y_n$. Since $(y_n)$ is bounded and $A$ is bounded
below, also $(x_n)$ is bounded. Since $K\in M(Y^\infty)$ it follows
that there exists an $x\in Y^\infty$ and a subsequence of $(x_n)$,
denoted again by $(x_n)$, such that $Kx_n\sto y-x$, so that
$x_n=y_n-Kx_n\sto y-(y-x)=x$. As $K\in S(Y^\infty)$ this implies
$Kx_n\sto Kx$ so that $Ax=x+Kx=y$. \Proofend

The above result has the following obvious corollary (phrased in the
spirit of Theorem \ref{prop_appl2}, which was the starting point of
this discussion).

\begin{Corollary} \label{cor_absrich_injinv}
If $A=I+K\in S(Y^\infty)$ is absolutely rich/periodic and $K\in M(Y^\infty)$
and if all limit operators of $A$ are injective, then all limit
operators of $A$ are invertible (with uniformly bounded inverses).
\end{Corollary}

Further down, in Theorem \ref{prop_NRBDO}, we will show that, in the
case when $A$ is also band-dominated (as defined in
\S\ref{sec_BDO}), this corollary holds more generally with
`absolutely rich' replaced by `norm-rich'; indeed in the
one-dimensional case $N=1$ we will show in Theorem \ref{prop_FC_JFA}
that this corollary holds even with `absolutely rich' replaced by
`rich'. We conclude the current section by a collection of results
for the general setting of all norm-rich/almost periodic operators
$A\in L(Y)$; a set that shall be denoted by
$\NR(Y)$\index{$Lndol$@$\NR(Y)$} for brevity. Our first result
follows from a slightly more general result which is Theorem 6.5.2
in Kurbatov \cite{Kurbatov}. We include a proof here for the
convenience of the reader.

\begin{Lemma} \label{lem_NRisBA}
$\NR(Y)$ is an inverse closed Banach subalgebra of $L(Y)$.
\end{Lemma}
\Proof Let $A,B\in\NR(Y)$ and take an arbitrary sequence
$h=(h(1),h(2),...)\subset\Z^N$. Pick a subsequence $g$ of $h$ such
that both $V_{-g(n)}AV_{g(n)}$ and $V_{-g(n)}BV_{g(n)}$ converge in
norm. Then, clearly, also $V_{-g(n)}(A+B)V_{g(n)}$ and
\[
V_{-g(n)}(AB)V_{g(n)}=(V_{-g(n)}AV_{g(n)})(V_{-g(n)}BV_{g(n)})
\]
converge in norm. To see that $\NR(Y)$ is closed in the operator
norm take $A_1,A_2,...\in\NR(Y)$ with $A_k\toto A\in L(Y)$ and an
arbitrary sequence $h=(h(n))_{n=1}^\infty\subset\Z^N$. Pick
subsequences $\cdots\subset h^{(2)}\subset h^{(1)}\subset h$ such
that, for every $k\in\N$,
\begin{equation} \label{eq_subseq}
\|V_{-h^{(k)}(m)}A_kV_{h^{(k)}(m)}-V_{-h^{(k)}(n)}A_kV_{h^{(k)}(n)}\|<
1/k,\qquad m,n>k,
\end{equation}
and put $g(n):=h^{(n)}(n)$ for all $n\in\N$. Then, for all $k\in\N$
and all $m,n>k$, noting that $g(n)=h^{(k)}(n')$ for some $n'\ge
n>k$,
\begin{eqnarray*}
\|V_{-g(m)}AV_{g(m)}-V_{-g(n)}AV_{g(n)}\| & \le &
\|V_{-g(m)}A_kV_{g(m)}-V_{-g(n)}A_kV_{g(n)}\|\\
& &\hspace{32.4mm} +\ 2\|A_k-A\|\\
&\le &1/k +\ 2\|A_k-A\|\ \to\ 0
\end{eqnarray*}
as $k\to\infty$. This shows that the sequence $(V_{-g(n)}AV_{g(n)})$
is Cauchy and therefore convergent in $L(Y)$. Since $g\subset h$, we
get that $A\in\NR(Y)$.

To see the inverse closedness suppose $A\in\NR(Y)$ is invertible in
$L(Y)$ and take an arbitrary sequence
$h=(h(1),h(2),...)\subset\Z^N$. Since $A\in\NR(Y)$, there is a
subsequence $g$ of $h$ such that $A_n:=V_{-g(n)}AV_{g(n)}\toto B$
for some $B\in L(Y)$. Since
$\|A_n^{-1}\|=\|V_{-g(n)}A^{-1}V_{g(n)}\|=\|A^{-1}\|$ is bounded
independently of $n$, it follows from a basic result on Banach
algebras (see, e.g. Lemma 1.3 of \cite{LiBook}) that $B$ is
invertible and
\[
A_n^{-1}\ =\ V_{-g(n)}A^{-1}V_{g(n)}\ \toto\ B^{-1},
\]
showing that $A^{-1}\in\NR(Y)$. \Proofend

\begin{Theorem} \label{prop_NR1}
For $A\in\NR(Y)$, the following holds.
\begin{itemize}
\item[(i)] If, for some sequence $h=(h(1),h(2),...)\subset\Z^N$ and $B\in
L(Y)$,
\[
V_{-h(n)}AV_{h(n)}\pto B\quad\textrm{holds,\ then}\quad
V_{-h(n)}AV_{h(n)}\toto B.
\]
\item[(ii)] $A\in\opsp(A)$ (i.e.\ $A$ is self-similar).
\item[(iii)] $\opsp(A)=\clos_{L(Y)}\Sh(A)$ is a compact subset
of $\NR(Y)$.
\item[(iv)] $A$ is invertible iff any one of its limit operators is
invertible.
\item[(v)] $\nu(A)=\nu(B)$ for all $B\in\opsp(A)$, so that $A$ is
bounded below iff $\opsp(A)$ is uniformly bounded below.
\item[(vi)] If $x$ is almost periodic, then $Ax$ is almost periodic.
\item[(vii)] If $A$ is invertible on $Y^\infty$, then it is invertible on
$\AP$.
\end{itemize}
\end{Theorem}
\Proof (i) Since $A\in\NR(Y)$, every subsequence of
$V_{-h(n)}AV_{h(n)}$ $(\pto B)$ has a norm-convergent subsequence
the limit of which must be $B$. But this proves norm convergence of
the whole sequence.

(ii) Let $h(n)=(n^2,0,...,0)\in\Z^N$ for every $n\in\N$. Since
$A\in\NR(Y)$, there is a subsequence $g$ of $h$ such that
$V_{-g(n)}AV_{g(n)}$ converges. But then
\[
\|V_{-(g(n+1)-g(n))}AV_{g(n+1)-g(n)}-A\|=\|V_{-g(n+1)}AV_{g(n+1)}-V_{-g(n)}AV_{g(n)}\|\to
0
\]
as $n\to\infty$, showing that $A=A_f\in\opsp(A)$ with
$f(n)=g(n+1)-g(n)\to\infty$.

(iii) The inclusion $\opsp(A)\subset\clos_{L(Y)}\Sh(A)$ follows from
(i). The reverse inclusion follows from (ii), from Theorem
\ref{lem_opsp} (ii) and the closedness of $\opsp(A)$ (see Theorem
\ref{lem_opsp} (iii) above or \cite[Corollary 3.96]{LiBook}). The
compactness of $\clos_{L(Y)}\Sh(A)$ follows from the relative
compactness of $\Sh(A)$ in $L(Y)$. By Lemma \ref{lem_NRisBA}, every
operator in $\clos_{L(Y)}\Sh(A)$ is norm-rich/almost periodic.

(iv) Take an arbitrary limit operator $A_h$ of $A$ and let
$h=(h(1),h(2),...)\subset\Z^N$ be such that
$A_n:=V_{-h(n)}AV_{h(n)}\pto A_h$ holds. By (i) we have that
$A_n\toto A_h$. If $A_h$ is invertible, then so is $A_n$ for every
large $n$, and therefore $A$ is invertible. Conversely, if $A$ is
invertible, then $A_h$ is invertible by a basic result on Banach
algebras (see e.g. \cite[Lemma 1.3]{LiBook}) since
$\|A_n^{-1}\|=\|A^{-1}\|$ is bounded.

(v) If $B\in\opsp(A)$, then, by (i), we have that
$V_{-h(n)}AV_{h(n)}\toto B$ for some sequence $h(1),h(2),...$ in
$\Z^N$. By \refeq{eq_lonoemineq} this implies that
$\nu(V_{-h(n)}AV_{h(n)})\to\nu(B)$ as $n\to\infty$. On the other
hand, since every $V_{h(n)}$ is an isometry, we have that
\[
\nu(V_{-h(n)}AV_{h(n)})\ =\ \inf_{\|x\|=1}\|V_{-h(n)}AV_{h(n)}x\| \
=\ \inf_{\|y\|=1}\|Ay\|\ =\ \nu(A)
\]
for every $n\in\N$, so that
$\nu(A)=\nu(V_{-h(n)}AV_{h(n)})\to\nu(B)$, i.e. $\nu(A)=\nu(B)$.

(vi) Let $h=(h(1),h(2),...)\subset\Z^N$ be arbitrary. If
$A\in\NR(Y^\infty)$ and $x\in\AP$ there is a subsequence $g$ of $h$
such that both $V_{g(n)}AV_{-g(n)}$ and $V_{g(n)}x$ converge in the
norm of $L(Y^\infty)$ and $Y^\infty$, respectively. But then also
\[
V_{g(n)}(Ax)\ =\ (V_{g(n)}AV_{-g(n)}) (V_{g(n)}x)
\]
converges in $Y^\infty$, which shows that $Ax\in\AP$.

(vii) If $A\in\NR(Y^\infty)$ is invertible on $Y^\infty$, then, by
Lemma \ref{lem_NRisBA}, also $A^{-1}\in\NR(Y^\infty)$. Now (vi)
shows that $x\in\AP$ iff $Ax\in\AP$. \Proofend

\section{Dual Space Arguments} \label{sec_dsa}
In the results we will present below dual space arguments will play
a role, in particular in the cases $p=1$ and $p=\infty$. Temporarily
set $Y=\ell^\infty(\Z^N,U)$ so that $Y_0 = c_0(\Z^N,U)$, and set
$Y_1 = \ell^1(\Z^N,U^*)$, where $U^*$\index{$U^*$} is the dual
space\index{dual space} of $U$. Note that $Y_1=Y_0^*$. Further,
$Y_1^*=Y$ if $U$ is reflexive, i.e.\ if $U=U^{**}$. In the general
case when $U$ is not reflexive we shall see that we can, in a
natural way, embed $Y$ as a closed subspace of $Y_1^*$.

For $x\in Y$, $y\in Y_1$, define the bilinear form $(\cdot,\cdot)$
on $(Y,Y_1)$ by
\begin{equation} \label{eq_bildef}
(x,y) := \sum_{j\in\Z^N} y_j(x_j), \quad \mbox{ for } x
=(x_j)_{j\in\Z^N}\in Y, \; y =(y_j)_{j\in\Z^N}\in Y_1,
\end{equation}
and note that, equipped with $(\cdot,\cdot)$, $(Y,Y_1)$ is a dual
system\index{dual system} in the sense e.g.\ of J\"orgens
\cite{Jorgens}. If $x\in Y_0$ and $y\in Y_1=Y_0^*$ then
$(x,y)=y(x)$. A similar equation holds if $U$ has a predual space
$U^\pa$, i.e.\ if there exists a Banach space $U^\pa$ such that
$(U^\pa)^* = U$. Then $Y$ is the dual space of $Y^\pa$, where $Y^\pa
:= \ell^1(\Z^N,U^\pa)$. Denote by $J_U$ the canonical embedding of
$U^\pa$ into its second dual $U^*$, given by $J_U u(v) = v(u)$,
$u\in U^\pa$, $v\in U$, and let $J^\pa:\ell^1(\Z^N,U^\pa)\to
\ell^1(\Z^N,U^*)$ be the natural embedding $J^\pa x = (J_U
x_j)_{j\in\Z^N}$. Note that both $J_U$ and $J^\pa$ are isometries.
Then
\begin{equation} \label{eq_blff}
(x,J^\pa y) = x(y), \quad x\in Y,\; y\in Y^\pa.
\end{equation}
An important observation is that, if $(x_n)\subset Y$, $x\in Y$, and
$y\in Y_1$, then
\begin{equation} \label{eq_cont_blf}
x_n\sto x\quad\Rightarrow\quad (x_n,y)\to (x,y).
\end{equation}

For every $A\in L_0(Y)$ (recall the definition (\ref{eq_L0}) of
$L_0(Y)$) let $A_0\in L(Y_0)$ be defined by $A_0:=
A|_{Y_0}$\index{$A_0$, restriction of $A$ to $Y_0$}. Then its
adjoint\index{adjoint operator}\index{operator!adjoint}\index{$A^*$,
adjoint of $A$} $A_0^*\in L(Y_0^*)=L(Y_1)$. From Corollary
\ref{cor_LYPinL0} we recall that, in particular, $A\in L_0(Y)$ if
$A\in L(Y,\P)\subset S(Y)$.

\begin{Lemma} \label{lem_trans}
If $A\in S(Y)\cap L_0(Y)$ then
$$
(Ax,y) = (x, A_0^*y), \quad x\in Y,\; y\in Y_1,
$$
i.e.\ $A$ is the {\em transpose}\index{transpose
operator}\index{operator!transpose} of $A_0^*$ with respect to the
dual system $(Y,Y_1)$.
\end{Lemma}
\Proof For $x\in Y_0$, $y\in Y_0^*=Y_1$,
$$
(Ax,y)=(A_0x,y) = (x,A_0^*y).
$$
Thus, for $x\in Y$, $y\in Y_1$,
$$
(AP_nx,y)=(P_nx,A_0^*y).
$$
Taking the limit $n\to\infty$, in view of (\ref{eq_cont_blf}) and
since $A\in S(Y)$, the result follows. \Proofend

It follows from Lemma \ref{lem_KK}(i) that if $A_0\in L(Y_0)$ is
Fredholm then it is invertible at infinity (cf.\ Remark
\ref{rem_Fred_invinf}). We shall see that under certain conditions
the same implication holds for $A\in L(Y)$.  Our tool to establish
this will be to relate Fredholmness of $A$ to that of $A_0$. To this
end a useful tool is to embed $Y$ as a closed subspace of
$Y_0^{**}=\ell^\infty(\Z^N,U^{**})$. Precisely, define $J:Y\to
Y_0^{**}$ by
$$
Jx(y) = (x,y),\quad x\in Y,\; y\in Y_0^*=Y_1.
$$
It is easy to check that $J$ is an isometry, so that $Y$ is
isometrically isomorphic to $\check Y:= J(Y)\subset Y_0^{**}$. For
$A\in L(Y)$ define $\check A\in L(\check Y)$ by $\check A :=
JAJ^{-1}$.

\begin{Lemma} \label{lem_hatA}
If $A\in S(Y)\cap L_0(Y)$, then $A_0^{**}(\check Y)\subset \check Y$
and $\check A = A_0^{**}|_{\check Y}$, so that
$$
\alpha(A_0)\leq \alpha(A)=\alpha(\check A) \leq \alpha(A_0^{**}).
$$
\end{Lemma}
\Proof For $x\in \check Y$, $y\in Y_1$, with $z:=J^{-1}x\in Y$,
$$
\check Ax(y) = (J(Az))(y)=(Az,y) = (z,A_0^*y),
$$
by Lemma \ref{lem_trans}, and
$$
A_0^{**}x(y) = x(A_0^*y) = (Jz)(A_0^*y)=(z,A_0^*y).
$$
 \Proofend

To make full use of the above observation, we need the following
characterisation, for a Banach space $Z$, of those operators $C\in
L(Z)$ whose range is closed, which is a standard corollary of the
open mapping theorem (applied to the injective operator $z+\ker
C\mapsto Cz$ from $Z/\ker C$ to $Z$, also see \cite[Theorem
XI.2.1]{GohGoldKash}): that
\begin{equation} \label{eq_closed}
C(Z) \mbox{ is closed }\quad \Leftrightarrow\quad \exists c>0 \mbox{
s.t. } \|Cz\|\geq c\inf_{y\in \ker C}\|z-y\|, \; \forall z\in Z.
\end{equation}
We also need the following consequence of the above
characterisation.

\begin{Lemma} \label{lem_bs} Suppose that $Z$ is a Banach space, $Z_0$ is a
closed subspace of
$Z$, $C(Z_0)\subset Z_0$, and set $C_0:= C|_{Z_0}$. If the range of $C$ is
closed and $\ker C=\ker C_0$ (i.e.\ $\ker C\subset Z_0$), then the
range of $C_0$ is also closed.
\end{Lemma}
\Proof If the conditions of the lemma are satisfied then, by
(\ref{eq_closed}), there exists $c>0$ such that $\|Cz\|\geq
c\inf_{y\in \ker C}\|z-y\|$, $z\in Z$. But, since $Z_0\subset Z$ and
$\ker C=\ker C_0$, this implies that $\|C_0z\|\geq c\inf_{y\in \ker
C_0}\|z-y\|$, $z\in Z_0$, so that the range of $C_0$ is closed.
\Proofend

\begin{Corollary}\label{cor_semi}
If $A\in S(Y)\cap L_0(Y)$ and $A_0$ is semi-Fredholm  with
$\alpha(A_0)< \infty$, then $A$ is semi-Fredholm and $\ker A=\ker
A_0$.
\end{Corollary}
\Proof If the conditions of the lemma are satisfied then, from
standard results on Fredholm operators (e.g.\ \cite{Jorgens}), we
have that $A_0^*$ and $A_0^{**}$ are also semi-Fredholm, and
$\alpha(A_0)=\beta(A_0^*)=\alpha(A_0^{**})$. Applying Lemma
\ref{lem_hatA}, it follows that
$\alpha(A_0)=\alpha(A)=\alpha(A_0^{**})$. Further, since
$\alpha(A_0)$ is finite and $\ker A_0\subset \ker A$, $\ker \check
A\subset \ker A_0^{**}$, it follows that $\ker A=\ker A_0$ and that
$\ker \check A = \ker A_0^{**}$. Applying Lemma \ref{lem_bs}, since
the range of $A_0^{**}$ is closed it follows that the range of
$\check A$ is closed and so $A(Y)$ is closed, and $A$ is
semi-Fredholm. \Proofend

We will prove the converse result only in the case when $U$ has a
predual space\index{predual space} $U^\pa$\index{$U^\pa$} and $A$
has a {\em preadjoint}\index{preadjoint
operator}\index{operator!preadjont} $A^\pa$\index{$A^\pa$,
preadjoint of $A$} on $Y^\pa := \ell^1(\Z^N,U^\pa)$. Recall that,
given a Banach space $X$, the Banach space $X^\pa$ is said to be a
{\em predual space} of $X$ if $X$ is isometrically isomorphic to
$(X^\pa)^*$. If $X$ has a predual space $X^\pa$ then $B^\pa\in
L(X^\pa)$ is said to be a preadjoint of $B\in L(X)$ if
$(B^\pa)^*=B$. If the Banach space $U$ has a predual space $U^\pa$
then a predual space of $Y= \ell^\infty(\Z^N,U)$ is $Y^\pa =
\ell^1(\Z^N,U^\pa)$. It is well-known that, if $X$ is a Banach space
which has a predual $X^\pa$ and $B\in L(X)$, then the following
statements are equivalent:
\begin{itemize}
\item[(i)] $B$ has a preadjoint $B^\pa\in L(X^\pa)$.\\[-3mm]
\item[(ii)] The adjoint $B^*$ maps $X^\pa$, understood as a subspace of
its second dual $(X^\pa)^{**}=X^*$, into itself.\\[-3mm]
\item[(iii)] $B$ is continuous in the weak$*$ topology\index{topology!weak$*$} on $X$.
\end{itemize}

Recalling the isometry $J^\pa:Y^\pa\to Y_1$ introduced above, let
$\check Y^\pa = J^\pa(Y^\pa)\subset Y_1=Y_0^*$, so that $\check
Y^\pa$ is isometrically isomorphic to $Y^\pa$. For $A^\pa\in
L(Y^\pa)$ let $\check A^\pa\in L(\check Y^\pa)$ be defined by
$\check A^\pa = J^\pa A^\pa (J^\pa)^{-1}$.

\begin{Lemma} \label{lem_firstembed}
If $A\in L_0(Y)$, $U$ has a predual $U^\pa$ and $A$ a preadjoint
$A^\pa\in L(Y^\pa)$, then $A_0^*(\check Y^\pa)\subset \check Y^\pa$
and $\check A^\pa = A_0^*|_{\check Y^\pa}$, so that $\ker \check
A^\pa \subset \ker A_0^*$.
\end{Lemma}
\Proof For $x\in \check Y^\pa$ and $y\in Y_0$, where $z :=
(J^\pa)^{-1}x\in Y^\pa$, using (\ref{eq_blff}),
$$
\check A^\pa x(y) = (y, \check A^\pa x)= y(A^\pa z)=Ay(z)=(Ay,x).
$$
Also,
$$
A_0^* x(y) = x(A_0y)=(A_0y,x)=(Ay,x).
$$
\Proofend

Let $J_1:Y_1\to Y^*$ be defined by
$$
J_1x(y):=Jy(x)=(y,x), \quad x\in Y_1,\; y\in Y.
$$

It is easy to check that $J_1$ is also an isometry. Let $\check
Y_1:= J_1(Y_1)\subset Y^*$, which is isometrically isomorphic to
$Y_1$. For $A_1\in L(Y_1)$ let $\check A_1 \in L(\check Y_1)$ be
defined by $\check A_1:= J_1 A_1 J_1^{-1}$.

\begin{Lemma} \label{lem_2ndembed}
If $A\in S(Y)\cap L_0(Y)$, then $A^*(\check Y_1)\subset \check Y_1$
and $\check A_0^* = A^*|_{\check Y_1}$, so that $\ker \check A_0^*
\subset \ker A^*$.
\end{Lemma}
\Proof For $x\in \check Y_1$, $y\in Y$, where $z=J_1^{-1}x\in Y_1$,
$$
\check A_0^*x(y)= J_1(A_0^*z)(y) = (y,A_0^*z)=(Ay,z),
$$
by Lemma \ref{lem_trans}. Also,
$$
A^*x(y) = x(Ay)=J_1z(Ay)=(Ay,z).
$$
\Proofend

We note that if the conditions of Lemmas \ref{lem_firstembed} and
\ref{lem_2ndembed} are satisfied, then
\begin{equation} \label{eq_ineq}
\alpha(A^\pa)\leq \alpha(A_0^*) \leq \alpha(A^*).
\end{equation}

\begin{Theorem} \label{Prop_Fred_equiv}
Suppose that $A\in S(Y)\cap L_0(Y)$, $U$ has a predual $U^\pa$ and
$A$ has a preadjoint $A^\pa\in L(Y^\pa)$. Then $A$ is Fredholm if
and only if $A_0$ is Fredholm and, if they are both Fredholm, then
$\alpha(A_0)=\alpha(A)$, $\beta(A_0)=\beta(A)$, and $\ker A = \ker
A_0$.
\end{Theorem}
\Proof Suppose first that $A_0$ is Fredholm . Then, by Corollary
\ref{cor_semi}, $A$ is semi-Fredholm and $\ker A=\ker A_0$. This
implies that $A^\pa$ and $A^*$ are also semi-Fredholm, and so, and
using (\ref{eq_ineq}),
$$
\beta(A) = \alpha(A^\pa) \leq \alpha(A_0^*)=\beta(A_0),
$$
so that $A$ is Fredholm. Moreover,
$$
\beta(A) = \alpha(A^*) \geq \alpha(A_0^*)=\beta(A_0)
$$
so $\beta(A)=\beta(A_0)$.

Conversely, if $A$ is Fredholm then so are $A^\pa$ and $A^*$ and
$\alpha(A^\pa)=\beta(A)=\alpha(A^*)$. Thus, by (\ref{eq_ineq}),
$\alpha(\check A_0^*)=\alpha(A_0^*)=\alpha(A^*)$ is finite and so it
follows from Lemma \ref{lem_2ndembed} that $\ker \check A_0^* =\ker
A^*$. Applying Lemma \ref{lem_bs} we see that the range of $\check
A_0^*$ is closed, so that the range of $A_0^*$ is closed and $A_0^*$
is semi-Fredholm. Thus $A_0$ is also semi-Fredholm , with
$\beta(A_0)=\alpha(A_0^*) =\alpha(A^*)< \infty$. But also
$\alpha(A_0)\leq \alpha(A)$ is finite, so $A_0$ is Fredholm.
\Proofend

Note that the above theorem and its proof simplifies greatly if the
Banach space $U$ is reflexive, in particular if $U$ is finite
dimensional. For then we can choose $U^\pa=U^*$ so that $Y^\pa=Y_1$
and $Y_0^{**} = Y$. Note also that, if the conditions of the above
theorem hold, in particular if $A$ has a preadjoint, then the above
theorem implies that $A$ is invertible if and only if $A_0$ is
invertible. But even without existence of a preadjoint, we can prove
this result in some cases; an observation which will be useful to us
later.

\begin{Lemma} \label{lem_half_invert}
If $A\in S(Y)\cap L_0(Y)$ and $A$ is invertible, then $A_0$ is
invertible.
\end{Lemma}
\Proof If $A$ is invertible then $A_0$ is injective and it follows
from Lemma \ref{lem_bs} that the range of $A_0$ is closed. Further,
since $A$ is the transpose of $A_0^*$ with respect to the dual
system $(Y,Y_1)$ it follows (see e.g.\ \cite{Jorgens}) that
$0=\beta(A)\geq \alpha(A_0^*)=\beta(A_0)$. Thus $A_0$ is surjective.
\Proofend

\begin{Lemma} \label{lem_other_half}
Suppose that $A\in L(Y,\P)$ or that $A=I+K$ with $K\in S(Y)\cap
M(Y)\cap L_0(Y)$, and suppose that $A_0$ is invertible. Then $A$ is
invertible.
\end{Lemma}
\Proof If the conditions of the lemma apply then, by Corollary
\ref{cor_semi}, $A$ is injective. In the case that $A\in L(Y,\P)$
then $A_0\in L(Y_0,\P)$ by Lemma \ref{lem3.13}, and since
$L(Y_0,\P)$ is inverse closed (Theorem 1.1.9 of \cite{RaRoSiBook}),
we have that $A_0^{-1}\in S(Y_0)$. This holds also by a modification
of the proof of Lemma \ref{lem_inv_closed} in the case that $A=I+K$
with $K\in S(Y)\cap M(Y)\cap L_0(Y)$. For if $(x_n)\subset Y_0$,
$x\in Y_0$, and $x_n\sto x$ then, defining $y_n : = A_0^{-1} x_n$,
\begin{equation} \label{eq_yncopy}
y_n + K y_n = x_n
\end{equation}
holds, and since $K\in M(Y)$ there exists a subsequence $(y_{n_m})$
and $y\in Y$ such that $x_{n_m} - K y_{n_m} \sto y$. From
(\ref{eq_yncopy}) it follows that $y_{n_m}\sto y$. Since $K\in
S(Y)$, it follows that $x_{n_m} - K y_{n_m} \sto x - Ky$. Thus $y =
x - Ky$, i.e. $Ay = x$. Note that, by injectivity of $A$, there is
only one $y\in Y$ with $Ay=x$ and that is $y=A_0^{-1}x\in Y_0$. We
have shown that $y_n = A_0^{-1}x_n$ has a subsequence strictly
converging to $y=A_0^{-1}x$. By the same argument, every subsequence
of $y_n$ has a subsequence strictly converging to $y$. Thus
$A_0^{-1}x_n \sto A_0^{-1} x$. So $A_0^{-1}\in S(Y_0)$.

Let $B\in S(Y)$ be the unique extension of $A_0^{-1}$ from $Y_0$ to
$Y$, which exists by Lemma \ref{lem3.13}. Then, for every $x\in Y$,
\[
BAx\ =\ s\!-\!\!\!\lim_{n\to\infty} BAP_n x \ =\ s\!-\!\lim
A_0^{-1}A_0P_n x\ =\ x
\]
and, similarly, $ABx=x$. So $A$ is invertible. \Proofend

\begin{Corollary} \label{cor_AA0}
For $A\in L(Y,\P)$ it holds that $A$ is invertible iff $A_0$ is
invertible. When both are invertible, then
$(A_0)^{-1}=(A^{-1})|_{Y_0}$.
\end{Corollary}
\Proof The first sentence follows immediately from the previous two
lemmas, and the equality concerning the two inverses is obvious if
both $A$ and $A_0$ are invertible. \Proofend

\section{Band-Dominated Operators} \label{sec_BDO}
Suppose $Y$ is one of the spaces $Y^p$ with $p\in\{0\}\cup
[1,\infty]$ and $\V$ is given by (\ref{eq_shifts}). For $m\in\Z^N$
let $E_m:U\to Y$\index{$E_m$} and $R_m:Y\to U$\index{$R_m$} be
extension\index{extension operator}\index{operator!extension} and
restriction operators\index{restriction
operator}\index{operator!restriction}, defined by
$E_my=(...,0,y,0,...)$, for $y\in U$, with the $y$ standing at the
$m$th place in the sequence, and by $R_mx=x(m)$, for
$x=(x(m))_{m\in\Z^N}\in Y$. To every operator $A\in L(Y)$ we can
associate a matrix $[A]=[a_{ij}]_{i,j\in\Z^N}$\index{$A_n$@$[A]$,
matrix representation of $A$} with entries $a_{ij}= R_iAE_j\in
L(U)$. For $x\in\tilde Y=\cup_{n\in\N} P_n(Y)$, we have that the
$i$th entry of $Ax$ is
\begin{equation} \label{eqn_matrix}
(Ax)(i)\ =\ \sum_{j\in\Z^N} a_{ij}\,x(j),\qquad i\in\Z^N.
\end{equation}
Since $Y_0$ is the norm closure of $\tilde Y$, clearly the entries
of $[A]$ determine $A|_{Y_0}$ uniquely, and so determine $A$ if
$Y=Y_0$, i.e. if $p\in\{0\}\cup [1,\infty)$. In the case $p=\infty$,
$[A]$ determines $A$ if $A\in S(Y)$, for then $Ax$ is given by
$\lim_{n\to\infty}AP_nx$, the limit taken in the $s$ sense.

Let $BO(Y)$\index{$BO(Y)$}\index{band operator}\index{operator!band}
denote the set of all operators $A\in L(Y,\P)$ such that $[A]$ is a
band matrix\index{band matrix}, that is, for some $w\in \N_0$ called
the {\em band-width}\index{band width} of $A$, $a_{ij}=0$ if
$|i-j|>w$.

Recall that, for $b=(b(m))_{m\in\Z^N}\in \ell^\infty(\Z^N,L(U))$,
the multiplication operator $M_b\in L(Y)$ is defined by
\refeq{eq_multop}. In terms of multiplication operators, an
alternative characterisation of $BO(Y)$ is the following
\cite{LiBook}: that $A\in L(Y)$ is a band operator of band-width $w$
iff
\begin{equation} \label{eq_BO}
 A = \sum_{|k|\le w} M_{b_k}V_k,
\end{equation}
with some $b_k\in \ell^\infty(\Z^N,L(U))$ for every $|k|\le w$. For
$A$ given by \refeq{eq_BO} the matrix representation $[A]$ has
$a_{ij}=b_{i-j}(i)$ for $|i-j|\le w$ and $a_{ij}=0$, otherwise.
Since $\V\subset L(Y,\P)$ and, clearly, $M_b\in L(Y,\P)$ for
$b\in\ell^\infty(\Z^N,L(U))$, every $A$ of the form (\ref{eq_BO}) is
in $L(Y,\P)$.

The linear space $BO(Y)$ is an algebra but is not closed with
respect to the norm in $L(Y)$. By taking the closure of $BO(Y)$ in
the operator norm of $L(Y)$ we obtain the Banach algebra $BDO(Y)$.
We refer to the elements of $BDO(Y)$\index{$BDO(Y)$} as {\em
band-dominated operators}\index{band-dominated
operator}\index{operator!band-dominated}. It should be noted that
$BDO(Y^p)$ depends on the exponent $p$ of the underlying space,
while $BO(Y^p)$ does not. Since $L(Y,\P)$ is closed in $L(Y)$, and
$BO(Y)\subset L(Y,\P)$, it follows that $BDO(Y)\subset L(Y,\P)$.

As a consequence of (\ref{eq_limops}), an operator $A\in BO(Y)$,
which has the form (\ref{eq_BO}), is a rich operator iff the
multiplication operators $M_{b_k}$ are rich for all $|k|\le w$,
i.e.\ iff (see Example \ref{ex_almostperiodic}) the set $\{V_i
b_k:i\in\Z^N\}$ is relatively sequentially compact in the strict
topology on $Y^\infty(L(U))$ for every $k$. Further, if $A\in
BDO(Y)$, in which case $A_n\toto A$ for some $(A_n)\subset BO(Y)$,
$A$ is rich if each $A_n$ is rich.

\begin{Lemma} \cite[Corollary 2.1.17]{RaRoSiBook} \label{lem_oprich}
If $U$ is finite-dimensional then every band-dominated operator is
rich.
\end{Lemma}

For band-dominated operators the notions of invertibility at
infinity (Definition \ref{inv_inf}) and $\P-$Fredholmness (Remark
\ref{rem_PFred}) coincide.

\begin{Lemma} \label{lem_invinf_PFred}
If $A\in BDO(Y)$, then the following statements are equivalent.
\begin{itemize}
\item[(i)] $A$ is invertible at infinity.
\item[(ii)] There exist $B\in BDO(Y)$ and $T_1,T_2\in K(Y,\P)$ such
that \refeq{eq_inv_inf} holds.
\item[(iii)] $A$ is $\P-$Fredholm.
\end{itemize}
\end{Lemma}
\Proof For $Y^p$ with $1\le p\le\infty$ this is precisely
Proposition 2.10 of \cite{LiBook} (the key idea, written down for
$Y^2$, is from \cite[Proposition 2.6]{RaRoSi2001}). The proof from
\cite{LiBook} literally transfers to $Y^0$. \Proofend

\noindent For band-dominated operators we have also the following
result; in fact the equivalence of (i), (ii) and (iii) in the next
lemma is true even for arbitrary $K\in L(Y)$.

\begin{Lemma} \label{lem_coll_comp_entries}
If $K\in BDO(Y)$, with $[K]=[\kappa_{ij}]_{i,j\in\Z^N}$ the matrix
representation of $K$, then the following statements are equivalent.
\begin{itemize}
\item[(i)] $\{V_{-k}KV_{k}:k\in\Z^N\}\cup\opsp(K)$ is uniformly Montel on
$(Y,s)$.
\item[(ii)] $\{V_{-k}KV_{k}:k\in\Z^N\}$ is uniformly Montel on $(Y,s)$.
\item[(iii)] $(V_{-k}KV_{k})_{k\in\Z^N}$ is asymptotically Montel on $(Y,s)$
and $K\in M(Y)$.
\item[(iv)] The set $\{\kappa_{ij}:i,j\in\Z^N\}\subset L(U)$ is collectively
compact.
\item[(v)] The set $\{\kappa_{ij}:i,j\in\Z^N,\,i-j=d\}\subset L(U)$ is
collectively
compact for every $d\in \Z^N$.
\end{itemize}
If $U$ is finite-dimensional, then (i)--(v) are also equivalent to:
\begin{itemize}
\item[(vi)] The set $\{\kappa_{ij}:i,j\in\Z^N,\,i-j=d\}\subset L(U)$ is bounded
for every $d\in \Z^N$.
\end{itemize}
\end{Lemma}
\Proof It is clear from the definitions that (i)$\Rightarrow$(ii),
(ii)$\Rightarrow$(iii) and that (iv)$\Rightarrow$(v). By Lemma
\ref{charact_uniformly_Montel}, (ii) implies (i).

Suppose now that (iii) holds and that
$h=(h(n))_{n=1}^\infty\subset\Z^N$ and that $(x_n)\subset Y$ is
bounded. If $h$ does not have a subsequence that tends to infinity,
then $h$ is bounded, and hence it has a subsequence that is
constant. In the case that $h$ has a subsequence that tends to
infinity, $(V_{-h(n)}KV_{h(n)}x_n)$ has a strictly convergent
subsequence since $(V_{-k}KV_{k})$ is asymptotically Montel. In the
case that $h$ has a constant subsequence, $(V_{-h(n)}KV_{h(n)}x_n)$
has a strictly convergent subsequence since $K\in M(Y)$. In either
case, we have shown that (ii) holds.

Next suppose that (ii) holds and that $i=(i(n))_{n=1}^\infty\subset
\Z^N$, $j=(j(n))_{n=1}^\infty\subset \Z^N$, and that $(u_n)\subset
U$ is bounded. For $n\in\Z^N$ define  $(x_n)\in Y$ by setting the
$i(n)-j(n)$ entry of $x_n$ equal to $u_n$ and setting the other
entries to zero. Then $(x_n)$ is bounded and the zeroth entry of
$(V_{-i(n)}KV_{i(n)}x_n)$ is $\kappa_{i(n),j(n)}u_n$. Since
$\{V_{-k}KV_{k}:k\in\Z^N\}$ is uniformly Montel,
$(V_{-i(n)}KV_{i(n)}x_n)(0)=\kappa_{i(n),j(n)}u_n$ has a convergent
subsequence.  Since $i$, $j$, and $(u_n)$ were arbitrary sequences,
we have shown that (iv) holds.

Finally, suppose that (v) holds. Then the set
$\{\kappa_{ij}:i,j\in\Z^N,\,|i-j|\le w\}$ is collectively compact
for every $w\in\N$. For every $M\in\N$, every
$h=(h(n))_{n=1}^\infty\subset \Z^N$, and every bounded sequence
$(x_n)\subset Y$, we have that the $i$th component of
$(V_{-h(n)}KV_{h(n)}P_Mx_n)$ is
$$
\sum_{|j|\le M} \kappa_{i+h(n),j+h(n)}x_n(j).
$$
Since $\{\kappa_{ij}:i,j\in\Z^N,\,|i-j|\le w\}$ is collectively
compact for each $w$, it follows that the $i$th component of
$(V_{-h(n)}KV_{h(n)}P_Mx_n)$ has a convergent subsequence for
every $M\in\Z$. Thus, by a diagonal argument,
$(V_{-h(n)}KV_{h(n)}P_Mx_n)$ has a strictly convergent
subsequence, for every $M\in\N$. Again by a diagonal argument, we
can find subsequences of $h$ and $(x_n)$, which we will denote
again by $h$ and $(x_n)$, such that $(V_{-h(n)}KV_{h(n)}P_nx_n)$
converges strictly to some $x\in Y$, so that
$P_mV_{-h(n)}KV_{h(n)}P_nx_n\to P_mx$ as $n\to\infty$, for each
$m$. Now \cite{RaRoSiBook}, since $K$ is band-dominated, it holds
for every $m\in\N$ that $P_mV_{-k}KV_{k}Q_n\toto0$ as
$n\to\infty$, uniformly in $k\in\Z^N$. Thus
$P_mV_{-h(n)}KV_{h(n)}x_n\to P_mx$, for every $m\in\N$, so that
$V_{-h(n)}KV_{h(n)}x_n\sto x$. We have shown that (ii) holds.

The equivalence of (vi) and (v) under the condition $\dim U<\infty$
is obvious since, in that case, a subset of $U$ is relatively
compact iff it is bounded. \Proofend

For brevity, and because we will frequently refer to this class of
operators in what follows, let us denote the set of all operators
$K\in BDO(Y)$ which are subject to the (equivalent) properties
(i)--(v) of Lemma \ref{lem_coll_comp_entries} by
$\UM(Y)$\index{$UM(Y)$@$\UM(Y)$}.

\begin{Lemma} \label{lem_7.23isBA}
The following statements hold.
\begin{itemize}
\item[(a)] The set $\UM(Y)$ is a Banach subspace of $BDO(Y)\cap M(Y)$.
\item[(b)] In particular, $\UM(Y)=BDO(Y)$ if $U$ is finite-dimensional.
\item[(c)] If $K\in\UM(Y)$ and $A\in BDO(Y)$, then $KA\in\UM(Y)$.
\item[(d)] If $M_b$ is rich and $K\in\UM(Y)$, then $M_bK\in\UM(Y)$.
\end{itemize}
\end{Lemma}
\Proof (a): By its definition, we have that $\UM(Y)\subset BDO(Y)$,
and from Lemma \ref{lem_coll_comp_entries} (iii) we get that
$\UM(Y)\subset M(Y)$. For the rest of this proof, we will use
property (ii) from Lemma \ref{lem_coll_comp_entries} to characterise
the set $\UM(Y)$.

If $S,T\in\UM(Y)$ and $\lambda\in\C$, then clearly $\lambda
S+T\in\UM(Y)$ since
\[
\{V_{-k}(\lambda S+T) V_k:k\in\Z^N\}\ \subset\ \lambda \{V_{-k}S
V_k:k\in\Z^N\}\ +\ \{V_{-k}T V_k:k\in\Z^N\}
\]
is uniformly Montel.

If $T_1,T_2,...\in\UM(Y)$ are such that $T_n\toto T$, then also
$T\in\UM(Y)$. To see this, take a sequence
$k=(k(1),k(2),...)\subset\Z^N$ and a sequence $(x_1,x_2,...)\subset
Y$ with $\mu:=\sup\|x_\ell\|<\infty$. By a simple diagonal argument,
we can pick a strictly monotonously increasing sequence
$s=(s(1),s(2),...)\subset\N$ such that
\[
V_{-k(s(\ell))}\,T_n\,V_{k(s(\ell))}\ x_{s(\ell)}
\]
converges strictly as $\ell\to\infty$ for every $n\in\N$. Let us
denote the strict limit by $y_n$, respectively. From
\[
\|y_{n_1}-y_{n_2}\|\ \le\
\sup_\ell\|V_{-k(s(\ell))}(T_{n_1}-T_{n_2})V_{k(s(\ell))}x_{s(\ell)}\|\
\le\ \|T_{n_1}-T_{n_2}\|\cdot\mu
\]
we see that $(y_n)$ is a Cauchy sequence in $Y$ and therefore
converges, to $y\in Y$, say. But then
$V_{-k(s(\ell))}\,T\,V_{k(s(\ell))}\ x_{s(\ell)}\sto y$ as
$\ell\to\infty$. Indeed, for all $M,n\in\N$,
\begin{eqnarray*}
\lefteqn{\|P_M(V_{-k(s(\ell))}\,T\,V_{k(s(\ell))}\,x_{s(\ell)}\ -\ y)\|}\\
&\le&\|P_M(V_{-k(s(\ell))}\,T_n\,V_{k(s(\ell))}\,x_{s(\ell)}\ -\
y_n)\|\ +\
\|P_M(y_n-y)\|\\
&&\hspace{47mm} +\
\|P_M(V_{-k(s(\ell))}\,(T-T_n)\,V_{k(s(\ell))}\,x_{s(\ell)})\|\\
&\le&\|P_M(V_{-k(s(\ell))}\,T_n\,V_{k(s(\ell))}\,x_{s(\ell)}\ -\
y_n)\|\ +\ \|y_n-y\|\ +\ \|T-T_n\|\cdot\mu
\end{eqnarray*}
holds. But, for every choice of $M,n\in\N$, the first term goes to
zero as $\ell\to\infty$, and the second and third term can be made
as small as desired by choosing $n$ sufficiently large.

(b): If $K\in BDO(Y)$ then property (vi) of Lemma
\ref{lem_coll_comp_entries} is automatically the case. Since this is
equivalent to properties (i)--(v) of the same lemma if $U$ is
finite-dimensional, we get that $K\in\UM(Y)$ then.

(c), (d): Let $K\in\UM(Y)$, $A\in BDO(Y)$ and $b\in Y^\infty(L(U))$
such that $M_b$ is rich. Take a sequence
$k=(k(1),k(2),...)\subset\Z^N$ and a bounded sequence
$(x_1,x_2,...)\subset Y$. Now, for every $\ell\in\N$, put
$y_\ell:=V_{-k(\ell)}AV_{k(\ell)}\,x_\ell$. Since $(y_\ell)$ is
bounded, $\{V_{-m}KV_m:m\in\Z^N\}$ is uniformly Montel and
$\{V_{-m}b:m\in\Z^N\}$ is relatively sequentially compact in the
strict topology on $Y^\infty(L(U))$ (since $M_b$ is rich, see
Example \ref{ex_almostperiodic}), we can pick a strictly
monotonously increasing sequence $s=(s(1),s(2),...)\subset\N$ such
that both $V_{-k(s(\ell))}KV_{k(s(\ell))}y_{s(\ell)}$ and
$V_{-k(s(\ell))}b$ converge strictly as $\ell\to\infty$. But then
$V_{-k(s(\ell))}\,(M_bKA)\,V_{k(s(\ell))}\ x_{s(\ell)}$ converges
strictly as $\ell\to\infty$ since, for every $m\in\N$,
\begin{eqnarray*}
\lefteqn{P_m\ V_{-k(s(\ell))}\,(M_bKA)\,V_{k(s(\ell))}\ x_{s(\ell)}}\\
&=& P_m\ (V_{-k(s(\ell))}M_bV_{k(s(\ell))})\ (V_{-k(s(\ell))}KV_{k(s(\ell))})\
(V_{-k(s(\ell))}AV_{k(s(\ell))}\ x_{s(\ell)})\\
&=& M_{P_m V_{-k(s(\ell))}b}\ P_m\,(V_{-k(s(\ell))}KV_{k(s(\ell))}\
y_{s(\ell)})
\end{eqnarray*}
converges in norm as $\ell\to\infty$, i.e. $M_bKA\in\UM(Y)$.
\Proofend


The purpose of the following two lemmas is to prove that, for every
$A\in L(Y^1(U))$, the operator spectra $\opsp(A)\subset L(Y^1(U))$
and $\opsp(A^*)\subset L((Y^1(U))^*)=L(Y^\infty(U^*))$ correspond
elementwise in terms of adjoints.

\begin{Lemma} \label{lem_lem2}
If $A\in L(Y^1(U))$, then
\[
\opsp(A^*)\ =\ \{B^* : B\in\opsp(A)\}.
\]
\end{Lemma}
\Proof It is a standard result that $B=A_h\in\opsp(A)$ implies
$B^*=(A_h)^*=(A^*)_h\in\opsp(A^*)$ (see, e.g. \cite[Proposition 3.4
e]{LiBook}).

For the reverse implication, suppose $C\in\opsp(A^*)\subset
L(Y^\infty(U^*))$. Then
\[
(V_{-h(m)}AV_{h(m)})^*=V_{-h(m)}A^*V_{h(m)}\pto C
\]
as $m\to\infty$ for some sequence $h(1),h(2),...\to\infty$ in
$\Z^N$. We will show in Lemma \ref{lem_lem1} that then $C=B^*$ and
$V_{-h_m}AV_{h_m}\pto B$, i.e. $B\in\opsp(A)$. \Proofend

\begin{Lemma} \label{lem_lem1}
The set of operators in $L(Y^\infty(U^*))$ that possess a preadjoint
in $L(Y^1(U))$ is sequentially closed under $\P-$convergence; that
is, if $A_1,A_2,...\in L(Y^1(U))$ and $A_m^*\pto C$ on
$Y^\infty(U^*)$, then there is a $B\in L(Y^1(U))$ such that $C=B^*$;
moreover, $A_m\pto B$ on $Y^1(U)$.
\end{Lemma}
\Proof From $A_m^*\pto C$ in $L(Y^\infty(U^*))$ and Lemma
\ref{lem_Pconv} we get that there is a $M>0$ such that
\begin{equation} \label{M}
\|A_m\|=\|A_m^*\|\le M,\qquad m\in\N.
\end{equation}
Moreover, for every $k\in\N$, it holds that
\begin{equation} \label{pconv}
P_k(A_m^*-C)\toto 0\qquad\textrm{ as }\qquad m\to\infty.
\end{equation}
So we get that $(P_k A_m^*)_{m=1}^\infty$ is a Cauchy sequence in
$L(Y^\infty(U^*))$ and therefore $(A_m P_k)_{m=1}^\infty$ is one in
$L(Y^1(U))$, for every fixed $k\in\N$. Denote the norm-limit of the
latter sequence by $B_k\in L(Y^1(U))$. As a consequence of (\ref{M})
we get that
\begin{equation} \label{C_knorm}
\|B_k\|=\|\lim_{m\to\infty} A_mP_k\|\le \sup_m \|A_mP_k\|\le
M,\qquad k\in\N.
\end{equation}
From $A_mP_k\toto B_k$ we get that $B_kP_k=B_k$ and, even more than
this, that
\begin{equation} \label{C_k}
B_r P_k\ =\ \lim_{m\to\infty} A_mP_rP_k\ =\ \lim_{m\to\infty}
A_mP_k\ =\ B_k,\qquad r\ge k.
\end{equation}
We will now show that the sequence $B_1,B_2,...$ strongly converges
in $Y^1(U)$. Therefore, take an arbitrary $x\in Y^1(U)$ and let us
verify that $(B_m x)$ is a Cauchy sequence in $Y^1(U)$. So choose
some $\varepsilon>0$. Since $Q_mx\to 0$ on $Y^1(U)$, there is an
$N\in\N$ such that
\begin{equation} \label{eps}
\|Q_N x\|<\frac\varepsilon{2M}.
\end{equation}
Now, for all $k,m\ge N$, the following holds
\begin{eqnarray*}
\|B_kx-B_mx\|&\le&\|(B_k-B_m)P_Nx\|+\|(B_k-B_m)Q_Nx\|\\
&\le&\|(B_kP_N-B_mP_N)x\|+\|B_k-B_m\|\cdot\|Q_Nx\|\\
&\le&\|(B_N-B_N)x\|+(\|B_k\|+\|B_m\|)\cdot\|Q_Nx\|\ <\ \varepsilon
\end{eqnarray*}
by (\ref{C_k}), (\ref{C_knorm}) and (\ref{eps}). Consequently, $(B_m
x)$ is a Cauchy sequence in $Y^1(U)$. Let us denote its limit in
$Y^1(U)$ by $Bx$, thereby defining an operator $B\in L(Y^1(U))$.
Passing to the strong limit as $r\to\infty$ in (\ref{C_k}), we get
\begin{equation} \label{C}
B P_k=B_k,\qquad k\in\N.
\end{equation}
Summing up, we have $A_mP_k\toto B_k=B P_k$, and hence
$(A_m-B)P_k\toto 0$ as $m\to\infty$, for all $k\in\N$. Passing to
adjoints in the latter gives $P_k(A_m^*-B^*)\toto 0$ in
$L(Y^\infty(U^*))$ as $m\to\infty$. If we subtract this from
(\ref{pconv}) we get $P_k(B^*-C)=0$ for all $k\in\N$, and
consequently $C=B^*$, by Lemma 1.30 a) in \cite{LiBook}. From
$A_m^*\pto C=B^*$ we then conclude
\[
\|(A_m-B)P_k\|=\|P_k(A_m^*-B^*)\|\to 0 \qquad\textrm{ as }\qquad
m\to\infty
\]
and
\[
\|P_k(A_m-B)\|=\|(A_m^*-B^*)P_k\|\to 0 \qquad\textrm{ as }\qquad
m\to\infty
\]
for every $k\in\N$, which, together with (\ref{M}) and again Lemma
\ref{lem_Pconv}, proves $A_m\pto B$. \Proofend

Our next statement is similar to Lemma \ref{lem_lem2}, but with
restriction from $Y$ to $Y_0$ instead of passing to the adjoint
operator.

\begin{Lemma} \label{lem_limopA0}
If $A\in L(Y,\P)$, then the limit operators of the restriction
$A_0:=A|_{Y_0}$ are the restrictions of the limit operators of $A$;
precisely,
\begin{equation} \label{eq_limopsY0}
\opsp(A_0)\ =\ \{ B|_{Y_0}:B\in\opsp(A)\}.
\end{equation}
In particular, the invertibility of all limit operators of $A_0$ in
$Y_0$ with uniform boundedness of their inverses is equivalent to
the same property for the limit operators of $A$ in $Y$.
\end{Lemma}
\Proof The proof of (\ref{eq_limopsY0}) consists of two
observations. The first one is that
$(V_{-\alpha}AV_\alpha)|_{Y_0}=V_{-\alpha}A_0V_\alpha$ for all
$\alpha\in\Z^N$, and the second one is that $A_m|_{Y_0}\pto A_0$ on
$Y_0$ iff $A_m\pto A$ on $Y$, for all $A_1,A_2,...\in L(Y,\P)$ since
\[
\|P_k(A_m|_{Y_0}-A_0)\|=\|\big(\,P_k(A_m-A)\,\big)|_{Y_0}\|=\|P_k(A_m-A)\|
\]
and its symmetric counterpart hold for all $k\in\N$ by the norm
equality in Lemma \ref{lem3.13}. The proof of the second sentence of
the lemma now follows from (\ref{eq_limopsY0}), Corollary
\ref{cor_AA0} and the norm equality in Lemma \ref{lem3.13} again.
\Proofend


For the choice of $Y$, $\P$ and $\V$ we are making in this section,
the operator spectrum of a rich band-dominated operator $A$ contains
enough information to characterise the invertibility at infinity of
$A$, which is the content of (iii) in the next theorem. In this
theorem, for an operator $A\in L(Y)$, we denote the {\em
spectrum}\index{spectrum} of $A$, i.e. the set of all $\lambda\in\C$
for which $\lambda I-A$ is not invertible, by
$\sp(A)$\index{$\psp$@$\sp(A)$}. We denote by
$\spess(A)$\index{$\pspe$@$\spess(A)$} the essential
spectrum\index{essential spectrum}\index{spectrum!essential} of $A$,
i.e.\ the set of $\lambda$ for which $\lambda I-A$ is not Fredholm .
We also define, for $\eps>0$, the {\em
$\eps$-pseudospectrum}\index{pseudospectrum}\index{spectrum!pseudo-}
of $A$, $\speps(A)$\index{$\pspeps$@$\speps(A)$}, by
$$
\speps(A) := \left\{\lambda\in\C\ :\ \lambda I-A \mbox{ is not
invertible or }||(\lambda I-A)^{-1}||\ge \eps^{-1}\right\}.
$$

\begin{Theorem} \label{prop_rich_invinf}
{\bf a) } Let $A$ be a rich band-dominated operator on $Y=Y^p(U)$
with a Banach space $U$ and some $p\in\{0\}\cup[1,\infty]$. Then the
following statements hold.
\begin{itemize}
\item[(i)] If $A$ is Fredholm and $p\ne\infty$ then $A$ is invertible at
infinity;
\item[(ii)] If $A$ is invertible at infinity and either $U$ is
finite-dimensional or $A=C+K$ with $C\in BDO(Y)$ invertible and
$K\in M(Y)$ then $A$ is Fredholm;
\item[(iii)] $A$ is invertible at infinity if and only if all limit operators
of $A$ are invertible and their inverses are uniformly bounded;
\item[(iv)] The condition on uniform boundedness in (iii) is
redundant if $p\in\{0,1,\infty\}$;
\item[(v)] It holds that $\sp(B) \subset \sp(A)$ for all
$B\in\opsp(A)$, indeed $\sp(B)\subset\spess(A)$, for $p\ne\infty$;
\item[(vi)] It holds that $\speps(B) \subset  \speps(A)$, for all
$B\in\opsp(A)$ and $\eps>0$.
\end{itemize}
{\bf b) } In the case $p=\infty$ if, in addition, it holds that $U$
has a predual $U^\pa$, and $A$ has a preadjoint, $A^\pa\in
L(Y^\pa)$, where $Y^\pa=Y^1(U^\pa)$, then (i) and (v) also apply for
$p=\infty$; that is, $A$ being Fredholm implies $A$ being invertible
at infinity, so that $\sp(B)\subset\spess(A)\subset\sp(A)$ for all
$B\in\opsp(A)$.
\end{Theorem}

\begin{Remark}
This theorem makes several additions and simplifications to
previously known results: (i), and therefore (v), is probably new
for $p=1$, and so is statement b). (ii) is a slight extension of
Proposition 2.15 (together with the first column of Figure 4) of
\cite{LiBook}. (iii) does not assume the existence of a preadjoint
operator (unlike Theorem 1 in \cite{LiDiss} and \cite{LiBook}) if
$p=\infty$. (iv) is probably new for $p=0$. Also note that our Lemma
\ref{lem_lem2} fills a gap in the proof of \cite[Proposition 3.6
a)]{LiBook} that is used in the proof of \cite[Theorem
 3.109]{LiBook} to deal with the case $p=1$. (vi) was only known when
$p=2$ and $U$ is a Hilbert space. For this setting, it follows from
Theorem 6.3.8 (b) of \cite{RaRoSiBook} which, in fact, states the
stronger result that the closure of the union of all $\speps(B)$
with $B\in\opsp(A)$ is equal to the $\eps$-pseudospectrum of the
coset of $A$ modulo $K(Y,\P)$.
\end{Remark}

{\it Proof of Theorem \ref{prop_rich_invinf}.} {\bf b) } Suppose the
predual $U^\pa$ and preadjoint $A^\pa$ exist and that $A$ is
Fredholm on $Y=Y^\infty(U)$. By Theorem \ref{Prop_Fred_equiv} (note
that $BDO(Y)\subset L(Y,\P)\subset S(Y)\cap L_0(Y)$), we have that
$A_0:=A|_{Y_0}$ is Fredholm on $Y_0=Y^0(U)$. From Lemma \ref{lem_KK}
(i) (also see Remark \ref{rem_Fred_invinf}) we get that then $A_0$
is invertible at infinity on $Y_0$, which, by Lemma
\ref{lem_invinf_PFred} means that we have $A_0B_0=I+S_0$ and
$B_0A_0=I+T_0$ for some $B_0\in L(Y_0,\P)$ and $S_0,T_0\in
K(Y_0,\P)$. If we use Lemma \ref{lem3.13} to extend both sides of
these two equalities to operators on $Y$, then we get that $A$ is
invertible at infinity on $Y$.

{\bf a)} (iii) For $p\in\{0\}\cup (1,\infty)$ we refer the reader to
\cite[Theorem 2.2.1]{RaRoSiBook}, and for $p=1$ (and also
$p\in(1,\infty)$) to \cite[Theorem 1]{LiBook}. It remains to study
the case $p=\infty$. The `if' part of statement (iii) is Proposition
3.16 in \cite{LiBook} (which does not use the existence of a
preadjoint). For the `only-if' part of (iii) we replace Proposition
3.12 from \cite{LiBook} (which needs the preadjoint) by the
following argument. Suppose $A$ is invertible at infinity on
$Y=Y^\infty$. By Lemma \ref{lem_invinf_PFred}, there are $B\in
BDO(Y)\subset L(Y,\P)$ and $S,T\in K(Y,\P)$ with $AB=I+S$ and
$BA=I+T$. Restricting both sides in both equalities to $Y_0$ we get
that, by Lemma \ref{lem3.13}, $A_0:=A|_{Y_0}$ is invertible at
infinity on $Y_0$, which, by our result (iii) for $p=0$, implies
that all limit operators of $A_0$ are invertible on $Y_0$ and their
inverses are uniformly bounded. From Lemma \ref{lem_limopA0} we now
get that also all limit operators of $A$ are invertible on $Y$ and
their inverses are uniformly bounded.

Statement (iv) for $p\in\{1,\infty\}$ is \cite[Theorem
3.109]{LiBook}. Precisely, the part for $p=\infty$ follows
immediately from \cite[Proposition 3.108]{LiBook}, and the $p=1$
part is a consequence of this and Lemma \ref{lem_lem2}. Indeed, if
all $B\in\opsp(A)$ are invertible on $Y^1(U)$ then also all their
adjoints $C=B^*$ are invertible on $Y^\infty(U^*)$, which, by Lemma
\ref{lem_lem2}, are all elements of $\opsp(A^*)$. Since $A^*\in
BDO(Y^\infty(U^*))$ is rich as well, we know from the results about
$p=\infty$ that
\begin{eqnarray*}
\sup_{B\in\opsp(A)}\|B^{-1}\|&=&\sup_{B\in\opsp(A)}\|(B^{-1})^*\|\
=\ \sup_{B\in\opsp(A)}\|(B^*)^{-1}\|\\
&=&\sup_{C=B^*\in\opsp(A^*)}\|C^{-1}\|\ <\ \infty
\end{eqnarray*}
since $B\in\opsp(A)$ iff $C=B^*\in\opsp(A^*)$, by Lemma
\ref{lem_lem2}. The statement (iv) for $p=0$ follows immediately
from Lemma \ref{lem_limopA0} (applied to the extension of $A$) and
the result for $p=\infty$.

(i) For $p\in\{0\}\cup(1,\infty)$ this follows immediately from
Lemma \ref{lem_KK} (i) (it was already pointed out in Remark
\ref{rem_Fred_invinf}). So let $p=1$ and suppose $A$ is Fredholm on
$Y=Y^1(U)$. Then $A^*$ is Fredholm on $Y^\infty(U^*)$. From part b)
of this theorem it follows that $A^*$ is invertible at infinity on
$Y^\infty(U^*)$. From (iii) we get that all limit operators of $A^*$
are invertible on $Y^\infty(U^*)$. By Lemma \ref{lem_lem2} this
implies that all limit operators of $A$ are invertible on $Y^1(U)$,
which, by (iii) and (iv), shows that $A$ is invertible at infinity.

(ii) Suppose $A$ is invertible at infinity. If $\dim U<\infty$ then
Lemma \ref{lem_KK} (ii) (also see Remark \ref{rem_Fred_invinf})
implies that $A$ is Fredholm. Alternatively, suppose that $A=C+K$
with $C\in BDO(Y)$ invertible and $K\in M(Y)$. From $C\in
BDO(Y)\subset L(Y,\P)$ we get, by \cite[Theorem 1.1.9]{RaRoSiBook},
that $C^{-1}\in L(Y,\P)\subset S(Y)$. Moreover, $K=A-C\in
BDO(Y)\subset S(Y)$ implies that $K\in S(Y)\cap M(Y)$ so that $A$ is
subject to the constraints in Theorem \ref{lem_Fredholm} which
proves that $A$ is Fredholm.

(v) For arbitrary $\lambda\in\C$, $\lambda I-B\in\opsp(\lambda I-A)$
iff $B\in\opsp(A)$. So it suffices to show that Fredholmness of a
rich band-dominated operator (for $p\ne\infty$) implies
invertibility of its limit operators. But this is an immediate
consequence of (i) and (iii).

(vi) From (iii) we know that, if $B\in\opsp(A)$ and $\lambda I-B$ is
not invertible then $\lambda I-A$ is not invertible (not even
invertible at infinity). So suppose $\lambda I-B$ is invertible. If
$\lambda I-A$ is not invertible then there is nothing to prove. If
also $\lambda I-A$ is invertible then, by Theorem \ref{lem_opsp}
(ix), which applies since $B\in\opsp(A)\subset L(Y,\P)$ as $A\in
BDO(Y)\subset L(Y,\P)$, it follows that
$$
\|(\lambda I-B)^{-1}\|=\frac 1{\nu(\lambda I-B)} \le \frac
1{\nu(\lambda I-A)}=\|(\lambda I-A)^{-1}\|.
$$
\Proofend

Now we will combine the results of Theorems \ref{prop_appl2} and
\ref{prop_rich_invinf}. Recall that the set $\UM(Y)$ was introduced
just before (and studied in) Lemma \ref{lem_7.23isBA}.

\begin{Corollary} \label{cor_combine}
Consider $Y=Y^\infty(U)$ where $U$ has a predual $U^\pa$, and
suppose $A=I-K\in BDO(Y)$ is rich, has a preadjoint $A^\pa\in
L(Y^\pa)$ where $Y^\pa=Y^1(U^\pa)$, and that $K\in\UM(Y)$. Then the
following statements are equivalent.
\begin{itemize}
\item[(a)] all limit operators of $A$ are injective ($\alpha(A_h)=0$
for all $A_h\in\opsp(A)$) and there is an $S-$dense subset,
$\sigma$, of $\opsp(A)$ such that $\beta(A_h)=0$ for all
$A_h\in\sigma$;
\item[(b)] all limit operators of $A$ are injective ($\alpha(A_h)=0$
for all $A_h\in\opsp(A)$) and there is an $S-$dense subset,
$\sigma$, of $\opsp(A)$ such that $\alpha(A_h^\pa)=0$ for all
$A_h\in\sigma$;
\item[(c)] $A$ is invertible at infinity;
\item[(d)] $A$ is Fredholm.
\end{itemize}
\end{Corollary}

\Proof Note first that, by Lemma \ref{lem_lem2}, each
$A_h\in\opsp(A)$ has a well-defined pre\-adjoint $A_h^\pa\in
L(Y^\pa)$ so that statement (b) is well-defined; in fact, by Lemma
\ref{lem_lem2}, $\{A_h^\pa : A_h\in\opsp(A)\}=\opsp(A^\pa)$. Since
always $\beta(A_h)\ge\alpha(A_h^\pa)$ \cite{Jorgens}, clearly
(a)$\Rightarrow$(b).

If (b) holds then, noting that property (i) of Lemma
\ref{lem_coll_comp_entries} implies that $\opsp(K)$ is uniformly
Montel on $(Y,s)$, applying Theorem \ref{prop_appl2}, $\opsp(A)$ is
uniformly bounded below, which implies that the range of each
$A_h\in\opsp(A)$ is closed. This implies that
$\beta(A_h)=\alpha(A_h^\pa)=0$ \cite{Jorgens} for each
$A_h\in\sigma$, so that (b)$\Rightarrow$(a) and each $A_h\in\sigma$
is surjective.

Applying Theorem \ref{prop_appl2} again, we see that all the
elements of $\opsp(A)$ are invertible and their inverses are
uniformly bounded. Applying Theorem \ref{prop_rich_invinf} we
conclude that (a)$\Leftrightarrow$(b)$\Leftrightarrow$(c).

The implication (c)$\Rightarrow$(d) follows from Theorem
\ref{prop_rich_invinf} (ii) with $C=I$ and $-K\in M(Y)$ by property
(iii) of Lemma \ref{lem_coll_comp_entries}. Finally,
(d)$\Rightarrow$(c) is Theorem \ref{prop_rich_invinf} b). \Proofend

We note that Corollary \ref{cor_combine}, for operators satisfying
the conditions of the corollary, reduces the problem of establishing
Fredholmness and/or invertibility at infinity on $Y^\infty(U)$ to
one of establishing injectivity of the elements of $\opsp(A)$ and a
subset of the elements of $\opsp(A^\pa)$. In applications in
mathematical physics this injectivity can sometimes be established
directly via energy or other arguments (e.g. \cite{CW_MathMeth97}),
this reminiscent of classical applications of boundary integral
equations in mathematical physics where $A=I+K$ with $K$ compact,
and injectivity of $A$ is established from equivalence with a
boundary value problem.

In the one-dimensional case $N=1$ a stronger version of Theorem
\ref{prop_appl1} can be shown, namely Theorem \ref{prop_FC_JFA}
below. This result is shown by establishing, in the case in which
$A=I-K\in BDO(Y)$ is rich, $K\in UM(Y)$, and all the limit operators
of $A$ are injective, the following three statements:

\begin{tabular}{rp{0.9\textwidth}}
a)& If $B\in\opsp(A)$ has a surjective limit operator then
$B$ is surjective itself.\\
b)& Every $B\in\opsp(A)$ has a self-similar limit operator.\\
c)& Self-similar limit operators (of $A$, including those of $B$)
are surjective.
\end{tabular}

That b) holds is Lemma \ref{prop_ss}. The proofs of a) and c) in
\cite{CWLi:JFA2007} are both examples of the application of Theorem
\ref{th4.1}.

\begin{Theorem} \cite{CWLi:JFA2007} \label{prop_FC_JFA}
Suppose that $Y=\ell^\infty(\Z,U)$, that $A=I-K\in BDO(Y)$ is rich,
that $K\in UM(Y)$, and that all the limit operators of $A$ are
injective.  Then all elements of $\opsp(A)$ are invertible and their
inverses are uniformly bounded.
\end{Theorem}

Combining this result with Corollary \ref{cor_combine} we have the following
simplified version of Corollary \ref{cor_combine} which holds in the
one-dimensional case.

\begin{Corollary} \label{cor_combine2}
Suppose $Y=\ell^\infty(\Z,U)$ where $U$ has a predual $U^\pa$, and
suppose $A=I-K\in BDO(Y)$ is rich, has a preadjoint $A^\pa\in
L(Y^\pa)$ where $Y^\pa=\ell^1(\Z,U^\pa)$, and that $K\in\UM(Y)$. Then the
following statements are equivalent:
\begin{itemize}
\item[(a)] all limit operators of $A$ are injective;
\item[(b)] all elements of $\opsp(A)$ are
invertible and their inverses are uniformly bounded;
\item[(c)] $A$ is invertible at infinity;
\item[(d)] $A$ is Fredholm.
\end{itemize}
\end{Corollary}

\section{Almost Periodic Band-Dominated Operators} \label{sec_NRBDO}
Recall the classes of all absolutely rich/periodic and all
norm-rich/almost periodic operators as introduced in Section
\ref{sec_NR+AR}. We will now look at operators that are
norm-rich/almost periodic and band-dominated at the same time. We
first show that every norm-rich/almost periodic band-dominated
operator can be approximated in the norm by norm-rich/almost
periodic band operators. The same statement holds with
`norm-rich/almost periodic' replaced by `rich', as was first pointed
out in \cite[Proposition 2.9]{LiDiss}. The proof of our lemma is
very similar to that of this related statement. An alternative, less
constructive method for approximating a norm rich/almost periodic
operator by a norm-rich/almost periodic band operator is described
in step I of the proof of \cite[Theorem 1]{Kurbatov1989}.

\begin{Remark} \label{rem_BOapprox}
Note that the norm-approximation of a given band-dominated operator
$A$ by a sequence $A_n$ of band operators is in general a more
involved problem than it might seem. In particular, it is not always
possible to achieve this approximation by letting $[A_n]$ be the
restriction of $[A]$ to $n$ of its diagonals (see Remarks 1.40 and
1.44 in \cite{LiBook}). Instead, for a given $A\in BDO(Y)$, in the
proof of \cite[Theorem 2.1.6]{RaRoSiBook} a sequence of band
operators
\begin{eqnarray}
\label{eq_opint1} A_n&=&\sum_{|k|\le n}c_{k,n}\ B_k,\qquad n\in\N,\\
\label{eq_opint2} \textrm{with}\quad
B_k&=&\int_{[0,2\pi]^N}M_{e_t}AM_{e_{-t}}\,e^{-i\,(t,k)}\ dt,\qquad
k\in\Z^N
\end{eqnarray}
is constructed, where $c_{k,n}\in\C$ and
$e_t(m)=e^{i(t_1m_1+...+t_Nm_N)}$ for all $m\in\Z^N$ and $t\in\R^N$.
This construction is such that each matrix $[B_k]$ is only supported
on the $k$th diagonal and $A_n\toto A$ as $n\to\infty$.
\end{Remark}

\begin{Lemma} \label{lem_BDOBONR}
For every band-dominated operator $A$ and the corresponding
approximating sequence $(A_n)$ of band operators \refeq{eq_opint1},
the following holds.
\begin{itemize}
\item[(i)] If $A$ is norm-rich/almost periodic, then each one of the band
operators $A_n$ is norm-rich/almost periodic.
\item[(ii)] If $A$ is rich and $\opsp(A)$ is uniformly Montel, then every
operator spectrum $\opsp(A_n)$ is uniformly Montel.
\end{itemize}
\end{Lemma}

\Proof Since the integrand in \refeq{eq_opint2} is continuous in
$t$, the integral can be understood in the Riemann sense and
therefore $B_k$ can be approximated in norm by the corresponding
Riemann sums
\begin{equation} \label{eq_riemannsums}
R_m^{(k)}\ = \left(\frac{2\pi}m\right)^N\
\sum_{j=1}^{m^N}M_{e_{t_{m,j}}}AM_{e_{-t_{m,j}}}\,e^{-i\,(t_{m,j},k)},\qquad
m\in\N
\end{equation}
as $m\to\infty$. Here $t_{m,j}\in T_{m,j}$ are arbitrary where
$\{T_{m,j}:j=1,...,m^N\}$ is a partition of $[0,2\pi]^N$ into
hyper-cubes of width $2\pi/m$ (also see the proof of \cite[Theorem
2.1.18]{RaRoSiBook}).

To prove (i) it suffices, by Lemma \ref{lem_NRisBA}, to show that
all Riemann sums $R_m^{(k)}$ are norm-rich/almost periodic. Since, as the
restriction of an almost periodic (even periodic) function
$\R^N\to\C$ to the integer grid $\Z^N$, the sequence $e_t$ is almost
periodic for every choice of $t\in\R^N$, we get that both
$M_{e_{t_{m,j}}}$ and $M_{e_{-t_{m,j}}}$ are norm-rich/almost periodic (see
Lemma
\ref{lem_BONR} below). By Lemma \ref{lem_NRisBA} again and
$A\in\NR(Y)$, it follows that then all of the Riemann sums
$R_m^{(k)}$ and consequently, all operators $B_k$ and $A_n$ are
norm-rich/almost periodic as well.

For the proof of (ii), let $A$ be rich and $\opsp(A)$ be uniformly
Montel. By Lemma \ref{charact_uniformly_Montel}, using that $A$ is
rich, we get that $A\in\UM(Y)$. Since every $M_{e_{t_{m,j}}}$ is
rich (even norm-rich/almost periodic), we get that
$M_{e_{t_{m,j}}}AM_{e_{-t_{m,j}}}\in\UM(Y)$ for all $m\in\N$ and
$j\in\{1,...,m^N\}$, by Lemma \ref{lem_7.23isBA} (c) and (d). This
fact, together with formulas
\refeq{eq_opint1}--\refeq{eq_riemannsums} and Lemma
\ref{lem_7.23isBA} (a), shows that all $R_m^{(k)}$, $B_k$ and $A_n$
are in $\UM(Y)$. But the latter implies that $\opsp(A_n)$ is
uniformly Montel, by Lemma \ref{charact_uniformly_Montel}. \Proofend

Having reduced the study of norm-rich/almost periodic band-dominated
operators to norm-rich/almost periodic band operators, we will now
have a closer look at the latter class.

\begin{Lemma} \label{lem_BONR}
A band operator is norm-rich/almost periodic iff all of its
diagonals are almost periodic; that means
\[
A=\sum_{k\in D}M_{b_k}V_k\ \in\ \NR(Y)\qquad\textrm{iff} \qquad
b_k\in\AP(L(U)),\quad \forall k\in D
\]
for all finite sets $D\subset\Z^N$ and $b_k\in Y^\infty(L(U))$,
$k\in D$.
\end{Lemma}
\Proof Let $D\subset\Z^N$ be finite, let $b_k\in Y^\infty(L(U))$ for
all $k\in D$, and put $A=\sum_{k\in D}M_{b_k}V_k$. Note that
\begin{equation} \label{eq_shiftsBO}
V_{-m}AV_m\ =\ \sum_{k\in D}M_{V_{-m}b_k}V_k
\end{equation}
for every $m\in\Z^N$. We show that a sequence of operators
\refeq{eq_shiftsBO} converges in the operator norm iff all of the
corresponding diagonals $V_{-m}b_k$ converge in the norm of
$Y^\infty(L(U))$.

Suppose $A\in\NR(Y)$ and take an arbitrary sequence
$h=(h(n))_{n\in\N}\subset\Z^N$. Then there exists a subsequence $g$
of $h$ such that $V_{-g(n)}AV_{g(n)}\toto C$ for some $C\in L(Y)$.
Then, for all $i,j\in\Z^N$, with $[C]=[c_{i,j}]$ and with the
restriction and extension operators $R_i$ and $E_j$ as introduced at
the beginning of Section \ref{sec_BDO}, we have that
\begin{eqnarray}
\nonumber\|V_{-g(n)}b_{i-j}(i)-c_{i,j}\|_{L(U)}&=&\|R_i(V_{-g(n)}AV_{g(n)}-C)E_j\|_{L(U)}\\
&\le&\|V_{-g(n)}AV_{g(n)}-C\|_{L(Y(U))}\to 0
\label{eq_matrixentriesconverge}
\end{eqnarray}
as $n\to\infty$. Now, for every $k\in\Z^N$, define $c_k\in
Y^\infty(L(U))$ by $c_k(i)=c_{i,i-k}$, so that $c_k$ is the $k-$th
diagonal of $C$. From \refeq{eq_matrixentriesconverge} we get that
$\|V_{-g(n)}b_k-c_k\|_{Y^\infty}\to 0$, so that $b_k\in\AP(L(U))$
for each $k\in\Z^N$.

Now, conversely, suppose that $b_k\in\AP(L(U))$ for all $k\in D$ and
take an arbitrary sequence $h=(h(n))_{n\in\N}$. Let
$\{k_1,k_2,...,k_m\}$ be an enumeration of $D\subset\Z^N$, and
choose a subsequence $h^{(1)}\subset h$ such that
$V_{-h^{(1)}(n)}b_{k_1}$ converges. From this choose a subsequence
$h^{(2)}\subset h^{(1)}$ such that also $V_{-h^{(2)}(n)}b_{k_2}$
converges, etc., until we arrive at a sequence $g:=h^{(m)}\subset h$
for which all $V_{-g(n)}b_k$ with $k\in D$ converge. Denote the
respective limits by $c_k\in Y^\infty(L(U))$. Then we have
$V_{-g(n)}AV_{g(n)}\toto \sum_{k\in D}M_{c_k}V_k=:C$ since
\begin{eqnarray*}
\|V_{-g(n)}AV_{g(n)}-C\|_{L(Y)}
&=&\left\|\sum_{k\in D}(M_{V_{-g(n)}b_k}-M_{c_k})V_k\right\|_{L(Y)}\\
&\le&\sum_{k\in D}\|V_{-g(n)}b_k-c_k\|_\infty\to 0
\end{eqnarray*}
as $n\to\infty$, showing that $A\in\NR(Y)$. \Proofend

To establish the main result of this section it is convenient to
introduce a further definition. For $1\leq r\leq N$ call $A\in L(Y)$
{\em $r$-partially periodic}\index{partially periodic
operator}\index{periodic
operator!partially}\index{operator!periodic!partially} if there
exist $m_1,\dots,m_r\in\N$ such that
\begin{equation} \label{r_absolute}
V_{m_je^{(j)}}A = A V_{m_je^{(j)}}, \quad j = 1,\dots, r.
\end{equation}
Then an $N$-partially periodic operator is precisely an absolutely
rich/periodic operator. By a $0$-partially periodic operator we
shall mean any operator in $L(Y)$. Similarly, we shall say that
$b\in Y^\infty(L(U))$ is {\em $r$-partially
periodic}\index{partially periodic sequence} if
\begin{equation} \label{r_periodic}
b(k+m_je^{(j)}) = b(k), \quad k\in\Z^N,\ j = 1,\dots, r.
\end{equation}
Then an $N$-partially periodic function is periodic. We shall say
that every $b\in Y^\infty(L(U))$ is $0$-partially periodic.

Arguing as in the proof of Lemma \ref{lem_BONR} we can show the
following result.

\begin{Lemma} \label{lem_BOPP}
A band operator is $r$-partially periodic iff all of its diagonals are
$r$-partially periodic; that means
\[
A=\sum_{k\in D}M_{b_k}V_k \mbox{ is $r$-partially periodic}\quad\textrm{iff}
\quad b_k \mbox{ is $r-$partially
periodic},\quad \forall k\in D,
\]
for all finite sets $D\subset\Z^N$ and $b_k\in Y^\infty(L(U))$,
$k\in D$.
\end{Lemma}

The proof of the following theorem on the invertibility of almost
periodic band operators depends in the first place on results on
invertibility of periodic band operators (Theorem
\ref{lem_A-inj-inv}) and, secondly, on the possibility of
approximation of almost periodic by periodic functions (cf.\ the
study of invertibility of elliptic, almost periodic differential and
pseudo-differential operators which dates back to
\cite{Muh1972,Shubin78}). This theorem is not new; it is Theorem 1
in \cite{Kurbatov1989} (and see \cite{Kurbatov}); in fact
\cite{Kurbatov1989} provides a proof of this result for the general
case $Y=Y^p$, $1\le p\le \infty$, by first reducing it to the case
$p=\infty$. However, it seems of interest to include our proof
which, while it has similarities with steps IV-VII in \cite[Theorem
1]{Kurbatov1989}, illustrates the application of Theorem
\ref{prop5.10} which carries a large part of the proof, and differs
also in applying an inductive argument.

\begin{Theorem} \label{prop_NRBO}
If $A=I+K$, $K\in BO(Y^\infty)$ is norm-rich/almost periodic,
$\opsp(K)$ is uniformly Montel, and $A$ is bounded below, then $A$
is invertible.
\end{Theorem}
\Proof Note first that, since $K\in\opsp(K)$ if $K$ is norm rich
(Theorem \ref{prop_NR1} (ii)), that $\opsp(K)$ uniformly Montel
implies $K\in M(Y^\infty)$.

We shall establish the theorem by proving, by induction, that, for
$r=0,1,\dots,N$,
\begin{eqnarray} \nonumber \label{eqn_ind}
\mbox{if $A$ satisfies the conditions of the theorem and, }\\
\mbox{additionally, $A$ is $r$-partially periodic, then $A$ is
invertible.}
\end{eqnarray}
Statement (\ref{eqn_ind}) for $r=0$ is precisely the theorem that we
wish to prove.

That (\ref{eqn_ind}) holds for $r=N$ follows from Theorem
\ref{lem_A-inj-inv}. Now suppose that (\ref{eqn_ind}) holds for
$r=s$, for some $s\in \{1,\dots,N\}$, and that $A$ satisfies the
conditions of the theorem and, additionally, is $(s-1)$-partially
periodic. We will show that this implies that $A$ is invertible, so
that (\ref{eqn_ind}) holds for $r=s-1$, proving the inductive step.
That $A$ is invertible will be proved by applying Theorem
\ref{prop5.10}.

First note that, by Lemmas \ref{lem_BONR} and \ref{lem_BOPP}, for
some finite set $D\subset\Z^N$ and $b_k\in Y^\infty(L(U))$, $k\in
D$, it holds that
$$
K=\sum_{k\in D}M_{b_k}V_k
$$
with each $b_k\in \AP(L(U))$ and $(s-1)$-partially periodic. For
$n\in\N$ let $h(n)\in \Z^N$ have $s$th component $n^2$ and all other
components zero. Since each $b_k\in \AP(L(U))$ we can choose a
subsequence $g$ of $h$ such that $V_{-g(n)}b_k$ is convergent as
$n\to\infty$ for each $k\in D$. Thus, defining $f(n)=g(n+1)-g(n)$,
$$
\|V_{-f(n)}b_k-b_k\|_\infty =
\|V_{-g(n+1)}b_k-V_{-g(n)}b_k\|_\infty\to 0
$$
as $n\to\infty$. Note that the $sth$ component of $f(n)$ is $\geq
2n+1$ and that all the other components of $f(n)$ are zero.

For $j\in \N$ define $\hat P_j\in L(Y^\infty(L(U)))$ by
$$
\hat P_j b(\ell) = b(\ell), \quad -\frac{j}{2}<\ell_s\leq
\frac{j}{2},
$$
for $\ell=(\ell_1,\dots,\ell_N)\in\Z^N$, and by the requirement that
$$
\hat P_jb(\ell+je^{(s)})=\hat P_jb(\ell), \quad \ell\in\Z^N.
$$
An important observation is that, for $j\in\N$, $b\in
Y^\infty(L(U))$,
\begin{equation} \label{eqn_shift}
\|P_j (\hat P_j b - b)\|_\infty \leq \|V_{je^{(s)}}b-b\|_\infty =
\|V_{-je^{(s)}}b-b\|_\infty.
\end{equation}

Using this notation, for each $k\in D$ let us define a sequence of
increasingly good approximations to $b_k$, each approximation being
$s$-partially periodic and almost periodic. In detail, for $n\in\N$
and $k\in D$, let
$$
b_k^{(n)} := \hat P_{f(n)}b_k.
$$
Further, define $K_n$, an approximation to $K$, by
$$
K_n=\sum_{k\in D}M_{b_k^{(n)}}V_k.
$$
Then it is clear that each $b_k^{(n)}$ is $s$-partially periodic. To
see that $b_k^{(n)}$ is also almost periodic, note that every
$\ell\in\Z^N$ can be written as $\ell=\hat\ell+\tilde\ell$, where
$\hat \ell=\ell_s e^{(s)}$ ($\ell_s$ the $s$th component of $\ell$).
Then $V_\ell=V_{\hat\ell}V_{\tilde \ell}$ and, for $j\in\N$,
$V_{\tilde\ell}$ commutes with $\hat P_j$. Thus, for every $n\in\N$
and $\ell\in\Z^N$,
\begin{equation} \label{eqn_swap}
V_\ell b_k^{(n)} = V_{\hat \ell} \hat P_{f(n)} V_{\tilde \ell}b_k =
V_{j e^{(s)}} \hat P_{f(n)} V_{\tilde \ell}b_k,
\end{equation}
for some $j\in\Z$ with $|j|\leq f(n)/2$. From this formula it is
clear that $\{V_\ell b_k^{(n)}:\ell\in\Z^N\}$ is relatively compact
if $\{V_\ell b_k:\ell\in\Z^N\}$ is relatively compact, so that
$b_k^{(n)}$ is almost periodic. Moreover it is clear that
$b_k^{(n)}\sto b_k$ as $n\to\infty$ for each $k\in D$, so that
$K_n\sto K$ and $A_n:=I+K_n\sto A$. Further, since each $b_k^{(n)}$
is $s$-partially periodic and almost periodic, it follows from
Lemmas \ref{lem_BONR} and \ref{lem_BOPP} that $K_n$ is
$s$-partially periodic and norm rich, and so therefore is $A_n$. Note
also that, by Lemma \ref{charact_uniformly_Montel} and the
equivalence of (iii) and (v) in Lemma \ref{lem_coll_comp_entries},
the set $\{b_k(\ell):\ell\in\Z^N\}$ is collectively compact for
every $k\in D$, so that
$\{b_k^{(n)}(\ell):\ell\in\Z^N,n\in\N\}\subset
\{b_k(\ell):\ell\in\Z^N\}$ is collectively compact. Arguing as in
the proof that (v) implies (iii) in Lemma
\ref{lem_coll_comp_entries}, we deduce that
$\{V_{-j}K_nV_j:j\in\Z^N,n\in\N\}$ is uniformly Montel, so in
particular, by the equivalence of (i) and (ii) in Lemma
\ref{lem_coll_comp_entries}, $\opsp(K_n)$ is uniformly Montel for
each $n$. By the inductive hypothesis it follows that if $A_n$ is
bounded below then $A_n$ is surjective, for each $n$.

We have shown that conditions (a)-(c) of Theorem \ref{prop5.10} are
satisfied by $A$ and $A_n$. To complete the proof we will show that
conditions (d) and (e) of Theorem \ref{prop5.10} are satisfied with
${\mathcal B} =\opsp(A)$. Since $A$ is bounded below which implies,
by (v) of Theorem \ref{prop_NR1}, that ${\mathcal B}$ is uniformly
bounded below, this choice of ${\mathcal B}$ satisfies (e) of
Theorem \ref{prop5.10}.

To see that (d) holds with ${\mathcal B} =\opsp(A)$, suppose that
$\ell(j)\in \Z^N$ and $n(j)\in \N$, for $j\in\N$, and that
$n(j)\to\infty$ as $j\to\infty$, and let
$$
B_j := V_{-\ell(j)}K_{n(j)}V_{\ell(j)} = \sum_{k\in D} M_{c_k^{(j)}}
V_k,
$$
where $c_k^{(j)} := V_{-\ell(j)} b_k^{n(j)}$. What we have to show
is that $(B_j)_{j\in\N}$ has a subsequence that $s$-converges to an
element of $\opsp(K)$. To see this, first note that, by
(\ref{eqn_swap}),
$$
c_k^{(j)}  = V_{p(j)e^{(s)}}\hat P_{f(n(j))} V_{-\tilde \ell(j)}
b_k,
$$
 for some $p(j)\in \Z$ with $|p(j)|\leq
f(n(j))/2$. Since each $b_k$ is almost periodic, setting
$\check{\ell}(j) := -p(j)e^{(s)}+\tilde \ell(j)$, we can find a
subsequence of $\check \ell$, which we will denote again by $\check
\ell$, such that, for each $k\in D$, there exists $\hat b_k\in
Y^\infty(L(U))$ such that
$$
\eps_j := \|V_{-\check \ell(j)} b_k - \hat b_k\|_\infty \to 0
$$
as $j\to\infty$. Then, for $m\in\N$ and all $j$ sufficiently large
such that $m\leq f(n(j))/2$,
\begin{eqnarray*}
\|P_m(c_k^{(j)} - \hat b_k)\|_\infty & \leq & \eps_j +
\|P_{m}(c_k^{(j)}- V_{-\check
\ell(j)}b_k)\|_\infty\\
& = & \eps_j + \|P_{m}V_{p(j)e^{(s)}}(\hat P_{f(n(j))}-I) V_{-\tilde
\ell(j)}b_k\|_\infty\\
& \leq & \eps_j +
\|P_{m+|p(j)|}(\hat P_{f(n(j))}-I) V_{-\tilde \ell(j)}b_k\|_\infty\\
& \leq & \eps_j + \|P_{f(n(j))}(\hat P_{f(n(j))}-I) V_{-\tilde
\ell(j)}b_k\|_\infty\\
& \leq & \eps_j + \|V_{-f(n(j))}b_k -b_k\|_\infty,
\end{eqnarray*}
by (\ref{eqn_shift}). Thus $\|P_m(c_k^{(j)} - \hat b_k)\|_\infty\to
0$ as $j\to\infty$, for every $m$, so that $c_k^{(j)} \sto \hat
b_k$.

Defining
$$
\hat K := \sum_{k\in D} M_{\hat b_k} V_k,
$$
we have that $\eps_j\to 0$ as $j\to\infty$ for each $k\in D$, so
that
$$
V_{-\check \ell(j)}K V_{\check \ell(j)} = \sum_{k\in D}
M_{V_{-\check \ell(j)} b_k} V_k \toto \sum_{k\in D} M_{\hat b_k} V_k
= \hat K.
$$
If $\check \ell$ has an unbounded subsequence, this implies that
$\hat K\in\opsp(K)$. If $\check \ell$ does not have an unbounded
subsequence then it has a constant subsequence, i.e.\
$\check\ell(j)=\tilde \ell$ for some $\tilde\ell\in\Z^N$ and
infinitely many $j\in\N$, so that $\hat K = V_{-\tilde \ell}K
V_{\tilde \ell}$. Since $K\in \opsp(K)$ implies (by Theorem
\ref{lem_opsp} (ii)) that $V_{-\tilde \ell}K V_{\tilde \ell}\in
\opsp(K)$, it follows also in this case that $\hat K \in\opsp(K)$.
But $c_k^{(j)} \sto \hat b_k$ for every $k\in D$ implies that
$$
B_j = \sum_{k\in D} M_{c_k^{(j)}} V_k \sto \sum_{k\in D} M_{\hat
b_k} V_k = \hat K,
$$
so that we have shown that $(B_j)_{j\in\N}$ has a subsequence that
$s$-converges to an element of $\opsp(K)$.
 \Proofend

Now we can finally prove the generalisation of Corollary
\ref{cor_absrich_injinv} to norm-rich/almost periodic band-dominated
operators. We note that results of this flavour in concrete cases,
in particular showing something close to equivalence of (i) and
(iii), date back at least to Shubin \cite{Shubin78} for scalar
elliptic differential operators with smooth almost periodic
coefficients, where the analogous statement to (iii) is termed the
{\em Favard condition}\index{Favard condition} (and see
\cite{CWLi:JFA2007}). Note also that we have already seen, in
Theorem \ref{prop_FC_JFA}, that in the one-dimensional case $N=1$
the equivalence of (iii) and (iv) holds even when $K$ is only rich
rather than norm rich/almost periodic.

\begin{Theorem} \label{prop_NRBDO}
If $A=I+K$ with $K\in BDO(Y^\infty)$ norm-rich/almost periodic and
$\opsp(K)$ is uniformly Montel, then the following statements are
equivalent.
\begin{itemize}
\item[(i)] $A$ is invertible;
\item[(ii)] $A$ is bounded below;
\item[(iii)] all limit operators of $A$ are injective;
\item[(iv)] all limit operators of $A$ are invertible with uniformly bounded
inverses.
\end{itemize}
\end{Theorem}
\Proof Let the conditions of the theorem be fulfilled. Implication
(iv)$\Rightarrow$(i) follows by Theorem \ref{prop_NR1} (ii);
(i)$\Rightarrow$(ii) is trivial; (ii)$\Rightarrow$(iii) follows by
Theorem \ref{prop_NR1} (v); and it remains to show
(iii)$\Rightarrow$(iv).

Suppose (iii) holds. By Theorem \ref{prop_appl2} it follows that
$\opsp(A)$ is uniformly bounded below, so that, by Theorem
\ref{prop_NR1} (v), $\nu(A_h)=\nu(A)>0$ for all $A_h\in\opsp(A)$.
Now, from Lemma \ref{lem_BDOBONR}, applied to $K$ in place of $A$,
we know that there exists a sequence $(K_n)$ of norm-rich/almost
periodic band operators with $\opsp(K_n)$ uniformly Montel for each
$n$ and $K_n\toto K$. From \refeq{eq_lonoemineq} it follows that
$\nu(I+K_n)\to\nu(I+K)=\nu(A)>0$, so that there exists a $n_0\in\N$
with
\[
\nu(I+K_n)\ >\ \frac 12 \nu(A)\ >\ 0,\quad  n\ge n_0.
\]
Invoking Theorem \ref{prop_NRBO}, we conclude that, for every $n\ge
n_0$, $I+K_n$ is invertible and
\[
\|(I+K_n)^{-1}\|\ =\ \frac 1{\nu(I+K_n)}\ <\ \frac 2{\nu(A)}\ <\
\infty.
\]
Since also $I+K_n \toto I+K = A$, we get, by \cite[Lemma
1.3]{LiBook}, that $A$ is invertible. From Theorem \ref{prop_NR1}
(iv) and (v) it now follows that each limit operator $A_h$ of $A$ is
invertible and
\[
\|(A_h)^{-1}\|\ =\ \frac 1{\nu(A_h)}\ =\ \frac 1{\nu(A)}\ <\ \infty.
\]
But this is (iv). \Proofend

\section{The Wiener Algebra} \label{subsection_wiener}
If we write of ``all spaces $Y$'' in this section, we always think
of the set of spaces $Y=Y^p(U)$ where $U$ is fixed, and only $p$
varies in $\{0\}\cup[1,\infty]$. For a band operator $A$ of the form
\refeq{eq_BO}, one clearly has
$$
\|A\|_{L(Y)}\ =\ \left\|\sum_{|k|\le w} M_{b_k}V_k\right\|_{L(Y)}\
\le\ \sum_{k\in\Z^N} \|b_k\|_{Y^\infty}\ =:\
\|A\|_\W\index{$\mida$@$\parallel\cdot\parallel_\W$}
$$
for all spaces $Y$, where we put $b_k=0$ if $|k|>w$. By
$\W=\W(U)$\index{$W$@$\W$}\index{$W(U)$@$\W(U)$} we denote the
closure of $BO(Y)$ with respect to $\|.\|_\W$. Equipped with
$\|.\|_\W$, the set $\W$ is a Banach algebra. Clearly $\W$ does not
depend on the underlying space $Y$ and we have, for all spaces $Y$,
$$
BO(Y)\ \subset\ \W\ \subset\ BDO(Y).
$$
Both inclusions are proper as the operators $A$ and $B$ with
$(Ax)(k)=x(-k)/|k|^2$ and $(Bx)(k)=x(-k)/|k|$ if $k\ne 0$,
$(Ax)(0)=0=(Bx)(0)$ show.

$\W$ is referred to as the {\em Wiener algebra}\index{Wiener
algebra}. Note that this is a natural (non-stationary) extension of
the classical algebra of all matrix operators with constant
diagonals and $\|A\|_W<\infty$ (which is isomorphic, via Fourier
transform, to the algebra of all periodic functions with absolutely
summable sequence of Fourier coefficients that is usually associated
with the name `Wiener algebra'). As a consequence of the definition
of $\W$ we get that
\begin{equation} \label{eq_Wienernorm}
\|A\|_{L(Y)}\ \le\ \|A\|_\W
\end{equation}
for all $A\in\W$ and all spaces $Y$. For operators in the Wiener
algebra $\W$ one has the following remarkable results (see Theorem
2.5.2 and Proposition 2.5.6 of \cite{RaRoSiBook}).

\begin{Lemma} \label{Wiener1}
For every $A\in\W$, the following holds.
\begin{itemize}
\item[(i)] If $A$ is invertible then $A^{-1}\in\W$, so that $\W$ is inverse
closed.
\item[(ii)] If $A$ is rich and $h=(h(n))\subset\Z^N$ tends to infinity
then there is a subsequence $g$ of $h$ such that the limit operator
$A_g$ exists with respect to all spaces $Y$. This limit operator
again belongs to $\W$, and $\|A_g\|_\W\le\|A\|_\W$.
\end{itemize}
\end{Lemma}

As a consequence of Lemma \ref{Wiener1} we get that, if $A\in\W$ is
invertible on one of the spaces $Y$ then it is invertible on all
spaces $Y$, and $\|A^{-1}\|_{L(Y)}\le\|A^{-1}\|_\W$. Another
consequence of this lemma is that the operator spectrum $\opsp(A)$
is contained in the Wiener algebra $\W$ and does not depend on the
underlying space $Y$ if $A\in\W$. So for $A\in\W$, the statement of
Theorem \ref{prop_rich_invinf} (iii) holds independently of the
underlying space $Y$. Moreover, by Theorem \ref{prop_rich_invinf}
(iv), also the uniform boundedness condition of the inverses is
redundant since this is true for $p\in\{0,1,\infty\}$ and
consequently, also for $p\in (1,\infty)$ by Riesz-Thorin
interpolation:

\begin{Theorem} \label{Wiener2}
If $A\in\W$ is rich then the following statements are equivalent.
\begin{itemize}
\item[(i)] $A$ is invertible at infinity on \textbf{one} of the spaces $Y$.
\item[(ii)] $A$ is invertible at infinity on \textbf{all} the spaces $Y$.
\item[(iii)] All limit operators of $A$ are invertible on \textbf{one} of the
spaces $Y$.
\item[(iv)] All limit operators of $A$ are invertible on \textbf{all} the
spaces
$Y$ and
\begin{equation} \label{eqn_bbound}
\sup_{p\in\{0\}\cup [1,\infty]}\;\sup_{A_h\in \opsp(A)}
\|A_h^{-1}\|_{L(Y^p)} <\infty.
\end{equation}
\end{itemize}
\end{Theorem}

\begin{Remark}
This theorem is a significant strengthening and simplification of
Theorem 2.5.7 in \cite{RaRoSiBook}. Theorem 2.5.7 requires that $U$
is reflexive, and, in the case that $U$ is reflexive, it implies
only a reduced version of our Theorem \ref{Wiener2} with the value
of $Y$ restricted to $Y^p$, $p\in\{0\}\cup (1,\infty)$, in
(i)--(iii).
\end{Remark}

{\it Proof of Theorem \ref{Wiener2}.} (i)$\Rightarrow$(iii) follows
from Theorem \ref{prop_rich_invinf} (iii).

(iii)$\Rightarrow$(iv): Suppose (iii) holds. We have observed
already that $\opsp(A)\subset\W$ is independent of the space $Y$ by
Lemma \ref{Wiener1} (ii). Applying Lemma \ref{Wiener1} (i) to the
limit operators of $A$, it follows that these limit operators are
invertible on all the spaces $Y$. By Theorem \ref{prop_rich_invinf}
(iv),
\[
s_p\ :=\ \sup_{A_h\in\opsp(A)}\|A_h^{-1}\|_{L(Y^p)}
\]
is finite for $p\in\{0,1,\infty\}$. Now, by Riesz-Thorin
interpolation (as demonstrated in the proof of \cite[Theorem
2.5.7]{RaRoSiBook}), we get that $ s_p\le s_1^{1/p}\
s_\infty^{1-1/p}<\infty$  for all $p\in(1,\infty)$, which proves
(iv).

(iv)$\Rightarrow$(ii) follows from Theorem \ref{prop_rich_invinf}
(iii).

Finally, (ii)$\Rightarrow$(i) is evident. \Proofend

From the above result and the relationship between invertibility at
infinity and Fredholmness, (see Remark \ref{rem_Fred_invinf} and
Theorem \ref{prop_rich_invinf} a)(i)+(ii) and b)), we can deduce
Corollary \ref{prop_W_selfadj} below, which relates Fredholmness to
invertibility of limit operators. In this corollary we require, for
the equivalence of (a)-(d) with (e), the existence of a predual
$U^\pa$ and, for $A$ considered as an operator on $Y^\infty$, the
existence of a preadjoint $A^\pa\in Y^1(U^\pa)$. The following
obvious lemma characterises existence of a preadjoint in terms of
existence of preadjoints of the elements of the matrix
representation of $A$.

\begin{Lemma} \label{lem_preadjoint} If $A\in \W(U)$, with $[A]=[a_{mn}]$, and $U$ has a
predual $U^\pa$, then $A$, considered as an operator on
$Y^\infty(U)$, has a preadjoint $A^\pa\in Y^1(U^\pa)$ iff each entry
$a_{mn}\in L(U)$ of the matrix representation of $A$ has a
preadjoint $a_{mn}^\pa\in L(U^\pa)$. If this latter condition holds
then a preadjoint is $A^\pa \in \W(U^\pa)$ with $[A^\pa]
=[a_{nm}^\pa]$. In particular, $A$ has a preadjoint if $U$ is
reflexive, given by $A^\pa\in \W(U^*)=\W(U^\pa)$, with $[A^\pa] =
[a_{nm}^*]$, where $a_{nm}^*\in L(U^*)=L(U^\pa)$ is the adjoint of
$a_{nm}$.
\end{Lemma}

\begin{Corollary} \label{prop_W_selfadj}
Suppose $A=I-K\in\W$ is rich, and $K\in\UM(Y)$. Then the following
statements are equivalent.
\begin{itemize}
\item[(a)] All limit operators of $A$ are injective on $Y^\infty$
and $\opsp(A)$ has an $S-$dense subset, $\sigma$, such that
$\beta(A_h)=0$ on $Y^\infty$ for all $A_h\in\sigma$;
\item[(b)] $A$ is invertible at infinity on \textbf{all} the spaces $Y$;
\item[(c)] $A$ is invertible at infinity on \textbf{one} of the spaces $Y$;
\item[(d)] $A$ is Fredholm on \textbf{all} the spaces $Y$.
\end{itemize}
In the case that $U$ has a predual $U^\pa$ and $A$, considered as an
operator on $Y^\infty(U)$, has a preadjoint $A^\pa\in Y^1(U^\pa)$,
then (a)--(d) are equivalent to
\begin{itemize}
\item[(e)] $A$ is Fredholm on \textbf{one} of the spaces $Y$.
\end{itemize}
\end{Corollary}

\Proof For the clarity of our argument we introduce two more
statements:
\begin{itemize}
\item[(f)] All limit operators of $A$ are invertible on $Y^\infty$;
\item[(g)] $A$ is invertible at infinity on $Y^2$.
\end{itemize}
Each of these will turn out to be equivalent to (a)--(d).

By Theorem \ref{prop_appl2}, statement (a) is equivalent to (f),
which, by Theorem \ref{Wiener2}, is equivalent to each of (b), (c)
and (g).

Since $K\in M(Y)$, the implication (b)$\Rightarrow$(d) follows from
Theorem \ref{prop_rich_invinf} (ii) (applied with $C=-I$).

Since, obviously, (d) implies Fredholmness of $A$ on $Y^2$, it also
implies (g), by Theorem \ref{prop_rich_invinf} (i). Another obvious
consequence of (d) is (e).

Finally, suppose $U^\pa$ and $A^\pa$ exist and (e) holds for
$Y=Y^p$. If $p=\infty$, then (c) follows by Theorem
\ref{prop_rich_invinf} b), and otherwise, if $p<\infty$, then (c)
follows by Theorem \ref{prop_rich_invinf} (i). \Proofend

The above corollary implies, for rich operators in the Wiener
algebra which are of the form $A=I-K\in\W$ with $K\in\UM(Y)$, i.e.
$K$ is subject to the (equivalent) properties (i)-(v) in Lemma
\ref{lem_coll_comp_entries},  and which possess a preadjoint, that
Fredholmness on one of the spaces $Y$ implies Fredholmness on all
spaces $Y$. The argument to show this is indirect: it depends on the
connection between Fredholmness and invertibility at infinity and on
the equivalence of (i) and (ii) in Theorem \ref{Wiener2}.

Recently, Lindner \cite{LiWiener} has studied directly the
invariance of the Fredholm property across the spaces $Y$ for
general operators in the Wiener algebra, but with a slight
restriction on the Banach space $U$, that it is either
finite-dimensional or possesses the {\em hyperplane
property}\index{hyperplane property}, meaning that it is isomorphic
to a subspace of $U$ of co-dimension 1. An equivalent
characterisation of the hyperplane property \cite{LiWiener} is that
there exists a $B\in L(U)$ which is Fredholm of index 1. This
characterisation suggests that infinite-dimensional Banach spaces
which do not have the hyperplane property are unusual. In
\cite{LiWienerPre} Lindner lists many sets of conditions on $U$
which ensure that $U$ has the hyperplane property, and recalls that
it was a long-standing open problem due to Banach -- the so-called
{\em hyperplane problem}\index{hyperplane problem} -- whether there
exist any infinite-dimensional Banach spaces which do not have the
hyperplane property. An example was finally constructed by Gowers in
\cite{Gowers94}, for which work, and the resolution of other
long-standing open questions posed by Banach, he received the Fields
medal in 1998.

The main result proved by Lindner \cite{LiWiener} is the following:

\begin{Theorem} \label{thm_lind_wie}
Suppose that $U$ is finite-dimensional or has the hyperplane
property and that $A\in\W(U)$. Then:
\begin{itemize}
\item[(a)] If $A$ is Fredholm on one of the spaces $Y^p$, with $p\in \{0\}\cup
[1,\infty)$, then $A$ is Fredholm on all spaces $Y$.
\item[(b)] If $U$ has a predual $U^\pa$ and $A$, considered as acting on
$Y^\infty$, has a preadjoint $A^\pa\in Y^1(U^\pa)$, then $A$ is
Fredholm on one of the spaces $Y$ iff $A$ is Fredholm  on all of the
spaces $Y$.
\end{itemize}
If $A$ is Fredholm on every space $Y$ then the index is the same on
each space.
\end{Theorem}

By combining Theorem \ref{Wiener2} and Corollary \ref{prop_W_selfadj} we
get the following result. In this result, again, by
$\spess(A)$\index{$\pspe$@$\spess(A)$} we denote the essential
spectrum\index{essential spectrum}\index{spectrum!essential} of $A$ (the
set of $\lambda$ for which $\lambda I-A$ is not Fredholm) and by
$\sp(A)$\index{$\psp$@$\sp(A)$} the ordinary spectrum\index{spectrum} of
$A$.

\begin{Corollary} \label{cor_spessformula}
Suppose $A=I-K\in\W$ is rich, $K\in\UM(Y)$, $U$ has a predual
$U^\pa$, and that $A$, considered as an operator on $Y^\infty(U)$,
has a preadjoint $A^\pa\in Y^1(U^\pa)$. Then statements (i)--(iv) of
Theorem \ref{Wiener2} and (a)--(e) of Corollary \ref{prop_W_selfadj}
are all equivalent. Further, on every space $Y$ it holds that
\begin{equation} \label{eq_spess}
\spess(A)\ =\ \bigcup_{B\in\opsp(A)} \sp(B).
\end{equation}
\end{Corollary}
\Proof It remains only to show that, for every $\lambda\in\C$,
$\lambda I-A = (\lambda-1)I+K$ is Fredholm iff $(\lambda-1) I +L$ is
invertible for every $L\in \opsp(K)$. For $\lambda\neq 1$ this
follows from the earlier part of the corollary. This is true also
for $\lambda=1$ when $U$ is finite-dimensional (see Corollary
\ref{cor7.25} below). When $U$ is infinite-dimensional and
$\lambda=1$, then $(\lambda-1)I+K=K\in UM(Y)\subset M(Y)$ and so
$(\lambda-1)I+L=L\in M(Y)$. This implies that all the entries of the
matrix representations of $(\lambda-1)I+K$ and $(\lambda-1)I+L$ are
compact (in $K(U)$). Since $U$ is infinite-dimensional, this implies
that $(\lambda-1)I+K$ and $(\lambda-1)I+L$ are not Fredholm.
\Proofend

In the particularly simple case of a finite-dimensional space $U$ we
have the following extended version of Corollary
\ref{cor_spessformula}.

\begin{Corollary} \label{cor7.25}
Suppose $A\in\W$ and $U$ is finite-dimensional. Then statements
(i)--(iv) of Theorem \ref{Wiener2} and (a)--(e) of Corollary
\ref{prop_W_selfadj} are all equivalent. Moreover, if $A$ is subject
to all these equivalent statements then the index of $A$ is the same
on each space $Y$. Further, on every space $Y$, \refeq{eq_spess}
holds.
\end{Corollary}
\Proof To see that the conditions of Corollary
\ref{cor_spessformula} are fulfilled, recall Lemma \ref{lem_oprich}
and remember that, by Lemma \ref{lem_7.23isBA} (b), $K=I-A\in\UM(Y)$
if $\dim U<\infty$. Also recall that finite-dimensional spaces $U$
are reflexive, so that $U^\pa=U^*$ is a predual and existence of a
preadjoint follows from Lemma \ref{lem_preadjoint}.

It remains to show that the index of $A$ is the same on all spaces
$Y$. This follows from Theorem \ref{thm_lind_wie} above. More
directly, let $A$ be Fredholm on one of the spaces $Y=Y^p$ with
$p\in\{0\}\cup[1,\infty]$ and denote its index by $\ind_p\, A$. From
the equivalence of (d) and (e) in Corollary \ref{prop_W_selfadj} we
get that $A$ is Fredholm on $Y^2$. Now proceed exactly as in the
proof of \cite[Lemma 2.1]{Roch_ellp} to arrive at
$\ind_p\,A=\ind_2\,A$ (see the remark below for some details).
\Proofend

\begin{Remark}
{\bf a) } In \cite[Lemma 2.1]{Roch_ellp} the independence of the
index just follows for the spaces $Y^p$ with $p\in\{0\}\cup
(1,\infty)$ because this is the setting for which \cite[Theorem
8]{RaRoSi1998} and \cite[Theorem 2.5.7]{RaRoSiBook} yield that the
Fredholm property is independent of the underlying space. Thanks to
the equivalence of (d) and (e) in Corollary \ref{prop_W_selfadj}, as
we proved earlier in Corollary \ref{cor_spessformula} (and
\ref{cor7.25}), we can now enlarge this framework to
$p\in\{0\}\cup[1,\infty]$.

{\bf b) } For simplicity, \cite[Lemma 2.1]{Roch_ellp} was only
stated for band operators but its proof actually applies to all
operators $A\in\W$. We make use of this fact in the proof of our
Corollary \ref{cor7.25}.

{\bf c)} Also, for simplicity, \cite[Lemma 2.1]{Roch_ellp} was only
stated for operators on $Y^p(\C)$ instead of $Y^p(U)$ with $n:=\dim
U<\infty$. This is not a loss of generality since these two spaces
are isomorphic and the discussed operator properties are preserved
under this isomorphism. Indeed, fix a basis in $U$ and let
$\ph:U\to\C^n$ refer to the isomorphism that maps $u\in U$ to its
coordinate vector $\ph(u)=:(\ph_1(u),...,\ph_n(u))\in\C^n$ with
respect to this basis. Then $\Phi:Y^p(U)\to Y^p(\C)$ with
\begin{eqnarray*}
(\Phi x)(n\cdot m_1+k,\,m_2,\,...\,,\,m_N)&=&
\ph_k\big(x(m_1,...,m_N)\big)\in\C,\\
& & k\in\{1,...,n\},\ m_1,...,m_N\in\Z
\end{eqnarray*}
for every $x\in Y^p(U)$ is an isomorphism. On the operator side,
roughly speaking, the matrix representation of an operator on
$Y^p(U)$ is an infinite matrix the entries of which are $n\times n$
matrices (operators on $U\cong\C^n$, via $\ph$). Via $\Phi$ this
matrix is identified, in a natural way, with an infinite matrix with
scalar entries, and this is the setting in which \cite[Lemma
2.1]{Roch_ellp} applies. Note that this identification preserves
membership of the Wiener algebra, Fredholmness, and the index of the
operator.
\end{Remark}

Finally, we note that in the one-dimensional case $N=1$ we have the
following refinement of Corollaries \ref{prop_W_selfadj},
\ref{cor_spessformula} and \ref{cor7.25}, as a consequence of
Theorem \ref{prop_FC_JFA}.

\begin{Corollary} \label{cor_oned_wie}
Suppose $N=1$ and that $A=I-K\in\W$ is rich and $K\in\UM(Y)$.  Then the
following statements are equivalent:
\begin{itemize}
\item[(a)] All limit operators of $A$ are injective on $Y^\infty$;
\item[(b)] All limit operators of $A$ are invertible on \textbf{one} of the
spaces $Y$;
\item[(c)] All limit operators of $A$ are invertible on \textbf{all} the spaces
$Y$ and (\ref{eqn_bbound}) holds;
\item[(d)] $A$ is invertible at infinity on \textbf{all} the spaces $Y$;
\item[(e)] $A$ is invertible at infinity on \textbf{one} of the spaces $Y$;
\item[(f)] $A$ is Fredholm on \textbf{all} the spaces $Y$.
\end{itemize}

In the case that $U$ has a predual $U^\pa$ and $A$, considered as an
operator on $Y^\infty(U)=\ell^\infty(\Z,U)$, has a preadjoint
$A^\pa$ on $Y^1(U^\pa)=\ell^1(\Z,U^\pa)$, then (a)--(f) are
equivalent to

\begin{itemize}
\item[(g)] $A$ is Fredholm on \textbf{one} of the spaces $Y$;
\end{itemize}

\noindent and on every space $Y$ it holds that
\begin{equation} \label{eq_specess_new}
\spess(A)\ =\ \bigcup_{B\in\opsp(A)} \sp(B)\ =\
\bigcup_{B\in\opsp(A)} \spp^\infty(B).
\end{equation}
\end{Corollary}

Here we denote by $\spp^\infty(B)$\index{$\pspecp$@$\spp^p(A)$} the
point spectrum\index{point spectrum}\index{spectrum!point} (set of
eigenvalues) of $B$, considered as an operator on $Y^\infty$.

\begin{Corollary} \label{cor7.25b}
Suppose $N=1$, $A\in\W$ and $U$ is finite-dimensional. Then
statements (a)--(g) of Corollary \ref{cor_oned_wie} are equivalent.
Moreover, if $A$ is subject to all these equivalent statements then
the index of $A$ is the same on each space $Y$. Further, on every
space $Y$, \refeq{eq_specess_new} holds.
\end{Corollary}

\chapter{Discrete Schr\"odinger Operators} \label{sec_Schroed} In this chapter
we illustrate the results of Chapter \ref{sec_ellp}, in particular
the results of Section \ref{subsection_wiener}, in the relatively
simple but practically relevant setting of
$Y=Y^p=Y^p(U)=\ell^p(\Z^N,U)$ with $p\in\{0\}\cup[1,\infty]$ and a
finite-dimensional space $U$. For applications to a class of
operators on $Y^p(U)$ with $U$ infinite-dimen\-sional, see Chapter
\ref{sec_intop}.

In this chapter we suppose that our operator $A$ is a \emph{discrete
Schr\"o\-dinger operator}\index{Schr\"odinger
operator}\index{operator!Schr\"odinger} on $Y$ in the sense e.g. of
\cite{Davies}. By this we mean that $A$ is of the form
\[
A=L+M_b
\]
with a translation invariant operator $L$, i.e.
$V_{-\alpha}LV_\alpha=L$ for all $\alpha\in\Z^N$, and with a
multiplication operator $M_b$, given by (\ref{eq_multop}), with
$b\in Y^\infty(L(U))$. A translation invariant operator $L$ on $Y$
is often referred to as a \emph{Laurent operator}\index{Laurent
operator}\index{operator!Laurent}, and the sequence $b$ is typically
called the \emph{potential}\index{potential} of $A$. The matrix
representation of $L$ is a {\em Laurent matrix}\index{Laurent
matrix} $[L]=[\lambda_{i-j}]_{i,j\in\Z^N}$ with $\lambda_k\in L(U)$
for all $k\in\Z^N$. To be able to apply the results of the previous
subsections we will suppose that $A=L+M_b\in L(Y^p)$, for $1\leq
p\leq \infty$, which is the case if $L\in\W$, i.e. if
\[
\|L\|_\W\ =\
\sum_{k\in\Z^N}\|\lambda_k\|<\infty.
\]

Discrete (or lattice) Schr\"odinger operators are widely studied in
mathematical physics (see e.g.\ \cite[\S XI.14]{ReedSimon3},
\cite{Simon82}, \cite{CarmonaLacroix91}, \cite{Kuchment1991},
\cite{PasturFigotinBook}, \cite{Teschl99}, \cite{AschJoye06}).

Particularly common and classical is the case where the Laurent
operator takes the form
\begin{equation} \label{eqn_schr_class}
L\ =\ \sum_{k=1}^N (V_{e^{(k)}}\,+\,V_{-e^{(k)}}),
\end{equation}
where $e^{(1)},...,e^{(N)}$\index{$\e$@$e^{(k)}$}\index{unit
coordinate vector} are the unit coordinate vectors in $\Z^N$. The
operator $A=L+M_b$ is then a discrete analogue of the second order
differential operator $-\Delta\,+\,M$ where $\Delta$ is the
Laplacian and $M$ is the operator of multiplication by a bounded
potential, both on $\R^N$.

Let $L^\prime$ be the Laurent operator with matrix representation
$[L^\prime]=[\lambda^*_{j-i}]_{i,j\in\Z^N}$. Then, identifying $U$
with $U^*$ (so that $U=U^*=U^\pa$), $L^\prime\in \W$ and $A^\prime :
= L^\prime + M_b\in \W$. Further, in the case $p=\infty$, when we
consider $A$ as an operator on $Y=Y^\infty$,
$A^\prime$\index{$A^\prime$, transpose of $A$}, considered as an
operator on $Y^1$, is the unique transpose\index{transpose
operator}\index{operator!transpose} of $A$ with respect to the dual
system $(Y^\infty,Y^1)$ of Section \ref{sec_dsa} and so the unique
preadjoint of $A$, for it holds for $x\in Y$, $y\in Y^1$, using
equations (\ref{eq_blff}), (\ref{eq_bildef}), and
(\ref{eqn_matrix}), that
$$
Ax(y) = (Ax,y)=(x,A^\prime y)=x(A^\prime y).
$$
We will say that $A$ is {\em  symmetric}\index{symmetric
operator}\index{operator!symmetric} if $A=A^\prime$. For example,
this is the case when $L$ is the classical operator
(\ref{eqn_schr_class}).

For $b\in Y^\infty(L(U))$ let $\Lim(b)$\index{$Lim(b)$@$\Lim(b)$}
denote the set of {\em limit functions}\index{limit function} of
$b$, by which we mean the set of all functions $b_h\in
Y^\infty(L(U))$ for which there exists a sequence $h:\N\to\Z^N$
tending to infinity such that
\begin{equation} \label{eqn_bhdef}
b_h(m) = \lim_{n\to\infty} b(m+h(n)), \quad m\in\Z^N.
\end{equation}
It follows from (\ref{eqn_Pconfmult}) that
$$
\opsp(A) = \{L+M_c:c\in \Lim(b)\}.
$$
Noting that Corollary \ref{cor7.25} applies to $A=L+M_b$, we have
the following result. In this result, for an operator $B\in\W$ we
denote by $\sp^p(B)$\index{$\pspa$@$\sp^p(A)$},
$\spess^p(B)$\index{$\pspea$@$\spess^p(A)$}, and
$\spp^p(B)$\index{$\pspecp$@$\spp^p(A)$}, respectively, the
spectrum\index{spectrum}, essential spectrum\index{essential
spectrum}\index{spectrum!essential}, and point spectrum\index{point
spectrum}\index{spectrum!point} (set of eigenvalues) of $B$
considered as an operator on $Y^p$.

\begin{Theorem} \label{prop_schr_first}
The following statements are equivalent:
\begin{itemize}
\item[(a)] $L+M_c$ is injective on $Y^\infty$ for all $c\in\Lim(b)$ and, for
some
$s$-dense subset, $\varsigma$, of $\Lim(b)$, $L^\prime+M_c$ is
injective on $Y^1$ for all $c\in\varsigma$.
\item[(b)] $L+M_c$ is invertible, for every $c\in\Lim(b)$, on one of the spaces
$Y^p$;
\item[(c)] $L+M_c$ is invertible on $Y^p$ for every $p$ and every
$c\in \Lim(b)$ and the inverses are uniformly bounded (in $p$ and
$c$);
\item[(d)] $A$ is Fredholm on one of the spaces $Y^p$;
\item[(e)] $A$ is Fredholm on all of the spaces $Y^p$ and the
index is the same on each space.
\end{itemize}
Thus, for every $p$ it holds that
\begin{eqnarray}
\label{eq_spess2a}\spess^p(A)&=& \bigcup_{c\in\Lim(b)} \sp^p(L+M_c)\\
\label{eq_spess2b}&=& \bigcup_{c\in\Lim(b)} [\spp^\infty(L+M_c)\
\cup\ \spp^1(L^\prime+M_c)]
\end{eqnarray}
and
\begin{eqnarray}
\label{eq_spess3a}\sp^p(A)&=&\spp^\infty(A)\ \cup\ \spp^1(A^\prime)\
\cup\
\spess^p(A)\\
\label{eq_spess3b}&=&\spp^1(A)\ \cup\ \spp^1(A^\prime)\ \cup\
\spess^p(A).
\end{eqnarray}
\end{Theorem}
\Proof From the equivalence of (a) and (b) in Corollary
\ref{cor_combine} we have that (a) is equivalent to the statement
that $\alpha(L+M_c)=0$ for all $c\in \Lim(b)$ and $\beta(L+M_c)=0$
for all $c\in\varsigma$. By Corollary \ref{cor7.25} this is
equivalent to (b)--(e). So it remains to prove
\refeq{eq_spess3a}+\refeq{eq_spess3b}. Equality \refeq{eq_spess3a}
follows since, as noted after Lemma \ref{Wiener1}, the spectrum of
$A$ does not depend on $p$, so that $\sp^p(A)=\sp^\infty(A)$.
Further, if $\lambda\in\sp^\infty(A)$ and $\lambda I-A$ is Fredholm,
then either $\alpha(\lambda I-A)\ne 0$ or $\beta(\lambda I-A)\ne 0$,
so that either $\lambda\in\spp^\infty(A)$ or
$\lambda\in\spp^1(A^\prime)$. To see equality \refeq{eq_spess3b},
note that $\spp^1(A)\subset\spp^\infty(A)$ since injectivity of
$\lambda I-A$ on $Y^\infty$ implies its injectivity on $Y^1\subset
Y^\infty$. Moreover, if $\lambda I-A$ is Fredholm then, by Theorem
\ref{Prop_Fred_equiv}, the kernel of $\lambda I - A$ is a subset of
$Y^0$. Since $(\lambda I - A)|_{Y^0}$ is Fredholm with the same
index on $Y^0$ and $Y^1\subset Y^0$ and since $Y^1$ is dense in
$Y^0$, it follows from a standard result on Fredholm operators
(e.g.\ \cite{ProssdorfSilbermann}) that the kernel of $\lambda I -
A$ is a subset of $Y^1$. Thus $\spp^\infty(A)\subset
\spess^\infty(A)\cup \spp^1(A)$.
 \Proofend

\begin{Remark}
We note that main parts of the above result, namely equality
(\ref{eq_spess2a}) and that the spectrum and essential spectrum do
not depend on $p\in [1,\infty]$, are well known (see e.g.\
\cite[Theorem 5.8.1]{RaRoSiBook}). The characterisation of the
essential spectrum by (\ref{eq_spess2b}) appears to be new.
\end{Remark}

Clearly, equations (\ref{eq_spess2a}) -- (\ref{eq_spess3b}) simplify
when $L$ is symmetric, for example if $L$ is given by
(\ref{eqn_schr_class}), since we then have that
$\spp^1(A^\prime)=\spp^1(A)$ and $\spp^1(L^\prime+M_c)\subset
\spp^\infty(L+M_c)$ for all $c\in \Lim(b)$. Simplifications also
occur when the potential $b$ is almost periodic, $b\in\AP(L(U))$, in
which case $\Lim(b)$ is precisely what is often called the {\em
hull}\index{hull} of $b$, the set $\clos_{Y^\infty(L(U))}\{V_k
b:k\in\Z^N\}$, the closure of the set of translates of $b$.

\begin{Theorem} \label{prop_SchroedAP}
If $b$ is almost periodic then, for all $p$ and all $\tilde b\in
\Lim(b)$,
\begin{equation} \label{eq_SchroeAP}
\spess^p(A)\ =\ \sp^p(A)\ =\ \sp^p(L+M_{\tilde b})\ =\
\bigcup_{c\in\Lim(b)} \spp^\infty(L+M_c).
\end{equation}
\end{Theorem}
\Proof $L$ is absolutely rich/periodic and so norm rich. Since $b$
is almost periodic, $M_b$ is norm rich by Lemma \ref{lem_BONR}. Thus
$A$ is norm rich. Further, $\opsp(A-I)$ is uniformly bounded by
Theorem \ref{lem_opsp} (i) and so uniformly Montel on $Y^\infty$ by
Corollary \ref{cor_UMUB}, since $\dim U <\infty$. The result thus
follows from Theorem \ref{prop_NRBDO}, Theorem \ref{prop_NR1} (iv),
and the equivalence of statements (b) and (d) in Theorem
\ref{prop_schr_first}. \Proofend

\begin{Remark} \label{rem7.4}
That $\spess^p(A) = \sp^p(A) = \sp^p(L+M_{\tilde b})$ for all
$\tilde b\in \Lim(b)$, the hull of $b$, is a classical result, see
e.g.\ \cite{Simon82,Teschl99,RaRoSiBook}. The result that
\[
\sp^p(A)\ =\ \bigcup_{c\in\Lim(b)} \spp^\infty(L+M_c)
\]
appears to be new in this generality. However, Avila and
Jitomirskaya have a comparable result for the case when $A$ is
self-adjoint (i.e. $L=L'$ and the potential $b$ is real-valued) on
$Y=\ell^2(\Z,\C)$ in their very recent paper \cite{AvilaJito} (see
Theorem 3.3 there). Moreover, analogous results for uniformly
elliptic differential operators on $\R^N$ with almost periodic
coefficients date back to Shubin \cite{Shubin78}.

Moreover, note that this result is well-known, as a part of
Floquet-Bloch theory \cite{Kuchment1991,Kuchment,Davies2007:Book},
in the case when $b$ is periodic; in fact one has the stronger
result in that case, at least when $L$ is given by
(\ref{eqn_schr_class}), that $\lambda$ is in the spectrum of $A$ iff
there exists a solution $x\in Y^\infty(U)$ of $\lambda x= Ax$ which
is quasi-periodic in the sense of \cite{Kuchment1991}. The latter
means that $x(m)=\exp(i k\cdot m) y(m)$ for all $m\in\Z^N$, where
$y\in Y^\infty(U)$ is periodic and $k\in\R^N$ is fixed, so that if
$x$ is quasi-periodic then it is certainly  almost periodic. Thus,
if $b$ is periodic then $\lambda$ is in the spectrum of $A$ iff
there exists a solution $x\in Y^\infty(U)$ of $\lambda x= Ax$ which
is almost periodic.

Natural questions are whether this statement still holds for the
case when $b$ is almost periodic, at least for $L$ given by
(\ref{eqn_schr_class}), or whether the weaker statement holds that
$\lambda$ is in the spectrum of $A$ iff, for some $c\in \Lim(b)$,
there exists an almost periodic solution $x\in Y^\infty(U)$ of
$\lambda x= (L+M_c)x$. The answer is, to our knowledge, unknown.
\end{Remark}

To illustrate the application of the above theorem in the 1D case
($N=1$) we consider a widely studied class of almost periodic
operators obtained by the following construction. For some $d\in\N$
let $B:\R^d\to L(U)$ be a continuous function satisfying
$$
B(s+m) = B(s), \quad s\in\R^d,\; m\in \Z^d.
$$
Let $\alpha = (\alpha_1,\dots,\alpha_d)\in\R^d$ and, for $s\in\R^d$
let $b_s:\Z\to L(U)$ be given by
\begin{equation} \label{eq_bs_def}
b_s(n) = B(\alpha n+s),\quad n\in\Z.
\end{equation}
If $\alpha_1,\dots ,\alpha_d$ are all rational, then $b_s$ is
periodic. Whatever the choice of $\alpha_1,\dots,\alpha_d$, $b_s$ is
almost periodic ($b_s\in\AP(L(U))$).

For $s\in\R^d$ let $[s]$ denote the coset $[s] = s+\Z^d$ in
$\R^d/\Z^d$. An interesting case is that in which
$1,\alpha_1,\alpha_2,\dots,\alpha_d$ are rationally independent, in
which case $\{[\alpha m]:m\in\Z\}$ is dense in $\R^d/\Z^d$. Then it
is a straightforward calculation to see that
\begin{equation} \label{eq_Limbs}
\Lim(b_s) = \{b_t:t\in\R^d\}.
\end{equation}
Thus, for this case, (\ref{eq_SchroeAP}) reads as
\begin{equation} \label{eq_SchroeAP2}
\spess^p(L+M_{b_s})\ =\ \sp^p(L+M_{b_s})\ = \bigcup_{t\in\R^d}
\spp^\infty(L+M_{b_t}).
\end{equation}

As a particular instance, this formula holds in the case when
$U=\C$, $d=1$, and $B(s) = \lambda\cos(2\pi s)$, $s\in\R$, for some
$\lambda\in\C$. Then
\begin{equation} \label{eq_mathieu}
b_s(n) = \lambda \cos(2\pi(\alpha n+s)), \quad n\in\Z,
\end{equation}
and (\ref{eq_SchroeAP2}) holds if $\alpha$ is irrational, in which
case $b_s$ is the so-called {\em almost Mathieu}\index{almost
Mathieu operator}\index{operator!almost Mathieu} potential.

We next modify the above example to illustrate the application of
Theorem \ref{prop_schr_first} in a particular 1D ($N=1$) case.

\begin{Example} \label{ex:LastSimon}
Define $b_s\in \AP(L(U))$ by (\ref{eq_bs_def}) and suppose that
$1,\alpha_1,\alpha_2,\dots,\alpha_d$ are rationally independent.
Suppose that $f:\Z\to\R^d$ satisfies
$$
\lim_{|n|\to\infty} |f(n+1)-f(n)|\to 0.
$$
Define $b\in Y^\infty(L(U))$ by
$$
b(n) = B(\alpha n + f(n)), \quad n\in\Z.
$$
Then it is straightforward to see that $\Lim(b) \subset
\{b_s:s\in\R^d\}$. Since \cite[Corollary 3.97]{LiBook}, $b_s\in
\Lim(b)$ implies that $\Lim(b_s)\subset \Lim(b)$, we have, by
(\ref{eq_Limbs}), that
$$
\Lim(b) = \{b_s:s\in\R^d\}.
$$
Thus, applying Theorem \ref{prop_schr_first} and
(\ref{eq_SchroeAP2}), we see that, for every $s\in\R^d$ and every
$p\in\{0\}\cup [1,\infty]$,
\begin{equation} \label{eq_SchroeAP3}
\spess^p(L+M_b) = \spess^p(L+M_{b_s})\ =\ \sp^p(L+M_{b_s})\ =
\bigcup_{t\in\R^d} \spp^\infty(L+M_{b_t}).
\end{equation}
We note that, in the special case that $L$ is given by
(\ref{eqn_schr_class}) (with $N=1$), $U=\C$, and $B$ is real-valued,
the statement that
$$
\spess^2(L+M_b) =\ \sp^2(L+M_{b_s})
$$
for all $s\in\R^d$ is Theorem 5.2 of Last and Simon
\cite{LastSimon06} (established by limit operator type arguments).
As a specific instance where (\ref{eq_SchroeAP3}) holds, let us take
$U=\C$, $d=1$, and $B(s) = \lambda\cos(2\pi s)$, $s\in\R$, for some
$\lambda\in\C$. Then $b_s$ is given by (\ref{eq_mathieu}) and,
taking (as one possible choice), $f(n) = |n|^{1/2}$, one has
$$
b(n) = \lambda \cos(2\pi(\alpha n+ |n|^{1/2}))
$$
(cf.\ \cite[Theorem 1.3]{LastSimon06}).
\end{Example}

As a further example we consider the case when $b$ is pseudo-ergodic
in the sense of Davies \cite{Davies}. Following \cite{Davies}, we
call the function $b\in Y^\infty(\Z^N,U)$ {\em
pseudo-ergodic}\index{pseudo-ergodic sequence}, if, for every
$\eps>0$, every finite set $S\subset\Z^N$ and every function $f:S\to
\Sigma:=\clos_U\,b(\Z^N)$, there is a $z\in\Z^N$ such that
$$
\|f(s)-b(z+s)\|_U<\eps,\qquad s\in S.
$$
One can show \cite[Corollary 3.70]{LiBook} that $b$ is
pseudo-ergodic iff $\Lim(b)$ is the set $\Sigma^{\Z^N}$ of all
functions $c:\Z^N\to \Sigma$. In particular, $b\in \Lim(b)$ if $b$
is pseudo-ergodic.

\begin{Theorem} \label{prop_pseudoerg}
If $b$ is pseudo-ergodic then, for all $p$,
\[
\spess^p(A)\ =\ \sp^p(A)\ =\ \bigcup_{c\in \Sigma^{\Z^N}}
\sp^p(L+M_c)\ =\ \bigcup_{c\in \Sigma^{\Z^N}} \spp^\infty(L+M_c).
\]
\end{Theorem}
\Proof The first two `$=$' signs follow from \refeq{eq_spess2a} and
the fact that $b\in\Lim(b)=\Sigma^{\Z^N}$. For the proof of the
remaining equality, we refer to the following $s$-dense subset of
$\Lim(b)=\Sigma^{\Z^N}$: Let $m_1=m_2=...=m_N=1$, and let
$\varsigma$ stand for the set of all periodic functions $x:\Z^N\to
\Sigma$, that is
\[
\varsigma\ :=\ \bigcup_{n\in\N}Y^\infty_n(\Sigma)
\]
with $Y^\infty_n(\Sigma)$ defined as in \refeq{eq_Yn} (with the
slight abuse of notation by writing $Y^\infty(\Sigma)$ for
$\Sigma^{\Z^N}$, i.e. the set of all functions $x:\Z^N\to \Sigma$).
Then $\varsigma$ is $s$-dense in $\Sigma^{\Z^N}$ as every $x\in
\Sigma^{\Z^N}$ can be strictly approximated by the sequence $(\tilde
P_nx)\subset\varsigma$ with $\tilde P_n$ as defined in
\refeq{eq_PnYn}. If $\lambda\in\C$ and all limit operators $\lambda
I-(L+M_c)$ of $\lambda I-A=\lambda I-(L+M_b)$, including those with
$c\in\varsigma$, are injective, then, by Theorem
\ref{lem_A-inj-inv}, we have that $\lambda I-(L+M_c)$ is surjective
for every $c\in\varsigma$. By the equivalence between (a) and (d) in
Theorem \ref{prop_schr_first}, this shows that $\lambda I-A=\lambda
I-(L+M_b)$ is Fredholm. \Proofend

\begin{Remark}
It is shown that
\[
\spess^2(A)\ =\ \sp^2(A)\ =\ \bigcup_{c\in \Sigma^{\Z^N}}
\sp^2(L+M_c)
\]
in \cite{Davies}. The result that $\sp^p(A)=\bigcup_{c\in
\Sigma^{\Z^N}} \spp^\infty(L+M_c)$ appears to be new.
\end{Remark}

The above theorems show that, in each of the cases $L$ symmetric,
$b$ almost periodic, and $b$ pseudo-ergodic, it holds that
\begin{equation} \label{eqn_conj}
\spess^p(A)\ =\ \bigcup_{c\in \Lim(c)} \spp^\infty(L+M_c).
\end{equation}
We conjecture that, in fact, this equation holds for all $c\in
Y^\infty(L(U))$. For $N=1$ this is no longer a conjecture, as we
showed in Corollary \ref{cor7.25b} (which follows from our more
general results in \cite{CWLi:JFA2007}, also see Theorem
\ref{prop_FC_JFA} above). For $N\ge 2$ however, this is an open
problem.

We finish this chapter with an example demonstrating how Theorem
\ref{prop_pseudoerg} can be used to compute spectra of Schr\"odinger
operators with random potential\index{random potential} $b$.

\begin{Example} \label{ex7.8}
Let $N=1$, $p\in[1,\infty]$, $U=\C$ and take a compact set $\Sigma$
in the complex plane. We compute the spectrum of $A=L+M_b$ as an
operator on $Y=Y^p(U)=\ell^p(\Z)$ where $L=V_{-1}$ is the backward
shift and the function values $b(k)$, $k\in\Z$, of the random
potential $b$ are chosen independently of each other from the set
$\Sigma$. We assume that, for every $\sigma\in\Sigma$ and $\eps>0$,
$\mathbb{P}(|b(k)-\sigma|<\eps)>0$. Then it is easy to see (the
argument is sometimes called `the Infinite Monkey Theorem' and it
follows from the Second Borel Cantelli Lemma, see \cite[Theorem
8.16]{CapinskiKopp} or \cite[Theorem 4.2.4]{Chung}) that, with
probability 1, $b$ is pseudo-ergodic.

For the calculation of the point spectra in Theorem
\ref{prop_pseudoerg}, let $c\in\Sigma^\Z$ and take $\lambda\in\C$.
If $x:\Z\to\C$ is a nontrivial solution of $(L+M_c)x=\lambda x$ then
$x(n_0)\ne 0$ for some $n_0\in\Z$ and $x(k+1)=(\lambda-c(k))\,x(k)$
for all $k\in\Z$. Note that $\lambda\ne c(k)$ for all $k<n_0$ since
otherwise $x(n_0)=0$ and, w.l.o.g., suppose that $x(n_0)=1$. As a
consequence we get that
\begin{equation} \label{eq_randprod}
x(n)\ =\ \left\{\begin{array}{cl} \prod\limits_{k=n_0}^{n-1}
(\lambda-c(k)),& n\ge n_0,\\
\prod\limits_{k=n}^{n_0-1}(\lambda-c(k))^{-1},&n<n_0
\end{array}\right.
\end{equation}
for every $n\in\Z$. Now put, for $r>0$,
\[
\Sigma_\cup^r\ :=\
\bigcup_{\sigma\in\Sigma}(\sigma+r\overline{\mathbb
D})\qquad\textrm{and}\qquad \Sigma_\cap^r\ :=\
\bigcap_{\sigma\in\Sigma}(\sigma+r{\mathbb D})
\]
with ${\mathbb D}$ denoting the open unit disk in $\C$ and
$\overline{\mathbb D}$ its closure.

Clearly, if $\lambda\not\in\Sigma_\cup^1$ then $|\lambda-\sigma|>1$
for all $\sigma\in\Sigma$ and hence, for every nontrivial
eigenfunction $x$ of $L+M_c$, we have that $|x(n)|\to\infty$ in
(\ref{eq_randprod}) as $n\to+\infty$, regardless of $c:\Z\to\Sigma$.

Similarly, if $\lambda\in\Sigma_\cap^1$ then $|\lambda-\sigma|<1$
for all $\sigma\in\Sigma$ and hence, for every nontrivial
eigenfunction $x$ of $L+M_c$, $|x(n)|\to\infty$ in
(\ref{eq_randprod}) as $n\to-\infty$, regardless of $c:\Z\to\Sigma$.
(Note that $n_0$ in (\ref{eq_randprod}) depends on $c$ and
$\lambda$.)

So in both cases, $(L+M_c)x=\lambda x$ has no nontrivial solution
$x\in Y^\infty$, so $\lambda\not\in\spp^\infty(L+M_c)$ for
all $c\in\Sigma^\Z$.
Now it remains to look at
$\lambda\in\Sigma_\cup^1\setminus\Sigma_\cap^1$. In this case, let
$\sigma,\tau\in\Sigma$ be such that $|\lambda-\sigma|\le
1\le|\lambda-\tau|$, which is possible by the choice of $\lambda$,
and put $c(k):=\tau$ for $k<0$ and $c(k):=\sigma$ for $k\ge 0$. Then
$c\in\Sigma^\Z$ and
\[
x\ =\
\big(\cdots\,,\,(\lambda-\tau)^{-2}\,,\,(\lambda-\tau)^{-2}\,,\,1\,,\,(\lambda-\sigma)^1\,,\,(\lambda-\sigma)^2\,,\,\cdots\big)\
\in\ Y^\infty,
\]
with the $1$ at position $n_0=0$, is an eigenvector of $L+M_c$
w.r.t. $\lambda$.

Summarising and using Theorem \ref{prop_pseudoerg}, we get that,
with probability 1,
\begin{equation} \label{eq_rand}
\sp^p A\ =\ \spess^p A\ = \bigcup_{c\in\Sigma^\Z}\spp^\infty(L+M_c)\
=\ \Sigma_\cup^1\setminus\Sigma_\cap^1,
\end{equation}
which confirms, in a simpler and more straightforward way, a result
of Trefethen, Contedini and Embree \cite[Theorem 8.1]{TrefContEmb}
(and see \cite[Section VIII]{TrefEmbBook}). Equation \eqref{eq_rand}
is illustrated, for two particular cases, in Figure \ref{fig:spec1}.
\begin{figure}[h]
\begin{center}
\includegraphics[width=0.9\textwidth]{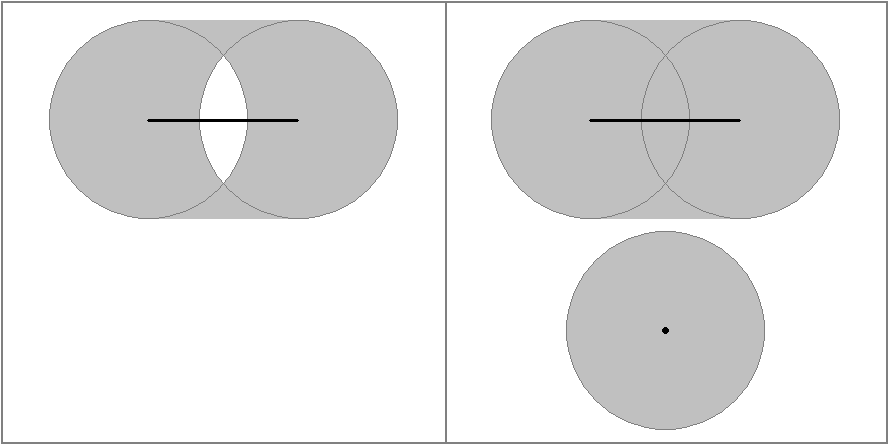}
\caption{\footnotesize The left image shows, as a gray shaded area,
$\sp^p A$ when $\Sigma$ is the black straight line of length $1.5$.
In the right image, one more point (the centre of the lower circle)
has been added to $\Sigma$ which results in
$\Sigma_\cap^1=\varnothing$.} \label{fig:spec1}
\end{center}
\end{figure}
By a simple scaling argument, it moreover follows from
(\ref{eq_rand}) that, with probability 1,
$$
\sp^p(\eps L+M_b)\ =\ \spess^p(\eps L+M_b)\ =\
\Sigma_\cup^\eps\setminus\Sigma_\cap^\eps.
$$
Thus
$$
\sp^p(\eps L+M_b)\ =\ \Sigma_\cup^\eps=\ \sp_\eps(M_b)
$$
for all $\eps>0$ small enough that $\Sigma_\cap^\eps=\varnothing$,
e.g. for $\eps\in (0,$ diam $\Sigma/2)$.
\end{Example}


\chapter{A Class of Integral Operators} \label{sec_intop}

In this chapter, we apply the results of Chapter \ref{sec_ellp} to
study a class of operators on
\[
Y\ =\ Y^p(U)\ =\ \left\{\begin{array}{cl}
\ell^p(\Z^N,U),&p\in[1,\infty],\\
c_0(\Z^N,U),&p=0
\end{array}\right.
\]
with $U=L^q([0,1]^N)$ for some fixed $q\in (1,\infty]$. In a natural
way (cf. \cite[\S 1.6.3]{Kurbatov}) we identify elements $x\in Y$
with equivalence classes of scalar-valued functions on $\R^N$ via
\begin{equation} \label{eq_Lpq}
\big(x(m)\big)(t)\ =\ f(m+t),\qquad m\in\Z^N,\ t\in[0,1]^N,
\end{equation}
and denote the set of all of these (equivalence classes of)
functions $f$ with
\[
\|f\|_{p,q}\ :=\ \|x\|_Y\ <\ \infty
\]
by $L^{p,q}(\R^N)$\index{$L^{p,q}(\R^N)$} or just $L^{p,q}$. Note
that $L^{q,q}(\R^N)=L^q(\R^N)$. Equipped with the norm
$\|.\|_{p,q}$, $L^{p,q}$ is a Banach space, and \refeq{eq_Lpq}
yields an isometric isomorphism between $L^{p,q}$ and $Y^p(U)$. We
will freely identify these two spaces and the notions of strict
convergence, limit operators, as well as the operators $P_m, V_k\in
L(Y^p(U))$ with the corresponding notions and operators on $L^{p,q}$
(cf. \cite[(1.3)]{LiBook}).

The operators we are going to study on $Y$ alias $L^{p,q}$ are
composed, via addition and composition, from two basic ingredients:
\begin{itemize}
\item For $b\in L^\infty:=L^\infty(\R^N)$, define the
{\it multiplication operator}\index{multiplication
operator}\index{operator!multiplication} $\M_b\in
L(L^{p,q})$\index{$M_ba$@$\M_b$} by
\[
\M_b f(t)\ =\ b(t)\,f(t),\qquad t\in\R^N
\]
for all $f\in L^{p,q}$. Via the identification \refeq{eq_Lpq}
between $L^{p,q}$ and $Y$, we can identify $\M_b$ with the
multiplication operator $M_c$ on $Y$ as defined in
\refeq{eq_multop}, where $c\in\ell^\infty(\Z^N,L(U))$ is such that
$(c(m)u)(t)=b(m+t)u(t)$ for all $m\in\Z^N$, $u\in U=L^q([0,1]^N)$
and $t\in [0,1]^N$. Recall from Example \ref{ex_almostperiodic} that
$\M_b$ is rich iff $\{b(\cdot+k)\}_{k\in\Z^N}$ is relatively
sequentially compact in the strict topology on $L^\infty$, in which
case we write $b\in\rich$. It is easy to check that $\rich$ is an
inverse closed Banach subalgebra of $L^\infty$.\\
\item For $\kappa\in L^1:=L^1(\R^N)$, define the {\it convolution
operator}\index{convolution operator}\index{operator!convolution}
$\Co_\kappa\in L(L^{p,q})$\index{$C_k$@$\Co_\kappa$} by
\[
\Co_\kappa f(t)\ =\ (\kappa\star f)(t)\ =\
\int_{\R^N}\kappa(s-t)\,f(s)\,ds,\qquad t\in\R^N
\]
for all $f\in L^{p,q}$ (cf. Example \ref{ex_conv}). As demonstrated
in \cite[Example 1.28]{LiBook}, the convolution operator
$\Co_\kappa$ on $L^{p,q}$ corresponds to a Laurent operator $L$ on
$Y$, where every entry $\lambda_k$ of
$[L]=[\lambda_{i-j}]_{i,j\in\Z^N}$ is the operator of convolution by
$\kappa(\cdot+k)$ on $U$. By Young's inequality, we get that
\[
\|\lambda_k\|\ \le\ \|\kappa|_{k+[-1,1]^N}\|_1
\]
for every $k\in\Z^N$, and by $\kappa\in L^1(\R^N)$ it follows that
$L\in\W=\W(U)$.
\end{itemize}

We denote by $\cA^o$\index{$A0a$@$\cA^o$} the smallest algebra in
$L(L^{p,q})$ containing all operators of these two types; that is
the set of all finite sum-products of operators of the form $\M_b$
and $\Co_\kappa$ with $b\in\rich$ and $\kappa\in L^1$. From the
above considerations it follows that every operator $A\in\cA^o$, if
identified with an operator on $Y$, is contained in the Wiener
algebra $\W$. By $\cA$\index{$A$@$\cA$} we denote the closure of
$\cA^o$ in the norm $\|.\|_\W$. Note that, by \refeq{eq_Wienernorm},
the closure of a set $S\subset\W$ in the $\W$-norm is always
contained in the closure of $S$ in the usual operator norm.

\begin{Lemma} \label{lem_preadjointexists}
The predual space $U^\pa$ exists and, if $p=\infty$, then every
$A\in\cA$ has a preadjoint operator $A^\pa$ on $L^{1,q'}$ with
$1/q+1/q'=1$.
\end{Lemma}
\Proof By the choice $q\in (1,\infty]$, it is clear that the predual
space $U^\pa$ of $U=L^q([0,1]^N)$ exists and can be identified with
$L^{q'}([0,1]^N)$, where $1/q+1/q'=1$, including the case $q'=1$ if
$q=\infty$. Now suppose $p=\infty$. Then the predual $Y^1(U^\pa)$ of
$Y=Y^\infty(U)$ exists and corresponds to $L^{1,q'}$ in the sense of
\refeq{eq_Lpq}. By \cite[Proposition 1.10]{LiBook} and the fact that
both multiplication and convolution operators have a preadjoint
operator (indeed, $\M_b^\pa=\M_b$ and
$\Co_\kappa^\pa=\Co_{\kappa(-\cdot)}$ for all $b\in L^\infty$ and
$\kappa\in L^1$), we see that indeed $A^\pa$ exists for every
$A\in\cA$. \Proofend

Now let
\[
\J^o\ :=\ \left\{\sum A_i\Co_{\kappa_i}B_i\ :\ A_i,B_i\in\cA^o,\
\kappa_i\in L^1\right\},\index{$J0$@$\J^o$}
\]
with the sum being finite, denote the smallest two-sided ideal of
$\cA^o$ containing all convolution operators $\Co_\kappa$ with
$\kappa\in L^1$, and let $\J$\index{$J$@$\J$} be its closure in the
norm $\|.\|_\W$, hence the smallest $\W$-closed two-sided ideal of
$\cA$ containing all $\Co_\kappa$.

\begin{Lemma} \label{lem_JisMontel}
It holds that $\J\subset\UM(L^{p,q})$. In particular, every operator
in $\J$ is Montel.
\end{Lemma}
\Proof The inclusion $\J\subset\UM(L^{p,q})$ follows from Lemma
\ref{lem_7.23isBA} (a),(c),(d) and the fact that
$\Co_\kappa\in\UM(L^{p,q})$ for all $\kappa\in L^1$ since the set
$\{V_{-k}KV_k:k\in\Z^N\}$ in Lemma \ref{lem_coll_comp_entries} (ii)
is just a singleton if $K=\Co_\kappa$. \Proofend

It can be shown that, in the same way as every $A\in\cA^o$ clearly
can be written as the sum of a multiplication operator and an
operator in $\J^o$, also every $A\in\cA$ can be uniquely written as
\begin{equation} \label{eq_decompA}
A\ =\ \M_b+K\qquad\textrm{with}\qquad
b\in\rich\quad\textrm{and}\quad K\in\J.
\end{equation}
This follows from \cite[Proposition 4.11]{LiBook} with $\cA$ and
$\J$ there replaced by the current meaning. As a consequence, we get
that the factor algebra $\cA/\J$ is isomorphic to $\rich$, and the
coset $A+\J$ of $A\in\cA$ is represented by the function $b\in\rich$
from \refeq{eq_decompA}.

\begin{Theorem} \label{prop_AFredholm1}
The operator \refeq{eq_decompA} is Fredholm iff it is invertible at
infinity and $b$ is invertible in $L^\infty$.
\end{Theorem}
\Proof If $A\in\cA$ is invertible at infinity and $b$ from its
representation \refeq{eq_decompA} is invertible in $L^\infty$, then
$A$ is Fredholm by Theorem \ref{prop_rich_invinf} (ii) and Lemma
\ref{lem_JisMontel}. Conversely, let $A\in\cA$ be Fredholm. By
Theorem \ref{prop_rich_invinf} (i) and b), together with Lemma
\ref{lem_preadjointexists}, we get that $A$ is invertible at
infinity. It remains to show that $b$ from \refeq{eq_decompA} is
invertible in $L^\infty$. To see this, take $B\in L(L^{p,q})$ and
$S,T\in K(L^{p,q})$ such that $AB=I+S$ and $BA=I+T$. Then, for every
$k\in\N$, we get that
\[
P_k\M_bBP_k+P_kKBP_k\ =\ P_kABP_k\ =\ P_k+P_kSP_k,
\]
and hence
\[
(P_k\M_bP_k)(P_kBP_k)\ =\ P_k\ +\ S'
\]
with $S'=P_kSP_k-P_kKBP_k\in K(L^{p,q}([-k,k]^N))$ by Lemmas
\ref{lem_JisMontel} and \ref{lem_charact_mon}. From the last
equality and its symmetric counter-part, we conclude that
$\M_{P_kb}=P_k\M_bP_k$ is Fredholm on $L^{p,q}([-k,k]^N)$, implying
that the function $P_kb$ is invertible in $L^\infty([-k,k]^N)$, by a
standard argument (see e.g. \cite[Lemma 2.42]{LiBook}). Since this
holds for every $k\in\N$, we get that $b$ is invertible in
$L^\infty$. \Proofend

As invertibility of $b$ in $L^\infty$ turned out to be necessary for
Fredholmness of \refeq{eq_decompA}, we will now, without loss of
generality, suppose that $b$ is invertible, therefore write
\[
A\ =\ \M_b\ +\ K\ =\ \M_b(I+K')\qquad\textrm{with}\qquad
K'=\M_{b^{-1}}K\in\J,
\]
and then merely study Fredholmness of $I+K'$. For this setting we
can show the analogous result of Theorem \ref{prop_schr_first}.

\begin{Theorem} \label{prop_AFredholm2}
For $A=I+K$ with $K\in\J$, the following statements are equivalent.
\begin{itemize}
\item[(a)] All limit operators of $A$ are injective on $L^{\infty,q}$
and $\opsp(A)$ has an $S-$dense subset of injective operators on
$L^{1,q}$;
\item[(b)] All limit operators of $A$ are invertible on one of the
spaces $L^{p,q}$ with $p\in\{0\}\cup[1,\infty]$;
\item[(c)] All limit operators $A_h$ of $A$ are invertible on all the
spaces $L^{p,q}$ with $p\in\{0\}\cup[1,\infty]$ and the inverses are
uniformly bounded (in $p$ and $h$);
\item[(d)] $A$ is Fredholm on one of the spaces $L^{p,q}$ with
$p\in\{0\}\cup[1,\infty]$;
\item[(e)] $A$ is Fredholm on all the spaces $L^{p,q}$ with
$p\in\{0\}\cup[1,\infty]$.
\end{itemize}
Thus, for every $p\in\{0\}\cup[1,\infty]$ it holds that
\begin{eqnarray}
\nonumber
\spess^p(A)&=& \bigcup_{A_h\in\opsp(A)} \sp^p(A_h)\\
\label{eq_spess4b}&=& \bigcup_{A_h\in\opsp(A)} [\spp^\infty(A_h)\
\cup\ \spp^1(A_h)]
\end{eqnarray}
and
\begin{eqnarray}
\nonumber
\sp^p(A)&=&\spp^\infty(A)\ \cup\ \spp^1(A^\prime)\ \cup\
\spess^p(A)\\
\nonumber
&=&\spp^1(A)\ \cup\ \spp^1(A^\prime)\ \cup\ \spess^p(A).
\end{eqnarray}
\end{Theorem}
\Proof We start by showing that $A$ is subject to the conditions in
Corollary \ref{cor_spessformula}. Clearly, $A=I+K$ is contained in
the Wiener algebra and it is rich, by (\ref{eq_limops}) and since
all generators $\M_b$ and $\Co_\kappa$ of $\cA$ are rich. Predual
$U^\pa$ and preadjoint $A^\pa$ exist by Lemma
\ref{lem_preadjointexists}, and $K\in\UM(L^{p,q})$ by Lemma
\ref{lem_JisMontel}.

The rest of this proof proceeds exactly as that of Theorem
\ref{prop_schr_first} with one difference: Unlike Theorem
\ref{prop_schr_first}, which rests on Corollary \ref{cor7.25}, we
here have an infinite-dimensional space $U$ and therefore we use
Corollary \ref{cor_spessformula}. \Proofend

\begin{Example}
The spectra of limit operators of $A=I+K$ can be written down
explicitly when $K\in\J$ is composed of convolution operators
$\Co_\kappa$ with $\kappa\in L^1$ and multiplication operators
$\M_b$ with a {\em slowly oscillating}\index{slowly oscillating
function} function $b\in L^\infty$. By the latter we mean that
\[
\esssup_{t\in [0,1]^N}\,|\,b(x+t)-b(x)\,|\ \to\
0\qquad\textrm{as}\qquad |x|\to\infty.
\]
In this case, the multiplication operator $\M_b$ is rich and all of
its limit operators are multiples of the identity. (In a sense, even
the reverse statement is true \cite{Li04}.) As a consequence, every
limit operator $A_h$ of $A=I+K$ is of the form $I+C_\kappa$ with
some $\kappa\in L^1$, in which case the set $\sp^p(A_h)$ is the
range of the function $1+F\kappa$ with $F$ being the Fourier
transform on $L^1$.
\end{Example}

One can now proceed similarly to Chapter \ref{sec_Schroed} to get
rid of the second injectivity condition in Theorem
\ref{prop_AFredholm2} (a) and the second point-spectrum in formula
\refeq{eq_spess4b}: Clearly, if $N=1$, then Corollary
\ref{cor_oned_wie} does exactly this for us. Otherwise, if $N\ge 2$,
one way to go is to restrict ourselves to symmetric operators in
$\cA$. We will illustrate another way, that is restricting the
generating multiplication operators $\M_b$ of $\cA$ to
norm-rich/almost periodic ones, i.e. to work with a type of almost
periodic functions $b\in L^\infty$.

Let $b\in L^\infty$ and put $c=(b|_{k+[0,1]^N})_{k\in\Z^N}\in
Y^\infty(L^\infty([0,1]^N))$. From Example \ref{ex_almostperiodic}
we know that the following are equivalent:
\begin{itemize}
\item The set $\{V_kb:k\in\Z^N\}$ is relatively compact in $L^\infty$.
\item The set $\{V_kc:k\in\Z^N\}$ is relatively compact in
$Y^\infty(L^\infty([0,1]^N))$.
\item $c$ is almost periodic, i.e. $c\in\AP(L^\infty([0,1]^N))$.
\item $M_c$ is norm-rich/almost periodic (and therefore rich) on $Y^p(U)$.
\end{itemize}
When this is the case, then we will say that $b\in L^\infty$ is {\it
$\Z$-almost periodic}\index{Z-almost@$\Z$-almost periodic function}
and write $b\in\LAP$\index{$L^z$@$\LAP(\R^N)$}. The set of
$\Z$-almost periodic functions on $\R^N$ is not to be confused with
the much smaller subset of almost periodic functions\index{almost
periodic function} on $\R^N$. Unlike almost periodic functions, the
functions in $\LAP$ do not need to be continuous (see \cite[\S
3.4.8]{LiBook} for more detail).

So let $\cAAP\subset\cA$\index{$AAPa$@$\cAAP$} denote the
$\W$-closure of the smallest algebra in $L(L^{p,q})$ that contains
all $\M_b$ with $b\in\LAP$ and all $\Co_\kappa$ with $\kappa\in
L^1$. Analogously, $\JAP$\index{$JAP$@$\JAP$} be the smallest
$\W$-closed two-sided ideal in $\cAAP$ containing all convolution
operators $\Co_\kappa$ with $\kappa\in L^1$. It is not hard to see
that $\JAP=\J\cap\cAAP$.

As seen before for operators in $\cA$, the study of Fredholmness in
$\cAAP$ can be reduced to studying operators of the form $I+K$ with
$K\in\JAP$. In this case, in full analogy to Theorem
\ref{prop_SchroedAP}, we have the following improved version of
formula \refeq{eq_spess4b}:

\begin{Theorem} \label{prop_intopAP}
If $A=I+K$ with $K\in\JAP$ then, for all $p\in\{0\}\cup [1,\infty]$
and all $A_h\in\opsp(A)$,
\[
\spess^p(A)\ =\ \sp^p(A)\ =\ \sp^p(A_h)\ =\ \bigcup_{B\in\opsp(A)}
\spp^\infty(B).
\]
\end{Theorem}
\Proof First of all, $A$ is norm-rich/almost periodic, by Lemma
\ref{lem_NRisBA} and since the generators of $\cAAP$ are
norm-rich/almost periodic. Secondly, $\opsp(K)$ is uniformly Montel,
by Lemmas \ref{lem_JisMontel} and \ref{charact_uniformly_Montel}.
Consequently, we may use Theorems \ref{prop_NRBDO} and
\ref{prop_NR1} (iv), which, together with the equivalence of
statements (b) and (d) in Theorem \ref{prop_AFredholm2}, prove this
formula. \Proofend


\chapter{Some Open Problems} \label{chap_openprobs}
We conclude with a small list of open problems the solutions of
which, we believe, could be crucial in extending the picture that we
have tried to draw in this text. We would like to see this list as
both a future agenda for ourselves as well as an invitation for the
interested reader.
\bigskip

{\it 1. Is a version of the limit operator theory possible with
$L(Y,\P)$ and $K(Y,\P)$ replaced by $S(Y)$ and $SN(Y)$?}\quad We
know from Lemma \ref{lem1}, Corollary \ref{cor4}, (\ref{KYP_def})
and Lemma \ref{LYP_iff} that $S(Y)$ and $SN(Y)$ are ``one-sided
versions'' of $L(Y,\P)$ and $K(Y,\P)$. Moreover, $SN(Y)$ is a closed
two-sided ideal in $S(Y)$ -- just like $K(Y,\P)$ is in $L(Y,\P)$.
The ideal $K(Y,\P)$ shapes the theory presented here in two ways: It
defines the notion of invertibility at infinity (see Definition
\ref{inv_inf}) and that of $\P$-convergence (see Definition
\ref{def_P_conv}). How do these properties change if one works with
the ideal $SN(Y)$ instead and what is the connection between the new
notions of `invertibility at infinity' and `limit operator'?

{\it 2.} A related but rather different question is the following: {\it To
what extent is a version of the limit operator theory possible with what
we termed `weak limit operators' in our discussion of the papers
\cite{Muh1981,Muh1985} on page \pageref{weaklim}? To be precise, to what
extent is a version of the limit operator theory possible replacing
$\P$-convergence ($\pto$) in Definition \ref{def_limop} of a limit
operator by the weaker $\sto$ convergence of operators introduced in
(\ref{eq_sto_def}) in Chapter \ref{chap_converg}?} This question is of
interest because weakening the concept of a limit operator widens the
class of operators that are rich, and so potentially widens the
applicability of the results in this text. Indications that versions of
at least some of our results of Chapter \ref{sec_ellp} may hold in this
context are: the results of Muhamadiev \cite{Muh1981,Muh1985} that we
discuss on page \pageref{weaklim}; and that it is the weaker $\sto$
convergence rather than $\pto$ which is required in the theorems in
Section \ref{sec:coll_comp_op_sp}.
\bigskip

{\it 3. Is $S(Y)$ inverse closed in $L(Y)$?}\quad We just mentioned
that $S(Y)$ is a ``one-sided version'' of $L(Y,\P)$ and we know from
\cite[Theorem 1.1.9]{RaRoSiBook} that $L(Y,\P)$ is inverse closed.
Looking back to problem 1, it seems a sensible question to ask
whether also $S(Y)$ is inverse closed. In Section
\ref{sec:algebraic} we were able to prove that if $A$ is invertible
in $S(Y)$ and $B=A+K\in S(Y)$ is a Montel perturbation that is
invertible in $L(Y)$ then $B^{-1}\in S(Y)$. This fact implies
inverse closedness of $S(Y)$ in cases where $Y=\hat Y$ and $P_n\in
K(Y)$ for all $n$.
\bigskip

{\it 4. Is the condition of existence of a predual $U^\pa$ and
preadjoint $A^\pa$ redundant in Theorem \ref{Prop_Fred_equiv}?}\quad
...and hence in Theorems \ref{prop_rich_invinf} and
\ref{thm_lind_wie} and Corollaries \ref{cor_combine2},
\ref{prop_W_selfadj}, \ref{cor_spessformula} and \ref{cor_oned_wie}?
Recall that we use the existence of predual and preadjoint to
conclude Fredholmness of $A_0=A|_{Y_0}$ from that of $A$ on
$Y^\infty$, including preservation of the index.
\bigskip

{\it 5. Is the hyperplane condition redundant in Theorem
\ref{thm_lind_wie}?}\quad  This possibility has been discussed in
the final section of \cite{LiWiener}. The answer is positive if the
following conjecture holds: $\ph(Y^p(U))=\ph(U)$, where $\ph(X)$
denotes the smallest positive integer $\kappa$ for which a Fredholm
operator of index $\kappa$ exists on the Banach space $X$ (with
$\ph(X)=0$ if all Fredholm operators on $X$ have index zero).
\bigskip

{\it 6. Is there a version of Theorem \ref{prop_FC_JFA} in arbitrary
dimensions $N\ge 1$?}\quad As we sketched briefly in the intro to
Theorem \ref{prop_FC_JFA}, the proof \cite{CWLi:JFA2007} of this
result consists of three steps. In two of these steps we use the
fact that $N=1$ and it is not clear to us whether and how this
condition could be removed. However, we are optimistic that the
result also holds for $N\ge 1$.
\bigskip

{\it 7. What does Theorem \ref{prop5.10} have to say about the
stability of classical approximation methods for operator equations
(e.g.\ finite section and discretisation methods)?} We have not
tackled this topic in this text, but some results in this
direction are in \cite{CW93,ChanZh97,MeierCW01,CW_NumMath} (and see
\cite{AnseloneSloan,RoSi,ProssdorfSilbermann,HaRoSi1,RaRoSi1998,RaRoSi2001,HaRoSi2,Li03,RaRoSiBook,LiRaRo04,LiBook,CWLi2,Roch_Memoir,SeidelSilbermann1,SeidelSilbermann2}).
\bigskip

\noindent And finally the classic:

{\it 8. Is the uniform boundedness condition in Theorem
\ref{prop_rich_invinf} (iii) also redundant for $p\in
(1,\infty)$?}\quad This question is as old as the first versions of
Theorem \ref{prop_rich_invinf} (iii). For $p\in\{0,1,\infty\}$ the
redundancy was shown in \cite{LiDiss,LiBook} and in this text. For
the remaining cases $p\in (1,\infty)$ the question remains open
(only for $A\in\W$ it is settled, by interpolation between $p=1$
and $p=\infty$). See
\cite[Section 3.3]{LiDiss} and \cite[Section 3.9]{LiBook} for a
little survey on this question and some humble attempts to tackle
it. Also note that the redundancy of this uniform boundedness condition
is sufficient (and in case of a normal operator $A$ also necessary)
for the set
\[
\bigcup_{B\in\opsp(A)} \sp(B)
\]
to be closed.


\backmatter

\newpage
\bibliographystyle{amsalpha}

\newpage
\printindex

\end{document}